\documentclass
[10pt,a4paper]
{article}
\usepackage{amsmath}
\usepackage{eucal}
\usepackage{amscd}
\usepackage{amssymb}
\usepackage{latexsym}
\usepackage{amsfonts}
\usepackage{makeidx}
\usepackage{graphicx}
\usepackage[greek,english,frenchb]{babel}

\makeindex

\newtheorem{rem}{Remark}[section]
 \makeatletter
    
    \@addtoreset{equation}{section}
  \makeatother

\allowdisplaybreaks[4]

\begin{document}
\title{From Euclid to Riemann and Beyond\footnote{Dedicated to the memory of Marcel Berger  (14 April 1927--15 October 2016).}\\
\vspace{0.5cm}
{\large -- How to describe the shape of the universe --}}
\author{Toshikazu Sunada}
%\dedicatory{Dedicated to the memory of Marcel Berger  (14 April 1927--15 October 2016).}
\date{}
\maketitle

%\selectlanguage{english}

The purpose of this essay is to trace the historical development of geometry while focusing on how we acquired mathematical tools for describing the ``shape of the universe." More specifically, our aim is to consider, without a claim to completeness, the origin of Riemannian geometry, which is indispensable to the description of the space of the universe as a ``generalized curved space."

But what is the meaning of ``shape of the universe"? The reader who has never encountered such an issue might say that this is a pointless question. It is surely hard to conceive of the universe as a geometric figure such as a plane or a sphere sitting in space for which we have vocabulary to describe its shape. For instance, we usually say that a plane is ``flat" and ``infinite," and a sphere is ``round" and ``finite." But in what way is it possible to make use of such phrases for the universe? Behind this inescapable question is the fact that the universe is not necessarily the ordinary 3D (3-dimensional) space where the traditional {\it synthetic geometry}---based on a property of parallels which turns out to underly the ``flatness" of space---is practiced. Indeed, as Einstein's\index{Einstein, Albert (1879--1955)} theory of {\it general relativity} (1915) claims, the universe is possibly ``curved" by gravitational effects. (To be exact, we need to handle 4D {\it curved space-times}; but for simplicity we do not take the ``time" into consideration, and hence treat the ``static" universe or the universe at any instant of time unless otherwise stated. We shall also disregard possible ``singularities" caused by ``black holes.") 

An obvious problem still remains to be grappled with, however. Even if we assent to the view that the universe is a sort of geometric figure, it is impossible for us to look out over the universe all at once because we are strictly confined in it. How can we tell the shape of the universe despite that? Before Albert Einstein\index{Einstein, Albert (1879--1955)} (1879--1955) created his theory, mathematics had already climbed such a height as to be capable to attack this issue. In this respect, Gauss\index{Gauss, Johann Carl Friedrich (1777--1855)} and Riemann\index{Riemann, Georg Friedrich Bernhard (1826--1866)} are the names we must, first and foremost, refer to as mathematicians who intensively investigated curved surfaces and spaces with the grand vision that their observations have opened up an entirely new horizon to {\it cosmology}. In particular, Riemann's work, which completely recast three thousand years of geometry executed in ``space as an a priori entity," played an absolutely decisive role when Einstein\index{Einstein, Albert (1879--1955)} established the theory of general relativity.

Gauss\index{Gauss, Johann Carl Friedrich (1777--1855)} and Riemann were, of course, not the first who were involved in cosmology. Throughout history, especially from ancient Greece to Renaissance Europe, mathematicians were, more often than not, astronomers at the same time, and hence the links between mathematics and cosmology are ancient, if not in the modern sense. Meanwhile, the venerable history of cosmology (and cosmogony) overlaps in large with the history of human thought, from a reflection on primitive religious concepts to an all-embracing understanding of the world order by dint of reason. It is therefore legitimate to lead off this essay with a rough sketch of philosophical and theological aspects of cosmology in the past, while especially focusing on the image of the universe held by scientists (see Koyr\'{e} \cite{Koy} for a detailed account).
In the course of our historical account, the reader will see how cosmology removed its religious guise through a long process of secularization and was finally established on a firm mathematical base. To be specific, Kepler\index{Kepler, Johannes (1571--1630)}, Galileo\index{Galileo Galilei (1564--1642)}, Fermat\index{Fermat, Pierre de (1601--1665)}, Descartes\index{Descartes, Ren\'{e} (1596--1650)}, Newton\index{Newton, Issac (1642--1726)}, Leibniz\index{Leibniz, Gottfried Wilhelm von (1646--1716)}, Euler\index{Euler, Leonhard (1707--1783)}, Lagrange\index{Lagrange, Joseph-Louis (1736--1813)}, and Laplace\index{Laplace, Pierre-Simon (1749--1827)} are on a short list of central figures who plowed directly or indirectly the way to the mathematization of cosmology.

The late 19th century occupies a special position in the history of mathematics. It was in this period that the autonomous progress of mathematics was getting apparent more than before. This is particularly the case after the notion of {\it set} was introduced by Cantor\index{Cantor, Georg (1845--1918)}. His theory---in concord with the theory of topological spaces\index{topological space}---allowed to bring in an entirely new concept of abstract space, with which one may talk not only about ({\it in}){\it finiteness} 
%or {\it infiniteness} 
of the universe in an {\it intrinsic manner}, 
an issue inherited since classical antiquity, but also about a global aspect of the universe despite that we human beings are confined to a very tiny and negligible planet in immeasurable space. In all these revolutions, Riemann's theory (in tandem with his embryonic research of {\it topology}\index{topology}) became encompassed in a broader, more adequate theory, and eventually led, with a wealth of new ideas and methodology, to modern geometry.

It is not too much to say that geometry as such (and mathematics in general) is   part of our cultural heritage because of the profound manifestations of the internal dynamism of human thought provoked by a ``sense of wonder" (the words Aristotle\index{Aristotle (384--322 BC)} used in the context of philosophy). Actually, its advances described in this essay may provide an indicator that reveals how human thought has been progressed over a long period of history.

\medskip
In this essay, I do not deliberately engage myself in the cutting-edge topics of differential geometry which, combined with topology, analysis, and algebraic geometry, have been highly cultivated since the latter half of the 20th century. I also do not touch on the {\it extrinsic} study of manifolds, i.e., the theory of configurations existing in space---no doubt an equally significant theme in modern differential geometry. In this sense, the bulk of my historical commentary might be quite a bit biased towards a narrow range of geometry. The reader interested in a history of geometry (and of mathematics overall) should consult Katz \cite{katz}.

\medskip

\noindent{\bf Acknowledgement}$~$ I wish to thank my colleague Jim Elwood whose suggestions were very helpful for the improvement of the first short draft of this essay. I am also grateful to Athanase Papadopoulos, Ken'ichi Ohshika, and Polly Wee Sy for the careful reading and for all their invaluable comments and suggestions to the enlarged version.

\section{Ancient models of the universe}\label{finiteor}

Since the inception of civilization that emerged in a number of far-flung places around the globe like ancient Egypt, Mesopotamia, ancient India, ancient China, and so on, mankind has struggled to understand the universe and especially how the world came to be as it is.\footnote{To know the birth and evolution of the universe is ``the problem of problems" at the present day as seen in the {\it Big Bang theory}, the most prevailing hypothesis of the birth.} Such attempts are seen in mythical tales about the birth of the world. ``Chaos", ``water", 
%``indistinguishable heaven and earth" 
and the like, were thought to be the fundamental entities in its beginning that was to grow gradually into the present state. For example, the epic of {\it Atrahasis} written about 1800 BCE contains a creation myth about the Mesopotamian gods {\it Enki} (god of water), {\it Anu} (god of sky), and {\it Enlil} (god of wind).
% who are depicted as superhumanly powerful humanoids. 
%For example, the {\it Old Testament}
%, compiled between the 10th century BCE and 400 BCE 
%says, ``In the beginning, God created the heavens and the earth." 
The Chinese myth in the {\it Three Five Historic Records} (the 3rd century CE) tells us that the universe in the beginning was like a big egg, inside of which was darkness chaos.\footnote{The term ``chaos" (\textgreek{q'aos}) is rooted in the poem {\it Theogony}, a major source on Greek mythology, by Hesiod ({\it ca.}~\!750 BCE--{\it ca.}~\!650 BCE). The epic poet describes Chaos as the primeval emptiness of the universe.} 

The question of the origin is tied to the question of future. Hindu cosmology has a unique feature in this respect. In contrast to the didactics in monotheistic Christianity describing the end of the world as a single event (with the Last Judgement) in history, the {\it Rig veda}, one of the oldest extant texts in Indo-European language composed between 1500 BCE--1200 BCE, alleges that our cosmos experiences a creation-destruction cycle almost endlessly (the view celebrated much later by Nietzsche). Furthermore, some literature (e.g. the {\it Bhagavata Purana} composed between the 8th and the 10th century CE or as early as the 6th century CE) mentions the ``multiverse" (infinitely many universes), which resuscitated as the modern astronomical theory of the {\it parallel universes}.

In the meantime, thinkers in colonial towns in Asia Minor, Magna Graecia, and mainland Greece, cultivated a love for systematizing phenomena on a rational basis, 
%(thereby christened {\it philosophers} (\textgreek{fil'osofos}), lovers of wisdom), 
as opposed to supernatural explanations typified by mythology  and folklore in which the cult of Olympian gods and goddesses are wrapped. Many of them were not only concerned with fundamental issues arising from everyday life, as represented by Socrates ({\it ca.}~\!470 BCE--399 BCE),
but also labored to mathematically understand multifarious phenomena and to construct an orderly system. They appreciated purity, universality, a certainty and an elegance of mathematics, the characteristics that all other forms of knowledge do not possess. Legend has it that Plato\index{Plato (c. 428--c. 348 BC)} ({\it ca.}~\!428 BCE--{\it ca.}~\!348 BCE)
%, Socrates' prized student, 
engraved the phrase ``Let no one ignorant of geometry enter here" at the entrance of the {\it Academeia}\index{Academeia} (\textgreek{>Akadhm'ia}) he founded in {\it ca}.~\!387 BCE in an outskirt of Athens. Whether or not this is historically real, Academeia\index{Academeia} indubitably put great emphasis on mathematics as 
%a branch or 
a prerequisite of philosophy.\footnote{Among Plato's\index{Plato (c. 428--c. 348 BC)} thirty-five extant dialogues, there are quite a few in which the characters (Socrates in particular) discuss mathematical knowledge in one form or another; say, {\it Hippias major}, {\it Meno}, {\it Parmenides}, {\it Theaetetus}, {\it Republic}, {\it Laws}, {\it Timaeus}\index{Timaeus@Timaeus (dialogue)}, {\it Philebus}. In the {\it Republic}, Book VII, Socrates says, ``Geometry is fully intellectual and a study leading to the good because it deals with absolute, unchanging truths, so that it should be part of education."}

Eventually, Greek philosophers began to speculate about the structure of the universe by deploying geometric apparatus.
% at their disposal. 
Exemplary is the {\it spherical model}\index{spherical model}, with the earth at the center, proposed by Plato\index{Plato (c. 428--c. 348 BC)} himself, and his two former students Eudoxus\index{Eudoxus (c. 408--c. 355 BC)} of Cnidus ({\it ca.}~\!408 BCE--{\it ca.}~\!355 BCE) and Aristotle\index{Aristotle (384--322 BC)} from Stagira (384 BCE--322 BCE).\footnote{Pythagoras\index{Pythagoras (c. 580--c. 500 BC)} of Samos ({\it ca.}~\!580 BCE--{\it ca.}~\!500 BCE) was the first to declare that the earth is a sphere, and that the universe has a soul and intelligence. Plato\index{Plato (c. 428--c. 348 BC)} was a devout {\it Pythagorean}\index{Pythagorean}, originally meaning the member of a mathematico-religious community created by Pythagoras\index{Pythagoras (c. 580--c. 500 BC)} ({\it ca}.~\!530 BCE) in Croton, Southern Italy.} In his dialogue {\it Timaeus}\index{Timaeus@Timaeus (dialogue)}, Plato\index{Plato (c. 428--c. 348 BC)} explored cosmogonical issues, and deliberated the nature of the physical world and human beings. He referred to the {\it demiurge} (\textgreek{dhmiourg'os}) as the creator of the world who chose a {\it round sphere} as the most appropriate shape that embraces within itself all the shapes there are. Meanwhile, Eudoxus\index{Eudoxus (c. 408--c. 355 BC)} proposed a sophisticated system of homocentric spheres rotating about different axes through the center of the Earth, as an answer to his mentor who set the question how to reduce the apparent motions of heavenly bodies to uniform circular motions. This was stated in the treatise {\it On velocities}---now lost, but Aristotle\index{Aristotle (384--322 BC)} knew about it. Eudoxus\index{Eudoxus (c. 408--c. 355 BC)} may have regarded his system simply as an abstract geometrical model; Aristotle\index{Aristotle (384--322 BC)} took it to be a description of the real world, and organized it into a kind of fixed hierarchy, conjoining with his metaphysical principle ({\it Metaphysics}, XII). 

\begin{figure}[htbp]
\vspace{-0.4cm}
\begin{center}
\hspace{1cm}\includegraphics[width=.58\linewidth]
{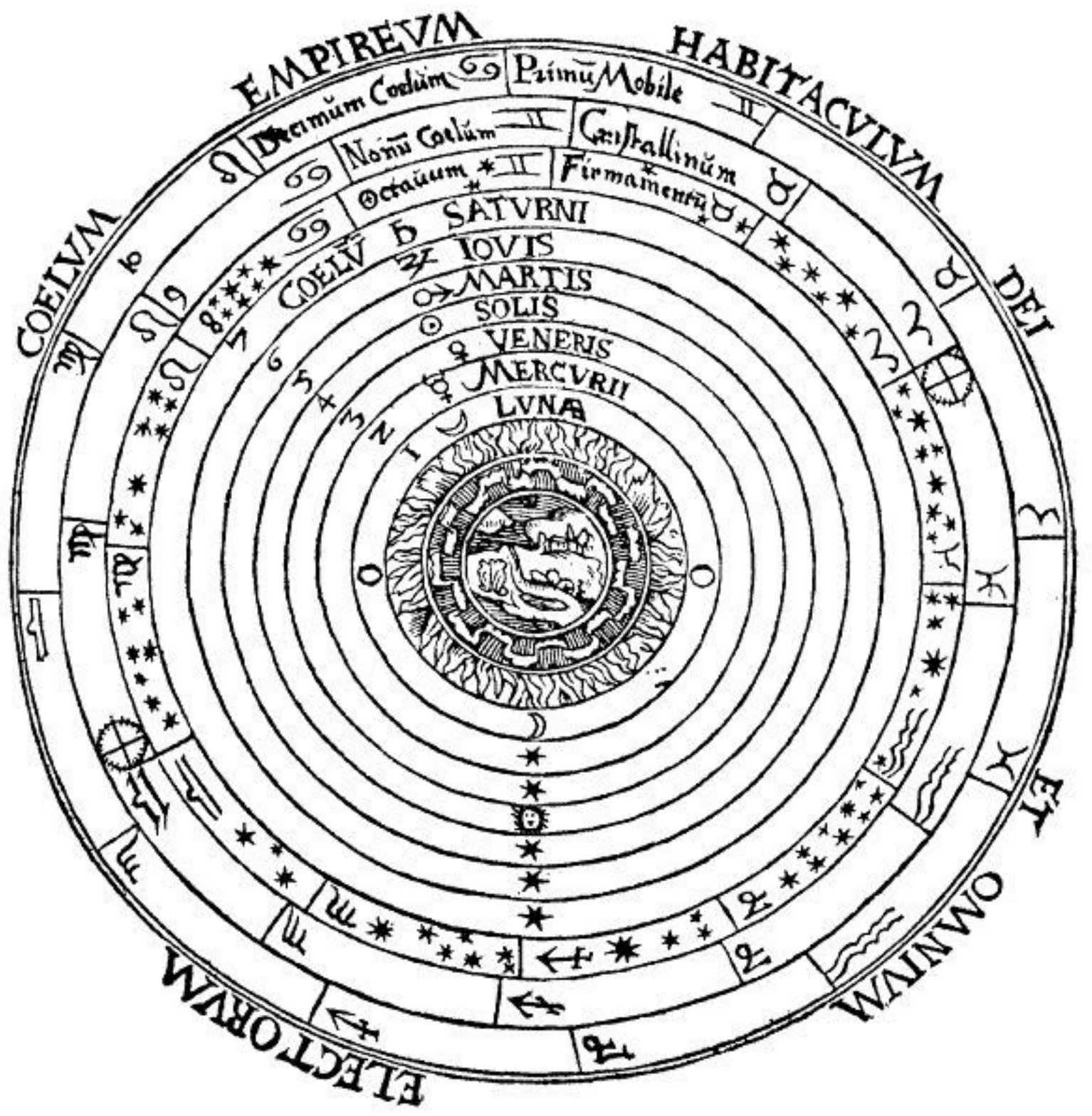}
\end{center}
\vspace{-0.8cm}
\caption{Aristotle's model (P.
Apianus, {\it Cosmographicus liber}, 1524)}
\label{fig:aristotle}
\vspace{-0.2cm}
\end{figure}

Aristotle's\index{Aristotle (384--322 BC)} spherical model\index{spherical model} (Fig.~\!\ref{fig:aristotle}) was refined later as the sophisticated geocentric theory\index{geocentric theory} by the Alexandrian astronomer Ptolemy\index{Ptolemy (c. 100--c. 170)} ({\it ca.}~\!100 CE--{\it ca.}~\!170 CE). His vision of the universe was set forth in his mathematico-astronomical treatise {\it Almagest}\index{almagest@Almagest}, and had been accepted for more than twelve hundred years in Western Europe and the Islamic world until Copernicus'\index{Copernicus, Nicolaus (1473--1543)} heliocentric theory\index{heliocentric theory} emerged (Sect.~\!\ref{sec:universe}) because his theory succeeded, up to a point, in describing the apparent motions of the sun, moon, and planets.\footnote{\label{fn:Alma}The title ``Almagest"\index{almagest@Almagest} was derived from the Arabic name meaning ``greatest." The original Greek title is {\it Mathematike Syntaxis} (\textgreek{Majhmatik`h S'untaxis}). It was rendered into Latin by Gerard\index{Gerard (c. 1114--1187)} of Cremona ({\it ca}.~\!1114--1187) from the Arabic version found in Toledo (1175), and became the most widely known in Western Europe before the Renaissance.}  What is notable in his system is the use of numerous {\it epicycles}\index{epicycle} (\textgreek{>ep'ikuklos}), where an epicycle is a small circle along which a planet is assumed to move, while each epicycle in turn moves along a larger circle ({\it deferent}).\footnote{The reasonable accuracy of Ptolemy's system\index{Ptolemy (c. 100--c. 170)} results from that an ``almost periodic" motion (in the sense of H. A. Bohr)---a presumable nature of planetary motions---can be represented to any desired degree of approximation by a superposition of circular motions.} This schema---a coinage of the Greek doctrine---was first set out by Apollonius of Perga\index{Apollonius (c. 262--c. 190 BC)} ({\it ca}.~\!262 BCE--{\it ca}.~\!190 BCE), and developed further by Hipparchus of Nicaea\index{Hipparchus (c. 190--c. 120 BC)} ({\it ca}.~\!190 BCE--{\it ca}.~\!120 BCE).

At all events, an intriguing (and arguable) feature of their model is the idea of unreachable ``outermost sphere." In Aristotle's\index{Aristotle (384--322 BC)} concept ({\it Physics}, VIII, 6), it is the domain of the Prime Mover (\textgreek{t`o pr'wth >ak'inhton}), a variant  of Plato's\index{Plato (c. 428--c. 348 BC)} notion in his cosmological argument unfolded in the {\it Timaeus}\index{Timaeus@Timaeus (dialogue)}, which caused the outermost sphere to rotate at a constant angular velocity.

We shall come back to the spherical model\index{spherical model} in Sect.~\!\ref{sec:universe} after discussing a relevant issue, and recount at some length how this peculiar model had an effect on philosophical and religious aspects of cosmology.

\section{What is infinity?}\label{sec:whatinfinity}
Besides the kinematical nature of celestial bodies, what seems to lie in the background of Aristotle's\index{Aristotle (384--322 BC)} view about the universe is his philosophical thought about {\it infinity}. Actually the issue of infinity was a favorite subject for Greek philosophers, dating back to the pre-Socratic period.\footnote{Anaximander of Miletus\index{Anaximander (c. 610--c. 546 BC)} ({\it ca.}~\!610 BCE--{\it ca.}~\!546 BCE) was the first who contemplated about infinity, and employed the word {\it apeiron} (\textgreek{>'apeiron} meaning ``unlimited") to explain all natural phenomena in the world. Anaxagoras of Clazomenae\index{Anaxagoras (c. 510--c. 428 BC)} ({\it ca.}~\!510 BCE--{\it ca.}~\!428 BCE) wrote the book {\it About Nature}, in which he says, ``All things were together, infinite in number." Assuming the infinite divisibility of matter, he avers, ``There is no smallest among the small and no largest among the large, but always something still smaller and something still larger."}

Deliberating over his predecessors'~\!vision, Aristotle\index{Aristotle (384--322 BC)} distinguished between {\it actual}\index{actual infinity} and {\it potential} infinity\index{potential infinity}. Briefly speaking, actual infinity is ``things" that are completed and definite, thereby being {\it transcendental} in nature, while potential infinity is ``things" that continue without terminating; more and more elements can be always added, but with no recognizable ending point, thus being what can be somehow corroborated within the scope of the capacity of human deed or thought. For a variety of possible reasons, Aristotle rejects actual infinity\index{actual infinity}, claiming that only potential infinity\index{potential infinity} exists, and captures space as something infinitely divisible into parts that are again infinitely divisible, and so on ({\it Physics}, III). His doctrine, which had satisfied nearly all scholars for a long time,\footnote{G. W. F. Hegel, a pivotal figure of German {\it idealism}, defended Aristotle's\index{Aristotle (384--322 BC)} perspective on the infinite in his {\it Wissenschaft der Logik} (1812--1816), though he used the terms ``true (absolute) infinity" and ``spurious infinity ({\it schlecht Unendlichkeit})" instead.} was employed by himself to find a way out of {\it Zeno's paradoxes}\index{Zeno's paradoxes}, especially the paradox (\textgreek{par'adoxos}) of Achilles and the Tortoise.\footnote{Zeno\index{Zeno  (c. 490--c. 430 BC)} of Elea ({\it ca.}~\!490 BCE--{\it ca.}~\!430 BCE) was a student of Parmenides\index{Parmenides (c. 515--c. 450 BC)} ({\it ca.}~\!515 BCE--{\it ca.}~\!450 BCE) who contended that the true reality is absolutely unitary, unchanging, eternal, ``the one."  To vindicate his teacher's tenet, Zeno offered the four arguable paradoxes ``The Dichotomy," ``Achilles and the Tortoise," ``The Arrow," and ``The Stadium."

During the Warring States period (476 BCE--221 BCE) in China, the {\it School of Names} cultivated a philosophy similar to Parmenides'\index{Parmenides (c. 515--c. 450 BC)}. Hui Shi ({\it ca}.~\!380 BCE--{\it ca}.~\!305 BCE) belonging to this school says, ``Ultimate greatness has no exterior, ultimate smallness has no interior."} This, in a rehashed form, says, ``{\it The fastest runner [Achilles] in a race can never overtake the slowest [tortoise], because the pursuer must first get to the point whence the pursued started, so that the slowest must always hold a lead}" ({\it Physics}, VI; 9, 239b15). To quote Aristotle's\index{Aristotle (384--322 BC)} reaction in response to the quibble, Zeno's argument exploits an ambiguity in the nature of `infinity' because Zeno seems to insist that Achilles cannot {\it complete} an infinite number of his actions (getting the point where the tortoise was); that is, `complete' is the word for actual infinity\index{actual infinity}, but `infinite number of actions' is the phrase for potential infinity\index{potential infinity}.

Aristotle's\index{Aristotle (384--322 BC)} resolute rejection of actual infinity\index{actual infinity} seems in part to come from the fact that mathematicians of those days had no adequate manner to treat continuous magnitude and could do, all in all, quite well without actual infinity\index{actual infinity}. This is distinctively  seen in the statement in Euclid's\index{Euclid ($4^{\rm th}$-$3^{\rm rd}$ c. BC)}
{\it Elements}\index{Elements}, Book IX, Prop.~\!20, about the infinitude of prime numbers, which deftly asserts, ``Prime numbers are more than any assigned multitude of prime numbers."\footnote{Euclid\index{Euclid ($4^{\rm th}$-$3^{\rm rd}$ c. BC)}
does not seem to entirely refrain from the use of non-potential infinity\index{potential infinity}. In the {\it Elements}\index{Elements}, Book X, Def.~\!3, he says, ``there exist straight lines {\it infinite in multitude} which are commensurable with a given one" (\cite{euclid}, Book X, p.10).}

On the other hand, it appears that Archimedes\index{Archimedes (c. 287--c. 212 BC)} of Syracuse ({\it ca.}~\!287 BCE--{\it ca.}~\!212 BCE) had a prescient view of infinity\index{actual infinity}, as adumbrated in the {\it Archimedes Palimpsest},\footnote{\label{fn:heiberg}In 1906, J. L. Heiberg, the leading authority on Archimedes, confirmed that the palimpsest, overwritten with a Christian religious text by 13th-century monks, included the {\it Method of Mechanical Theorems}, one of Archimedes' lost works by that time. Our knowledge of Archimedes was greatly enriched by this fabulous discovery.} a 10th-century Byzantine Greek copy housed at the Metochion of the Holy Sepulcher in Jerusalem. To our astonishment, in the 174-pages text (specifically in the {\it Method},  Prop.~\!14), he duly compared two infinite collections of certain geometric objects by means of a {\it one-to-one correspondence}\index{one-to-one correspondence} (henceforth OTOC); see \cite{netz}. This is surely related to the concept of actual infinity\index{actual infinity} that has been revived in 19th century (Sect.~\!\ref{subsec:transfinite}).

Archimedes\index{Archimedes (c. 287--c. 212 BC)} was also a master of the {\it method of exhaustion}\index{method of exhaustion} (or the method of double contradiction) which originated, however incomplete, with Antiphon the Sophist\index{Antiphon (c. 480--c. 411 BC)} ({\it ca.}~\!480 BCE--{\it ca}.~\!411 BCE) and Bryson of Heraclea\index{Bryson (born in c. 450 BC)} (born in {\it ca}.~\!450 BCE), and exploited by Eudoxus\index{Eudoxus (c. 408--c. 355 BC)} to avoid flaws that may happen when we treat infinity in a naive way. Such a flaw is found in the claim by the two originators; they contended that it is possible to construct, with compass and straightedge, a square with the same area as a given circle $C$. Their argument (criticized roundly by Aristotle\index{Aristotle (384--322 BC)} in the {\it Posterior Analytics} I,~\!9,~\!75b40) is as follows: From the correct fact that such a construction is possible for a given polygon (Euclid's\index{Euclid ($4^{\rm th}$-$3^{\rm rd}$ c. BC)}
{\it Elements}\index{Elements}, Book II, Prop.~\!14), they elicited the incorrect consequence that the same is true for $C$ on the grounds that the regular polygon with $2^n$ edges inscribed in $C$ {\it eventually coincides with} $C$ as we let $n$ increase endlessly (needless to say, a passage to a limit does not necessarily preserve the given properties).\footnote{The approach taken by Antiphon\index{Antiphon (c. 480--c. 411 BC)} and Bryson\index{Bryson (born in c. 450 BC)}, however incorrect, was appropriately adopted by Archimedes\index{Archimedes (c. 287--c. 212 BC)}, who obtained $223/71<\pi<22/7$ 
using two regular polygons of 96 sides inscribed and circumscribed to a circle
({\it Measurement of the Circle}).}

Incidentally, the problem Antiphon\index{Antiphon (c. 480--c. 411 BC)} and Bryson\index{Bryson (born in c. 450 BC)} challenged is designated ``squaring the circle"\index{squaring the circle}, one of the three big problems on constructions in Greek geometry; the other problems are ``doubling the cube"\index{doubling the cube} and ``trisecting the angle." Anaxagoras\index{Anaxagoras (c. 510--c. 428 BC)} is the first who worked on squaring the circle, while Hippocrates of Chios\index{Hippocrates (c. 470--c. 410 BC)} ({\it ca.}~\!470 BCE--{\it ca.}~\!410 BCE) squared a {\it lune}. This held out hope to square the circle for a time, but all attempts met with failure. It was in 1882 that this problem came to a conclusion; Ferdinand von Lindemann\index{Lindemann, Ferdinand von (1852--1939)} (1852--1939) proved the impossibility of such a construction by showing that $\pi$ $=$ $3.141592\cdots$ is a {\it transcendental number}; i.e., $\pi$ is not a solution of an algebraic equation with integral coefficients (1882). 
Doubling the cube and trisecting the angle are also impossible as  P. L. Wantzel showed (1837).

In the background of the method of exhaustion\index{method of exhaustion} is the premise that a given quantity $\alpha$ (not necessarily numerical) can be made smaller than another one given beforehand by successively halving $A$. This, in modern terms, says that for any quantity $\beta$, there exists a natural number $n$ with $\alpha/2^n<\beta$; thus the premise goes along well with Anaxagoras' view of ``limitless smallness"\footnote{The word ``exhaustion" for this reasoning was first used by Gr\'{e}goire de Saint-Vincent (1584--1667); {\it Opus Geometricum Quadraturae Circuli et Sectionum Conti}, 1647, p.~\!739.} and is, though restricted to very special situations, regarded as a harbinger of the {\it predicate calculus} in modern logic---the branch of logic that deliberately deals with quantified statements such as ``there exists an $x$ such that $\cdots$" or ``for any $x$, $\cdots$"--- and particularly the $\epsilon$-$\delta$ argument\index{epsilon@$\epsilon$-$\delta$ argument} invented for the rigorous treatment of {\it limits} in the 19th century which adequately avoids endless processes (Remark \ref{rem:function}). A variant of this premise is\!: ``Given two quantities, one can find a multiple of either which will exceed the other." This is what we call the {\it Axiom of Archimedes}\index{Axiom of Archimedes} since it was explicitly formulated by Archimedes\index{Archimedes (c. 287--c. 212 BC)} in his work {\it On the Sphere and Cylinder}.\footnote{The {\it Elements}\index{Elements}, Book V, Def.~\!4 states, ``Magnitudes (\textgreek{m'egejos}) are said to have a ratio to one another {\it which can, when multiplied, exceed one another}." A similar statement appears in Aristotle's\index{Aristotle (384--322 BC)} {\it Physics} VIII, 10, 266b 2. The name ``Axiom of Archimedes" was given by O. Stolz in 1883, and was adopt by David Hilbert\index{Hilbert, David (1862--1943)} (1862--1943) in a modern treatment of Euclidean geometry\index{Euclidean geometry} (Remark \ref{rem:list} (2)).} Eudoxus\index{Eudoxus (c. 408--c. 355 BC)} and Archimedes\index{Archimedes (c. 287--c. 212 BC)} combined these premises with {\it reductio ad absurdum}\index{reductio ad absurdum} (proof by contradiction\index{proof by contradiction}) in a judicious manner,\footnote{Reductio ad absurdum is a reasoning in the Eleatics philosophy initiated by Parmenides\index{Parmenides (c. 515--c. 450 BC)}. Incidentally, the Chinese word for ``contradiction" is ``m\'{a}od\`{u}n," literally ``spear-shield," stemming from an anecdote in the {\it Han Feizi}, an ancient Chinese text attributed to the political philosopher Han Fei ({\it ca}.~\!280 BCE--{\it ca}.~\!233 BCE). The story goes as follows. A dealer of spears and shields advertised that one of his spears could pierce any shield, and at the same time said that one of his shields could defend from all spear attacks. Then one customer queried the dealer what happens when the shield and spear he mentioned would be used in a fight.} and established various results on area and volume by highly sophisticated arguments. Noteworthy is that, in his computations of ratios of areas or volumes of two figures, Archimedes\index{Archimedes (c. 287--c. 212 BC)} made use of {\it infinitesimals} in a way similar to the one in integral calculus at the early stage (e.g. {\it Quadrature of the Parabola}).

Some 1500 years later, the issue of infinity was examined by two scholars. Thomas Bradwardine\index{Bradwardine, Thomas (c. 1290--c. 1349)} ({\it ca}.~\!1290--{\it ca}.~\!1349)---a key figure in the {\it Oxford Calculators}\index{Oxford Calculators} (a group of mathematics-oriented thinkers associated with Merton College)---employed the principle of OTOC\index{one-to-one correspondence} in discussing the aspects of infinite. His observation is construed, in modern terms, as ``an infinite subset could {\it be equal} to its proper subset" (\!{\it Tractatus de continuo}, 1328--1335). This is a polemic work directed against atomistic thinkers in his time such as  N. d'Autr\'{e}court who, following the classical atomic concept, considered that matter and space were all made up of indivisible atoms, as opposed to Aristotle's\index{Aristotle (384--322 BC)} doctrine. After a while, Nicole Oresme\index{Oresme, Nicole (1320--1382)} (1320--1382), a significant scholar of the later Middle Ages, elaborated the thought that the collection of odd numbers is {\it not smaller than} that of natural numbers because it is possible to count the odd numbers by the natural numbers (\!{\it Physics Commentary}, around 1345).\footnote{Oresme was in favor of the spherical model\index{spherical model} of the universe, but took a noncommittal attitude regarding the Aristotelian theory of the stationary Earth and a rotating sphere of the fixed stars ({\it Livre du ciel et du monde} and {\it Questiones super De celo, Questiones de spera}).}

Meanwhile, the magic of infinity has bewildered Galileo Galilei\index{Galileo Galilei (1564--1642)} (1564--1642). In his {\it Discorsi e Dimostrazioni Matematiche Intorno a Due Nuove Scienze} (1638), written during house arrest as the result of the Inquisition, he says, in a similar vein to Oresme's\index{Oresme, Nicole (1320--1382)} claim, ``Even though the {\it number} of squares should be less than that of all natural numbers, there appears to exist as many squares as natural numbers because any natural number corresponds to its square, and any square corresponds to its square root" ({\it Galileo's paradox}\index{Galileo's paradox}).

\begin{rem}\label{rem:naturalnumber}{\rm {\small 
(1) The principle of OTOC\index{one-to-one correspondence} seems to have its roots deep in a fundamental faculty of human brain. Thinking  back on how the species acquired ``natural numbers" when they did not yet have any clue about numerals to count things, we are led to the speculation that they relied on the principle of OTOC. More specifically, ancient people are supposed to check whether their cattle put out to pasture returned safely to their shed, without any loss, by drawing on OTOC\index{one-to-one correspondence} with some identical things such as sticks of twigs prepared in advance. They extended the same manner in counting things in various aggregations. This experience over many generations presumably made people aware that there is ``something in common" behind, not depending on things, and they eventually got a way of identifying as an entity all aggregations among which there are OTOCs\index{one-to-one correspondence}. This entity is nothing but a natural number.\footnote{A. N. Whitehead said, ``The first man who noticed the analogy between a group of seven fishes and a group of seven days made a notable advance in the history of thought" ({\it Science and the Modern World}, 1929).}

Once reached this stage, it did not take much time for people to come up with assigning symbols 
%(as I, II, III,$\ldots$ in Roman numerals) 
and names, 
%(as one, two, three,$\ldots$ in English), 
and evolving primitive arithmetic through empirical fumbling, especially the {\it division algorithm}\index{division algorithm}, if not in the way we express it today. In fact, numerical notation systems allowing us to represent numbers with a few symbols evolved from the division algorithm\index{division algorithm}. For instance, lining up rod-like symbols as $|,~||,~|||,~||||,~|||||,~||||||,\cdots$ is the most primitive way to represent numerals, as seen in the first few of the Babylonian 
%Roman 
and Chinese numerals.
In the {\it quinary system}, 
%(the numerical system with base $5$), 
we split up a given $|||\cdots |$ into several groups each of which consists of five rods $|||||$ together with a group consisting of rods less than five (this is an infantile form of the division algorithm\index{division algorithm}). We then replace each $|||||$ by a new symbol, e.g. $\top$; thus, for example, $|||||||||||| =|||||~~ |||||~~ ||$ is exhibited as $\top ~\top~||$. We do the same for the symbol $\top$, and continue this procedure
(this story is, of course, considerably simplified).

\smallskip

(2) The {\it decimal system}, contrived in India between the 1st and 4th centuries, made it possible to express every number by finitely many symbols $0,1,\ldots,9$.\footnote{The explicit use of zero as a symbol more than a placeholder was made much later. In India, zero was initially represented by a point. The first record of the use of the symbol ``0" is dated in 876 (inscribed on a stone at the Chaturbhuja Temple). The historical process leading to the term "zero" is as follows: ``sunya" in Hindu meaning {\it emptiness} $\rightarrow$ ``sifr" (cipher) in Arabic $\rightarrow$ ``zephirum" in Latin used in 1202 by Fibonacci\index{Fibonacci, Leonardo (c. 1170 --c. 1250)} $\rightarrow$ ``zero" in Italian ({\it ca}.~\!1600).} This epochal format
%, which was pushed further by a more effective use of the division algorithm\index{division algorithm}, 
amplified arithmetic considerably, and was brought through the Islamic world to Europe in the 10th century. The names to be mentioned  are Muhammad ibn M$\overline{\rm u}$s$\overline{\rm a}$ al-Khw$\overline{\rm a}$rizm\=\i\index{alk@alk@al-Khw$\overline{\rm a}$rizm\=\i}~({\it ca.}~\!780 CE--{\it ca.}~\!850 CE), Gerbert of Aurillac\index{Gerbert (c. 946--1003)} ({\it ca}.~\!946 CE--1003),  and Leonardo Fibonacci\index{Fibonacci, Leonardo (c. 1170 --c. 1250)} 
%(Leonardo of Pisa, 
({\it ca}.~\!1170 --{\it ca}~\!1250). In {\it ca}.~\!825 CE, al-Khw$\overline{\rm a}$rizm\=\i~wrote a treatise on the decimal system. His Arabic text (now lost) was  translated into Latin with the title {\it Algoritmi de numero Indorum}, 
%in the 12th century, 
most likely by Adelard of Bath\index{Adelard (c. 1080--c. 1152)} ({\it ca}.~\!1080--{\it ca}.~\!1152). Gerbert\index{Gerbert (c. 946--1003)} is said to be the first to introduce the decimal system in Europe (probably without the numeral zero). Fibonacci\index{Fibonacci, Leonardo (c. 1170 --c. 1250)}, best known by the sequence with his name, travelled extensively around the Mediterranean coast, and assimilated plenty of knowledge including Islamic mathematics. His {\it Liber Abaci} (1202) popularized the decimal system in Europe. 

The decimal notation---practically convenient and theoretically being of avail---is firmly planted in our brain as a mental image of numbers.
\hfill$\Box$

} 
}
\end{rem}

\section{Is the universe infinite or finite?}\label{sec:universe}

Let us turn to the issue of the universe. Were Aristotle's\index{Aristotle (384--322 BC)} spherical model\index{spherical model} correct, there would be its ``outside." So the question at once arises as to what the outside means after all. Is it something substantial or just speculative fabrication? Defending this ridiculous image of the universe, Aristotle\index{Aristotle (384--322 BC)} explained away by saying, ``there is neither place, nor void outside the heavens by reason that the heaven does not exist inside another thing" ({\it On the Heavens}, I, 9). 

Aristotle's\index{Aristotle (384--322 BC)} model continued to be employed even in the time when the Christian tenet dominated the scholarly world in Europe. Especially, 
influenced by Thomas Aquinas (1225--1274), who incorporated extensive Aristotelian philosophy throughout his own theology, the scientific substratum in Christianity was synthesized with the Aristotelian physics.\footnote{In his masterwork {\it Summa Theologiae}, Aquinas asserts that God is perfect, complete, and embodies actual infinity\index{actual infinity}. To reconcile his assertion with Aristotle's\index{Aristotle (384--322 BC)} doctrine, he says, ``God is an infinite that has bounds."} After a while, the Florentine poet Dante Alighieri\index{Dante, Alighieri (1265--1321)} (1265--1321) reinforced Aristotle's\index{Aristotle (384--322 BC)} idea of ``outside" in his unfinished work {\it Convivio}. He says, ``in the supreme edifice of the universe, all the world is included, and beyond which is nothing; and it is not in space, but was formed solely in the Primal Mind." He further offered an imaginal vision of an intriguing macrocosm in order to turn down the queer consequence of the Biblical concept of ascent to Heaven and descent to Hell which connotes that Hell is the center of the spherical universe (Sect.~\!\ref{subsec:topodes}).

An exception is Nicolaus Cusanus\index{Cusanus, Nicolaus (1401--1464)} (Nicholas of Cusa, 1401--1464),
a first-class scientist of his time. He alleges, refuting the prevalent outlook of the world, ``The universe is not finite in the sense of physically unboundedness since there is no special center in the universe, and hence the outermost sphere cannot be a boundary." At the same time, he argued finiteness of the universe, by which he meant to say ``privatively infinite" because ``the world cannot be conceived of as finite, albeit it is not infinite" ({\it De Docta ignorantia}; 1440). This rather contradictory dictum elaborated in the ground-breaking book is a consequence of his analysis of conceivability and a parallelism between the universe and God. 

A hundred years later, a dramatic turnabout took place. The heliocentric theory\index{heliocentric theory} proposed by Nicolaus Copernicus\index{Copernicus, Nicolaus (1473--1543)} (1473--1543) in his {\it De revolutionibus orbium coelestium} brings out the apparent retrograde motion of planets better than Ptolemy's\index{Ptolemy (c. 100--c. 170)} theory.
% though he still called for 34 epicycles\index{epicycle}. 
He already got his ideas---relying largely on Arabic astronomy typified by al-Batt$\overline{\rm a}$n\=\i
%\index{alb@
\index{al-Batt$\overline{\rm a}$n\=\i~(c. 858--929)}\footnote{\label{fn:al-Batt}His work {\it Kit$\overline{\rm a}$b az-Z\=\i j} (Book of Astronomical Tables), translated into Latin as {\it De Motu Stellarum} by Plato Tiburtinus in 1116, was quoted by Tycho Brahe\index{Tycho Brahe (1546--1601)}, Kepler\index{Kepler, Johannes (1571--1630)} and Galileo\index{Galileo Galilei (1564--1642)}.} ({\it ca}.~\!858 CE--929 CE)---some time before 1541 ({\it Commentariolus}), but resisted, in spite of his friends' persuasion, to make his theory public since he was afraid of being a target of contempt. It was in 1543, just before his death, that his work (with a preface which puts the accent on the hypothetical nature of the contents) was brought out.

Now, what did Copernicus\index{Copernicus, Nicolaus (1473--1543)} think about the size of the universe? His model of the universe is spherical with the outermost consisting of motionless, fixed stars; thus being not much different from Aristotle's\index{Aristotle (384--322 BC)} model in this respect. Meanwhile, Giordano Bruno\index{Bruno, Giordano (1548--1600)} (1548--1600) argued against the outermost sphere (while accepting Copernicus' theory), reasoning that the infinite power of God would not have produced a finite creation. He was quite explicit in his belief that there is no special place as the center of the universe ({\it De l'infinito universo et mondi}, 1584); thus indicating that the universe is endless, limitless, and {\it homogeneous}\index{homogeneous}. His outspoken views, containing  Hermetic elements, scandalized Catholics and Protestants alike.\footnote{Some modern scholars consider him as a magus, not a pioneer of science because {\it Hermeticism} is an ancient spiritual and magical tradition. Others, however, have argued that the magical world-view in Hermeticism was a necessary precursor to the {\it Scientific Revolution}.} Consequently, he was condemned as an impenitent heretic and eventually burned alive at the stake after 8 years' solitary confinement. 

Bruno's view was partly shared by his contemporary Thomas Digges\index{Digges, Thomas (1546--1595)} (1546--1595). He was the first to expound the Copernican system in English. In an appendix to a new edition of his father's book {\it A Prognostication everlasting}, Digges discarded Copernicus'\index{Copernicus, Nicolaus (1473--1543)} notion of a fixed shell of immovable stars, presuming infinitely many stars at varying distances (1576)---a more tenable reasoning than Bruno's. What should be thought over here is that the religious atmosphere in England in his time was different from that in the Continent because of the English Reformation that started in the reign of Henry VIII.. 

\medskip

Interrupting the chronological account, we shall go back to ancient Greece, where we find forerunners of Cusanus\index{Cusanus, Nicolaus (1401--1464)}, Copernicus\index{Copernicus, Nicolaus (1473--1543)}, and Bruno\index{Bruno, Giordano (1548--1600)}. Among them, Archytas\index{Archytas (428--347 BC)} of Tarentum (428 BCE--347 BCE) is a precursor of Bruno\index{Bruno, Giordano (1548--1600)}. He inferred that any place in our space looks the same, and hence no boundary can exist; otherwise we have a completely different sight at a boundary point. He thus concluded infiniteness of the universe.\footnote{Anaximander\index{Anaximander (c. 610--c. 546 BC)}, belonging to the previous generation, conceived a mechanical model of the world based on a non-mythological explanatory hypothesis, and alleged that the Earth has a disk shape and is floating very still in the center of the infinite, not supported by anything, while Anaxagoras\index{Anaxagoras (c. 510--c. 428 BC)} portrayed the sun as a {\it mass of blazing metal}, and postulated that Mind (\textgreek{no\~us}, Nous) was the initiating and governing principle of the cosmos (\textgreek{k'osmos}).} The Pythagorean\index{Pythagorean} Philolaus\index{Philolaus (c. 470--c. 385 BC)} ({\it ca.}~\!470 BCE--{\it ca.}~\!385 BCE) relinquished the geocentric model, saying that the Earth, Sun, and stars revolve around an {\it unseen central fire}. As alluded to in Archimedes'\index{Archimedes (c. 287--c. 212 BC)} work {\it Sand Reckoner}, Aristarchus of Samos\index{Aristarchus (c. 310--c. 230 BC)} ({\it ca.}~\!310 BCE--{\it ca.}~\!230 BCE)---apparently influenced by Philolaus\index{Philolaus (c. 470--c. 385 BC)}---had speculated that the sun is at the center of the solar system for the reason that the geocentric theory\index{geocentric theory} is against his conclusion that the sun is much bigger than the moon. In his work {\it On the Sizes and Distances}, he figured out the distances and sizes of the sun and the moon, under the assumption that the moon receives light from the sun. In the {\it Sand Reckoner}, Archimedes\index{Archimedes (c. 287--c. 212 BC)} himself harbored the ambition to estimate the size of the universe in the wake of Aristarchus'\index{Aristarchus (c. 310--c. 230 BC)} heliocentric spherical model\index{spherical model}. To this end, he proposed a peculiar number system to remedy the inadequacies of the Greek one and expressed the number of sand grains filling a cosmological sphere. The presupposition he made is that the ratio of the diameter of the universe to the diameter of the orbit of the earth around the sun equals the ratio of the diameter of the orbit of the earth around the sun to the diameter of the earth. 

Aristarchus'\index{Aristarchus (c. 310--c. 230 BC)} heliocentric model was espoused by Seleucus\index{Seleucus (born c. 190 BC)} of Seleucia, a Hellenistic astronomer (born {\it ca}.~\!190 BCE), who developed a method to compute planetary positions. In the end, however, the Greek heliocentrism\index{heliocentric theory} had been long forgotten. Even Hipparchus\index{Hipparchus (c. 190--c. 120 BC)} who undertook to find a more accurate distance between the sun and the moon took a step backwards.

Now in passing, we shall make special mention of Alexandria founded at the mouth of the Nile in 332 BCE by Alexander the Great, after his conquest of Egypt. It was Ptolemy~\!I, the founder of the Ptolemaic dynasty, who raised the city to a center of {\it Hellenistic culture}, an offspring of Greek culture that flourished around the Mediterranean after the decline of Athens. He and his son Ptolemy~\!II---both held academic activities in high esteem---established the Great Museum\index{Great Museum} (\textgreek{Mouse\~ion}), where
%, where as many as 700,000 papyrus scrolls collected by Ptolemies' avid hunt were housed. 
many thinkers and scientists from the Mediterranean world studied and collaborated with each other. Archimedes\index{Archimedes (c. 287--c. 212 BC)} stayed there sometime in adolescence.
%; indeed, the {\it cosmopolitanism} is often taken to be a remarkable hallmark of the Hellenistic period. 

\begin{figure}[htbp]
%\vspace{-1.7cm}
\vspace{-0.3cm}
\begin{center}
%\special{pdf: minorversion=6}
\includegraphics[width=.45\linewidth]{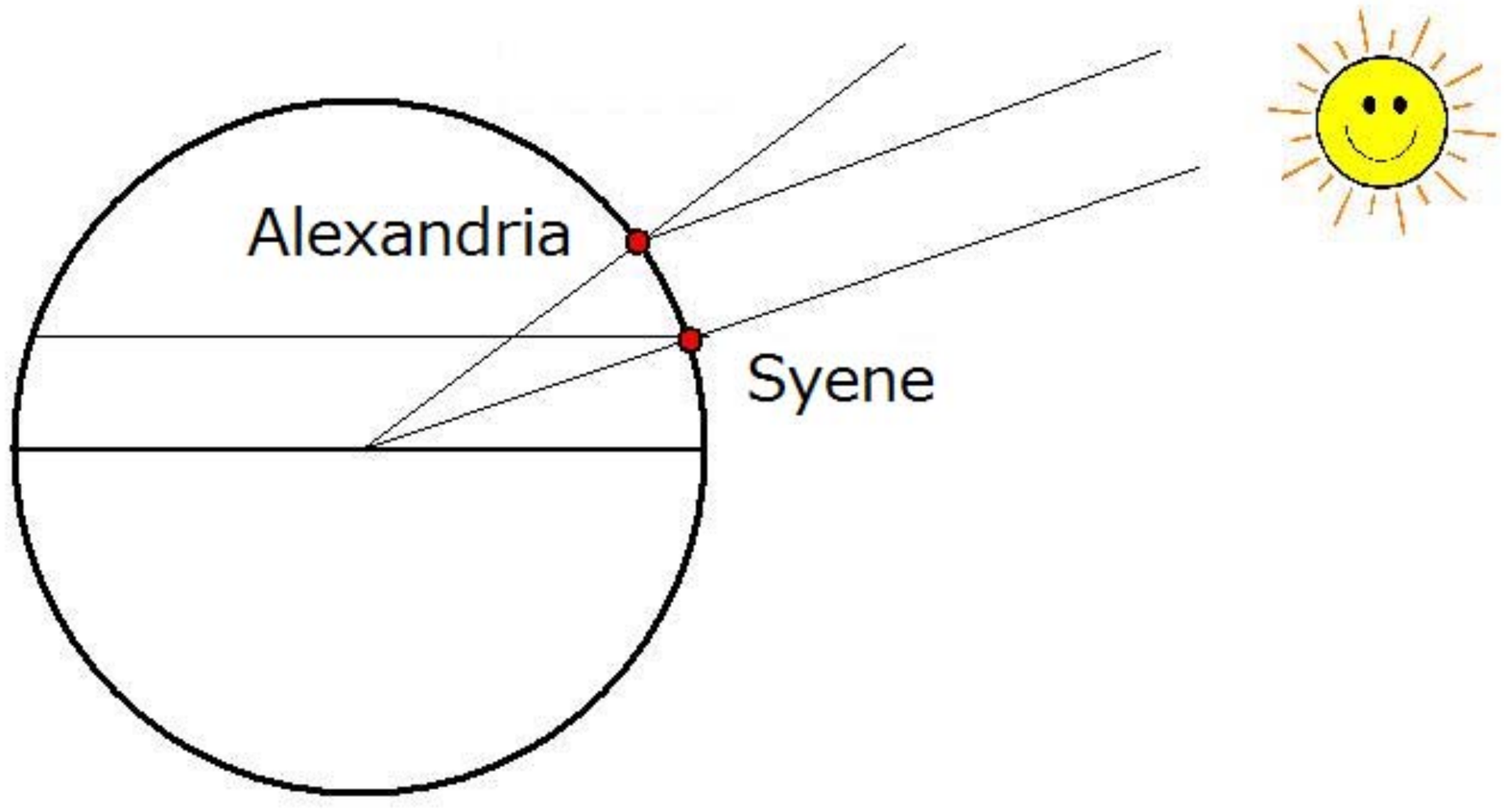}
\end{center}
\vspace{-1.5cm}
\caption{The circumference of the earth}\label{fig:sun}
\vspace{-0.1cm}
\end{figure}

Prominent among them is Eratosthenes of Cyrene\index{Eratosthenes (c. 276--c. 195 BC)} ({\it ca.}~\!276 BCE--{\it ca.}~\!195 BCE), the third chief librarian appointed by Ptolemy II and famous for an algorithm for finding all prime numbers.\footnote{Archimedes'\index{Archimedes (c. 287--c. 212 BC)} {\it Method} mentioned in fn.~\!\ref{fn:heiberg} takes the form of a letter to Eratosthenes.} Pushing further the view about the spherical earth, he adroitly calculated the meridian of the earth. The outcome was 46,620 km, about 16\% greater than the actual value. His work {\it On the measurement of the earth} was lost, but the book {\it On the circular motions of the celestial bodies} by Cleomedes (died {\it ca}.~\!489 CE) explains Eratosthenes'\index{Eratosthenes (c. 276--c. 195 BC)} deduction relying on a property of parallels and the observation that at noon on the summer solstice, the sun casts no shadow in Syene (now Aswan), while it casts a shadow on one fiftieth of a circle ($7.2$ degree) in Alexandria on the same degree of longitude as Syene (Fig.~\!\ref{fig:sun}). His resulting value is deduced from the distance between the two cities (5040 stades$=$925 km), which he inferred from the number of days that caravans require for travelling between them. 

The Ptolemaic dynasty lasted until the Roman conquest of Egypt in 30 BCE, but even afterwards Alexandria maintained its position as the center of scientific activity (see Sect.~\!\ref{subsec:discovery}).

\medskip

Finishing the excursion into ancient times, we shall turn back to the later stage of the Renaissance, the age that science was progressively separated from theology, and scientists came to direct their attention to scientific evidence rather than to theological accounts (thus science was becoming an occasional annoyance to the religious authorities). Representative in this period (religiously tumultuous times) are the Catholic Galileo\index{Galileo Galilei (1564--1642)} and the Protestant Johannes Kepler\index{Kepler, Johannes (1571--1630)} (1571--1630); both gave a final death blow to Aristotelian/Ptolemaic theory.

Kepler\index{Kepler, Johannes (1571--1630)} discovered the {\it three laws of planetary motion}, based on the data of Mars' motion recorded by Tycho Brahe\index{Tycho Brahe (1546--1601)} (1546--1601).\footnote{Tycho Brahe, the last of the major naked-eye astronomers, rejected the heliocentricism\index{heliocentricism}. His model of the universe is a combination of the Ptolemaic and the Copernican systems, in which the planets revolved around the Sun, which in turn moved around the stationary Earth.} His first and second laws were elaborated in the {\it Astronomia nova}  (1609). The first law asserts that the orbit of a planet is elliptical in shape with the center of the sun being located at one focus. The second law says that a line segment joining the sun and a planet sweeps out equal areas in equal intervals of time.  The third law stated in the {\it Epitome Astronomiae Copernicanae} (1617--1621) and the {\it Harmonices Mundi} (1619) maintains that the square of the orbital period of a planet is proportional to the cube of the semi-major axis of its orbit.

\begin{figure}[htbp]
%\vspace{-1.7cm}
\vspace{-0.2cm}
\begin{center}
\includegraphics[width=.28\linewidth]{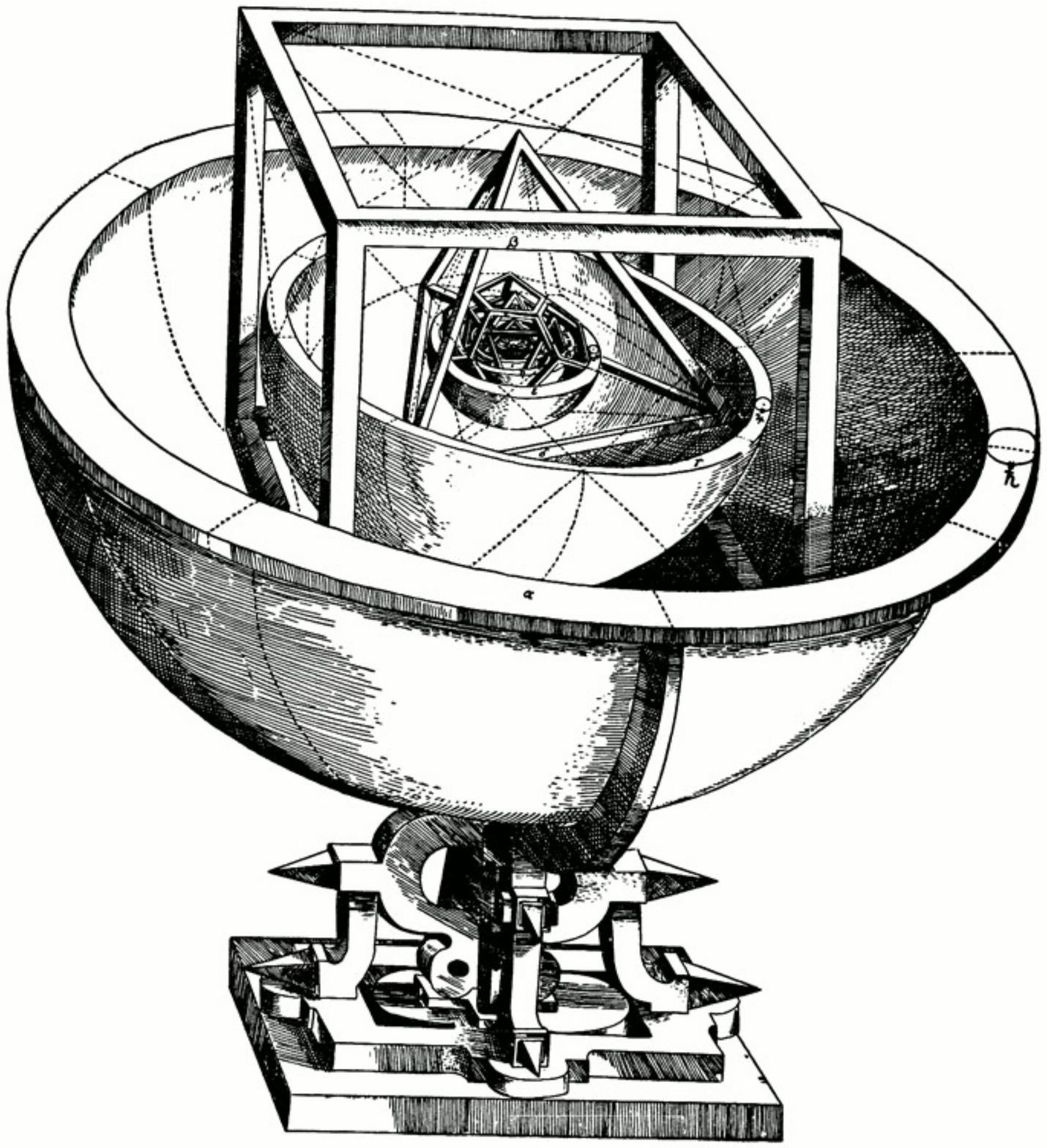}
\end{center}
\vspace{-1.8cm}
\caption{Kepler's model of the Solar system}\label{fig:kepler}
\vspace{-0.2cm}
\end{figure}

Before this epoch-making discovery, however, he attempted to explain the distances in the Solar system by means of ``regular convex polyhedra" inscribing and circumscribed by spheres, where the six spheres separating those solids correspond to Saturn, Jupiter, Mars, Earth, Venus, and Mercury ({\it Mysterium Cosmographicum}, 1596; Fig.~\!\ref{fig:kepler}). He also rejected the infiniteness of the universe, on the basis of an astronomical speculation on the one hand, and on the traditional scholastic doctrine on the other (for this reason, Kepler\index{Kepler, Johannes (1571--1630)} is portrayed as the last astronomer of the Renaissance, and not the first of the new age).

Galileo\index{Galileo Galilei (1564--1642)} was not affected by Aristotelian\index{Aristotle (384--322 BC)} prejudice that an experiment was an interference with the natural course of Nature. Without performing any experiment at all, medieval successors of Aristotle\index{Aristotle (384--322 BC)} insisted, for instance, that a projectile is pushed along by the force they called ``impetus," and if the impetus is expended, the object should fall straight to the ground. At the apogee of his scientific career, Galileo\index{Galileo Galilei (1564--1642)} had conducted experiments on projectile motion, and observed that their trajectories are always parabolic (1604--1608). He further investigated the relation between the distance $d$ an object (say, a ball) falls and the time $t$ that passes during the fall, and found the formula $d=\frac{1}{2}g t^2$ (in modern terms), where $g=9.8~\!{\rm m}/ {\rm s}^2$, the {\it gravitational acceleration} ({\it Discorsi}). The crux of this formula is that the drop distance does not depend on the mass of a falling object, contrary to Aristotle's\index{Aristotle (384--322 BC)} prediction ({\it Physics}, IV, 8, 215a25).\footnote{\label{fn:pisa}The {\it mean speed theorem} due to the Oxford Calculators\index{Oxford Calculators} (nothing but the formula for the area of a {\it trapezoid} from today's view) is a precursor of ``the law of falling bodies."}

Galileo\index{Galileo Galilei (1564--1642)} turned his eyes on the universe. With his handmade telescope, he made various astronomical observations and confirmed that four moons orbited Jupiter. This prompted him to defend the heliocentricism\index{heliocentricism}, for provided that Aristotle\index{Aristotle (384--322 BC)} were right about all things orbiting earth, these moons could not exist ({\it Sidereus Nuncius}, 1610). As an inevitable consequence, his Copernican view led to the condemnation at the Inquisition (1616, 1633). After being forced to racant, he was thrown into a more thorny position than Kepler\index{Kepler, Johannes (1571--1630)}, and had to avoid dangerous issues that may provoke the Church. Regarding the size of the universe, he denied the existence of the celestial sphere as the limit, while he was reluctant to say definitely that the universe is infinite, because of censorship by the Church. In a letter to F. Liceti on February 10, 1640, Galileo\index{Galileo Galilei (1564--1642)} says, ``the question about the size of the universe is beyond human knowledge; it can be only answered by the {\it Bible} and a divine revelation."

\medskip

Next to Kepler\index{Kepler, Johannes (1571--1630)} and Galileo\index{Galileo Galilei (1564--1642)} are Descartes\index{Descartes, Ren\'{e} (1596--1650)} and Pascal\index{Pascal, Blaise (1623--1662)}, the intellectual heroes of 17th-century France, who laid the starting point of the Enlightenment with their inquiries into truth and the limits of reason.   

Ren\'{e} Descartes\index{Descartes, Ren\'{e} (1596--1650)} (1596--1650), whose natural philosophy agrees in broad lines with Galileo's\index{Galileo Galilei (1564--1642)} one,  defended the heliocentrism, by saying that it is much simpler and distinct, but hesitated to make his opinion public upon the news of the Inquisition's conviction of Galileo\index{Galileo Galilei (1564--1642)} (1633). As regards the size of the universe, he maintained at first that the universe is finite; but later became ambivalent. In a letter in 1649 to the rationalist theologian H. More of Cambridge, he conceded, after a long dispute, that the universe must have infinite expanse because one cannot think of the limit; for if the limit exists, one cannot help but think of the outside of space which must be the same as our space.

While, for Descartes\index{Descartes, Ren\'{e} (1596--1650)}, the {\it self} is prior and independent of any knowledge of the world, Blaise Pascal\index{Pascal, Blaise (1623--1662)} (1623--1662) pronounced that the order is reversed; that is, knowledge of the world is a prerequisite to knowledge of the self. In the {\it Pens\'{e}es} (1670), in comparison with the disproportion between our justice and God's, he writes, ``Unity added to infinity does not increase it at all [....]: the finite is annihilated in the presence of the infinite and becomes pure nothingness" (fragment 418).  Further, he says, ``Nature is an infinite sphere of which the center is everywhere and the circumference nowhere" (fragment 199).\footnote{This sentence, reminiscent of Cusanus'\index{Cusanus, Nicolaus (1401--1464)} view, may be from the {\it Liber XXIV philosophorum} (Book of the 24 Philosophers), an influential philosophical and theological medieval text, usually attributed to Hermes Trismegistus, the purported author of the {\it Corpus Hermeticum}.}

\medskip

As described hitherto, the posture to pay attention to the universe is a steady tradition of European culture. In its background, philosophy was nurtured in the bosom of cosmogony, and could not be separated from religion because of the intimacy between them. Especially after the rise of Christianity, Western scholars had to be confronted, at the risk of their life in the worst case, with God as a creator. Even after the tension between science and religion was eased, scientists could not be entirely free from God, whatever His image is. Such milieu led quite naturally to probing questions about our universe in return.

In summary, ``whether the universe is finite or infinite," whatever it means, is an esoteric issue in religion, metaphysics, and astronomy in Europe that had intrigued and baffled mankind since dim antiquity. What matters most of all is whether the outside is necessary when we talk about finiteness of the universe.

\section{From Descartes to Newton and Leibniz}\label{sec:descartes}

Intuition for space surrounding us was the driving force for people to puff up their image of the universe. Needless to say, from the ancient times to the Middle Ages and even the early modern times, the mental image that nearly all people had is, if not being particularly conscious about it, the one described by {\it Euclidean space}\index{Euclidean space}, the model of our space named after the Alexandrian Euclid\index{Euclid ($4^{\rm th}$-$3^{\rm rd}$ c. BC)}
 (\textgreek{E>ukle'idhs}; {\it ca}.~300 BCE; see Sect.~\!\ref{subsec:elements} for the details). Even the advocates of the spherical model\index{spherical model} of the universe imagined Euclidean space\index{Euclidean space} as the entity embracing all. Synthetic geometry executed on this model is what we call {\it Euclidean geometry}\index{Euclidean geometry}. It started with a collection of geometrical results acquired in Egypt and Mesopotamia by empirical investigations or experience of land surveys and constructions of magnificent and imposing structures, and had been systemized, as the search of universals, through the efforts of Greek geometers.

Putting it briefly, Euclidean space\index{Euclidean space} 
is {\it homogeneous}\index{homogeneous} in the sense that there is no special place, and it is {\it isotropic}\index{isotropic}; i.e., there is no special direction.\footnote{Euclidean space holds one more significant aspect expressed by a property of parallels, but its true meaning had not been comprehended for quite a while; see Sect.~\!\ref{subsec:fifth}.} These features 
%(reminiscent of Bruno's\index{Bruno, Giordano (1548--1600)} vision of the universe) 
are not expressly indicated in Greek geometry, but are guaranteed by a property of {\it congruence}\index{congruence};\footnote{
% For example, 
Prop.~\!4 in the {\it Elements}\index{Elements}, Book I, is the first of the congruence\index{congruence} propositions.}
 namely any geometric figure can be rotated and moved to an arbitrary place while keeping its shape and size.
What should be pinpointed here is that it was not until the 19th century that people explicitly conceptualized Euclidean space\index{Euclidean space} (or space as a mathematical object which fits in with our spatial intuition). Until then, {\it geometry} meant only Euclidean geometry\index{Euclidean geometry}, and the space where geometry is performed was considered as an a priori entity, or what amounts to the same thing, our space---a place of storage in which objects are recognized, and a place of manufacture in which objects are constructed---was not an entity for which we investigate whether our understanding of it is right or wrong, thereby all the propositions of geometry being considered ``absolute truth."

As stated above, Euclidean geometry\index{Euclidean geometry} had been a lofty edifice for nearly 2000 years that nobody could break down. Only the appearance changed when Descartes\index{Descartes, Ren\'{e} (1596--1650)} invented the so-called {\it algebraic method}. This epoch-making method is elucidated in his {\it La G\'{e}om\'{e}trie}\index{Lageom@La G\'{e}om\'{e}trie} \cite{des}, one of the three essays attached to the philosophical and autobiographical treatise {\it Discours de la m\'{e}thode pour bien conduire sa raison, et chercher la v\'{e}rit\'{e} dans les sciences} (1637).\footnote{{\it La G\'{e}om\'{e}trie} consists of Book I ({\it Des probl\`{e}mes qu'on peut construire sans y employer que des cercles et des lignes droites}), Book II ({\it De la nature des lignes courbes}), and Book III ({\it De la construction des probl\`{e}mes solides ou plus que solides}). The other two essays are {\it La Dioptrique} and {\it Les M\'{e}t\'{e}ores}.}

Descartes'\index{Descartes, Ren\'{e} (1596--1650)} prime concern was, though he was trained in religion with the still-authoritarian nature, to find principles that one can know as true without any scruples. In the self-imposed search for {\it certainty}, he linked philosophy with science, and had the confidence that certainty could be found in the mathematical proofs having the apodictic character. Upon his emphasis on lucid methodology and dissatisfaction with the ancient arcane method, as already glimpsed in his {\it Regulae ad Directionem Ingenii} (1628), he attempted to create ``universal mathematics" by bridging the gap between arithmetic and geometry that used to be thought of as different terrains. Indeed, the Greeks definitely distinguished geometric quantities from numerical values, and even thought that length, area, and volume belong to different categories, thus lumping them together makes no sense. Under such shackles (and being devoid of symbolic algebra), the ``equality," ``addition/subtraction," and the ``large/small relation" for two figures in the same category were defined by means of geometric operations.\footnote{The Greeks had difficulty to handle irrationals in their arithmetic and were forced to replace {\it algebraic} manipulations by geometric ones (\cite{neuge}). Actually, geometry had been thought of as far more general than {\it arithmetic} as seen in Aristotle's\index{Aristotle (384--322 BC)} words, ``We cannot prove geometric truths by arithmetic" ({\it Posterior Analytics}, I, 7; see also Plato's\index{Plato (c. 428--c. 348 BC)} {\it Philebus}, 56d).

The word ``algebra," stemmed from the Arabic title {\it al-Kit$\overline{a}b$ al-mukhtasar fi his$\overline{a}$b al-jabr walmuq$\overline{a}$bala} (The Compendious Book on Calculation by Completion and Balancing) of the book written approximately 830 CE by al-Khw$\overline{\rm a}$rizm\=\i~({\it ca}.~\!780--{\it ca}.~\!850)\index{alk@al-Khw$\overline{\rm a}$rizm\=\i~(c. 780--c. 850)} (Remark \ref{rem:naturalnumber} (2)), was first imported into Europe in the early Middle Ages as a medical term, meaning ``the joining together of what is broken."
%; for instance, in the novel {\it Don Quixote} (1605), Miguel de Cervantes calls a bone-setter an {\it algebrista}. 
In addition, {\it algorithm}\index{algorithm}, meaning a process to be followed in calculations, is a transliteration of his surname al-Khw$\overline{\rm a}$rizm\=\i.} Specifically, they considered that two polygons (resp. polyhedra) are ``equal" if they are {\it scissors-congruent}\index{scissors-congruent}; i.e., if the first can be cut into finitely many polygonal (resp. polyhedral) pieces that can be reassembled to yield the second.
% (Fig.~\!\ref{fig:geoalg} illustrates a way to verify that a triangle is equal to one-half of a rectangle with the same base and height; cf. \index{Euclid ($4^{\rm th}$-$3^{\rm rd}$ c. BC)} Euclid's {\it Elements}\index{Elements}, Book I, Prop.~\!41).

%\begin{figure}[htbp]
%\vspace{-16.2cm}
%\vspace{-0.7cm}
%\begin{center}
%\includegraphics[width=0.9\linewidth]
%{geometricalgebranew.pdf}
%\end{center}
%\vspace{-6.2cm}
%\caption{Scissors-congruence}\label{fig:geoalg}
%\vspace{-0.4cm}
%\end{figure}

\begin{rem}{\rm {\small
Any two polygons with the same numerical area are scissors-congruent as shown independently by W. Wallace in 1807, Farkas Bolyai in 1832, and P. Gerwien in 1833. Gauss\index{Gauss, Johann Carl Friedrich (1777--1855)} questioned whether this is the case for polyhedra  in two letters to  his former student C. L. Gerling dated 8 and 17 April, 1844 ({\it Werke}, VIII, 241--42). In 1900, Hilbert\index{Hilbert, David (1862--1943)} put Gauss's\index{Gauss, Johann Carl Friedrich (1777--1855)} question as the third problem in his list of the 23 open problems at the second ICM. M. W. Dehn, a student of Hilbert\index{Hilbert, David (1862--1943)}, found two tetrahedra with the same volume, but non-scissors-congruent (1901); see \cite{sushi}.
\hfill$\Box$

}
}
\end{rem}

\begin{figure}[htbp]
%\vspace{-13.5cm}
\vspace{-0.8cm}
\begin{center}
\includegraphics[width=.86\linewidth]
{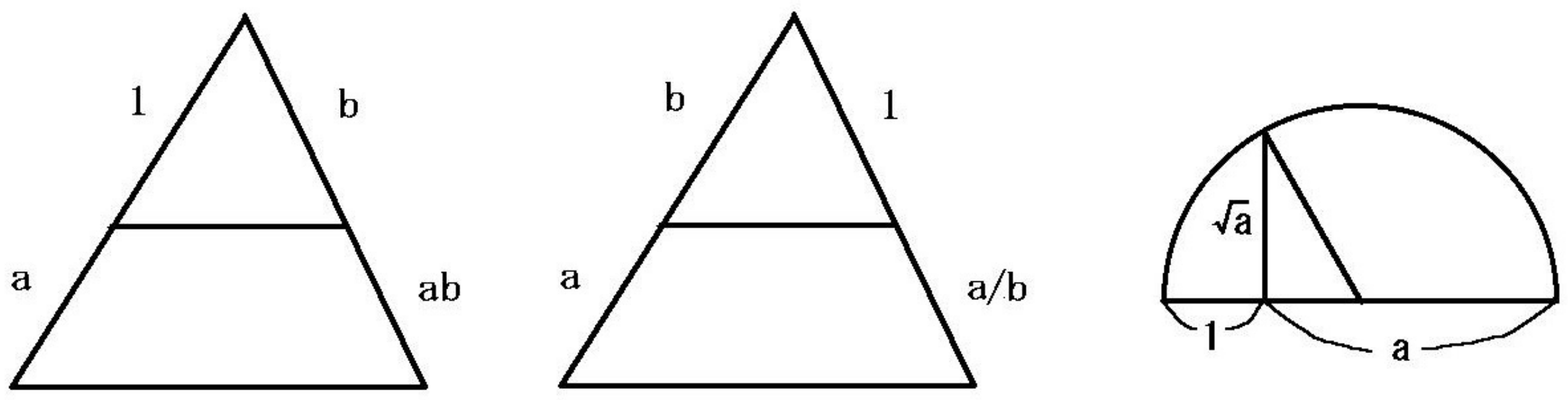}
\end{center}
\vspace{-5.2cm}
\caption{Operations of segments}\label{fig:decart}
\vspace{-0.2cm}
\end{figure}

In his essay {\it La G\'{e}om\'{e}trie}\index{Lageom@La G\'{e}om\'{e}trie}, Descartes\index{Descartes, Ren\'{e} (1596--1650)} made public that all kinds of geometric magnitude can be unified by representing them as line segments, and introduced the fundamental rules of calculation in this framework.\footnote{A letter dated March 26, 1619 to I. Beeckman indicates that Descartes\index{Descartes, Ren\'{e} (1596--1650)} already held a rough scheme at the age of 22. Pierre de Fermat\index{Fermat, Pierre de (1601--1665)} (1601--1665) argued a similar idea in his {\it Ad locos plano et solidos isagoge} (1636), though not published in his life time.} His idea is to fix a unit length, with which he defines addition, subtraction, multiplication, division, and the extraction of roots for segments, appealing to the proportional relations of similar triangles and the Pythagorean Theorem\index{Pythagorean Theorem}, as indicated in Fig.~\!\ref{fig:decart}. Then he adopts a {\it single axis} to represent these operations; thereby replacing the Greek {\it geometric algebra}\index{geometric algebra} by ``numerical" algebra. On this basic format, he handles various problems for a class of algebraic curves. Priding himself on his invention, he says, in a fond familiar letter dated November 1643 to Elisabeth (Princess Palatine of Bohemia), that extra effort to make up the constructions and proofs by means of Euclid's theorems is no more than a poor excuse for self-congratulation of petty geometers because such effort does not require any scientific mind.

Afterwards, Franciscus van Schooten\index{Schooten, Franciscus van (1615--1660)} (1615--1660), who met Descartes\index{Descartes, Ren\'{e} (1596--1650)} in Leiden (1632) and read the still-unpublished {\it La G\'{e}om\'{e}trie}\index{Lageom@La G\'{e}om\'{e}trie}, translated the French text into Latin, and endeavored to disseminate the algebraic method to the scientific community (1649). In the second edition (1659--61), he added annotations  and transformed Descartes'\index{Descartes, Ren\'{e} (1596--1650)} approach into a systematic theory, which made it far more accessible to a large number of audience.

Descartes'\index{Descartes, Ren\'{e} (1596--1650)} method, along with algebraic 
notations dating back to Fran\c{c}ois Vi\`{e}te (1540--1603)\index{Vi\`{e}te, Fran\c{c}ois (1540--1603)} who made a momentous step towards modern algebra,\footnote{Vi\`{e}te, {\it In artem analyticem isagoge} (1591).} was inherited as {\it analytic geometry}\index{analytic geometry} afterwards. With this progression, geometric figures were to be transplanted to algebraic objects in the coordinate plane $\mathbb{R}^2$ or the coordinate space $\mathbb{R}^3$.\footnote{It was Leibniz\index{Leibniz, Gottfried Wilhelm von (1646--1716)} who first introduced the ``coordinate system" in its present sense ({\it De linea ex lineis numero infinitis ordinatim ductis inter se concurrentibus formata, easque omnes tangente, ac de novo in ea re Analysis infinitorum usu}, Acta Eruditorum, {\bf 11} (1692), 168--171). Incidentally, Oresme\index{Oresme, Nicole (1320--1382)} used constructions similar to rectangular coordinates in the {\it Questiones super geometriam Euclidis} and {\it Tractatus de configurationibus qualitatum et motuum} ({\it ca}.~\!1370), but there is no evidence of linking algebra and geometry.} What is more, this new discipline was integrally connected with the calculus initiated by Issac Newton\index{Newton, Issac (1642--1726)} (1642--1726) and Gottfried Wilhelm von Leibniz\index{Leibniz, Gottfried Wilhelm von (1646--1716)} (1646--1716), independently and almost simultaneously, which is to be a vital necessity when we talk about the shape of our universe.

Newton\index{Newton, Issac (1642--1726)} learnt much from the {\it Exercitationes mathematicae libri quinque} (1657) by van Schooten\index{Schooten, Franciscus van (1615--1660)} and {\it Clavis Mathematicae} (1631) by William Oughtred\index{Oughtred, William (1574--1660)} (1574--1660) when he was a student of Trinity College, Cambridge. ``Fluxion" is Newton's underlying term in his differential calculus, meaning the instantaneous rate of change of a {\it fluent}, a varying (flowing) quantity. His idea, to which Newton was led by personal communication with Issac Barrow\index{Barrow, Issac (1630--1677)} (1630--1677) and also by Wallis' book {\it Arithmetica infinitorum} (1656),\footnote{John Wallis\index{Wallis, John (1616--1703)} (1616--1703) propagated Descartes'\index{Descartes, Ren\'{e} (1596--1650)} idea in Great Britain through his work.} is stated in two manuscripts; one is of October 1666 written when he was evacuated in a neighborhood of his family home at Woolsthorpe during the Great Plague, and another is the {\it Tractatus De Methodis Serierum et Fluxionum} (1671). In his calculus, Newton made use of power series in a systematic way.

Meanwhile, Leibniz\index{Leibniz, Gottfried Wilhelm von (1646--1716)} almost completed calculus while staying in Paris (1673--1676) as a diplomat of the Electorate of Mainz in order to get Germany back on its feet from the exhaustion caused by the Thirty Years' War. Although the aspired end of diplomacy was not attained, he had an opportunity to meet Christiaan Huygens\index{Huygens, Christiaan (1629--1695)} (1629--1695), a leading scientist of his time, who happened to be invited by Louis XIV as a founding member of the Acad\'{e}mie des Sciences. Very helpful for him was the suggestion by Huygens\index{Huygens, Christiaan (1629--1695)} to read Pascal's\index{Pascal, Blaise (1623--1662)} {\it Lettre de Monsieur Dettonville$\cdots$} (1658). Further, his perusal of Descartes'\index{Descartes, Ren\'{e} (1596--1650)} work was an assistance in strengthening the basis of his thought. Around the time when Leibniz\index{Leibniz, Gottfried Wilhelm von (1646--1716)} returned to Germany, he improved the presentation of calculus, and brought it out as two papers; {\it Nova methodus pro maximis et minimis}, and {\it De geometria recondita et analysi indivisibilium atque infinitorum}.\footnote{His papers were published in Acta Eruditorum, {\bf 3} (1684), 467--473 and {\bf 5} (1686), 292--300. In the second paper, he states, ``With this idea, geometry will make a far more greater strides than Vi\`{e}te's
%\index{Vi\`{e}te, Fran\c{c}ois (1540--1603)} 
and  Descartes'\index{Descartes, Ren\'{e} (1596--1650)}."}

Their calculus, though still something of mystery (Remark \ref{rem:infinixx} (2)), provided a powerful tool which enabled to handle more complicated geometric figures than the ones that Greek geometers treated in an {\it ad hoc} manner.\footnote{In the {\it Elements}\index{Elements}, Book III, Def.~\!2, Euclid\index{Euclid ($4^{\rm th}$-$3^{\rm rd}$ c. BC)}
says, ``a straight line is said to touch a circle which, meeting the circle and being produced, does not cut the circle." Curiously, in the {\it Metaphysics}, III, 998a, Aristotle\index{Aristotle (384--322 BC)} reports that Protagoras ({\it ca}.~\!490 BCE--{\it ca}.~\!420 BCE) argued against geometry that a straight line cannot be perceived to touch a circle at only one point.} Indeed, differential calculus provides a unified recipe to find tangents of a general curve (say, transcendental curves Descartes\index{Descartes, Ren\'{e} (1596--1650)} did not deal with), and integral calculus (called the ``inverse tangent problem" by Leibniz\index{Leibniz, Gottfried Wilhelm von (1646--1716)}) allows more latitude in calculating the area of a general figure without any ingenious trick. Paramountly important is the discovery that tangent (a local concept) and area (a global concept) are linked through {\it the fundamental theorem of calculus}\index{fundamental theorem of calculus} (FTC).\footnote{I. Barrow\index{Barrow, Issac (1630--1677)}, who strongly opposed to Descartes'\index{Descartes, Ren\'{e} (1596--1650)} method in contrast with Wallis\index{Wallis, John (1616--1703)}, knew the FTC in a geometric form (Prop.~\!11, Lecture 10 of his {\it Lectiones Geometricae} delivered in 1664--1666). A transparent proof was given by Leibniz\index{Leibniz, Gottfried Wilhelm von (1646--1716)} in his {\it De geometria}, 1686.}

There are other lines of evolution of Descartes'\index{Descartes, Ren\'{e} (1596--1650)} method and its offspring. Analytic geometry\index{analytic geometry} paved the way for {\it higher-dimensional} geometry (Sect.~\!\ref{sec:dimension}), which was to be incorporated into Riemann's\index{Riemann, Georg Friedrich Bernhard (1826--1866)} theory of curved spaces (Sect.~\!\ref{sec:riemann}). {\it Algebraic geometry} is regarded as the ultimate incarnation of Cartesian geometry, which came of age from the late 19th century to the early 20th century. Descartes'\index{Descartes, Ren\'{e} (1596--1650)} method also prompted to set up mathematical theories  in a strictly logical manner, not relying on geometry, but relying on arithmetic wherein concepts involved are entirely framed in the language of real numbers or systems of such (Sect.~\!\ref{sec:after1}). This signifies in some sense that mathematicians eventually become aware of {\it mathematical reality}---a sort of Plato's\index{Plato (c. 428--c. 348 BC)} world of {\it Forms} (\textgreek{{e'idos}} or its cognates)---distinguished from {\it physical reality} or the empirical world.

\begin{rem}\label{rem:analyticalsoc} {\rm {\small 
(1) Symbols and symbolization has played a significant role in the history of mathematics. Prior to Vi\`{e}te\index{Vi\`{e}te, Fran\c{c}ois (1540--1603)}, 
the symbols $+, - $ were used by H. Grammateus in {\it Ayn new Kunstlich Buech} (1518). Actually, these symbols already appeared in 1489 to indicate surplus and deficit for the mercantile purpose. In 1525, the symbol $\sqrt{\cdot}$ was invented by C. Rudolff, a student of  Grammateus. Subsequently, R. Recorde adopted the equal sign $ = $ in 1557. After Vi\`{e}te, Oughtred\index{Oughtred, William (1574--1660)} innovated the symbol $\times$ for multiplication in his {\it Clavis Mathematicae} (1631). In the same year, the book {\it Artis analyticae praxis ad aequationes algebraicas resolvendas} by Thomas Harriot\index{Harriot, Thomas (c. 1560 --1621)} ({\it ca}.~\!1560 --1621), which left a great mark on the history of symbolic algebras, was brought out posthumously, wherein the symbols ``$<, >$" representing the magnitude relation appear. For all of this, the usage of symbols in Descartes'\index{Descartes, Ren\'{e} (1596--1650)} time was almost in agreement with that of today. 

\smallskip

(2) 
 When integral calculus was still in its nascent stage, Kepler\index{Kepler, Johannes (1571--1630)} computed the volumes of solids of revolution.\footnote{{\it Nova stereometria doliorum vinariorum} (1615).} After a while, Bonaventura Cavalieri\index{Cavalieri, Bonaventura (1598--1647)} (1598--1647), a disciple of Galileo\index{Galileo Galilei (1564--1642)}, proposed a naive but a systematic method, known as {\it Cavalieri's principle}, to find areas and volumes of general figures, which is thought of as intermediating between the Greek quadrature and integral calculus.\footnote{{\it Geometria indivisibilibus continuorum nova quadam ratione promota} (1635, 1653).} The linchpin of his innovation is the notion of {\it indivisibles}\index{indivisible}; he regards, for instance, a plane region as being composed of an infinite number of parallel lines, each considered to be an infinitesimally thin rectangle.

\smallskip

(3) 
Although calculus initiated by Newton\index{Newton, Issac (1642--1726)} and Leibniz\index{Leibniz, Gottfried Wilhelm von (1646--1716)} was the same in essence, their styles differed in a crucial way. The difference was reflected in the notation they used. For instance, Leibniz\index{Leibniz, Gottfried Wilhelm von (1646--1716)} invented the symbol $dx/dt$, while Newton employed the dot notation $\dot{x}$. This difference arose, as often said, that the British mathematics under a baneful influence of Newton's\index{Newton, Issac (1642--1726)} authority had fallen seriously behind the Continental counterpart. They were jolted by the {\it M\'{e}canique C\'{e}leste} by Pierre-Simon Laplace\index{Laplace, Pierre-Simon (1749--1827)} (1749--1827) which was affiliated with the ``Leibnizian school." John Playfair\index{Playfair, John (1748--1819)} (1748--1819) says, ``We will venture to say that the number of those in this island who can read the {\it M\'{e}canique C\'{e}leste} with any tolerable facility is small indeed" (1808).

\smallskip

(4)
Until the 19th century through a pre-classical theory implicit in the work of Newton\index{Newton, Issac (1642--1726)}, Leibniz\index{Leibniz, Gottfried Wilhelm von (1646--1716)} and their successors, real numbers had been grasped as points on a straight line---the vision that dates back to Descartes\index{Descartes, Ren\'{e} (1596--1650)} and is credited to J. Wallis\index{Wallis, John (1616--1703)}---or ``something" (without telling what they are) approximated by rational numbers.\footnote{Simon Stevin\index{Stevin, Simon (1548-1620)} (1548-1620) renovated the notation for decimal fractions to make all computations easier, which contributed to the apprehension of the nature of real numbers. In the 35-page booklet {\it De Thiende} (1585), he expressed $27.847$ by 27$\textcircled{\scriptsize0}$8$\textcircled{\scriptsize1}$4$\textcircled{\scriptsize2}$7$\textcircled{\scriptsize3}$ for instance. } 

Thus calculus had been built on a fragile base. In a letter to his benefactor dated March 29, 1826, Niels Henrik Abel\index{Abel, Niels Henrik (1802--1829)} (1802--1829)  writes, ``the tremendous obscurity one undoubtedly finds in analysis today. It lacks all plan and system [$\cdots$] The worst of it is, it has never been treated with rigor. There are very few theorems in advanced analysis that have been demonstrated with complete rigor" (Remark \ref{rem:function}).

An entrenched foundation to calculus was furnished by Richard Dedekind\index{Dedekind, Richard (1831--1916)} (1831--1916). He  constructed the real number system by appealing to what we now call the method of {\it Dedekind cuts}\index{Dedekind cuts}\index{Dedekind, Richard (1831--1916)}, partitions of the rational numbers into two sets such that all numbers in one set are smaller than all numbers in the other (1872).\footnote{
%{\it Stetigkeit und irrationale Zahlen}, 1872. 
In the same year, Cantor\index{Cantor, Georg (1845--1918)} arrived at another definition of the real numbers (fn.~\!\ref{fn:realnumber}).} To say the least, this understanding 
%of real numbers 
heralds the emancipation of analysis from geometry.

Long ago, Eudoxus\index{Eudoxus (c. 408--c. 355 BC)} developed an idea similar to Dedekind's\index{Dedekind, Richard (1831--1916)} in essence. Before that, Greek geometers had tacitly assumed that two magnitudes $\alpha$, $\beta$ are always commensurable\index{commensurable} (\textgreek{s'ummetros}); i.e., there exist two natural numbers $m,n$ such that $m\alpha=n\beta$. Thus, with the discovery of incommensurable magnitudes  which results from the Pythagorean Theorem\index{Pythagorean Theorem}, the issue arose as to how to define equality of two ratios. The impeccable definition attributed to Eudoxus\index{Eudoxus (c. 408--c. 355 BC)} is : ``Magnitudes are said to be in the same ratio, the first to the second and the third to the fourth, when, if any equimultiples whatever are taken of the first and third, and any equimultiples whatever of the second and fourth, the former equimultiples alike exceed, are alike equal to, or alike fall short of, the latter equimultiples respectively taken in corresponding order" (\cite{euclid}, Book V, p.114). In modern terminology, this is expressed as ``Given four magnitudes $\alpha,\beta,\gamma,\delta$, the two ratios $\alpha:\beta$ and $\gamma:\delta$ are said to be the same if for all natural numbers $m,n$, it be the case that according as $m\alpha\gtreqqless n\beta$, so also is $m\gamma \gtreqqless n\delta$."

Related to the theory of proportions is the {\it anthyphairesis}\index{anthyphairesis} (\textgreek{>anjufa'iresis}), a method to find the ratio of two magnitudes. A special case is {\it Euclid's algorithm}\index{Euclid ($4^{\rm th}$-$3^{\rm rd}$ c. BC)}
used to compute the greatest common divisors of two numbers ({\it Elements}\index{Elements}, Book VII, Prop.~\!2; see \cite{negre}). Specifically, for two magnitudes $\alpha>\beta$, it proceeds as $\alpha=n_1\beta+\gamma_1 ~(0<\gamma_1<\beta),~\beta=n_2\gamma_1+\gamma_2 ~(0<\gamma_2<\gamma_1),~\gamma_1=n_3\gamma_2+\gamma_3 ~(0<\gamma_3<\gamma_2),~\gamma_2=n_4\gamma_3+\gamma_4 ~(0<\gamma_4<\gamma_3),~\cdots$ (Book X, Prop.~\!2). Here, $\alpha$ and $\beta$ are commensurable if and only if this process terminates after a finite number of steps (i.e., $\gamma_n=0$ for some $n$). This idea was handed down to us as the technique to obtain {\it continued-fraction expansions}\index{continued fraction} afterward. Namely, for positive real numbers $\alpha,\beta$, eliminating $\gamma_i$ in the above, we have

\medskip
\qquad$
\displaystyle\frac{\alpha}{\beta}=n_1+\frac{1}{n_2+}~\frac{1}{n_3+}~\frac{1}{n_4+}~\cdots~~~~~(=[n_1,n_2,n_3,n_4\ldots])
$.\hfill$\Box$

}
}
\end{rem}

\section{A new approach in classical geometry}\label{sec:proj}

\begin{figure}[htbp]
%\vspace{-11.8cm}
\vspace{-0.6cm}
\begin{center}
\includegraphics[width=.7\linewidth]
{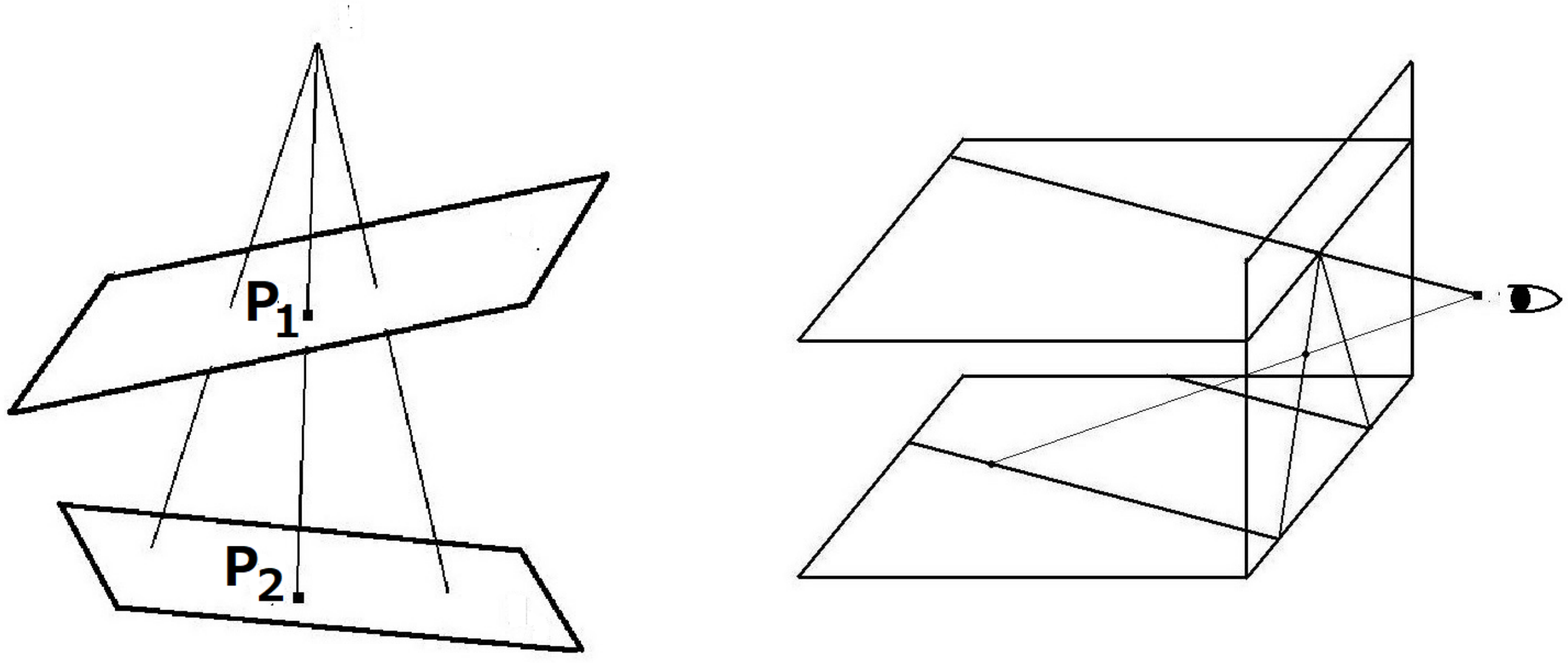}
\end{center}
\vspace{-3.1cm}
\caption{Perspective}\label{fig:perspective}
\vspace{-0.1cm}
\end{figure}

Leaving Euclidean geometry\index{Euclidean geometry} aside for a moment, we shall touch on {\it projective geometry}\index{projective geometry}, the embryonic subject that came up during the Renaissance period and had matured in the 19th century. This new approach is exclusively concerned with quantity-independent properties such as ``three points are on a line" ({\it collinearity}) and ``three lines intersect at a point" ({\it concurrency}) that are invariant under {\it projective transformations}\index{projective transformation} (the left of Fig.~\!\ref{fig:perspective}).\footnote{\label{fn:erlangen} Imitating this wording, one may say that Euclidean geometry\index{Euclidean geometry} is a geometry that treats invariant properties under congruence\index{congruence} transformations. Such a perspective is ascribed to Christian Felix Klein\index{Klein, Christian Felix (1849--1925)} (1849--1925), who promulgated his scheme in the booklet {\it Vergleichende Betrachtungen \"{u}ber neuere geometrische Forschungen} (``Erlangen Program"\index{Erlangen Program} in short) in 1872.}

Projective geometry\index{projective geometry} has a historical link to (linear) {\it perspective}\index{perspective}, an innovative skill in drawing contrived by Filippo Brunelleschi(1377--1446), which allowed Renaissance artists to portray a thing in space and landscapes as someone actually might see them.\footnote{One of the earliest who have used perspective is Masaccio. He limned the {\it San Giovenale Triptych} (1422) based on the principles he learned from Brunelleschi, which emblematizes the transition from medieval mysticism to the Renaissance spirit.} The theoretical aspect of perspective\index{perspective} was investigated by Leon Battista Alberti  (1404--1472) and Piero della Francesca (1416--1492). In the preface of his book {\it Della Pittura} (1435) dedicated to Brunelleschi, Alberti emphatically writes, ``I used to marvel and at the same time to grieve that so many excellent and superior arts and sciences from our most vigorous antique past could now seem lacking and almost wholly lost. $\cdots$ Since this work [due to Brunelleschi] seems impossible of execution in our time, if I judge rightly, it was probably unknown and unthought of among the Ancients" (translated by J. R. Spencer). Piero's book {\it De Prospectiva Pingendi} (1475) exemplifies the effective symbiosis of geometry and art; he actually had profound knowledge of Greek geometry and made a transcription of a Latin translation of Archimedes' work\index{Archimedes (c. 287--c. 212 BC)}.

\begin{figure}[htbp]
%\vspace{-12.5cm}
\vspace{-0.5cm}
\begin{center}
\includegraphics[width=1.02\linewidth]
{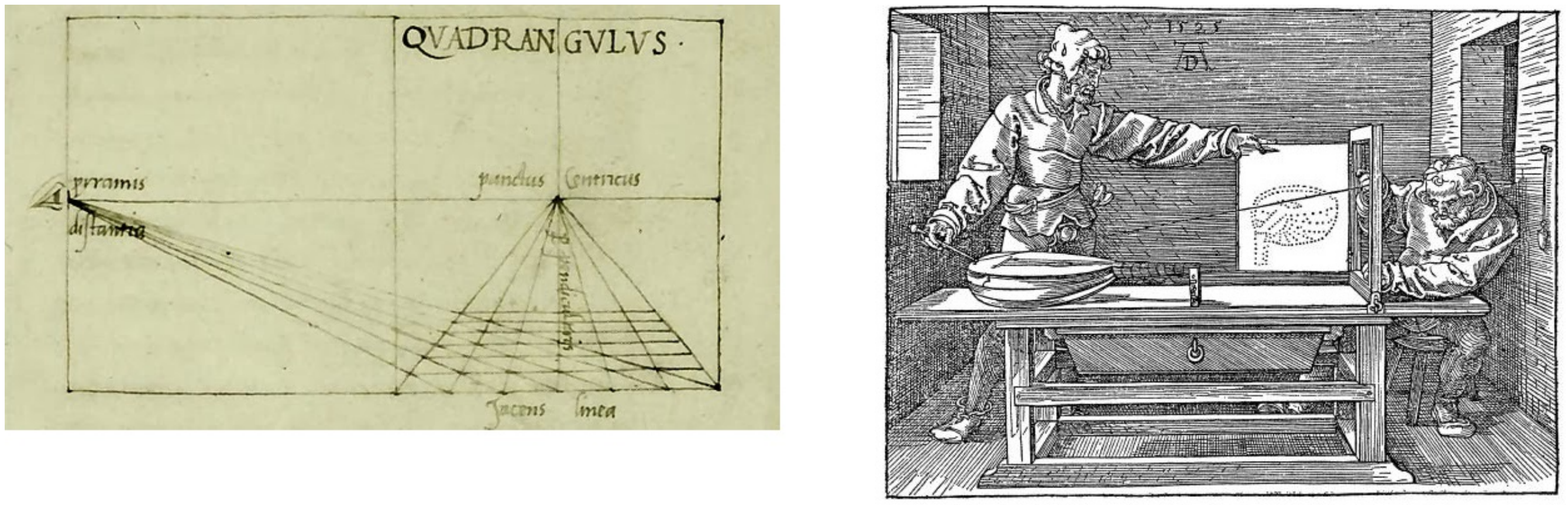}
\end{center}
\vspace{-4.9cm}
\caption{Alberti's drawing in {\it Della Pictura} and D\"{u}rer's illustration}\label{fig:perspective0}
\vspace{-0.1cm}
\end{figure}

The ``Renaissance Man" Leonardo da Vinci\index{da Vinci, Leonardo (1452--1519)} (1452--1519)---stimulated by Alberti's book---fully deserves his reputation as a true master of perspective\index{perspective}. His technique is particularly seen in his study for {\it A Magis adoratur} ({\it ca}.~\!1481). He says, ``Perspective is nothing else than seeing a place or objects behind a plane of glass, quite transparent, on the surface of which the objects behind the glass are to be drawn" (\cite{richter}). Albrecht D\"{u}rer (1471--1528), a weighty figure of the Northern Renaissance who shared da Vinci's pursuit of art, discussed in his work {\it Underweysung der Messung} (1525) an assortment of mechanisms for drawing in perspective from models (the right of Fig.~\!\ref{fig:perspective0}). This work, in which he touched on ``Doubling the cube," is the first on advanced mathematics in German.

Subsequently, Federico Commandino\index{Commandino, Federico (1506--1575)} (1506--1575) published the work entitled {\it Commentarius in planisphaerium Ptolemaei} (1558). This, a commentary on  Ptolemy's\index{Ptolemy (c. 100--c. 170)} work {\it Planisphaerium}, includes an account of the {\it stereographic projection} of the celestial sphere, and is a work on perspective\index{perspective} from a mathematical viewpoint (concurrently, he translated several work of ancient scholars).

\begin{rem}\label{fn:witch} {\rm {\small 
Da Vinci\index{da Vinci, Leonardo (1452--1519)} drew the illustrations of the regular polyhedra in the book {\it De divina proportione} (1509) by Fra Luca Bartolomeo de Pacioli\index{Pacioli, Fra Luca Bartolomeo de (c. 1447--1517)} ({\it ca}.~\!1447--1517), a pupil of Piero. The theme of the book is mathematical and artistic proportions, especially the {\it golden ratio}\index{golden ratio}---the positive solution $(1+\sqrt{5})/2=[1,1,1,\ldots]$ of the equation $x^2=x+1$---and its application in architecture. The golden ratio appears in the {\it Elements}\index{Elements}, Book II, Prop.~\!11; Book IV, Prop.~\!10--11; Book XIII, Prop.~\!1-–6, 8-–11, 16-–18,\footnote{In Book VI, Def.~\!3, Euclid\index{Euclid ($4^{\rm th}$-$3^{\rm rd}$ c. BC)}
says, ``A straight line is said to have been cut in {\it extreme and mean ratio} (\textgreek{>'akros ka`i m'esos l'ogos}) when, as the whole line is to the greater segment, so is the greater to the less." In today's terms, ``the segment $AB$ is cut at $C$  in extreme and mean ratio when $AB : AC\!=\!AC : CB$." Note that $AB/AC$ equals the golden ratio. Kepler\index{Kepler, Johannes (1571--1630)}, who proved that $\phi$ is the limit of the ratio of consecutive numbers in the Fibonacci\index{Fibonacci, Leonardo (c. 1170 --c. 1250)} sequence, described it as ``a fundamental tool for God in creating the universe" ({\it Mysterium Cosmographicum}).} and  is commonly represented by the Greek letter $\phi$ after the sculptor Pheidias (\textgreek{Feid'ias}, {\it ca}.~\!480 BCE--{\it ca}.~\!430 BCE) who is said to have employed it in his work.\hfill $\Box$

}
}
\end{rem}

The origin of projective geometry\index{projective geometry} can be traced back to the work of Apollonius of Perga\index{Apollonius (c. 262--c. 190 BC)} on {\it conic sections}\index{conic section} and Pappus of Alexandria\index{Pappus (c. 290--c. 350)} ({\it ca.}~\!290 CE--{\it ca.}~\!350 CE).\footnote{The first who studied conic sections\index{conic section} is Menaechmus ({\it ca}.\!380 BCE--{\it ca}.\!320 BCE), who used them to solve ``doubling the cube\index{doubling the cube}." The names {\it ellipse} (\textgreek{>'elleiyis}), {\it parabola} (\textgreek{parabolh}) and {\it hyperbola} (\textgreek{'uperbolh}) for conic curves were introduced by Apollonius\index{Apollonius (c. 262--c. 190 BC)}.}  In his {\it Collection}, Book VII, where Pappus\index{Pappus (c. 290--c. 350)} enunciated his {\it hexagon theorem}\index{Pappus' hexagon theorem},  he made use of a concept equivalent to what we call now the {\it cross-ratio}\index{cross-ratio}, which was to enrich projective geometry\index{projective geometry} from the quantitative side since it is invariant under projective transformations\index{projective transformation}.  Here the cross-ratio\index{cross-ratio} of collinear points $A$, $B$, $C$ and $D$ is defined as $[A,B,C,D]$ $=\displaystyle (AC\cdot BD)/(BC\cdot AD)$, where each of the distances is signed according to a consistent orientation of the line.\footnote{The cross-ratio\index{cross-ratio} of real numbers $x_1,x_2,x_3,x_4$ is defined by $[x_1,x_2,x_3,x_4]\!:=(x_3-x_1)(x_4-x_2)/(x_3-x_2)(x_4-x_1)$. The projective invariance of cross-ratios amounts to the identity $[T(x_1),T(x_2),T(x_3),T(x_4)]=[x_1,x_2,x_3,x_4]$, where $T(x)=(ax+b)/(cx+d)$ $(ad-bc\neq 0)$.} 

Projective geometry as a solid discipline was essentially initiated by Girard Desargues\index{Desargues, Girard (1591--1661)} (1591--1661), a coeval of Descartes\index{Descartes, Ren\'{e} (1596--1650)}.
In 1639, motivated by a practical purpose pertinent to perspective\index{perspective}, he published the {\it Brouillon project d'une atteinte aux \'{e}v\'{e}nemens des rencontres du Cone avec un Plan}.\footnote{The celebrated Desargues' theorem\index{Desargues' theorem} showed up in an appendix of the book {\it Exemple de l'une des mani\`{e}res universelles du S.G.D.L. touchant la pratique de la perspective} published in 1648 by his friend A. Bosse.} Pascal\index{Pascal, Blaise (1623--1662)} took a strong interest in Desargues'\index{Desargues, Girard (1591--1661)} work, and studied conic sections\index{conic section} in depth  at the age of 16. Philippe de La Hire\index{La Hire, Philippe de (1640--1718)} (1640--1718)  was affected by Desargues' work, too. He is most noted for the work {\it Sectiones Conicae in novem libros distributatae} (1685). However, partly because Descartes'\index{Descartes, Ren\'{e} (1596--1650)} method spread among the mathematical community in those days and because of his peculiar style of writing, his work remained unrecognized 
%for a long period.\footnote{
until his work was republished in 1864 by 
No\"{e}l-Germinal Poudra. Without being aware of Desargues' work,
%Almost 200 years later, 
Jean-Victor Poncelet\index{Poncelet, Jean-Victor (1788--1867)} (1788--1867) laid a firm foundation for projective geometry\index{projective geometry} in his masterpiece {\it Trait\'{e} des propri\'{e}t\'{e}s projectives des figures} (1822).

In any case, projective geometry\index{projective geometry} falls within classical geometry at this moment; but it has an interesting hallmark in view of infiniteness. Guidobaldi del Monte (1505--1607) in the {\it Perspectivae Libri VI} (1600) and Kepler\index{Kepler, Johannes (1571--1630)} in his {\it Astronomiae pars optica} (1604) proposed the idea of {\it points at infinity}\index{point at infinity}. Subsequently, having in mind the ``horizon" in perspective\index{perspective} drawing (Fig.~\!\ref{fig:perspective}), Desargues\index{Desargues, Girard (1591--1661)} appends points at infinity to the plane as idealized limiting points at the ``end," and argues that two parallel lines intersect at a point at infinity. Later, the ``plane {\it plus} points at infinity\index{point at infinity}" came to be termed the {\it projective plane}\index{projective plane} (Sect.~\!\ref{sec:after1}).

\section{Euclid's legacy in physics and philosophy}\label{sec:legacy}

Not surprisingly, Euclidean geometry\index{Euclidean geometry} predominated in physics for a long time. Salient is Newton\index{Newton, Issac (1642--1726)} who adopted Euclidean geometry as the base of his grand work {\it Philosophiae Naturalis Principia Mathematica} (1687, 1713, 1726; briefly called the {\it Principia}\index{Principia@Principia}), in which he formulated the {\it law of inertia},\footnote{The law of inertia was essentially discovered by Galileo\index{Galileo Galilei (1564--1642)} during the first decade of the 17th century though he did not understand the law in the general way (the term ``inertia" was first introduced by Kepler\index{Kepler, Johannes (1571--1630)} in his {\it Epitome Astronomiae Copernicanae}). The general formulation of the law was devised by Galileo's\index{Galileo Galilei (1564--1642)} pupils and Descartes\index{Descartes, Ren\'{e} (1596--1650)} (the {\it Principles of Philosophy}, 1644).} the {\it law of motion}, the {\it law of action-reaction} and the {\it law of universal gravitation}\index{law of universal gravitation}.

It was in 1666, exactly the same year when he worked out calculus, that Newton\index{Newton, Issac (1642--1726)} found a clue leading to the law of universal gravitation. Hagiography has it that he was inspired, in a stroke of genius, to formulate the law while watching the fall of an apple from a tree. Leaving aside whether this is a fact or not, all we can definitely say is, he concluded that the force acting on objects on the ground (say, apples) acts on the moon as well; he thus found the extraordinarily significant connection between the terrestrial and celestial which had been thought of as being independent of each other.   

To go further with his subsequent observation, let $O$ be the center of the earth, and let $P$ be the position of the moon. We denote by $v$ and $a$ the speed and acceleration of $P$, respectively. The acceleration is directed towards $O$, and $a=v^2/R$, where $R$ is the distance between $O$ and $P$. If the orbital period of the moon is $T$, then $vT=2\pi R$, so that $a=4\pi^2R/T^2$. Next, applying Kepler's\index{Kepler, Johannes (1571--1630)} third law to a circular motion, we have $R^3=cT^2$ $(c>0)$. Hence $a=4\pi^2cR^{-2}$ and $ma=4\pi^2mcR^{-2}$, which is the force acting on $P$ in view of the law of motion, and the law of universal gravitation follows. Conversely, based on his laws, Newton\index{Newton, Issac (1642--1726)} derived Kepler's\index{Kepler, Johannes (1571--1630)} laws through far-reaching deductions.\footnote{Newton\index{Newton, Issac (1642--1726)} tried to reconfirm his theory of gravitation by explaining the motion of the moon observed by J. Flamsteed, but the 3-body problem for the Moon, Earth and Sun turned out be too much complicated to accomplish his goal (he planned to carry a desired result as a centerpiece in the new edition of the {\it Principia}\index{Principia@Principia}).}
% Indeed, motions in the many-body problem could be {\it chaotic}.}

\begin{rem}\label{rem:gravitationnewton}{\rm {\small We let $\boldsymbol{x}(t)=\big(x(t),y(t),z(t)\big)\in \mathbb{R}^3$ be the position of a point mass (particle) in motion at time $t$. The law of motion\index{Newton's law of motion} is expressed as $m\ddot{\boldsymbol{x}}=\boldsymbol{F}$, where $m$ is the {\it inertial} mass\index{inertial mass} of the particle, and $\boldsymbol{F}$ stands for a force. The law of universal gravitation\index{law of universal gravitation} says that a point mass fixed at the origin $O$ attracts the point mass at $\boldsymbol{x}$ by the force $\boldsymbol{F}$ $=-\displaystyle GMm\|\boldsymbol{x}\|^{-3}\boldsymbol{x}$,\footnote{Throughout, $\|\cdot\|$ and $\langle\cdot,\cdot\rangle$ stand for {\it norm} and {\it inner product}, respectively.} where $G=6.67408 \times 10^{-11} {\rm m}^3 {\rm kg}^{-1} {\rm s}^{-2}$ is the {\it gravitational constant}, and $M$ is the {\it gravitational mass} at the origin. What should be stressed is,  as stated in the opening paragraph of the {\it Principia}\index{Principia@Principia}, that the inertial mass\index{inertial mass} coincides with the gravitational mass\index{gravitational mass} under a suitable system of units (this is by no means self-evident), so that the Newtonian equation is expressed as $\ddot{\boldsymbol{x}}=\displaystyle -GM\|\boldsymbol{x}\|^{-3}\boldsymbol{x}$, namely, the acceleration does not depend on the mass $m$ as Galileo\index{Galileo Galilei (1564--1642)} observed for a falling object. Kepler's\index{Kepler, Johannes (1571--1630)} laws of planetary motion\index{Kepler's laws of planetary motion} are deduced from this equation.\hfill$\Box$

}
}
\end{rem}

The Newton's\index{Newton, Issac (1642--1726)} cardinal importance in  the history of cosmology lies in his understanding of physical space. According to him, Euclidean space\index{Euclidean space} is an ``absolute space,"\index{absolute space} with which one can say that an object must be either in a state of absolute rest or moving at some absolute speed. To justify this, he assumed 
%to be able to  take up an ` proposing 
that the fixed stars can be a basis of ``inertial frame\index{inertial frame}." His {\it bucket experiment},\footnote{\label{fn:bucket}When a bucket of water hung by a long cord is twisted and released, the surface of water, initially being flat, is eventually distorted into a paraboloid-like shape by the effect of centrifugal force. This shape shows that the water is rotating as well, despite the fact that the water is at rest relative to the bucket. From this, Newton\index{Newton, Issac (1642--1726)} concluded that the force applied to water does not depend on the relative motion between the bucket and the water, but it results from the absolute rotation of water in the stationary absolute space\index{absolute space} (fn.~\!\ref{fn:Mach}).} elucidated in the {\it Scholium} to Book 1 of the {\it Principia}\index{Principia@Principia}, is an attempt to support the existence of absolute motion.

Setting aside the theoretical matters, what is peculiar about the {\it Principia}\index{Principia@Principia} is that Newton\index{Newton, Issac (1642--1726)} adhered to the Euclidean\index{Euclid ($4^{\rm th}$-$3^{\rm rd}$ c. BC)}
style, and used laboriously the traditional theory of proportions in Book V of the {\it Elements}\index{Elements} and the results in Apollonius'\index{Apollonius (c. 262--c. 190 BC)} {\it Conics} even though he was well acquainted with Descartes'\index{Descartes, Ren\'{e} (1596--1650)} work, which was, as a matter of course, more appropriate for the presentation. He testified, writing about himself in the third person, ``By the help of the {\it Analysis}, Mr.~\!Newton found out most of Propositions of his {\it Principia Philosophiae}\index{Principia@Principia}: but because the Ancients for making things certain admitted nothing into Geometry before it was demonstrated synthetically, he demonstrated the Propositions synthetically, that the System of Heavens might be founded upon good Geometry. And this makes it now difficult for unskillful Men to see the Analysis by which those Propositions were found out" (Phil.~\!Trans.~\!R.~\!Soc., {\bf 29} (1715), 173--224).

Apart from the style of the {\it Principia}\index{Principia@Principia}, this epoch-making work (and its offspring) was ``the theory of everything" in the centuries to follow as far as the classical description of the world is concerned. Indeed, Newton's\index{Newton, Issac (1642--1726)} laws seemed to tap all the secrets of nature, and hence to be the last word in physics. It was in the 19th century that physical phenomena inexplicable by the Newtonian mechanics were discovered one after another. A notable example is electromagnetic phenomena which are in discord with his mechanics in the fundamental level; see Remark \ref{rem:spacialrelativity}~\!(4). Moreover various physicochemical phenomena called for entirely new explanations as well, and eventually led to quantum physics. 

\medskip

Newton\index{Newton, Issac (1642--1726)} was a pious Unitarian. His theological thought had him say, ``Space is God's boundless uniform sensorium."\footnote{The medieval tradition had so much effect on Newton\index{Newton, Issac (1642--1726)} that he was deeply involved in alchemy and regarded the universe as a cryptogram set by God. As J. M. Keynes says, he was not the first of the age of reason, but the last of the magicians ({\it Newton the Man}, 1947).}  This gratuitous comment to his own comprehensive theory provoked a criticism from Leibniz\index{Leibniz, Gottfried Wilhelm von (1646--1716)}, who said that God does not need a {\it sense organ} to perceive objects. Moreover, on the basis of ``the principle of sufficient reason" and ``the principle of the identity of indiscernibles," Leibniz\index{Leibniz, Gottfried Wilhelm von (1646--1716)} claimed that space is merely relations between objects, thereby no absolute location in space, and that time is order of succession. His thought was unfolded in a series of long letters between 1715 and 1716 to a friend of Newton\index{Newton, Issac (1642--1726)}, S. Clarke.\footnote{\label{fn:clarke}Clarke (1717), {\it A Collection of Papers, which passed between the late Learned Mr. Leibniz, and Dr. Clarke, In the Years 1715 and 1716}.} Having said that, however, the existence of God is an issue which Leibniz\index{Leibniz, Gottfried Wilhelm von (1646--1716)} could not sidestep. To be specific, his principle of sufficient reason made him assert that nothing happens without a reason, and that all reasons are {\it ex hypothesi} God's reasons. One may ask, for instance, ``Why would God have created the universe here, rather than somewhere else?" That is, when God created the universe, He had an infinite number of choices. According to Leibniz\index{Leibniz, Gottfried Wilhelm von (1646--1716)}, He would choose the best one among different possible worlds. As will be explained later (Sect.~\!\ref{sec:geodesic}), his insistence, though having a strong theological inclination, is relevant to a fundamental physical principle. 

At all events, the period that begins with Newton\index{Newton, Issac (1642--1726)} and Leibniz\index{Leibniz, Gottfried Wilhelm von (1646--1716)} corresponds to the commencement of the close relationship between mathematics and physics. From then on, both disciplines have securely influenced each other.

\medskip

Immanuel Kant\index{Kant, Immanuel (1724--1804)} (1724--1804), who stimulated the birth of German idealism, was influenced by the rationalist philosophy represented by Descartes\index{Descartes, Ren\'{e} (1596--1650)} and Leibniz\index{Leibniz, Gottfried Wilhelm von (1646--1716)} on the one hand, and troubled by Hume's thoroughgoing skepticism on the other.\footnote{The 
%avowed 
empiricist David Hume is a successor of F. Bacon, T. Hobbes, J. Locke, and G. Berkeley. He says that an orderly universe does not necessarily prove the existence of God ({\it An Enquiry Concerning Human Understanding}, 1748). Hobbes, an Aristotle's critic\index{Aristotle (384--322 BC)}, 
described the world as mere ``matter in motion," maybe the most colorless depiction of the universe since the ancient atomists, but he did not abandon God in his cosmology ({\it Leviathan}, 1651).} He felt the need to rebuild metaphysics to argue against Hume's view. Being also dissatisfied with the state of affairs surrounding metaphysics, in contrast to the scientific model cultivated in his days, Kant\index{Kant, Immanuel (1724--1804)} strove to lay his philosophical foundation of reason and judgement on secure grounds. 

He was 
an enthusiast of Newtonian physics framed on Euclidean geometry\index{Euclidean geometry}, and asked himself, ``Are space and time real existences? Are they only determinations or relations of things as Leibniz\index{Leibniz, Gottfried Wilhelm von (1646--1716)} insists?" His inquiry brought him to the conclusion that Euclidean geometry\index{Euclidean geometry} is the inevitable necessity of thought and {\it inherent in nature} because our space (as an ``absolute" entity) is a ``sacrosanct" framework for all and any experience. He also held that space (as a ``relative" entity) is the ``subjective constitution of the mind" (see Sect.~\!\ref{sec:right}), He stressed further that mathematical propositions are formal descriptions of the a priori structures of space and time.\footnote{Affected by the science in his time, Kant\index{Kant, Immanuel (1724--1804)} developed further the {\it nebular hypothesis} proposed by E. Swedenborg in 1734 (this hypothesis that the solar system condensed from a cloud of rotating gas was discussed later by Laplace\index{Laplace, Pierre-Simon (1749--1827)} in 1796).}
 
In the work {\it Allgemeine Naturgeschichte und Theorie des Himmels} (1755), Kant\index{Kant, Immanuel (1724--1804)} discussed infiniteness of the universe. He believed that the universe is infinite, because God is the ``infinite being," and creates the universe in proportion to his power (recall Bruno's\index{Bruno, Giordano (1548--1600)} view). In turn, he argues that it is not philosophically possible to decide whether the universe is infinite or finite; that is, we are incapable of perceiving so large distance, because the mind is finite.\footnote{Hegel dismissed the claim that the universe extends infinitely.}

\medskip

As Kant\index{Kant, Immanuel (1724--1804)} emphasized in the {\it Kritik der reinen Vernunft} (1781), 
abstract speculation must be pursued without losing touch with reality grasped by intuition; otherwise the outcome could be empty. This is not least the case for scientific knowledge. Yet intuition about our space turned out to be a clinging constraint on us; it was not easy to free ourselves from it.

\section{Gauss --Intrinsic description of the universe}\label{sec:gauss}

In Christianity, God is portrayed to be omnipresent, yet His whereabouts is unknown. Nonetheless, the laity usually personify God (if not a person literally), and believe that He ``dwells" outside the universe (or {\it empyrean}).
%, and can take a view of the whole universe from the outside at once.
Imagining such a deity is not altogether extravagant from the view that we will take up in investigating a model of the universe though we do not posit the existence of God. What is available in place of God is, of course, {\it mathematics}.\footnote{Robert Grosseteste\index{Grosseteste, Robert (c. 1175--1253)} ({\it ca}.~\!1175--1253) in medieval Oxford says, ``mathematics is the most supreme of all sciences since every natural science ultimately depends on it" ({\it De Luce}, 1225). He maintained, by the way, that stars and planets in a set of nested spheres around the earth were formed by crystallization of matter after the birth of the universe in an explosion.}  Thus, we pose the question, ``Is it mathematically possible to tell how the universe is curved without mentioning any outside of the universe?"

The answer is ``Yes." To explain why, let us take a look at smooth {\it surfaces} as 2D toy models of the universe.\footnote{Throughout, ``smooth" means ``infinite differentiable."} That is, we human beings are supposed to be confined in a surface. Here space in which the surface is located is thought of as the ``outside," and is called the {\it ambient space}\index{ambient space}. A person in the outside plays the role of God incarnate. 
%since he is capable of telling the shape of the whole surface. 
In our story, we allow him to exist for the time being; that is, we investigate surfaces in an {\it extrinsic}\index{extrinsic} manner.
At the final stage, we are disposed to remove him (and even the outside) from the scene, and to let the inhabitants in the surface find a way to understand his universe. In other words, what we shall do is, starting from an extrinsic study of surfaces, to look for some geometric concepts with which one can talk about how surfaces are curved without the assistance of the outside. Such concepts will be called {\it intrinsic}\index{intrinsic}, one of the most vital key terms in this essay.

A few words about intrinsicness\index{intrinsicness}.
Imagine that an inhabitant confined in a surface wants to examine the structure of his universe. The only method he holds is, like us, measuring the distance between two points and the angle between two directions. In our universe, we use the {\it light ray} to this end, which in a uniform medium propagates along shortest curves (straight lines) in space;  a very special case of {\it Fermat's principle of least time}\index{Fermat's principle of least time} saying that the path taken between two points by a light ray is the path which can be traversed in the least time.\footnote{This principle was stated in a letter (January 1, 1662) to M. C. de la Chambre. Employing it, he deduced the {\it law of refraction}. The first discoverer of the law is  Harriot\index{Harriot, Thomas (c. 1560 --1621)} (July 1601), but he died before publishing. W. Snellius rediscovered the law (1621). Descartes\index{Descartes, Ren\'{e} (1596--1650)} derived the law in {\it La Dioptrique} under a wrong setup; actually he believed that light's speed was instantaneous. The finite propagation of light was demonstrated by O. R{\o}mer in 1676.} Hence shortest curves on a surface are regarded as something similar to paths traced by light rays (see Sect.~\!\ref{sec:geodesic}), and it is reasonable to say that a quantity (or a concept) attached to a surface is ``intrinsic"\index{intrinsic} if it is derived from distance and angle measured by using shortest curves.

The progenitor of the intrinsic\index{intrinsic} theory of surfaces is Johann Carl Friedrich Gauss\index{Gauss, Johann Carl Friedrich (1777--1855)} (1777--1855). In 1828, he brought out the memoir {\it Disquisitiones generales circa superficies curvas}, where he formulated {\it bona fide} intrinsic\index{intrinsic} curvature, which came to be called {\it Gaussian curvature}\index{Gaussian curvature} later.

Just for the comparison purpose, let us primarily look at the curvature of plane curves. Subsequent to the pioneering studies by Kepler\index{Kepler, Johannes (1571--1630)} ({\it Opera}, vol.~\!2, p.~\!175) and Huygens\index{Huygens, Christiaan (1629--1695)} (1653-4), Newton\index{Newton, Issac (1642--1726)} used his calculus to compute curvature.\footnote{{\it Tractus de methodis serierum et fluxionum} (1670-1671).} His approach is to compare a curve with {\it uniformly curved figures}; i.e., circles or straight lines. To be exact, given a point $p$ on a smooth curve $C$, he singles out the circle (or the straight line) $C_0$ that has the highest possible order of contact with $C$ at $p$, and then defines the curvature\index{curvature} of $C$ at $p$ to be the reciprocal of the radius of $C_0$ (when $C_0$ is a straight line, the curvature is defined to be $0$); thus the smaller the circle $C_0$ is, the more curved $C$ is. To simplify his computation, take a coordinate system $(x,y)$ such that the $x$-axis is the tangent line of $C$ at $p$, and express $C$ around $p$ as the graph associated with a smooth function $y\!=\!f(x)$ with $f(0)\!=\!f'(0)\!=\!0$ and $f''(0)\!\geq\!0$. Compare $f(x)$ $=$ $\frac{1}{2}f''(0)x^2+\frac{1}{3!}f'''(\theta x)x^3$ $~(0\!<\!\theta\!<\!1)$ (Taylor's theorem) with the function $f_0(x)$ $=\!R-\sqrt{R^2-x^2}\!=\!\frac{1}{2}R^{-1}x^2+\cdots$ corresponding to the circle $C_0$ of radius $R$, tangent to the $x$-axis at $x=0$.
%\footnote{Taylor's theorem\index{Taylor's theorem} is named after B. Taylor who stated the result (without the remainder term) in 1712 and published it in the {\it Methodus Incrementorum Directa et Inversa}, 1715, though J. Gregory had obtained the same result nearly 40 years earlier which is stated in a letter to the publisher John Collins on 15 February 1671. The remainder term was provided much later by Joseph-Louis Lagrange\index{Lagrange, Joseph-Louis (1736--1813)} (1736--1813).} 
Thus $C_0$ has the second order of contact with $C$ at $(0,0)$ if and only if $R=1/f''(0)$. Therefore the curvature\index{curvature} is equal to $f''(0)$.

\begin{rem}\label{rem:signedcurvature}{\rm {\small 
For an {\it oriented curve} $C$ (a curve with a consistent direction defined along the curve), we select a coordinate system $(x,y)$ such that the direction of $x$-axis coincides with that of $C$ at $p$. Then $f''(0)$ is possibly negative. We call $\kappa(p):$ $=$ $f''(0)$ the {\it signed curvature}\index{signed curvature} of $C$ at $p$. If $c:[a,b]\longrightarrow \mathbb{R}^2$ is a parameterization of $C$ with $\|\dot{c}\|\equiv 1$, then $\ddot{c}$ is perpendicular to $\dot{c}$, and $\ddot{c}(s)=\kappa(c(s))\boldsymbol{n}(s)$, where $\boldsymbol{n}(s)$ is the unit vector obtained by the $90^{\circ}$ counterclockwise rotation of $\dot{c}(s)$. Writing $c(s)=(x(s),y(s))$, we have $\boldsymbol{n}(s)=(-\dot{y}(s),\dot{x}(s))$, and $\kappa(c(s))=\langle\ddot{c}(s),\boldsymbol{n}(s)\rangle=\dot{x}(s)\ddot{y}(s)-\dot{y}(s)\ddot{x}(s)$.\hfill$\Box$

}
}
\end{rem}   

In the 18th century, Continental mathematicians polished calculus as an effectual instrument to handle diverse problems in sciences. 
Leonhard Euler\index{Euler, Leonhard (1707--1783)} (1707--1783) was among those who helped to brings calculus to completion and employed it as a tool for the extrinsic\index{extrinsic} geometry of surfaces. His particular interest was the curvature of the curve obtained by cutting a surface with the plane spanned by a tangent line and the normal line. This study brought him to the notion of {\it principal curvature}\index{principal curvatures}.\footnote{{\it Recherches sur la courbure des surfaces}, M\'{e}m.~\!Acad.~\!R.~\!Sci.~\!Berlin, {\bf 16} (1767), 119-143.} Then in 1795, Gaspard Monge\index{Monge, Gaspard (1746--1818)} (1746--1818) published the book {\it Application de l'analyse \`{a} la g\'{e}om\'{e}trie}, an early influential work on differential geometry containing a refinement of Euler's\index{Euler, Leonhard (1707--1783)} work. However neither Euler\index{Euler, Leonhard (1707--1783)} nor Monge\index{Monge, Gaspard (1746--1818)} arrived at intrinsic\index{intrinsic} curvature. Here came Gauss\index{Gauss, Johann Carl Friedrich (1777--1855)}, who found how to formulate it with absolute confidence in its significance.

\begin{figure}[htbp]
%\vspace{-7.2cm}
\vspace{-0.5cm}
\begin{center}
\includegraphics[width=.58\linewidth]
{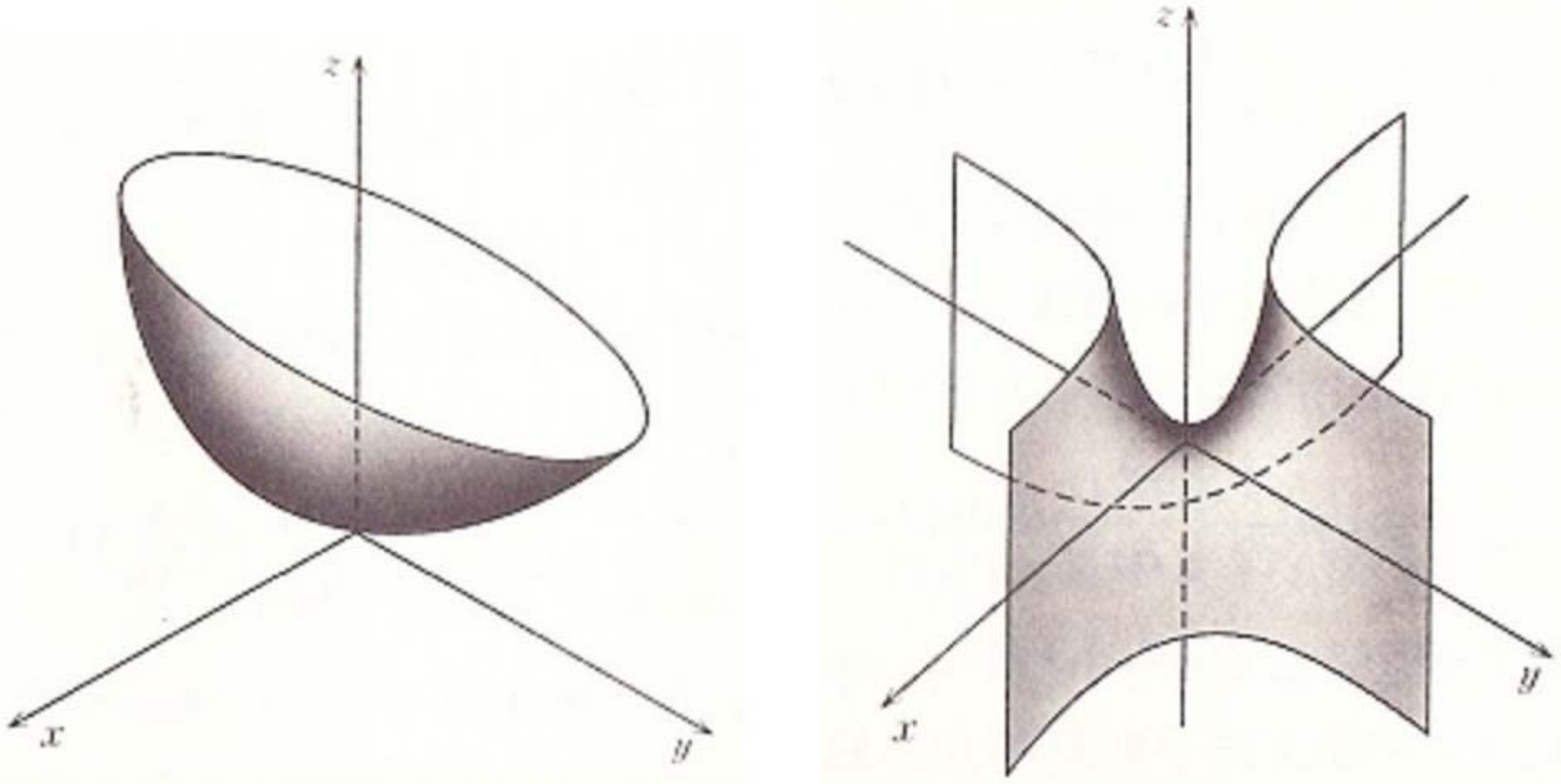}
\end{center}
\vspace{-1.8cm}
\caption{Paraboloid and Hyperboloid}\label{fig:surfacexx}
\vspace{-0.2cm}
\end{figure}

What comes to mind when defining the curvature\index{curvature} of a surface $S$ at $p$ is the use of a coordinate system $(x,y,z)$ with origin $p$ such that the $xy$-plane is tangent to $S$ at $p$,\footnote{\label{fn:gaussmap}Gauss's\index{Gauss, Johann Carl Friedrich (1777--1855)} original definition uses the map ({\it Gauss map}) from $S$ to the unit sphere which associates to each point on $S$ its oriented unit normal vector ({\it Disquisitiones generales}, Art.~\!6).} with which $S$ is locally expressed as a graph of a smooth function $f(x,y)$ defined around $(x,y)\!=\!(0,0)$ with $f(0,0)\!=\!f_x(0,0)\!=\!f_y(0,0)\!=\!0$, where $\displaystyle f_x\!=\!\partial f/\partial x$ and $\displaystyle f_y\!=$ $\partial f/\partial y$. Then we obtain $f(x,y)\!=\!\frac{1}{2}(ax^2+2bxy+cy^2)$ $+r(x,y)$, where $a\!=\!f_{xx}(0,0)$, $b\!=\!f_{xy}(0,0)$, $c\!=\!f_{yy}(0,0)$, and $r(x,y)$ is of higher degree than the second. The shape of the graph (and hence the shape of $S$ around $p$) is roughly determined by the quantity $ac-b^2$; it is paraboloid-like (resp. hyperboloid-like) if $ac-b^2\!>\!0$ (resp. $ac-b^2\!<\!0$) (Fig.~\!\ref{fig:surfacexx}). Given all this, we define the curvature $K_S(p)$ to be $ac-b^2$, which is independent of the choice of a coordinate system. For example, $K_S\!\equiv\!R^{-2}$ for the sphere of radius $R$.

To account the exact meaning of intrinsicness\index{intrinsicness} of the Gaussian curvature\index{Gaussian curvature}, we consider another surface $S'$, and let $\varPhi$ be a one-to-one correspondence\index{one-to-one correspondence} from $S$ to $S'$ preserving the distance of two points and the angle between two directions (such $\varPhi$ is called an {\it isometry}\index{isometry} in the present-day terms; see Sect.~\!\ref{sect:reimenniangeometry}). If we would adopt the approach mentioned above on intrinsicness\index{intrinsicness}, Gauss's\index{Gauss, Johann Carl Friedrich (1777--1855)} outcome could be rephrased as ``$K_{S'}(\varPhi(p))=K_S(p)$." In particular, $K_S\equiv 0$ for a cylindrical surface $S$, since $S$, which is ostensibly curved, is deformed without distortion to a strip in the plane. This example articulately signifies that the Gaussian curvature\index{Gaussian curvature} could be possibly defined without the aid of the outside.

Executing lengthy computations, Gauss\index{Gauss, Johann Carl Friedrich (1777--1855)} confirmed the intrinsicness\index{intrinsicness} of $K_S$ (see Sect.~\!\ref{sec:after2}).\footnote{Crucial 
in his argument is the commutativity of partial differentiation.
%; say $F_{xxy}=F_{xyx}$. 
Euler\index{Euler, Leonhard (1707--1783)} and A. C. Clairaut knew it, but the first correct proof was given by K. H. A. Schwarz in 1873.} Since the goal is more than a pleasing astonishment, Gauss\index{Gauss, Johann Carl Friedrich (1777--1855)} named his outcome ``Theorema Egregium\index{Theorema Egregium} (remarkable theorem)" ({\it Disquisitiones generales}, Art.~\!12). Actually, it is no exaggeration to say that this theorem is a breakthrough not only in geometry, but also in cosmology. The following is his own words (1825) unfolded in a letter to the P. A. Hansen:

\medskip 

{\it This research is deeply entwined with much else, I would like to say, with metaphysics of space and I find it difficult to shake off the consequences of this, such as for instance the true metaphysics of negative or imaginary quantities} ({\it Werke}, XII, p.~\!8. See \cite{bos}).\footnote{Even in the 18th century, there were a few 
%unrepentant people 
who did not accept negative number. For example, F. Maseres said, ``[negative numbers] darken the very whole doctrines of the equations and make dark of the things which are in their nature excessively obvious and simple" 
%({\it Dissertation on the use of the negative sign in algebra}, 
(1758). Although limited to the practical use, negative numbers already appeared in the {\it Jiuzhang Suanshu} (The Nine Chapters on the Mathematical Art) written in the period of the Han Dynasty (202 BCE--8 CE). The reason why the Chinese accepted negative numbers is, in all likelihood, that they had the concept of {\it yin} (negative or dark) and {\it yang} (positive or bright).} 

\medskip

Although the Theorema Egregium\index{Theorema Egregium} has a deep metaphysical nature as Gauss\index{Gauss, Johann Carl Friedrich (1777--1855)} put it, his incipient motivation came from the practical activity pertaining to the {\it geodetic survey} of Hanover to link up with the existing Danish grid that started in 1818 and continued for 8 years. He traveled across the state in order to cover the whole country. Fueled by this arduous activity, he brought out many papers on the geometry of curves and surfaces from about 1820 onwards. The Theorema Egregium\index{Theorema Egregium} is a culmination of his intensive study in this field.

\begin{rem}\label{fn:principal} {\rm {\small 
 In her study of elasticity (1831), M.-S. Germain introduced the {\it mean curvature} $H_S(p):=\frac{1}{2}(a+c)$, which appears in a characterization of {\it minimal surfaces}\index{minimal surface}; that is, the surface area of $S$ is locally minimized if and only if  $H_S\!\equiv\!0$. Putting $A\!=\!\begin{pmatrix} a& b\\ b &c\end{pmatrix}$, we have $K_S(p)\!=\!\det A$ and $H_S(p)\!=\!\frac{1}{2}{\rm tr}~\!A$. The eigenvalues $\kappa_1(p),~\!\kappa_2(p)$ of $A$ are called the {\it principal curvatures}\index{principal curvatures} of $S$. Obviously $K_S(p)\!=\!\kappa_1(p)\kappa_2(p)$ ({\it Disquisitiones generales}, Art.~\!8) and $H_S(p)\!=\!\frac{1}{2}(\kappa_1(p)+\kappa_2(p))$.\hfill$\Box$

}
}
\end{rem}

\section{Number systems}
It is off the subject, but while looking at the background behind Gauss's\index{Gauss, Johann Carl Friedrich (1777--1855)} opinion about complex numbers quoted above, we shall prepare some material related to what we shall describe later (Sect.~\!\ref{subsec:topodes}).

Gerolamo Cardano\index{Cardano, Gerolamo (1501--1576)} (1501--1576), who stood out above the rest at that time as an algebraicist, acknowledged the existence of complex numbers in his {\it Ars Magna} (1545) that contains the first printed solutions to cubic and quartic equations. He observed that his formula may possibly involve complex numbers even when applied to a cubic equation possessing only real solutions.\footnote{For example, his formula applied to the cubic equation $x^3-15x-4~\!(=(x-4)(x+2+\sqrt{3})(x+2-\sqrt{3}))=0$ yields a solution $x$ $=$ $\sqrt[3]{2+\sqrt{-124}}+\sqrt[3]{2-\sqrt{-124}}$.} Yet, he did not understand the nature of complex numbers and even thought that they are useless. Subsequently, Rafael Bombelli\index{Bombelli, Rafael (1526--1572)} (1526--1572) took down the rules for multiplication of complex numbers ({\it binomio}) in his Book I of {\it L'Algebra}, while Descartes\index{Descartes, Ren\'{e} (1596--1650)} used the term ``imaginary number" ({\it nombre imaginaire}) with the meaning of contempt. Wallis\index{Wallis, John (1616--1703)} insisted that imaginary numbers are not unuseful and absurd when properly understood by using a geometric model just like negative numbers ({\it Algebra}, Vol. II, Chap. LXVI, 1673). Meanwhile, Euler\index{Euler, Leonhard (1707--1783)} persuaded himself that there is an advantage for the use of imaginary numbers (1751),
%\footnote{{\it Recherches sur les racines imaginaires des \'{e}quations}, M\'{e}m.~\!Acad.~\!R.~\!Sci.~\!Berlin, {\bf 5} (1751),  222--288.} 
and introduced the symbol $i$ for the imaginary unit (1755). It was around 1740 that he found the earth-shaking formula $e^{\sqrt{-1}\theta}=\cos\theta+\sqrt{-1}\sin\theta$.

In 1799, 
Gauss\index{Gauss, Johann Carl Friedrich (1777--1855)} gained his doctorate in absentia from the University of Helmstedt. The subject of his thesis is {\it the fundamental theorem of algebra}\index{fundamental theorem of algebra},\footnote{{\it Demonstratio nova theorematis omnem functionem algebraicam rationalem integram unius variabilis in factores reales primi vel secundi gradus resolvi posse}. However, his proof had a gap. He gave three other proofs (1816, 1849) (see Sect.~\!\ref{subsec:topodes} for a ``topological" proof).} which, with a long history of inspiring many mathematicians, says that every polynomial equation with complex coefficients has at least one complex root.\footnote{Albert Girard\index{Girard, Albert (1595--1632)} (1595--1632) was the first to suggest the fundamental theorem of algebra\index{fundamental theorem of algebra} (FTA) ({\it Invention Nouvelle en l'Alg\`{e}bre}, 1629). Meanwhile, Leibniz\index{Leibniz, Gottfried Wilhelm von (1646--1716)} claimed that $x^4+1=0$ affords a counterexample for the FTA (1702).
He assumed that the square root of $i$ should be a more complicated ``imaginary" entity, not able to be expressed as $a+bi$; but the truth is that $x=\pm \frac{1}{\sqrt{2}}(1\pm i)$ (any double sign). Euler\index{Euler, Leonhard (1707--1783)} (1742) and D'Alembert\index{d'Alembert, Jean-Baptiste le Rond (1717--1783)} (1746) counted on the FTA when they dealt with indefinite integrals of rational functions.} As the title of his thesis indicates, he shied away from the use of imaginary numbers, notwithstanding that Gauss\index{Gauss, Johann Carl Friedrich (1777--1855)} already had a tangible image. In a letter to F. W. Bessel, dated December 18, 1811, he communicated, ``first I would like to ask anyone who wishes to introduce a new function into analysis to explain whether he wishes it to be applied merely to real quantities, and regard imaginary values of the argument only as an appendage, or whether he agrees with my thesis that in the realm of quantities the imaginaries $a+b\sqrt{-1}$ have to be accorded equal rights with the reals. Here it is not a question of practical value; analysis is for me an independent science, which would suffer serious loss of beauty and completeness, and would have constantly to impose very tiresome restrictions on truths which would hold generally otherwise, if these imaginary quantities were to be neglected..." ({\it Werke}, X, p.~\!366--367; see \cite{koch}).

In 1833, William Rowan Hamilton\index{Hamilton, William Rowan (1805-1865)} (1805-1865) conceived the idea to express a complex number $a+b\sqrt{-1}$ by the point $(a,b)$ in the coordinate plane $\mathbb{R}^2$,\footnote{Gauss\index{Gauss, Johann Carl Friedrich (1777--1855)} had been in possession of the geometric representation of complex numbers since 1796. In 1797, C. Wessel presented a memoir to the Copenhagen Academy of Sciences in which he announced the same idea, 
% and interpreted addition as translation and multiplication as rotation, 
but it did not attract attention. In 1806, J. R. Argand went public with the same formulation.} which motivated him to construct 3D ``numbers" with arithmetical operations similar to complex numbers, but the quest was a dead end. Alternatively, he discovered {\it quaternions}\index{quaternion}, the 4D numbers at the expense of commutativity of multiplication (1843). Specifically, the multiplication is defined by 

\vspace{-0.5cm}
{\small 
\begin{eqnarray}
&&( a_1, b_1, c_1, d_1)\cdot ( a_2, b_2, c_2, d_2  )
= ( a_1a_2- b_1b_2 -c_1c_2-d_1   d_2, 
~a_1 b_2 + b_1a_2 \nonumber\\
&& + c_1 d_2 - d_1 c_2,
~ a_1 c_2 - b_1 d_2   + c_1 a_2   + d_1b_2,~  a_1d_2  + b_1c_2- c_1b_2   + d_1a_2).\label{fn:quat}
\end{eqnarray}

}
\noindent In Hamilton's\index{Hamilton, William Rowan (1805-1865)} notations, $i=(0,1,0,0)$, $j=(0,0,1,0)$, $k=(0,0,0,1)$, so that quaternions\index{quaternion} are represented in the form: $a + bi + cj + dk$ $~(a,b,c,d\in \mathbb{R})$. A fact deserving special mention is that putting $\|\alpha\|=\sqrt{a^2+b^2+c^2+d^2}$ for $\alpha=(a,b,c,d)$, we have $\|\alpha\cdot\beta\|=\|\alpha\|\|\beta\|$, the identity similar to $|zw|=|z||w|$ for complex numbers, which agrees with the  {\it four-square identity} discovered by Euler\index{Euler, Leonhard (1707--1783)} in 1748 during his investigation on sums of four squares (see Remark \ref{rem:Grave}).

It may be said, incidentally, that Gauss\index{Gauss, Johann Carl Friedrich (1777--1855)} had already discovered quaternions\index{quaternion} in 1818 ({\it Mutationen des Raumes}, {\it Werke}, VIII, 357--361). What he observed is, though his original expression is slightly different, that the rotation group ${\rm SO}(3)$ is parameterized by $(a,b,c,d)$ with $a^2\!+\!b^2\!+\!c^2\!+\!d^2\!=\!1$ as 
$$ 
A(a,b,c,d)=
\begin{pmatrix}
a^2-b^2-c^2+d^2 & -2(ab+cd) & 2(bd-ac)\\
2(ab-cd) & a^2-b^2+c^2-d^2& -2(ad+bc)\\
2(ac+bd) & 2(ad-bc) & a^2+b^2-c^2-d^2
\end{pmatrix},
$$
and 
$
A(a_1,b_1,c_1,d_1)A(a_2,b_2,c_2,d_2)\!=\!A\big((a_1,b_1,c_1,d_1)\cdot(a_2,b_2,c_2,$ $d_2)\big)
$.
Note here that quaternions\index{quaternion} $\alpha=a+bi+cj+dk$ with $\|\alpha\|=1$ form a group, which is identified with the {\it spin group}\index{spin group} ${\rm Spin}(3)$, pertinent to the quantum version of {\it angular momentum} that describes {\it internal degrees of freedom} of electrons. Additionally, the kernel of the homomorphism $\alpha\mapsto A(\alpha)$ is a homomorphism of ${\rm Spin}(3)$ onto ${\rm SO}(3)$ whose kernel is $\{\pm 1\}$ (Remark \ref{spinxxx}). 

\begin{rem}\label{fn:complexline} {\rm {\small 
Euler\index{Euler, Leonhard (1707--1783)} discovered his formula while studying differential equations with constant coefficients (1748). Previously, de Moivre derived a formula that later brought on what we now call {\it de Moivre's theorem} (1707), a precursor of Euler's formula. Subsequently, R. Cotes had given the formula $ix\!=\!\log (\cos x\!+\!i\sin x)$ in 1714,
but overlooked that the {\it complex logarithm}\index{complex logarithm} assumes infinitely many values, differing by multiples of $2\pi i$, as found by Euler\index{Euler, Leonhard (1707--1783)} (1746).\footnote{The {\it logarithmic function}\index{logarithmic function} originated from the computational demands of the late 16th century in observational astronomy, long-distance navigation, and geodesy. The main contributors are John Napier\index{Napier, John (1550--1617)} (1550--1617), Joost B\"{u}rgi, Henry Briggs, Oughtred\index{Oughtred, William (1574--1660)}, and Kepler\index{Kepler, Johannes (1571--1630)}. The term {\it logarithm}, literally ``ratio-number" from \textgreek{l'ogos} and \textgreek{`arijm'os}, was coined by Napier\index{Napier, John (1550--1617)} in his pamphlet {\it Mirifici Logarithmorum Canonis Descriptio} (1614). Kepler\index{Kepler, Johannes (1571--1630)} studied Napier's pamphlet in 1619, and published the logarithmic tables {\it Chilias Logarithmoria} (1624), which he used in his calculations of the {\it Tabulae Rudolphinae Astronomicae} (1627), a star catalog and planetary tables based on the observations of Tycho Brahe\index{Tycho Brahe (1546--1601)}.}

Before the problem of the complex logarithm\index{complex logarithm} was settled by Euler\index{Euler, Leonhard (1707--1783)}, there were arguments between Leibniz\index{Leibniz, Gottfried Wilhelm von (1646--1716)} and Johann Bernoulli\index{Bernoulli, Johann (1667--1748)} (1667--1748) in 1712 and then between Euler and Jean-Baptiste le Rond d'Alembert\index{d'Alembert, Jean-Baptiste le Rond (1717--1783)} (1717--1783) around 1746 about $\log x$ for a negative $x$. Leibniz says that $\log(-1)$ does not exist, while Bernoulli alleges that $\log (-x)=\log x$. In the above-mentioned letter to Bessel in 1811,  Gauss\index{Gauss, Johann Carl Friedrich (1777--1855)} renewed the question about the {\it multivaluedness} of the complex logarithm\index{complex logarithm}, representing $\log z$ as the complex line integral $\int_c 1/z~\!dz$, where $c$ is a curve in $\mathbb{C}\backslash \{0\}$ joining $1$ and $z$ (see Sect.~\!\ref{subsec:topodes}). He further referred to what we now call {\it Cauchy's integral theorem}.\footnote{Augustin-Louis Cauchy\index{Cauchy, Augustin-Louis (1789--1857)} defined complex line integrals in 1825 and established his theorem in 1851.}
\hfill$\Box$

}
}
\end{rem}

\section{Euclid's Elements}\label{subsec:elements}
Gauss\index{Gauss, Johann Carl Friedrich (1777--1855)} is one of the discoverers of {\it non-Euclidean geometry}\index{non-Euclidean geometry} (a geometry consistently built under the assumption that there are more than one line through a point which do not meet a given line), but he did not publish his work and only disclosed it in letters to his friends.\footnote{For instance, his discovery was mentioned in a letter to F. A. Taurinus dated November 8, 1824 ({\it Werke}, VIII, p.~\!186).} His hesitation in part came from the atmosphere in those days; that is, his new geometry was entirely against the predominant Kantianism. He wanted to avoid controversial issues,\footnote{In a letter to Bessel dated January 27, 1829, Gauss\index{Gauss, Johann Carl Friedrich (1777--1855)} writes, ``I fear the cry of the Boetians [known as vulgarians] if I were to voice my views" ({\it Werke}, VIII, p.~\!200).} and requested his friends not to make it public.

The other discoverers of the new geometry are Nikolai Ivanovich Lobachevsky\index{Lobachevsky, Nikolai Ivanovich (1792--1856)} (1792--1856) and J\'{a}nos Bolyai\index{Bolyai, J\'{a}nos (1802--1860)} (1802--1860).\footnote{Strictly speaking, both Lobachevsky and Bolyai discovered the 3D geometry, while Gauss\index{Gauss, Johann Carl Friedrich (1777--1855)} dealt only with the 2D case (Sect.~\!\ref{subsec:discovery}).} In 1826, Lobachevsky stated publicly that the new geometry exists.\footnote{Lobachevsky,  {\it Exposition succincte des principes de la g\'{e}om\'{e}trie avec une d\'{e}monstration rigoureuse du th\'{e}or\`{e}me des parall\`{e}les};{\it O nachalakh geometrii}, Kazanskii Vestnik, (1829--1830).} Fourteen years later,  as this research did not draw attention, he issued a small book in German containing a summary of his work to appeal to the mathematical community,\footnote{{\it Geometrischen Untersuchungen zur Theorie der Parallellinien} (1840).} and then put out his final work {\it Pang\'{e}om\'{e}trie ou pr\'{e}cis de g\'{e}om\'{e}trie fond\'{e}e sur une th\'{e}orie g\'{e}n\'{e}rale et rigoureuse des parall\`{e}les}, a year before his death (1855); see \cite{loba}. Meanwhile, in a break with tradition, Bolyai independently began his head-on approach to the non-Euclidean properties as early as 1823. In a letter to his father Farkas, on November 3, 1823, J\'{a}nos\index{Bolyai, J\'{a}nos (1802--1860)} says with confidence, ``I discovered a whole new world out of nothing," and made public his discovery in a 26-pages appendix {\it Scientiam spatii absolute veram exhibens} to his father's book.\footnote{ {\it Tentamen juventutem studiosam in elementa Matheseos Purae} (1832).}

Gauss\index{Gauss, Johann Carl Friedrich (1777--1855)} was a perfectionist by all accounts. He did not publish many outcomes, fearing that they were never perfect enough; his motto was ``{\it pauca sed matura}---few, but ripe." After his death it was discovered that many results credited to others had been already worked out by him (remember Cauchy's\index{Cauchy, Augustin-Louis (1789--1857)} integral theorem; Remark \ref{fn:complexline}). The extraordinariness of Gauss\index{Gauss, Johann Carl Friedrich (1777--1855)} may be narrated by the joke: {\it Suppose you discovered something new. If Gauss would be still alive and would browse your result, then he would say, ``Ah, that's in my paper," and surely would take an unpublished article out of a drawer in his desk.} This was no joke for J\'{a}nos Bolyai\index{Bolyai, J\'{a}nos (1802--1860)}. Gauss\index{Gauss, Johann Carl Friedrich (1777--1855)} received a copy of his father's book in 1832, and endorsed the discovery in the response. To J\'{a}nos' disappointment, his cool reply reads, ``The entire contents of your son's work coincides almost exactly with my own meditation, which has occupied my mind for from thirty to thirty-five years" ({\it Werke}, VIII, p.~\!220). So much disappointed, J\'{a}nos wholly withdrew from scientific activity, though Gauss\index{Gauss, Johann Carl Friedrich (1777--1855)} wrote to Gerling (February 14, 1832), ``I consider the young geometer Bolyai a genius of the first rank" (see \cite{cox}). As regards Lobachevsky\index{Lobachevsky, Nikolai Ivanovich (1792--1856)}, it was not until 1841 that his work was recognized by Gauss\index{Gauss, Johann Carl Friedrich (1777--1855)} who looked over the {\it Geometrischen Untersuchungen} by chance. Gauss was very impressed with it, and praised Lobachevsky\index{Lobachevsky, Nikolai Ivanovich (1792--1856)} as a clever mathematician in a letter to his friend J. F. Encke (February 1841; {\it Werke}, VIII, p.~\!232).\footnote{F. K. Schweikart, a professor of law and an uncle of Taurinus, had a germinal idea of non-Euclidean geometry\index{non-Euclidean geometry} (``astral geometry" in his term), and asked for Gauss\index{Gauss, Johann Carl Friedrich (1777--1855)}, through Gerling, to comment on the idea in 1818. In response to the request (March 16, 1819), Gauss\index{Gauss, Johann Carl Friedrich (1777--1855)} communicated, with compliments, that he concurred in Schweikart's observation ({\it Werke}, VIII, p.~\!181). However Schweikart neither elaborated his idea nor published any result.}

\medskip

Now, what is known about Euclid\index{Euclid ($4^{\rm th}$-$3^{\rm rd}$ c. BC)}
 and his {\it Elements}\index{Elements}? It is popularly thought that Euclid was attached to Plato's Academeia\index{Academeia}\index{Plato (c. 428--c. 348 BC)} or at least was affiliated with it, and that he was temporarily a member of the Great Museum \index{Great Museum} during the reign of Ptolemy I. Truth to tell, very little is known about his life. Archimedes\index{Archimedes (c. 287--c. 212 BC)}, whose life, too, overlapped the reign of Ptolemy I, mentioned briefly Euclid\index{Euclid ($4^{\rm th}$-$3^{\rm rd}$ c. BC)}
 in his {\it On the Sphere and the Cylinder I.2},  though being a fiction made by the posterity. In the {\it Collection VII}, Pappus\index{Pappus (c. 290--c. 350)} records  that Apollonius\index{Apollonius (c. 262--c. 190 BC)} spent a long time with the disciples of Euclid\index{Euclid ($4^{\rm th}$-$3^{\rm rd}$ c. BC)}
in Alexandria. Apollonius\index{Apollonius (c. 262--c. 190 BC)} himself referred to Euclid's achievement with the statement that his method surpasses Euclid's\index{Euclid ($4^{\rm th}$-$3^{\rm rd}$ c. BC)}
{\it Conics}. Proclus Diadochus\index{Proclus (412--485)} (412 CE--485 CE) and Joannes Stobabeing (the 5th century CE) recorded the oft-told (but not trustworthy) anecdotes about Euclid. The sure thing is that he is the author of the definitive mathematical treatise {\it Elements}\index{Elements} (\textgreek{Stoiqe\~ia}),\footnote{Besides the {\it Elements}\index{Elements}, at least five works of Euclid\index{Euclid ($4^{\rm th}$-$3^{\rm rd}$ c. BC)}
have survived to today; {\it Data}, {\it On Divisions of Figures}, {\it Catoptrics}, {\it Phaenomena}, {\it Optics}. There are a few other works that are attributed to Euclid but have been lost: {\it Conics}, {\it Porisms}, {\it Pseudaria}, {\it Surface Loci}, {\it Mechanics}.} which is a crystallization of the plane and space geometry known in his day ({\it ca}.~\!300 BCE) and was so predominant that all earlier texts were driven out. It is solely to the {\it Elements}\index{Elements}, one of the most influential books in long history, which Euclid\index{Euclid ($4^{\rm th}$-$3^{\rm rd}$ c. BC)}
owes his abiding fame in spite of his shadowy profile.

The {\it Elements}\index{Elements}, consisting of thirteen books, starts with the construction of equilateral triangles and ends up with the classification of the regular solids; {\it tetrahedron}, {\it cube}, {\it octahedron}, {\it dodecahedron}, and {\it icosahedron} (Book XIII, Prop.~\!465). This suggests the influence of Pythagorean\index{Pythagorean} mysticism that put great emphasis on aesthetic issues, which, nowadays, are described in terms of symmetry\index{symmetry} (\textgreek{summetr'ia}).\footnote{``Symmetry" is explicated by {\it group actions}, with which both spatial homogeneity and isotropy are described. Klein's\index{Klein, Christian Felix (1849--1925)} Erlangen Program\index{Erlangen Program} is also explained in this framework.

Aetius (1st- or 2nd-century CE) says that Pythagoras\index{Pythagoras (c. 580--c. 500 BC)} discovered the five regular solids ({\it De Placitis}), while Proclus\index{Proclus (412--485)} argues that it was the historical Theaetetus\index{Theaetetus (c. 417--c. 369 BC)} ({\it ca.}~\!417 BCE--{\it ca.}~\!369 BCE) who theoretically constructed the five solids. Plato's\index{Plato (c. 428--c. 348 BC)} {\it Timaeus}\index{Timaeus@Timaeus (dialogue)} (53b5–6) provides the earliest known description of these solids as a group. Nowadays, regular polyhedra are discussed often in relation to finite subgroups ({\it polyhedral groups}) of ${\rm SO}(3)$ or ${\rm O}(3)$. For example, the {\it icosahedral group} is the rotational symmetry group of the icosahedron (see Remark \ref{spinxxx}).} Moreover, it has an orderly organization consisting of 23 fundamental definitions, 5 postulates,  5 axioms, and 465 propositions. The ``reductio ad absurdum"\index{reductio ad absurdum} is made full use of as a powerful gambit (Prop.~\!6 in Book I is the first one proved by reductio ad absurdum). 

As often said, only a few theorems are thought of having been discovered by Euclid\index{Euclid ($4^{\rm th}$-$3^{\rm rd}$ c. BC)}
himself. The theorem of angle-sum\index{theorem of angle-sum} (Book I, Prop.~\!32), the Pythagorean Theorem\index{Pythagorean Theorem} (Book I, Prop.~\!47), and perhaps the {\it triangle inequality}\index{triangle inequality} (Book I, Prop.~\!20) are attributed to the Pythagoreans\index{Pythagorean}. The theory of proportions (Book V) and the method of exhaustion\index{method of exhaustion} (Book X, Prop.~\!1, Book XII) are, as we mentioned before, greatly indebted to Eudoxus\index{Eudoxus (c. 408--c. 355 BC)}.

\medskip

We shall talk about a checkered history of the {\it Elements}\index{Elements}, originally written on fragile papyrus scrolls.\footnote{In the tiny fragment ``Papyrus Oxyrhynchus 29"
% discovered in 1897, 
housed in the library of the University of Pennsylvania, one can see a diagram related to Prop.~\!5 of Book II, which can be construed in modern terms as $ab+\big(\frac{a+b}{2}-b\big)^2=\big(\frac{a+b}{2}\big)^2$.} The book survived in Alexandria and Byzantium under the control of the Roman Empire. In contrast, Roman people---the potentates of the world during this period---by and large expressed little interest in pure sciences as testified by Marcus Tullius Cicero (BCE 106--BCE 43), who confesses, ``Among the Greeks, nothing was more glorious than mathematics. We, however, have limited the usability of this art to measuring and calculating" ({\it Tusculanae Disputationes}, {\it ca}.~\!45 BCE). Luckily, Alexandria is far away from the central government of the Empire, disintegrating in scandal and corruption. Situated as it was and thanks to the Hellenistic tradition, pure science flourished there. The people involved were Heron ({\it ca.}~\!10 CE--{\it ca.}~\!70 CE), Menelaus\index{Menelaus (c. 70--c. 140)} ({\it ca.}~\!70 CE--{\it ca.}~\!140 CE), Ptolemy\index{Ptolemy (c. 100--c. 170)}, Diophantus ({\it ca.}~\!215 CE--{\it ca.}~\!285 CE), Pappus\index{Pappus (c. 290--c. 350)}, Theon Alexandricus ({\it ca.}~\!335 CE--{\it ca.}~\!405 CE), and Hypatia ({\it ca.}~\!350 CE--415 CE).\footnote{Diophantus stands out as unique because he applied himself to some algebraic problems in the time when geometry was still in the saddle. Theon is known for editing the {\it Elements}\index{Elements} with commentaries. His version was the only Greek text known, until a hand-written copy of the {\it Elements}\index{Elements} in the 10th century was discovered in the Vatican library (1808). Hypatia, a daughter of Theon, is the earliest female mathematician in history, and is far better known for her awful death; she was dragged to her death by fanatical Christian mobs.}

But even in the glorious city, scientific attitude declined with the lapse of time, as the academic priority had been given, if anything, to annotations on predecessor's work, and eventually fell into a state of decadence. This was paralleled by the steady decay of the city itself that all too often suffered overwhelming natural and manmade disasters. After a long tumultuous times, the great city was finally turned over to Muslim hands (641 CE). 

In the Dark Ages, Europeans could no longer understand the {\it Elements}\index{Elements}. In this circumstance, Anicius Manlius Boethius ({\it ca}.~\!480 CE--524 CE) and Flavius Magnus Aurelius Cassiodorus Senator ({\it ca}.~\!485 CE--{\it ca}.~\!585 CE) are very rare typical scholars in this period who were largely influential during the Middle Ages. 
%The former made an adapted translation of the {\it Introductionis Arithmeticae} by Nicomachus ({\it ca}.~\!60 CE--{\it ca}.~\!120 CE). The latter wrote two books {\it De Institutione Divinarum Litterarum} and {\it De Artibus ac Disciplinis Liberalium Litterarum}, discussing the quadrivium, as a guide for his monks. 
%Other people, if we dare to name, are St.~\!Isidore of Seville ({\it ca}.~\!560 CE--636 CE), St.~\!Bede ({\it ca}.~\!673 CE--735 CE), Alcuin of York ({\it ca}.~\!735 CE--804 CE), and Gerbert\index{Gerbert} (Remark \ref{rem:naturalnumber}~\!(2)). All of them were more or less connected to Church, and only contributed to elementary calculation methods or number theory that lapses into {\it numerology}.
Meanwhile, the {\it Elements}\index{Elements} was brought from Byzantium to the Islamic world in {\it ca}.~\!760 CE, and was taught in Bagdad, Cordova, and Toledo  where there were large-scale libraries comparable with the Great Museum\index{Great Museum} in Alexandria.\footnote{Worthy to mention is the {\it House of Wisdom} (Bayt al-Hikma) in Bagdad.
%, a melting pot of different ethnic groups in the Islamic world. 
It was founded by Caliph Harun al-Rashid ({\it ca}.~\!763 CE--809 CE) as an academic institute devoted to translations, research, and education, and was culminated under his son Caliph al-Mamoon. Al-Khw$\overline{\rm a}$rizm\=\i\index{alk@al-Khw$\overline{\rm a}$rizm\=\i~(c. 780--c. 850)} was a scholar at this institute.}

After approximately 400 years of blank period in Europe, Adelard of Bath\index{Adelard (c. 1080--c. 1152)} (Remark \ref{rem:naturalnumber}~\!(2)), who undertook a journey to Cordova around 1120 to delve into Arabic texts of Greek classics, obtained a copy of the book and translated it into Latin. Later on, Gerard \index{Gerard (c. 1114--1187)} (fn.~\!\ref{fn:Alma}) translated it into Latin from another Arabic version procured in Toledo. His translation is considered being close to the Greek original text. Then, in the mid-renaissance, only 27 years later from the type printing of the {\it Bible} by Johannes Gutenberg in Mainz, the first printed edition---based on the Latin version from the Arabic by Campanus of Novara who probably had access to Adelard's translation---was published in Venice. Thereafter, versions of the {\it Elements}\index{Elements} in various languages had been printed. The first English-language edition was printed in 1570 by H. Billingsley.
%, a merchant and Lord Mayor of London. 
Then in 1607, Matteo Ricci and Xu Guangqi translated into Chinese the first six volumes of the Latin version published in 1574 by Christopher Clavius.

One of the first multicolored books {\it The First Six Books of the Elements\index{Elements} of Euclid} by O. Byrne was printed in 1847 in London. A. De Morgan was very critical of its non-traditional style embellished by pictorial proofs, and judged this inventive attempt to be nonsense ({\it A Budget of Paradoxes}, 1872). At that time, 
%(the Victorian Era)
there was considerable debate about how to teach geometry; say, the pros and cons of whether to adopt a new approach to Euclidean geometry\index{Euclidean geometry} in teaching. C. Dodgson in Oxford, known as Lewis Carroll, assumed a critical attitude to his colleague Playfair\index{Playfair, John (1748--1819)} who tried to simplify the Euclid's\index{Euclid ($4^{\rm th}$-$3^{\rm rd}$ c. BC)}
proofs by introducing algebraic notations. In the little book {\it Euclid and his Modern Rivals} printed in 1879, he argues, ``no sufficient reasons have yet been shown for abandoning [Euclid's {\it Elements}\index{Elements}] in favour of any one of the modern Manuals which have been offered as substitutes."  G. B. Halsted says, in the translator's introduction to J. Bolyai's\index{Bolyai, J\'{a}nos (1802--1860)} {\it Scientiam spatii},  ``Even today (1895), in the vast system of examinations carried out by the British Government, by Oxford, and by Cambridge, no proof of a theorem in geometry will be accepted which infringes Euclid's sequence of propositions." As a matter of course, it is now seldom to make direct use of Euclid's\index{Euclid ($4^{\rm th}$-$3^{\rm rd}$ c. BC)}
propositions for university examinations, but the book went through more than 2000 editions to date, and has been (and still is) the encouragement and guide of scientific thought. Furthermore, modern mathematics inherits much of its style from the {\it Elements}\index{Elements}; in particular, many mathematical theories, even if not all of them, begin with axiomatic systems (Sect.~\!\ref{subsec:transfinite}). 

What about the original Greek text? It was Heiberg (fn.~\!\ref{fn:heiberg}) who tracked down all extant manuscripts---the aforementioned Greek Vatican manuscript is one of them---to produce a definitive Greek text together with its Latin translation and prolegomena ({\it Euclidis Opera Omnia}, 8 vols, 1883--1916). This is a trailblazing achievement that has become the basis of later researches on Euclid\index{Euclid ($4^{\rm th}$-$3^{\rm rd}$ c. BC)}.

\section{The Fifth Postulate}\label{subsec:fifth}
Euclid's\index{Euclid ($4^{\rm th}$-$3^{\rm rd}$ c. BC)}
greatest contribution is his daring selection of a few postulates as major premises (if not completely adequate from today's view), of which Greek geometers thought as self-evident truth requiring no proof. He built up a deductive system in which every theorem is derived from these postulates. Among his five postulates, the last one, called the {\it Fifth Postulate}\index{Fifth Postulate} (FP henceforward), constitutes the core of our story on non-Euclidean geometry\index{non-Euclidean geometry}. It says: 

\medskip

{\it If a straight line falls on two straight lines in such a manner that the interior angles on the same side are together less than two right angles, then the straight lines, if produced indefinitely, meet on that side on which are the angles less than the two right angles} (\cite{euclid}, Book I).

\medskip

The theorem of angle-sum\index{theorem of angle-sum} and the Pythagorean Theorem\index{Pythagorean Theorem} rely on this postulate. For example, what is required for the former is the fact that ``if a straight line fall upon two parallel straight lines it will make the alternate angles equal to one another" (Prop.~\!29 in Book I), which is a consequence of the FP.\footnote{In contrast to these two theorems, the triangle inequality\index{triangle inequality} is proved  without invoking the FP (Book I, Prop.~\!20).} 

The first four postulates have fairly simple forms; say, the first one declares, in today's terms, that given two points, there exists a straight line through these points, whereas the necessity of the last postulate is by no means overt, owing not only to the intricacy of its formulation, but also to the fact that the converse of Prop.~\!29 is proved without an appeal to the last postulate (Prop.~\!27). Actually the last one had been believed to be redundant and hence to be a ``proposition" all through modern times until the second decade of the 19th century (\cite{bon}).

\begin{figure}[htbp]
%\vspace{-4.7cm}
\vspace{-0.4cm}
\begin{center}
\includegraphics[width=.34\linewidth]
{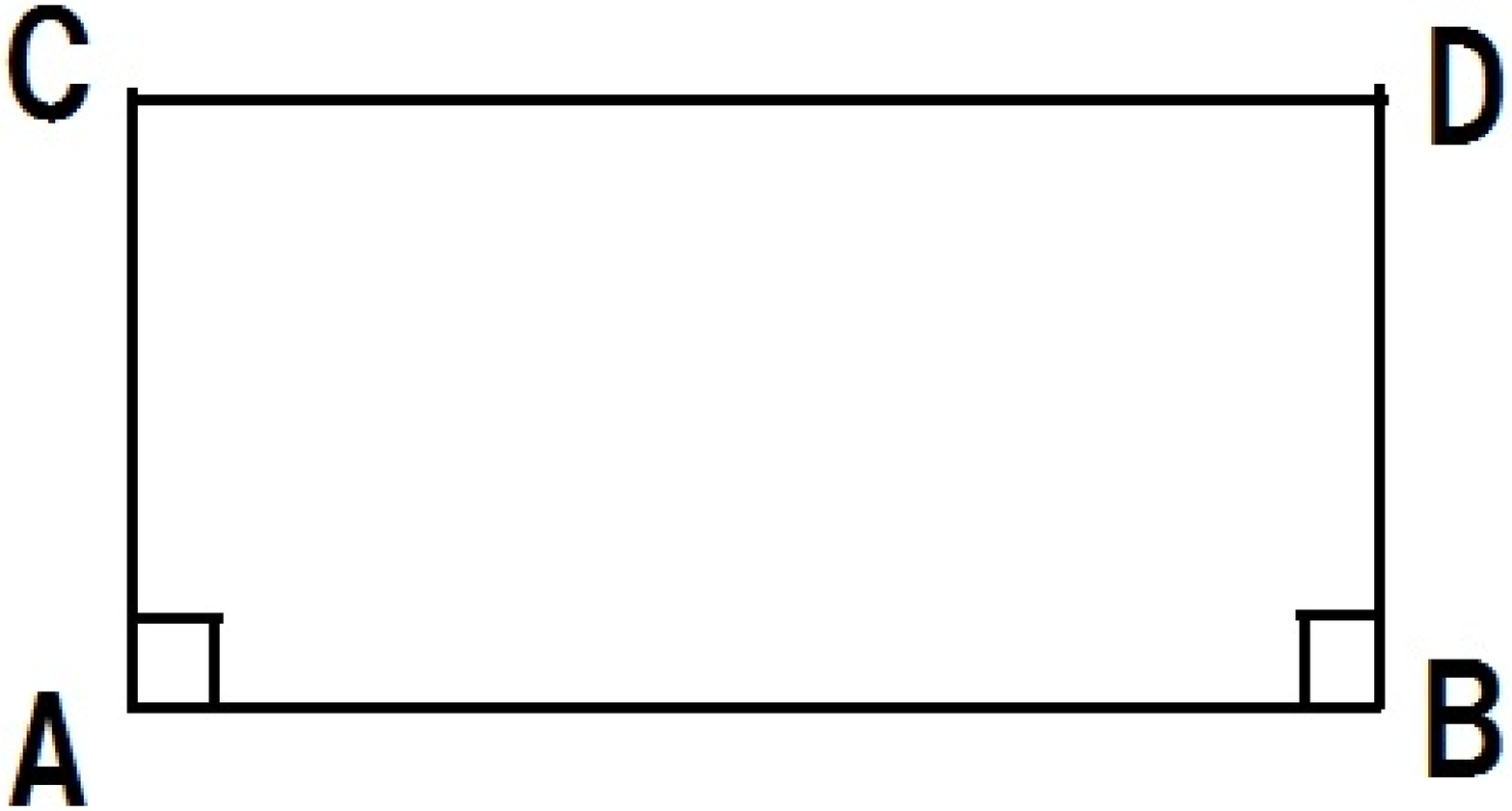}
\end{center}
\vspace{-1.3cm}
\caption{Saccheri's quadrilateral}\label{fig:sacc}
\vspace{-0.1cm}
\end{figure}

An elaborate and stalwart attempt to vindicate Euclidean geometry\index{Euclidean geometry} by ``proving" the FP\index{Fifth Postulate} was done by Girolamo Saccheri\index{Saccheri, Girolamo (1667-1733)} (1667-1733), an Italian Jesuit priest, who explored the consequences under the negation of the postulate, hoping to reach a contradiction. For this sake, he made use of quadrilaterals of which two opposite sides are equal to each other and perpendicular to the base, and set up the following three hypotheses (Fig.~\!\ref{fig:sacc});\footnote{Omar Khayyam\index{Khayyam, Omar (1048--1131)} (1048--1131) dealt with the FP\index{Fifth Postulate} in his {\it Explanation of the Difficulties in the Postulates of Euclid}. He is considered a predecessor of Saccheri\index{Saccheri, Girolamo (1667-1733)}.} 

\smallskip
(i) {\it the hypothesis of the right angle}: $\angle C=\angle D=\angle R$, 

\smallskip

(ii) {\it the hypothesis of the obtuse angle}: $\angle C=\angle D>\angle R$, and 

\smallskip
(iii) {\it the hypothesis of the acute angle}: $\angle C=\angle D<\angle R$, 

\smallskip

\noindent where the first hypothesis is equivalent to the FP\index{Fifth Postulate}. Examining these hypotheses meticulously, he concluded that hypothesis (ii) 
contradicts the infinite extent of plane, and is false, but could not exclude hypothesis (iii). At the final stage of his argument, Saccheri\index{Saccheri, Girolamo (1667-1733)} appeals to intuition about our space, and says, ``{\it Hypothesis anguli acuti est absolute falsa; quia repugnans naturae liniae rectae}" [``the hypothesis of the acute angle is absolutely false; because it is repugnant to the nature of the straight line”] (\cite{bon}). In summary, he virtually proved that, under the negation of the FP\index{Fifth Postulate}, the angle-sum of any triangle is less than $\pi$ (thus the existence of a triangle with angle-sum of $\pi$ implies the FP\index{Fifth Postulate}).\footnote{Saccheri, {\it Euclides ab omni naevo vindicatus sive conatus geometricus quo stabiliuntur prima ipsa geometriae principia} 
(1733). See \cite{sacc}.}

Johann Heinrich Lambert\index{Lambert, Johann Heinrich (1728--1777)} (1728--1777)---conceivably familiar with Saccheri's\index{Saccheri, Girolamo (1667-1733)} work because, in his {\it Theorie der Parallellinien} (1766), he quoted the work of G. S. Kl\"{u}gel who listed nearly 30 attempts to prove the FP\index{Fifth Postulate} including Saccheri's\index{Saccheri, Girolamo (1667-1733)} work\footnote{{\it Conatuum praecipuorum theoriam parallelarum demonstrandi recensio} (1763).}---investigated the FP\index{Fifth Postulate} in alignment with an idea resembling that of Saccheri\index{Saccheri, Girolamo (1667-1733)} (\cite{bon}, p.~\!44).\footnote{Lambert proved that the transfer time from  $\boldsymbol{x}_1$ to $\boldsymbol{x}_2$ along a Keplerian\index{Kepler, Johannes (1571--1630)} orbit is independent of the orbit's eccentricity and depends only on $\|\boldsymbol{x}_1\|+\|\boldsymbol{x}_2\|$ and $\|\boldsymbol{x}_1-\boldsymbol{x}_2\|$ ({\it \"{U}ber die Eigenschaften der Kometenbewegung}, 1761). He was the first to express Newton's\index{Newton, Issac (1642--1726)} second law of motion in the notation of differential calculus ({\it Vis Viva}, 1783).}  His work, published after the author's death, bristles with far-sighted observations and highly speculative remarks (see \cite{papa2}). One of his revelations is that, under the negation of the FP\index{Fifth Postulate}, the {\it angular defect}\index{angular defect} $\pi-(\angle A+\angle B+\angle C)$ for a triangle $\triangle ABC$ is proportional to the area of the triangle. This is inferred from the fact that if we define $m(P)$ for a polygon $P$ with $n$ sides by setting
$m(P)=(n-2)\pi-(\text{the sum of inner angles of}~P)$, then $m(\cdot)$ and the area functional ${\rm Area}(\cdot)$ share the ``additive property"; i.e., $m(P)=m(P_1)+m(P_2)$ for a polygon $P$ composed of two polygons $P_1$, $P_2$. With this observation in hand, he ratiocinated that non-Euclidean plane-geometry (if it exists) has a close resemblance to {\it spherical geometry}.

Let us briefly touch on spherical geometry\index{spherical geometry}, which is nearly as old as Euclidean geometry\index{Euclidean geometry}, and was developed in connection with geography and astronomy. In his {\it Sphaerica}, extant only in an Arabic translation, Menelaus\index{Menelaus (c. 70--c. 140)} introduced the notion of {\it spherical triangle}\index{spherical triangle}, a figure formed on a sphere by three great circular arcs intersecting pairwise in three vertices. He observed that the angle-sum of a spherical triangle is greater than $\pi$ (\cite{rashed}).  Greek spherical geometry gave birth to {\it spherical trigonometry}\index{spherical trigonometry} that deals with the relations between trigonometric functions of the sides and angles of spherical triangles. At a later time, it was treated in Islamic mathematics and then by Vi\`{e}te\index{Vi\`{e}te, Fran\c{c}ois (1540--1603)}, Napier\index{Napier, John (1550--1617)}, and others. A culmination is the following fundamental formulas given by Euler\index{Euler, Leonhard (1707--1783)}.\footnote{{\it Trigonometria sphaerica universa, ex primis principiis breviter et dilucide derivata}, Acta Academiae Scientarum Imperialis Petropolitinae, {\bf 3} (1782), 72--86 (see \cite{papa1}). The spherical cosine law (\ref{eq:bbb}) is seen in  al-Batt$\overline{\rm a}$n\=\i's
%\index{alb@
\index{al-Batt$\overline{\rm a}$n\=\i~(c. 858--929)} work {\it Kit$\overline{\rm a}$b az-Z\=\i j} (fn.~\!\ref{fn:al-Batt}).}  
\begin{align}
\cos\angle A&=-\cos\angle B\cos \angle C+\sin \angle B\sin\angle C\cos \frac{a}{R},\label{eq:aaa}\\
 \cos\frac{a}{R}&=\cos\frac{b}{R}\cos\frac{c}{R}+\sin\frac{b}{R}\sin\frac{c}{R}\cos\angle A,\label{eq:bbb}
\end{align}
where $\triangle ABC$ is a spherical triangle on the sphere of radius $R$, and $a, b, c$ are the arc-length of the edges corresponding to $A, B, C$, respectively.

\begin{figure}[htbp]
\vspace{-0.3cm}
\begin{center}
\includegraphics[width=.3\linewidth]
{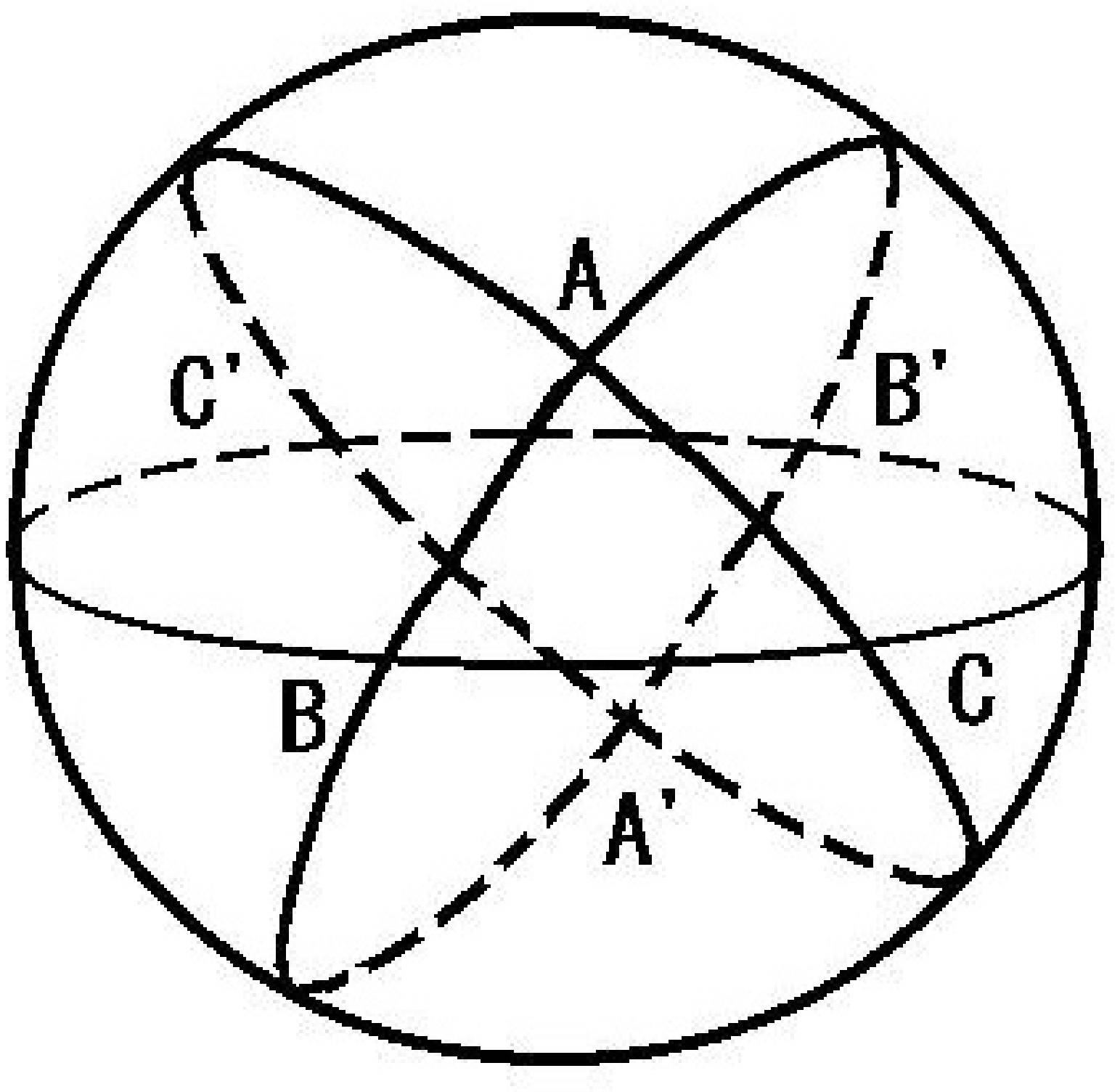}
\end{center}
\vspace{-2.8cm}
\caption{A spherical triangle}\label{fig:sphere}
\vspace{-0.2cm}
\end{figure}

Now, what Lambert\index{Lambert, Johann Heinrich (1728--1777)} paid attention to is the formula $(\angle A+\angle B+\angle C)-\pi=R^{-2}{\rm Area}(\triangle ABC)$, a refinement of Menelaus'\index{Menelaus (c. 70--c. 140)} observation on spherical triangles discovered by Harriot\index{Harriot, Thomas (c. 1560 --1621)} (1603) and Girard\index{Girard, Albert (1595--1632)} ({\it Invention nouvelle en algebra}, 1629), and probably known by Regiomontanus\index{Regiomontanus (1436--1476)}.\footnote{The Latin name of Johannes M\"{u}ller von K\"{o}nigsberg (1436--1476), the most important astronomer of the 15th century who may have arrived at the heliocentricism\index{heliocentricism}.} Replacing the radius $R$ by $\sqrt{-1}R$ in this formula,  we obtain $\pi-(\angle A+\angle B+\angle C)=R^{-2}{\rm Area}(\triangle ABC)$, which fits in with his outcome on the angular defect\index{angular defect}, and brought him to the speculation that an ``imaginary sphere"\index{imaginary sphere} (under a suitable justification) provides a model of ``non-Euclidean plane," but he did not  pursue this idea any further.

\begin{rem}\label{rem:list}{\rm {\small 
(1) Here is a list (by no means exhaustive) of people who tried to deduce the FP\index{Fifth Postulate} from other propositions: Posidonius of Apameia ({\it ca.}~\!135 BCE--{\it ca.}~\!51 BCE), Geminus of Rhodes (the 1st century BCE), Ptolemy\index{Ptolemy (c. 100--c. 170)}, Proclus\index{Proclus (412--485)}, Nas\=\i r al-D\=\i n T$\overline{\rm u}$s\=\i (1201--1274), C. Clavius (1538--1612), P. A. Cataldi (1548--1626), G. A. Borelli (1608--1679), G. Vitale (1633--1711), Leibniz\index{Leibniz, Gottfried Wilhelm von (1646--1716)}, Wallis\index{Wallis, John (1616--1703)}, J. S. K\"{o}nig (1712--1757), Adrien-Marie Legendre\index{Legendre, Adrien-Marie (1752--1833)} (1752--1833), B. F. Thibaut (1775--1832), and F. Bolyai. Gauss\index{Gauss, Johann Carl Friedrich (1777--1855)}, J. Bolyai\index{Bolyai, J\'{a}nos (1802--1860)}, and Lobachevsky\index{Lobachevsky, Nikolai Ivanovich (1792--1856)} were, for that matter, no exception at the beginning of their studies of parallels. What is in common to their arguments (except for the great trio) is that they beg the question; that is, they assume what they intend to prove. 

\smallskip

(2) The system of postulates in the {\it Elements}\index{Elements} is not complete. An airtight system was provided by Hilbert\index{Hilbert, David (1862--1943)} in the text {\it Grundlagen der Geometrie} (1899). His system consists of five main axioms;  I. {\it Axioms of Incidence}, II. {\it Axioms of Order}, III. {\it Axioms of Congruence}, IV. {\it Axiom of Parallels}, V. {\it Axioms of Continuity}.\footnote{ Moritz Pasch (1843--1930), Giuseppe Peano\index{Peano, Giuseppe (1858--1932)} (1858--1932), and  Mario Pieri (1860--1913) contributed to the axiomatic foundation of Euclidean geometry\index{Euclidean geometry}.} The Axiom of Parallels\index{Axiom of Parallels} says, ``Given a line and a point outside it, there is exactly one line through the given point which lies in the plane of the given line and point so that the two lines do not meet." This statement, clearly given by Proclus\index{Proclus (412--485)} in his commentary, was adopted by Playfair\index{Playfair, John (1748--1819)} in his edition of the {\it Elements}\index{Elements} (1795) as an alternative for the FP\index{Fifth Postulate}. Hilbert's\index{Hilbert, David (1862--1943)} Axioms of Continuity consists of the axiom of completeness related to Dedekind cuts\index{Dedekind cuts} (Remark \ref{rem:analyticalsoc} (4)) and the Axiom of Archimedes\index{Axiom of Archimedes}.\hfill$\Box$

}
}
\end{rem}

\section{Discovery of non-Euclidean geometry}\label{subsec:discovery}

We now go back to Gauss's\index{Gauss, Johann Carl Friedrich (1777--1855)} discovery of non-Euclidean geometry\index{non-Euclidean geometry}. In his letter to Taurinus in 1824, Gauss\index{Gauss, Johann Carl Friedrich (1777--1855)} mentions, ``The assumption that angle-sum [of triangle] is less than $180^{\circ}$ leads to a curious geometry, quite different from ours [the euclidean] but thoroughly consistent, which I have developed to my entire satisfaction. The theorems of this geometry appear to be paradoxical, and, to the uninitiated, absurd, but calm, steady reflection reveals that they contain nothing at all impossible" (see \cite{blu}).

Now, when did Gauss\index{Gauss, Johann Carl Friedrich (1777--1855)} become aware of the new geometry\index{non-Euclidean geometry}? A few tantalizing pieces of evidence are in Gauss's {\it Mathematisches Tagebuch}, a record of his discoveries from 1796 to 1814 ({\it Werke}, X).\footnote{Gauss's\index{Gauss, Johann Carl Friedrich (1777--1855)} diary---which contains 146 entries written down from time to time, and most of which consist of brief and somewhat cryptical statements---was kept by his bereaved until 1899. The 80th entry (October, 1797) announces the discovery of the proof of the FTA\index{fundamental theorem of algebra}.} The 72nd item written on July 28, 1797 says, ``Plani possibilitatem demostravi" [``I have demonstrated the {\it possibility} of a plane"]. This seems to indicate Gauss's\index{Gauss, Johann Carl Friedrich (1777--1855)} interest in the foundation of Euclidean geometry\index{Euclidean geometry}, of which he never lost sight from that time onward. Following that, in the 99th entry dated September 1799, he wrote, ``In principiis Geometriae egregios progressus fecimus" [``On the foundation of geometry, I could make a remarkable progress"]. It is difficult to imagine what this oblique statement means, but it could have something to do with non-Euclidean geometry\index{non-Euclidean geometry}.\footnote{\label{fn:antiKant}In a letter to F. Bolyai, March 6, 1832 ({\it Werke} VIII, p.~\!224), Gauss\index{Gauss, Johann Carl Friedrich (1777--1855)} alluded that Kant\index{Kant, Immanuel (1724--1804)} was wrong in asserting that space is only the form of our intuition.} In the same year, Gauss\index{Gauss, Johann Carl Friedrich (1777--1855)} confided to F. Bolyai to the effect that  he had doubt about the truth of geometry.\footnote{In the letter, he also said, ``If one could prove the existence of a triangle of arbitrarily great area, all Euclidean geometry would be validated" (see \cite{cox}).} However, his letter to F. Bolyai in 1804 witnessed that he was still suspicious of the new geometry ({\it Werke}, VIII, p.~\!160). All of the available witnesses point to the possibility that it is in 1817 at the latest that Gauss\index{Gauss, Johann Carl Friedrich (1777--1855)} disengaged himself from his preconceptions, and began to deal seriously with the new geometry. In a letter to his friend H. Olbers (April 28, 1817), he elusively says, ``I am becoming more and more convinced that the necessity of {\it our geometry} cannot be proved" ({\it Werke}, VIII, p.~\!177).

In all likelihood, Gauss\index{Gauss, Johann Carl Friedrich (1777--1855)} had a project to bring out an exposition of the ``Nicht-Euklidische Geometrie". His unpublished notes about the new theory of parallels might be part of this exposition ({\it Werke}, VIII, p.~\!202--209).  A letter addressed to Schumacher on May 17, 1831 says, ``In the last few weeks, I have begun to write up a few of my own meditations, which in part are already about 40 years old."\footnote{In a letter to Schumacher of July 12, 1831, Gauss\index{Gauss, Johann Carl Friedrich (1777--1855)} said that the circumference of a non-Euclidean circle of radius $r$ is $\pi k(e^{r/k}-e^{-r/k})$ where $k$ is a positive constant ({\it Werke}, VIII, p.~\!215), which turns out to coincide with the radius $R$ of the imaginary sphere.} His project was, however, suspended upon the receipt of a copy of J. Bolyai's work. Accordingly, the detailed study of the new geometry\index{non-Euclidean geometry} was put in the hands of Bolyai\index{Bolyai, J\'{a}nos (1802--1860)} and Lobachevsky\index{Lobachevsky, Nikolai Ivanovich (1792--1856)}, who established, under the negation of the FP\index{Fifth Postulate}, an analogue of trigonometry, with which they could be confident of the consistency of the new geometry\index{non-Euclidean geometry} (\cite{loba} and \cite{papa}).

The fundamental formulas in non-Euclidean trigonometry\index{non-Euclidean trigonometry} are given by
\begin{eqnarray}
&&\cos\angle A=-\cos\angle B\cos \angle C+\sin \angle B\sin\angle C\cosh \frac{a}{R},\label{eq:nonsphe1}\\
&& \cosh\frac{a}{R}=\cosh\frac{b}{R}\cosh\frac{c}{R}-\sinh\frac{b}{R}\sinh\frac{c}{R}\cos\angle A.\label{eq:nonsphe2}
\end{eqnarray}
These are obtained from the trigonometric formulas established by Bolyai\index{Bolyai, J\'{a}nos (1802--1860)} and Lobachevsky\index{Lobachevsky, Nikolai Ivanovich (1792--1856)}.\footnote{The non-Euclidean formulas were found by Taurinus ({\it Geometriae Prima Elementa}, 1826), but the outcome did not wholly convince him of the possibility of a new geometry.} Notice here that, thanks to Euler's\index{Euler, Leonhard (1707--1783)} formula, we get (\ref{eq:nonsphe1}) and (\ref{eq:nonsphe2}) by transforming $R$ into $\sqrt{-1}R$ in (\ref{eq:aaa}) and (\ref{eq:bbb}), thereby justifying Lambert's\index{Lambert, Johann Heinrich (1728--1777)} anticipation of a new geometry on the imaginary sphere.

We should underline here that non-Euclidean geometry\index{non-Euclidean geometry} involves a positive parameter $R$ just like spherical geometry involving the radius of a sphere. As Lambert\index{Lambert, Johann Heinrich (1728--1777)} took notice of it, this implies that, once $R$ is selected, the {\it absolute unit} of length is determined.\footnote{In the ordinary geometry, only a relative meaning is attached to the choice of a particular unit in the measurement. The concept of ``similarity" arises exactly from this fact.} Actually, the non-existence of an absolute unit of length is equivalent to the FP\index{Fifth Postulate}. Failing to notice this fact, Legendre\index{Legendre, Adrien-Marie (1752--1833)} believed that he succeeded in proving the FP\index{Fifth Postulate}.\footnote{{\it L\'{e}flexions sur diff\'{e}rentes mani\`{e}res de d\'{e}montrer la th\'{e}orie des parall\'{e}les ou le th\'{e}or\`{e}me sur la somme des trois angles du triangle}, 1833. In his futile attempt to prove the FP\index{Fifth Postulate}, Legendre\index{Legendre, Adrien-Marie (1752--1833)} rediscovered Saccheri's\index{Saccheri, Girolamo (1667-1733)} result. The error in his ``proof" was pointed out by Gauss\index{Gauss, Johann Carl Friedrich (1777--1855)} in a letter to Gerling (1816; {\it Werke}, VIII, p.168, 175). (By 1808, Gauss\index{Gauss, Johann Carl Friedrich (1777--1855)} realized that the new geometry must have an absolute unit of length; {\it Werke} VIII, p.165.)}

\medskip
Thus the issue of the FP\index{Fifth Postulate}, which had perplexed many able scholars for more than 2000 years, was almost settled down (the reason why we say ``almost" will become clear later). More than that, mathematicians grew to realize that postulates in general are {\it not} statements which are regarded as being established, accepted, or {\it self-evidently true}, but a sort of predefined rules. 

On reflection, Euclid's\index{Euclid ($4^{\rm th}$-$3^{\rm rd}$ c. BC)}
choice of the FP\index{Fifth Postulate} shows his keen insight and deserves great veneration. Thomas Little Heath\index{Heath, Thomas Little (1861--1940)} (1861--1940), a leading authority of Greek mathematics, remarked that ``when we consider the countless successive attempts made through more than twenty centuries to prove the Postulate, many of them by geometers of ability, we cannot but admire the genius of the man who concluded that such a hypothesis, which he found necessary to the validity of his whole system of geometry, was really indemonstrable" (\cite{euclid}, p.~\!202).

\section{Geodesics}\label{sec:geodesic}
In Sect.~\!\ref{sec:gauss}, we stressed the role of shortest curves in the intrinsic\index{intrinsic} measurement in a surface. In order to outline another work of Gauss\index{Gauss, Johann Carl Friedrich (1777--1855)} related to non-Euclidean geometry\index{non-Euclidean geometry}, we shall take a look at such curves from a general point of view. 

The length of a curve $c(t)=(x(t),y(t),z(t))$, $a\le t\leq b$, on a surface $S$  is given by $\ell(c)=\int_a^b\sqrt{\dot{x}(t)^2+\dot{y}(t)^2+\dot{z}(t)^2}~\!dt$. Thus our issue boils down to finding a curve $c$ that minimizes $\ell(c)$ among all curves on $S$ satisfying the boundary condition $c(a)=p$, $c(b)=q$. Johann Bernoulli\index{Bernoulli, Johann (1667--1748)} raised the problem of the shortest distance between two points on a convex surface (1697). Following his senescent teacher Johann, Euler\index{Euler, Leonhard (1707--1783)} obtained a differential equation for shortest curves on a surface given by an equation of the form $F(x,y,z)=0$ (1732).

Finding shortest curves is a quintessential problem in the {\it calculus of variations}\index{calculus of variations}, a central plank of mathematical analysis originated by Euler\index{Euler, Leonhard (1707--1783)}, which deals with maximizing or minimizing a general {\it functional} (i.e., a function whose ``variable" is functions or maps). A bud of this newly-born field is glimpsed in Fermat's principle of least time\index{Fermat's principle of least time}, but the derivation of Snell's law\index{Snell's law} is a straightforward application of the usual extreme value problem. The starting point in effect was the {\it brachistochrone\index{brachistochrone} curve problem} raised by Johann Bernoulli\index{Bernoulli, Johann (1667--1748)} (1696);\footnote{``Brachistochrone" is from the Greek words \textgreek{br'aqistos} (shortest) and \textgreek{qr'onos} (time). The first to consider the brachistochrone problem is Galileo\index{Galileo Galilei (1564--1642)} ({\it Two New Sciences}).} the problem of finding the shape of the curve down which a sliding particle, starting at rest from $P_1=(x_1,y_1)$ and accelerated by gravity, will slip (without friction) to a given point $P_2=(x_2,y_2)$ in the least time. To be exact, the problem is to find a function $y=y(x)$ which minimizes the integral $\int_{x_1}^{x_2}\sqrt{\frac{1+y'{}^2}{2gy}}~\!dx$.

In the paper {\it Elementa Calculi Variationum} read to the Berlin Academy on September 16, 1756, Euler\index{Euler, Leonhard (1707--1783)} unified various variational problems\index{variational problem}. At about the same time, Joseph-Louis Lagrange\index{Lagrange, Joseph-Louis (1736--1813)} (1736--1813) simplified Euler's earlier analysis and derived what we now call the {\it Euler-Lagrange equation}\index{Euler-Lagrange equation} (E-L equation). Especially, his letter to Euler dated August 12, 1755  (when he was only 19 years old!) contains the technique used today. Here in general (however restricted to the case of one-variable), the E-L equation\index{Euler-Lagrange equation} associated with the functional of the form  
\begin{equation}\label{eq:E-L}
\mathcal{F}[\boldsymbol{q}]=\int_a^bL\big(t,\boldsymbol{q}(t),\dot{\boldsymbol{q}}(t)\big)~\!dt,\quad \boldsymbol{q}(t)=(q_1(t),\ldots,q_n(t))
\end{equation}
is given by 
\begin{equation}\label{eq:E-L1}
\frac{\partial L}{\partial q_i}-\frac{d}{dt}\frac{\partial L}{\partial \dot{q}_i}=0~~\quad~(i=1,\ldots,n).
\end{equation}
What Lagrange\index{Lagrange, Joseph-Louis (1736--1813)} proved is that $\boldsymbol{q}(t)$ minimizing $\mathcal{F}$ is a solution of (\ref{eq:E-L1}). 

Applying the E-L equation\index{Euler-Lagrange equation} to length-minimizing curves on a surface, we obtain an ordinary differential equation for $x(t), y(t), z(t)$. What should be aware of is that a solution of this equation is only a {\it locally} length-minimizing curve; not necessarily shortest. Such a curve is called a {\it geodesic}\index{geodesic}, the term (coming from geodesy) coined by Joseph Liouville\index{Liouville, Joseph (1809--1882)} (1809--1882).\footnote{{\it De la ligne g\'{e}od\'{e}sique sur un ellipsoide quelconque}, J.~\!Math.~\!Pures Appl., {\bf 9} (1844), 401--408. Strictly speaking, as a parameter of a geodesic\index{geodesic}, we adopt the one proportional to arc-length. We should note that $\ell(c)$ is left invariant under a parameter change.} 

A geodesic\index{geodesic} may be defined as a trajectory of {\it force-free motion}. To explain this, we regard a curve $c(t)$ in a surface $S$ as the trajectory of a particle in motion. The velocity vector $\dot{c}(t)$ is a vector tangent to $S$ at $c(t)$, but the acceleration vector $\ddot{c}(t)$ is not necessarily tangent to $S$. If the {\it constraint force} confining the particle to $S$ is the only force, then $\ddot{c}(t)$ must be perpendicular to the tangent plane of $S$ at $c(t)$ because the constraint force always acts vertically on $S$. Such a curve $c$ turns out to be a geodesic\index{geodesic}. Hence, if we denote by $\frac{D}{dt}\frac{dc}{dt}$ the tangential part of $\ddot{c}$, then $c$ is a geodesic\index{geodesic} if and only if $\frac{D}{dt}\frac{dc}{dt}=0$. Since this is a second-order ordinary differential equation (see (\ref{eq:geodesic}) and Sect.~\!\ref{sect:reimenniangeometry}), we have a unique geodesic\index{geodesic} $c$ defined on an interval containing $0$ such that $c(0)=p$, $\dot{c}(0)=\boldsymbol{\xi}$ for given a point $p$ and a vector $\boldsymbol{\xi}$ tangent to $S$ at $p$.

There is an alternative derivation of the equation $\frac{D}{dt}\frac{dc}{dt}=0$, based on the {\it principle of least action}\index{principle of least action} that was enunciated by Pierre Louis Moreau de Maupertuis\index{Maupertuis, Pierre Louis Moreau de (1698--1759)} (1698--1759) and occupies its position in the physicist's pantheon,\footnote{\label{fn:cassini}{\it Accord de diff\'{e}rentes loix de la nature qui avoient jusqu'ici paru incompatibles}, read to the French Academy on April 15, 1744. His paper 
{\it Loi du repos des corps}, Histoire Acad.~\!R.~\!Sci.~\!Paris (1740), 170--176, is its forerunner.} according to which a particle travels by the {\it action}-minimizing path between two points. The action he defined is the  ``sum of $mvs$" along the entire trajectory of a particle with mass $m$, where $s$ is the length of a tiny segment on which the velocity of the particle is $v$.\footnote{On reflection, the principle of least action\index{principle of least action} seems to tell us that present events are dependent on later events in a certain manner; so, though a bit exaggerated, a natural occurrence seems to be founded on an {\it intention} directed to a certain end. For such a teleological and purposive flavor (dating back to Aristotle\index{Aristotle (384--322 BC)}), Maupertuis\index{Maupertuis, Pierre Louis Moreau de (1698--1759)} thought that his cardinal principle proves a rational order in nature, and buttresses the existence of a rational God.}

The obscurity of Maupertuis' formulation was removed by Euler\index{Euler, Leonhard (1707--1783)} (1744)and Lagrange\index{Lagrange, Joseph-Louis (1736--1813)} (1760).
%\footnote{Euler, {\it Additamentum} 2 to his {\it Methodus inveniendi lineas curvas maximi minimive proprietate gaudentes} (1744). Lagrange, {\it Application de la M\'{e}thode pr\'{e}c\'{e}dente \`{a} la solution de diff\'{e}rens Probl\`{e}mes de Dynamique}, M\'{e}langes de philosophie et de math\'{e}matique de la soci\'{e}t\'{e} royale de Turin, (1760-1761), 196--298.} 
Specifically the {\it action integral}\index{action integral} $\mathcal{S}$ applied to a particle in motion under a given {\it potential energy} $V$ (the concept conceived by Lagrange\index{Lagrange, Joseph-Louis (1736--1813)} in 1773) is\!: $\mathcal{S}\!=\!\int_a^b \Big(\frac{1}{2}m\|\dot{c}(t)\|^2-u(c(t))\Big)dt$. The E-L equation\index{Euler-Lagrange equation} for $\mathcal{S}$ reduces to the Newtonian equation\!: $\ddot{c}\!=\!-{\rm grad}~\!u$.
%\footnote{Lagrange, {\it M\'{e}canique Analytique}, 1788.} 
(For example, the gravitational potential\index{gravitational potential} energy generated by the mass distribution $\mu$ is $u(\boldsymbol{x})\!=\!-G\int\|\boldsymbol{x}-\boldsymbol{y}\|^{-1}$ $d\mu(\boldsymbol{y})$.) This suggests us to define the action integral\index{action integral} for a curve $c$ in $S$ by $E(c)\!=\!\int_a^b \|\dot{c}(t)\|^2~\!dt$. Indeed, the E-L equation\index{Euler-Lagrange equation} for this functional coincide with $\frac{D}{dt}\frac{dc}{dt}=0$.

\medskip

The minimum (maximum) problems come up in various scenes in geometry. For instance, the (classical) {\it isoperimetric theorem} whose origin goes back to ancient Greece states, ``Among all planar regions with a given perimeter, the circle encloses the greatest area." Needless to say, minimal surfaces\index{minimal surface} are solutions of a minimum problem.
%\footnote{The study of minimal surfaces\index{minimal surface} dates back to Euler\index{Euler} (1744).}  
In particular, finding a surface of the smallest area with a fixed bounday curve is known as {\it Plateau's problem}, the problem raised by Lagrange\index{Lagrange, Joseph-Louis (1736--1813)} in 1760 and investigated in depth in the 20th century onward.\footnote{
%Lagrange, {\it Essai d'une nouvelle m\'{e}thode pour d\'{e}terminer les maxima et les minima des formules int\'{e}grales ind\'{e}finies}, Misc. Philos.-Math. Soc. Priv. Taurinensis, {\bf 2} (1760-1762), 173--195. 
Plateau conducted experiments with soap films to study their configurations (1849).} Minimum/maximum principles hold emphatic importance even in different fields; say, in {\it potential theory} (Remark \ref{rem:infinixx} (1)), in {\it statistical mechanics}, in {\it information science}, in the design of crystal structures (\cite{sunaxx}), and much else besides. Even more surprisingly, so too does in the characterization of the ``shape" of the universe as will be explained later.

\begin{rem}\label{rem:momentum} {\rm {\small
(1) Hamilton\index{Hamilton, William Rowan (1805-1865)} reformulated the E-L equation\index{Euler-Lagrange equation} in the context of mechanics as follows (1835).
%\footnote{{\it On a General Method in Dynamics}, Phil.~\!Trans.~\!R.~\!Soc., (1834), 247--308; (1835), 95--144.} 
For simplicity, suppose that the integrand $L$ in (\ref{eq:E-L}) (called a {\it Lagrangian}) does not involve the time variable $t$. Making the change of variables $(q_1,\ldots,q_n,$ $\dot{q}_1,$ $\ldots,\dot{q}_n)\mapsto (q_1,\ldots,q_n,$ $p_1,\ldots,$ $p_n)$ where $p_i=\displaystyle\partial L/\partial \dot{q}_i$, and defining the function $H$ (the {\it Hamiltonian}) with the new variables $(\boldsymbol{p},\boldsymbol{q})$ by $H(\boldsymbol{p},\boldsymbol{q})=\sum_ip_i\dot{q}_i-L(\boldsymbol{q},\dot{\boldsymbol{q}})$, we may transform (\ref{eq:E-L1}) into what we call the {\it Hamiltonian equation}\index{Hamiltonian equation}:
\begin{equation*}%\label{eq:hamiltoneq}
\frac{dp_i}{dt}  =  -\frac{\partial H}{\partial q_i},\quad
\frac{dq_i}{dt}  =  \frac{\partial H}{\partial p_i},\quad (i=1,\ldots,n).
\end{equation*}
A general form of the {\it law of the conservation of energy} is then embodied by the fact that $H(\boldsymbol{p}(t),\boldsymbol{q}(t))$ is constant for any solution $(\boldsymbol{q}(t),\boldsymbol{p}(t))$. In the case where $L\!=\!\frac{m}{2}\sum_{i=1}^n\dot{q}_i{}^2-V(\boldsymbol{q})$, we have $\boldsymbol{p}\!=\!m\dot{\boldsymbol{q}}$ and $H\!=\!\frac{1}{2m}\sum_{i=1}^np_i{}^2+V(\boldsymbol{q})$ $(=\!\frac{m}{2}\sum_{i=1}^n\dot{q}_i{}^2\!+ V(\boldsymbol{q}))$. Hence $\boldsymbol{p}$ in this case is the {\it momentum}, and $H$ is the {\it total mechanical energy}; that is, the sum of the {\it kinetic energy} and the potential energy. Hamilton's\index{Hamilton, William Rowan (1805-1865)} formulation gained further prominence with the advent of statistical mechanics and quantum mechanics.

\smallskip

(2) To outline a maximum principle in statistical mechanics,  consider a simple quantum mechanical system with possible microscopic states $1,\ldots,n$ whose energies are $E_1,\ldots,E_n$.\footnote{A quantum mechanical system is described by a self-adjoint operator $\widehat{H}$ acting on a Hilbert\index{Hilbert, David (1862--1943)} space. A unit eigenvector of $\widehat{H}$ with eigenvalue $E$ represents a state with energy $E$.}
A (macroscopic) {\it state} of the system is described by a probability distribution $P=\{p_i\}_{i=1}^n$, where $p_i$ is the probability that the microscopic state $i$ occurs during the system's fluctuations. In this setting, the {\it internal energy} $U(P)$ and the {\it entropy} $S(P)$ are defined by $U(P)=\sum_{i=1}^nE_ip_i,~S(P)=-k_B \sum_{i=1}p_i\log p_i$, respectively, where $k_B$ is {\it Boltzmann's constant}. Among all states $P$ with $U=U(P)$, there is a unique $P^0=\{p_i^0\}_{i=1}^n$ with the maximum entropy; i.e. $S(P)\leq S(P^0)$. It is given by the {\it Gibbs distribution}\footnote{The term ``entropy," from Greek \textgreek{en} (``in") + \textgreek{trop'h} (``a turning"), was coined in the early 1850s by R. J. E. Clausius. In the 1870s, L. Boltzmann gave its statistical definition.} representing the equilibrium state
$
p_i^0=\exp(-E_i/k_BT)/Z(1/k_BT);$ $Z(\beta):=\sum_{i=1}^n\exp(-\beta E_i)
$,
where  $T$ is the {\it temperature}, determined by the equation
$
U=-d/d\beta\big|_{\beta=1/k_BT}\log Z(\beta)
$.

The {\it second law of thermodynamics} (formulated by Clausius) says that entropy of an isolated system invariably increases because the system evolves towards equilibrium. Therefore the total entropy of the universe is continually increasing.  
\hfill$\Box$

}
}
\end{rem}

\section{Curvature and non-Euclidean geometry}\label{sect:curvnon}
In a nutshell, what Gauss\index{Gauss, Johann Carl Friedrich (1777--1855)} discovered during his study of the new 
geometry is that homogeneity and isotropy are not enough to characterize Euclidean space; namely, his discovery suggests that there is another ``space" with these properties, which is to be called ``non-Euclidean space"\index{non-Euclidean space}. However, we have to stress that the way Gauss\index{Gauss, Johann Carl Friedrich (1777--1855)} and others searched for their new geometry was synthetic. At this stage, therefore, the ``new space" has not yet been explicitly described. 

\begin{figure}[htbp]
%\vspace{-5.8cm}
\vspace{-0.3cm}
\begin{center}
\includegraphics[width=.41\linewidth]
{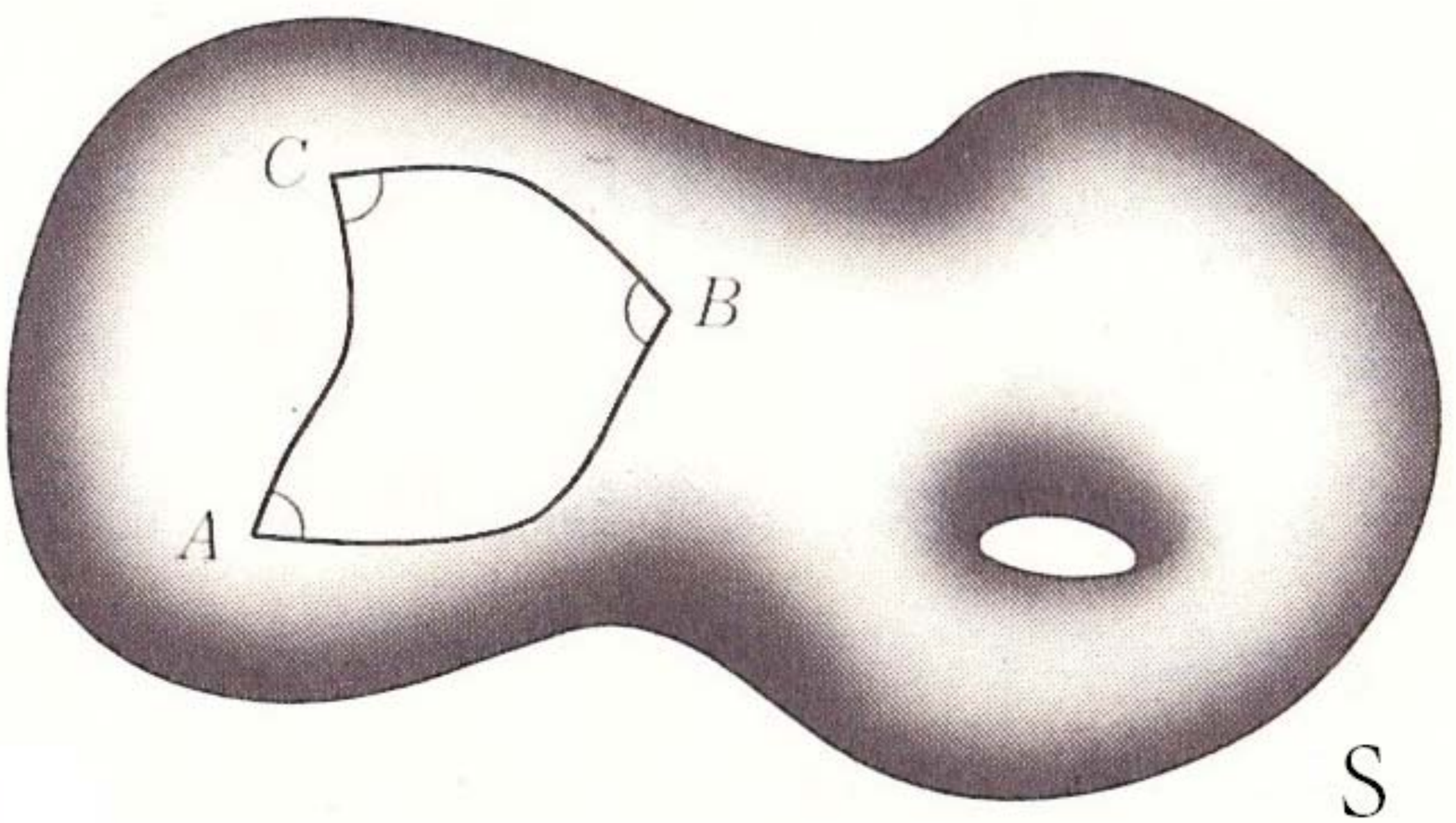}
\end{center}
\vspace{-1.2cm}
\caption{A triangle on a surface}\label{fig:surfacetri}
\vspace{-0.2cm}
\end{figure}

Gauss\index{Gauss, Johann Carl Friedrich (1777--1855)} did not refer to the non-Euclidean space\index{non-Euclidean space}, but perhaps he was aware that the ``new space" may yield an alternative model of the universe. In truth, the moment when he perceived this is by far the most significant turning point in cosmology.\footnote{In the {\it Disquisitiones generales} (Art.~\!28), Gauss\index{Gauss, Johann Carl Friedrich (1777--1855)} refers to the numerical data for the triangle formed by the three mountain peaks Brocken, Hohenhagen, and Inselsberg. Some say that with these data Gauss\index{Gauss, Johann Carl Friedrich (1777--1855)} wished to determine whether the universe is Euclidean. What we can say at the very least is that he compared a geodesic triangle on the earth surface with a plane rectilinear triangle.  Meanwhile,  Lobachevsky\index{Lobachevsky, Nikolai Ivanovich (1792--1856)} attempted to measure the angle sum of a triangle by analyzing data on the parallax of stars  ({\it  Exposition succincte}; see \cite{mona}).} What is more, Gauss\index{Gauss, Johann Carl Friedrich (1777--1855)} seemed to know an intimate relation between his intrinsic\index{intrinsic} theory of surfaces and non-Euclidean (plane) geometry because he obtaind the arresting formula ({\it Disquisitiones generales}, Art.~\!20)\!: 
\begin{equation}\label{eq:gaussformula}
\iint_{\triangle ABC}K_S~\!d\sigma=(\angle A+\angle B+\angle C)-\pi,
\end{equation}
where $\triangle ABC$ is a {\it geodesic triangle}\index{geodesic triangle} (a triangle on $S$ formed by the arcs of three geodesics\index{geodesic} segments; Fig.~\!\ref{fig:surfacetri}), and $d\sigma$ is the {\it surface element}\index{surface element} of $S$ (Sect.~\!\ref{sec:riemann}). Note that this is a generalization of Harriot-Girard's\index{Harriot, Thomas (c. 1560 --1621)}\index{Girard, Albert (1595--1632)} formula.

\begin{figure}[htbp]
%\vspace{-6cm}
\vspace{-0.3cm}
\begin{center}
\includegraphics[width=.42\linewidth]
{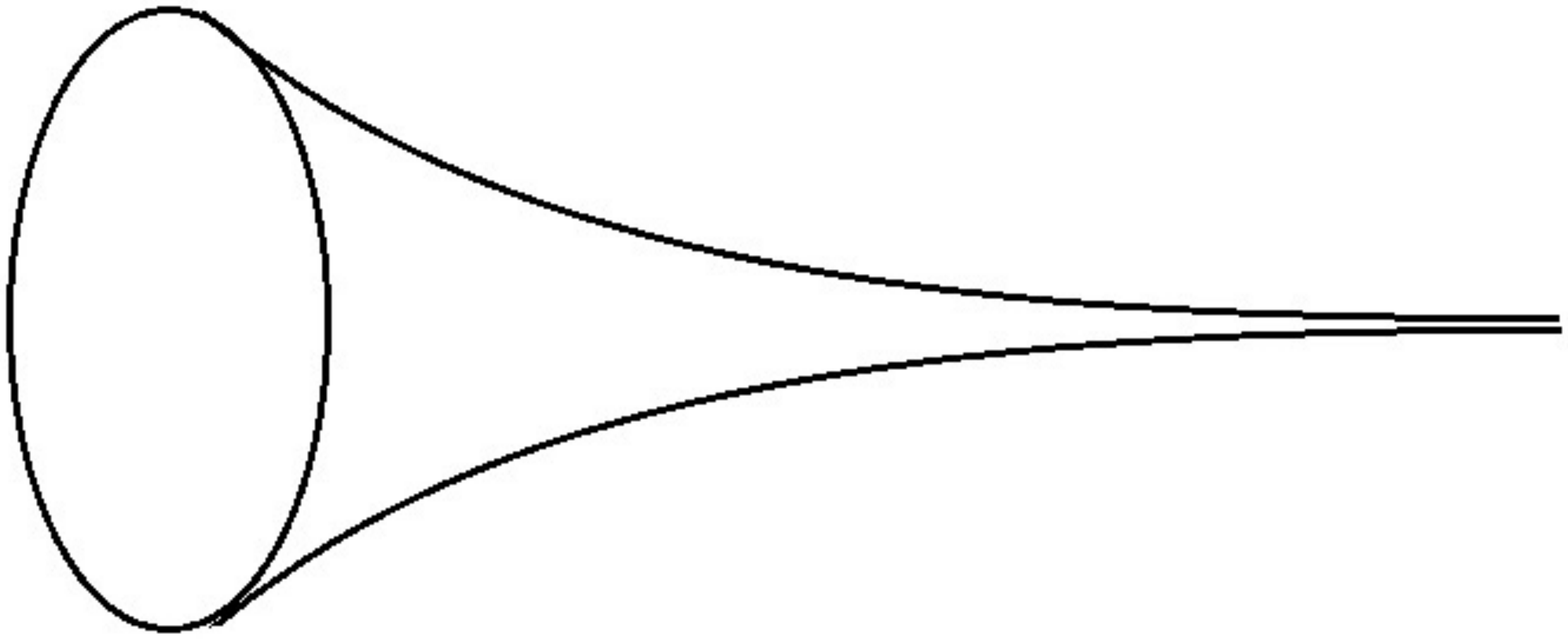}
\end{center}
\vspace{-1.7cm}
\caption{Das Gegenst\"{u}ck der Kugel}\label{fig:constantcurv}
\vspace{-0.2cm}
\end{figure}

From (\ref{eq:gaussformula}), a ``plane-like surface" of constant negative curvature (if exists) is expected to be a non-Euclidean plane-model, and hence finding such a surface became an urgent issue. In a note dated about 1827, Gauss\index{Gauss, Johann Carl Friedrich (1777--1855)} jotted down his study about the surface of revolution generated by the {\it tractrix}\index{tractrix} (Fig.~\!\ref{fig:constantcurv}),\footnote{The tractrix \index{tractrix} is the graph of the function $y=\displaystyle R\log ((R+\sqrt{R^2-x^2})/x)$ $+\sqrt{R^2-x^2}$, which shows up as a locus of an object obtained from dragging it along a line with an inextensible string. The architect C. Perrault proposed a problem related to the tractrix (1670) for the first time. The appellation ``tractrix" was coinded by Huygens\index{Huygens, Christiaan (1629--1695)} in 1692.} which he called ``das Gegenst\"{u}ck der Kugel" (the opposite of the sphere).\footnote{\label{fn:pseudo}This was to be called the {\it pseudosphere}\index{pseudosphere} by Eugenio Beltrami\index{Beltrami, Eugenio (1835--1900)} (1835--1900); see fn.~\!\ref{fn:Saggio}.} This is an example of a surface with constant negative curvature, but he did not mention clearly that it has to do with the new geometry ({\it Werke}, VIII, p.~\!265).\footnote{F. Minding found that the trigonometry for geodesic triangles\index{geodesic triangle} on this surface enjoys analogy to spherical trigonometry\index{spherical trigonometry} (1840). However he did not notice that his finding is equivalent to the non-Euclidean trigonometric relations.} The reason might be because he could not construct such a surface with a planar shape. The truth is that the non-Euclidean plane cannot be realized as an ordinary surface in space as Hilbert\index{Hilbert, David (1862--1943)} showed in 1901. In order to apprehend its true identity, we need the notion of manifold (Sect.~\!\ref{sec:after1}).

Gauss\index{Gauss, Johann Carl Friedrich (1777--1855)} may have known (\ref{eq:gaussformula}) as early as 1816, so it may have served as a forerunner to the Theorema Egregium\index{Theorema Egregium}. In his {\it Disquisitiones generales}, Art.~\!20, he says, ``This formula, if I am not wrong, should be counted among the most elegant ones of the theory of surfaces."

Behind (\ref{eq:gaussformula}), there seems to be his commitment to electromagnetism;\footnote{The rapid progress of electromagnetism\index{electromagnetism} originated with H. C. Oersted who verified by experiment that electricity and magnetism are not independent (1820).} in 1813, he proved a special case of what we call now {\it Gauss's\index{Gauss, Johann Carl Friedrich (1777--1855)} flux law} that correlates the distribution of electric charge with the resulting electric field. The law is expressed by a formula connecting an integral inside a closed surface with a surface integral over the boundary ({\it divergence theorem}\index{divergence theorem}; {\it Werke}, V, p.~\!5--7).

\begin{rem}{\rm {\small 
The divergence theorem\index{divergence theorem} asserts that   
$$
\iiint_D {\rm div}~\!X~\!d\boldsymbol{x}=\iint_S\langle \boldsymbol{n},X\rangle~\!d\sigma
$$
for a vector field $X$ defined on a domain $D\subset \mathbb{R}^3$ with boundary $S$, where $\boldsymbol{n}$ is the outer unit normal of the surface $S$. A special case was treated by Lagrange\index{Lagrange, Joseph-Louis (1736--1813)} (1762), and by Gauss\index{Gauss, Johann Carl Friedrich (1777--1855)} (1813). It was M. V. Ostrogradsky who gave the first proof for the general case (1826). The 2D version is known as Green's theorem\index{Green's theorem} (1828): 
\begin{equation}\label{eq:lineintegral0}
\iint_D\Big(\frac{\partial g}{\partial x}-\frac{\partial f}{\partial y}
\Big)~\!dxdy=\int_C fdx+gdy,
\end{equation}
where $D$ is the domain surrounded by a piecewise smooth, simple closed counterclockwise oriented curve $C$.
%\footnote{George Green\index{Green, George (1793--1841)} (1793--1841), {\it The Application of Mathematical Analysis to the Theories of Electricity and Magnetism} (1828).} 
Another important formula is {\it Stokes' formula}\index{Stokes' formula} for a bordered surface $S$ with a normal unit vector field $\boldsymbol{n}$: 
$$
\iint_S \langle \boldsymbol{n},{\rm rot}~\!X\rangle~\!d\sigma=\int_C\langle \dot{c}(s),X(c(s))\rangle~\!ds,
$$  
where we represent the boundary curve $C$ of $S$ by an arc-length parameterization $c(s)$ in such a way that $\dot{c}(s)\times \boldsymbol{n}(c(s))$ points towards the interior of $S$.\footnote{The statement of Stokes' formula appeared as a postscript to a letter (July 2, 1850) from William Thomson (Lord Kelvin\index{Kelvin (William Thomson; 1824--1907)}, 1824--1907) to G. Stokes. Stokes set the theorem as a question on the 1854 Smith's Prize  Examination at  Cambridge.}
\hfill$\Box$

}
}
\end{rem}

A formula that deserves to be mentioned is derived from (\ref{eq:gaussformula}).% with a modicum of effort. 
For a closed surface $S$ with a division by a collection of geodesic triangles\index{geodesic triangle} such that each triangle side is entirely shared by two adjacent triangles, we have
\begin{equation}\label{eq:gaussbonnet}
\int_SK_S~\!d\sigma=2\pi (v-e+f),
\end{equation}
where $v, e, f$ are the number of vertices, edges, and faces, respectively. This formula, due to Walther Franz Anton von Dyck\index{von Dyck, Walther Franz Anton (1856--1934)} (1856--1934), is called the {\it Gauss-Bonnet formula}\index{Gauss-Bonnet formula}, though they never referred to (\ref{eq:gaussbonnet}) in their work.\footnote{von Dyck, {\it Beitr\"{a}ge zur Analysis situs I. Aufsatz. Ein- und zweidimensionale Mannigfaltigkeiten}, Math.~\!Ann., {\bf 32} (1888), 457--512. Pierre Ossian Bonnet\index{Bonnet, Pierre Ossian (1819--1892)} (1819--1892) generalized (\ref{eq:gaussformula}) to triangles with non-geodesic sides (1848).}
%({\it M\'{e}moire sur la th\'{e}orie g\'{e}n\'{e}rale des surfaces}, J. \'{E}cole Polytech,. {\bf 19} (1848), 1-–146).} 
It tells us that $\chi(S):=v-e+f$ (the {\it Euler characteristic}\index{Euler characteristic} of $S$) is independent of the choice of a triangulation. In particular, applying to the unit sphere $S^2$, we get the {\it polyhedron formula}\index{Euler's polyhedron formula} $\chi(S^2)=2$ discovered by Euler\index{Euler, Leonhard (1707--1783)} in 1750 (Sect.~\!\ref{subsec:topodes}).

The proof of (\ref{eq:gaussbonnet}) goes as follows. The sum of the interior angles gathering at each vertex is $2\pi$, so that 
{\small 
$$
\int_SK_Sd\sigma=\sum_{\triangle ABC}\iint_{\triangle ABC}K_S~\!d\sigma=\sum_{\triangle ABC}\big\{(\angle A+\angle B+\angle C)-\pi\big\}=2\pi v-f\pi,
$$
}
\hspace{-0.15cm}while $3f=2e$ because each edge is counted twice when we count edges of all triangles. Consequently, the left-hand side of (\ref{eq:gaussbonnet}) is equal to $\pi(2v-f)=$ $\pi$ $(2v$ $-2e+2f)=2\pi(v-e+f)$, as desired.

 The Gauss-Bonnet formula\index{Gauss-Bonnet formula} opened the doorway to {\it global analysis} that knits together three fields; differential geometry, topology\index{topology}, and analysis (Sect.~\!\ref{subsec:topodes}).

\medskip

After Gauss's\index{Gauss, Johann Carl Friedrich (1777--1855)} monumental feat, cosmology began to remove itself from religion and even philosophy, and to snuggle up to higher mathematics which is unconstrained by perceptual experience; thus God walked away from the center stage. However the outside of space still remains as a backdrop. It was Riemann\index{Riemann, Georg Friedrich Bernhard (1826--1866)}, the next protagonist in our story, who at last erased the outside.

\section{The dimension of space}\label{sec:dimension}
Deferring the account of Riemann's achievement to the next section, we shall pause for a moment to discuss {\it higher-dimensional spaces}.

{\it Dimension}\index{dimension} is, in a naive sense, a measurable extent of a particular kind, such as length, breadth and depth, or the maximum number of independent directions in a refined sense. Physically it is the number of degrees of freedom available for movement in a space.  Aristotle\index{Aristotle (384--322 BC)} was among the earliest who contemplated dimension; he says, ``magnitude which is continuous in one dimension is length; in two breadth, in three depth, and there is no magnitude besides these" ({\it Metaphysics}, V, 11-14). Ptolemy\index{Ptolemy (c. 100--c. 170)} likewise asserted in his book {\it On Distance} that one cannot consider ``space of dimension more than three," and even gave its proof as testified by Simplicius of Cilicia ({\it ca.}~\!490 CE--{\it ca.}~\!560 CE). 

In Sect.~\!\ref{sec:gauss}, we speculated what the inhabitants in a surface would do when they want to know the shape of their universe. Surfaces are 2D, whilst our universe is 3D. Why are we in the 3D space, not in another?\footnote{Closely related to this question is the {\it anthropic principle} (Brandon Carter, 1973), the philosophical deliberation that the universe is inseparable from mankind's very existence.} However difficult it may have been to furnish a definite answer to this question, there is no barrier to thinking of higher-dimensional models of space. Actually we have a clue in analytic geometry\index{analytic geometry} to formulate coordinate spaces of general dimension.\footnote{Analytic geometry\index{analytic geometry} acquired its whole range in the 18th century with A. Parent, Clairaut, Euler\index{Euler, Leonhard (1707--1783)}, G. Cramer, Lagrange\index{Lagrange, Joseph-Louis (1736--1813)} and many others. La Hire's\index{La Hire, Philippe de (1640--1718)} work {\it Nouveaux \'{e}l\'{e}mens des sections coniques, les lieu g\'{e}om\'{e}triques} (1679) summarized the progress of analytic geometry\index{analytic geometry} during a few decades since Descartes\index{Descartes, Ren\'{e} (1596--1650)}. His work contains some ideas leading to the extension of space to more than three dimensions.}

It was in the 18th century that mathematicians began to contemplate the possibility of the use of $\mathbb{R}^d$ as a gadget to develop mechanics. As Lagrange\index{Lagrange, Joseph-Louis (1736--1813)} and D'Alembert\index{d'Alembert, Jean-Baptiste le Rond (1717--1783)} suggested, it is natural to add the time parameter as the fourth dimension. Moreover, it is technically convenient to describe a mass system of $N$ particles located at $(x_{1}, y_1,z_1), (x_2,y_2,z_2),$ $\ldots,$ $(x_N,y_N,z_N)$ as the single point $(x_1,y_1,z_1,x_2,y_2,z_2,\ldots,x_N,y_N,z_N)$ in $\mathbb{R}^{3N}$. Then, around the 1830s, the theory of $n$-ple integrals was evolved. For instance, Carl Gustav Jacob Jacobi (1804--1851) computed the volume of a higher-dimensional sphere in connection with quadratic forms though he deliberately avoided geometrical language (1834).

Certain it is that, by the early 19th century, the time was ripe for developing higher-dimensional geometry. In 1843, Arthur Cayley\index{Cayley, Arthur (1821--1895)} published a small memoir on $\mathbb{R}^n$, however the paper is characteristically algebraic. In the same year, Hamilton\index{Hamilton, William Rowan (1805-1865)} conceived the 4D number system (Sect.~\!\ref{sec:gauss}). A year later, {\it Die Lineale Ausdehnungslehre, ein neuer Zweig der Mathematik}, a pioneer work on vector spaces of general dimension, was published by Hermann Grassmann\index{Grassmann, Hermann (1809--1877)} (1809--1877).\footnote{In this book, he created what we now call the {\it exterior algebra}\index{exterior algebra} (or {\it Grassmann algebra}), which, denigrated at the time, was to exert a profound influence on the formation of algebra and geometry; see Sect.~\!\ref{sect:reimenniangeometry}. In this sense, he was quite ahead of his time.}   

A point at issue is how to imagine $\mathbb{R}^d$ as a higher-dimensional analogue of plane and space. A lighthearted way is to execute a geometry in $\mathbb{R}^d$ by mimicking geometry in $\mathbb{R}^2$ and $\mathbb{R}^3$ (thus geometric terminology may be of psychological assistance when we imagine $\mathbb{R}^d$).
It is to the credit of Ludwig Schl\"{a}fli\index{Schl\"{a}fli, Ludwig (1814--1895)} (1814--1895) to have inaugurated the study of the geometry of $\mathbb{R}^d$ along this line. In his magnum opus, {\it Theorie der vielfachen Kontinuit\"{a}t}, worked out from 1850 to 1852,\footnote{It was in 1901, after his death, that the entire manuscript was published in Denkschriften Der Schweizerischen Naturforschenden Gesellschaft, {\bf 38} (1901), 1--237.} Schl\"{a}fli defined the distance between $A=(x_1,\ldots,x_d)$ and $B=(y_1,\ldots,y_d)$ in $\mathbb{R}^d$ by
%\begin{equation}\label{eq:eucdist}
$d(A,B):=\sqrt{(x_1-y_1)^2+\cdots+(x_d-y_d)^2}$.
%\end{equation}
%\begin{equation}\label{eq:angle1}
%\cos \angle AOB:=\frac{x_1y_1+\cdots+x_dy_d}{\sqrt{x_1{}^2+\cdots+x_d{}^2}\sqrt{y_1{}^2+\cdots+y_d{}^2}}.
%\end{equation}

As expected, many results in $\mathbb{R}^d$ are just direct generalizations of the facts that hold in plane and space, but there are several results depending heavily on dimension. Schl\"{a}fli's\index{Schl\"{a}fli, Ludwig (1814--1895)} classification of higher-dimensional convex regular polytopes is such an example. Specifically, there are six regular polytopes in $\mathbb{R}^4$, while there are only three regular polytopes in $\mathbb{R}^d$, $d> 4$ (the analogues of the tetrahedron, the cube and the octahedron). Another dimension-dependent example---rather recent one originating with Newton's\index{Newton, Issac (1642--1726)} curiosity regarding celestial bodies---is the maximum possible kissing number $k(d)$ for $\mathbb{R}^d$, where the {\it kissing number}\index{kissing number} is the number of non-overlapping unit spheres $S^d$ $=\{(x_1,x_2,\ldots,x_{d+1})|~(x_1-a_1)^2+\cdots+(x_{d+1}-a_{d+1})^2=1\}$ that can be arranged in such a way that each of them touches another given unit sphere. For example, $k(2)$ $=6$, $k(3)=12$ and $k(4)=24$, but curiously $k(d)$ is unknown for $d>4$ except for $k(8)=240$ and $k(24)=196,560$.\footnote{The fact $k(3)=12$ was conjectured by Newton\index{Newton, Issac (1642--1726)}, while D. Gregory, a nephew of James Gregory, insisted that $k(3)=13$ in a discussion taking place at Cambridge in 1694. A proof in favor of Newton\index{Newton, Issac (1642--1726)} was given by K. Sch\"{u}tte and van der Waerden (1953).} Such phenomena enrich geometry; thereby, as opposed to Kant's\index{Kant, Immanuel (1724--1804)} dictum, higher-dimensional geometry turns out to be non-empty and fruitful.

\begin{rem}\label{rem:Grave}{\rm {\small 
As mentioned in Sect.~\!\ref{sec:gauss}, there exists no 3D number system that inherits, in part, the arithmetic properties of real numbers, while a 4D system exists. Hence a dimension-dependence phenomenon occurs in number systems. In 1843, J. T. Graves discovered the {\it octonions}, a non-associative 8D number system, inspired by Hamilton's\index{Hamilton, William Rowan (1805-1865)} discovery of quaternions\index{quaternion}.\footnote{Slightly before the publication of his paper, Cayley\index{Cayley, Arthur (1821--1895)} discovered the octonions\index{octonion}. Octonions are often referred to as {\it Cayley numbers}.} We thus have four number systems. These systems share the following properties; writing $\boldsymbol{a}\cdot\boldsymbol{b}$ for multiplication of $\boldsymbol{a}=(a_1,$ $\ldots,$ $a_d), \boldsymbol{b}=(b_1,\ldots,b_d)$ $~(d=1,2,4,8)$, we have  $(\boldsymbol{a}+\boldsymbol{b})\cdot \boldsymbol{c}\!=\!\boldsymbol{a}\cdot\boldsymbol{c}+\boldsymbol{b}\cdot \boldsymbol{c},~\boldsymbol{a}\cdot (\boldsymbol{b}+\boldsymbol{c})\!=\boldsymbol{a}\cdot\boldsymbol{b}+\boldsymbol{a}\cdot\boldsymbol{c}$,~$k(\boldsymbol{a}\cdot\boldsymbol{b})\!=\!(k\boldsymbol{a})\cdot\boldsymbol{b}\!=\!\boldsymbol{a}\cdot(k\boldsymbol{b})$ $(k\!\in\!\mathbb{R})$, and $\|\boldsymbol{a}\cdot\boldsymbol{b}\|=\|\boldsymbol{a}\|\|\boldsymbol{b}\|$. Here $\|\boldsymbol{a}\|^2\!=\!a_1{}^2\!+\!\cdots\!+\!a_d{}^2$. Interestingly, if $\mathbb{R}^d$ has multiplication having these properties, then $d$ must be 1,2,4, or 8 as A. Hurwitz showed in 1898.

Quaternions and octonions---the products of pure thought--- have applications in modern physics; thus being more than mathematical fabrications.
\hfill$\Box$

}
}

\end{rem}

\section{Riemann --The universe as a manifold}\label{sec:riemann}
Gauss\index{Gauss, Johann Carl Friedrich (1777--1855)} brought about a revolution in geometry as he envisaged. Indeed, his study of surfaces turned out to be the germ of a ``science of generalized spaces" that was essentially embodied by Georg Friedrich Bernhard Riemann\index{Riemann, Georg Friedrich Bernhard (1826--1866)} (1826--1866) in his Habilitationsschrift {\it \"{U}ber die Hypothesen welche der Geometrie zu Grunde liegen} presented to the Council of G\"{o}ttingen University.\footnote{\label{fn footnotexx}{\it Habilitationsschrift} is a postdoctoral thesis required for qualification as a lecturer. This is his second stage of the procedure (see fn.~\!\ref{fn:firststage} for his thesis at the first stage). The title of Riemann's thesis is an implicit criticism of Kant\index{Kant, Immanuel (1724--1804)}.}

It was Gauss\index{Gauss, Johann Carl Friedrich (1777--1855)}---only one year before the death at the age of 78---who asked his brilliant disciple to prepare a thesis on the foundations of geometry, out of the three possible topics that Riemann proposed. The two other topics were on electrical theory, and unlike the one chosen by his mentor, he was well-prepared for them. So Riemann had to work out the topic Gauss\index{Gauss, Johann Carl Friedrich (1777--1855)} recommended, and in a matter of months he could bring it to completion without qualm. 

The aged Gauss\index{Gauss, Johann Carl Friedrich (1777--1855)}, albeit being sick at that time, attended Riemann's probationary lecture delivered on 10 June 1854 in anticipation of hearing something new related to his old work. ``Among Riemann's audience, only Gauss\index{Gauss, Johann Carl Friedrich (1777--1855)} was able to appreciate the depth of Riemann's thoughts. ... The lecture exceeded all his expectations and greatly surprised him. Returning to the faculty meeting, he spoke with the greatest praise and rare enthusiasm to W. Weber about the depth of the thoughts that Riemann had presented" (\cite{mona}).

In his lecture aimed to a largely non-mathematical audience, Riemann\index{Riemann, Georg Friedrich Bernhard (1826--1866)} posed profound questions about how geometry is connected to the world we live in. The thrust being as such, the content was much more philosophical than mathematical though he humbly denied this (Part I).

At the outset, he says that the traditional geometry assumes the notion of space and some premises  which are merely nominal, while the relation of these assumptions remains in darkness because the general notion of {\it multiply extended magnitude} remained entirely unworked. He thus sets himself the task of constructing the notion of {\it manifold}\index{manifold} of general dimension.

There is something of note here. In his setup, Riemann\index{Riemann, Georg Friedrich Bernhard (1826--1866)} does not set a limit to the dimension of space. Remember that, when he inaugurated his study, higher-dimensional geometry as an abstract theory was just formulated. This might be a precipitating factor for his commitment to higher-dimensional curved spaces.  His general treatment of space was, however, not a mere idea or pipe dream, and had truly a decisive influence on modern cosmology that employs manifolds of dimension more than four (in 
superstring theory, space-time is 
10D).

Riemann\index{Riemann, Georg Friedrich Bernhard (1826--1866)} continues his deliberation on manifolds, which is bolstered by his immense insight into the principles behind the proper understanding of our space; but to go straight to the heart of his scheme, let us take a shortcut at the expense of his true motive. 

To single out a reasonable 
model of the universe, we must take into account the fact that our universe looks Euclidean in so far as we examine our vicinity. In the light of this observation, we define a $d$-dimensional manifold to be a ``generalized space" on which one can set up a curvilinear coordinate system $(x_1,$ $\ldots,x_d)$ around each point $p_0$, so that we may indicate each point $p$ near $p_0$ by its coordinates $(x_1,$ $\ldots, x_d)$. In addition, the higher-dimensional Pythagorean Theorem\index{Pythagorean Theorem} is supposed to hold in the infinitesimal sense.
To explore what this means, we recall that the distance $\Delta s$ between $(y_1,\ldots,y_d)$ and its displacement $(y_1+\Delta y_1,\ldots,y_d+\Delta y_d)$ in $\mathbb{R}^d$ is given by $(\Delta s)^2=(\Delta y_1)^2+\cdots+(\Delta y_d)^2$. When a skew coordinate system $(x_1,\ldots,x_d)$ is taken instead of the Cartesian one, this expression is modified as $(\Delta s)^2=\sum_{i,j=1}^dg_{ij}\Delta x_i\Delta x_j$. Here if the coordinate transformation between $(y_1,\ldots,y_d)$ and $(x_1,\ldots,x_d)$ is given by $y_i=\sum_{j=1}^dr_{ij}x_j$, then  $\Delta y_i=\sum_{j=1}^dr_{ij}\Delta x_j$, and hence $g_{ij}=\sum_{h=1}^d r_{hi}r_{hj}$. The square matrix $(g_{ij})$ is a positive symmetric matrix.  Since in the infinitesimally small scale, the curvilinear coordinate system is regarded as a skew coordinate system, it is natural to express the ``infinitesimal" Pythagorean Theorem as 
\begin{equation}\label{eq:infinitesimalform}
ds^2=\sum_{i,j=1}^dg_{ij}dx_i dx_j,
\end{equation}
where the ``line element" $ds$ stands for the ``distance between infinitesimally nearby points" $(x_1,\ldots,x_d)$ and $(x_1+dx_1,\ldots,x_d+dx_d)$, and  $(g_{ij})$ is, in turn, a function of the variables $x_1,\ldots,x_d$ with values in positive symmetric matrices. The right-hand side is called {\it the first fundamental form}\index{first fundamental form} (or {\it Riemannian metric}\index{Riemannian metric} in today's term; Sect.~\!\ref{sect:reimenniangeometry}), which is literally fundamental in Riemann's theory of manifolds. With it, one can calculate the length of a curve $c(t)$ $=(u_1(t),\ldots,u_d(t))$ $~(a\leq t\leq b)$ by the integral
\begin{equation}\label{eq:length}
\int_a^b\Big\{\sum_{i,j=1}^dg_{ij}\frac{du_i}{dt}\frac{du_j}{dt}\Big\}^{1/2}~\!dt.
\end{equation}
In imitation of the case of surfaces, we say that $c$ is a geodesic\index{geodesic} if it is locally length-minimizing. The distance of two points $p,q$ is defined to be the infimum of the length of curves joining $p,q$. 
%In addition, the angle $\theta$ between two directions towards $(x_1+du_1,\ldots,x_d+du_d)$ and $(x_1+dv_1,\ldots,x_d+dv_d)$ is defined by 
%$$
%\cos\theta =\Big(\displaystyle\sum_{i,j=1}^dg_{ij}du_idv_j\Big)\Big/\Big\{\displaystyle\sum_{i,j=1}^dg_{ij}du_idu_j\Big\}^{1/2}\Big\{\displaystyle\sum_{i,j=1}^dg_{ij}dv_idv_j\Big\}^{1/2}.
%$$

The next task---a highlight of his thesis---is to define ``curvature" in his setting. To this end, Riemann\index{Riemann, Georg Friedrich Bernhard (1826--1866)} takes a union of geodesics passing  through a given point $p$ to form a surface $S$, and then defines what we now call the ``sectional curvature"\index{sectional curvature} to be $K_S(p)$ (Part II, \S 3; see Sect.~\!\ref{sect:reimenniangeometry} for a modern definition). 

The manifold is said to have constant curvature if $K_S(p)$ does not depend on the choice of $p$ and $S$. He observes that a manifold with constant curvature $\alpha$ has a coordinate system $(x_1,\ldots,$ $x_d)$ around each point such that\begin{equation}\label{eq:constcurve}
ds^2=\displaystyle\sum_{i=1}^ddx_i{}^2\Big/\displaystyle\Big(1+\frac{\alpha}{4}\sum_{i=1}^dx_i{}^2\Big)^2.
\end{equation}
In particular, the curvature\index{curvature} of a manifold covered by a single coordinate system $(x_1,\ldots,x_d)$ with the fundamental form\index{first fundamental form} $\sum_{i=1}^ddx_i{}^2$ vanishes everywhere. Obviously this manifold, for which he used the term ``flat," is the $d$-dimensional coordinate space $\mathbb{R}^d$ with the distance 
%(\ref{eq:eucdist})
 introduced by Schl\"{a}fli\index{Schl\"{a}fli, Ludwig (1814--1895)}.

Related to the curvature is a question about {\it unboundedness} (Unbegrenzheit) and {\it infinite extent} (Unendlichkeit) of the universe (Part III, \S 2). He said that this kind of inquiry is possible even when our empirical determination is beyond the limit of observation. For example, {\it if the independence of bodies from positions could be assumed}, then the universe is of constant curvature, and hence must be {\it finite} provided that the curvature is positive (see Sect.~\!\ref{subsec:transfinite} for the exact meaning of finiteness). He justifies his argument by saying that a manifold with positive constant curvature has a sphere-like shape.\footnote{\label{fn:spehere}Riemann's claim should be phrased as ``a $d$-dimensional manifold with positive constant curvature is locally a portion of the sphere $S^d$." A typical example besides $S^d$ is the {\it projective space}\index{projective space} $P^d(\mathbb{R})$ obtained from $S^d$ by identifying $(x_1,\ldots,x_{d+1})$ with $(-x_1,\ldots, -x_{d+1})$, on which projective geometry\index{projective geometry} is developed (see Sect.~\!\ref{sec:after1}). In the 3D case, the classification of such manifolds is related to polyhedral groups. We will come across an example in Sect.~\!\ref{subsec:topodes}.}

Riemann's\index{Riemann, Georg Friedrich Bernhard (1826--1866)} consideration is not limited to finite dimensions; he suggests the possibility to contrive the theory of infinite-dimensional manifolds which may allow to handle spaces of functions or mappings {\it en bloc} from a view similar to the finite-dimensional case (Part I, \S 3). Thus, at this point of time, he might already had an inspiration that variational problems\index{variational problem} are regarded as extreme value problems on infinite-dimensional manifolds which a little later led to his use of what he called {\it Dirichlet's principle}\index{Dirichlet's principle} in complex analysis (1857).\footnote{\label{fn:Abelcschen}{\it Theorie der Abel'schen Funktionen}, J.~\!f\"{u}r die Reine und Angew.~\!Math., {\bf 54}, 115--155.} 

At the end of his lecture, Riemann ponders on the question of the validity of the hypotheses of geometry in the infinitely small realm. He says that the question of the validity of the hypotheses of geometry in the infinitely small is bound up with the question of the ground of the metric relations of space, and concluded the lecture with the following confident statement.

\medskip

{\it Researches starting from general notions, like the investigation we have just made, can only be useful in preventing this work from being hampered by too narrow views, and progress in knowledge of the interdependence of things from being checked by traditional prejudices. This leads us into the domain of another science, of physics, into which the object of this work does not allow us to go to-day} (translated by W. K. Clifford, {\it Nature}, Vol. VIII, 1873; 
see also \cite{sp}).

\begin{rem}\label{rem:infinixx}{\rm {\small 
(1) Dirichlet's principle\index{Dirichlet's principle} was already conceived by Gauss\index{Gauss, Johann Carl Friedrich (1777--1855)} in 1839 to solve a boundary value problem for harmonic functions by means of a minimum principle.\footnote{Kelvin\index{Kelvin (William Thomson; 1824--1907)} used the principle in 1847 in his study of electric fields. The most complete formulation was given by Johann Peter Gustav Lejeune Dirichlet\index{Dirichlet, Johann Peter Gustav Lejeune (1805--1859)} (1805--1859) around 1847 in his lectures on potential theory at Berlin.} Specifically, it says that a solution $u$ of the {\it Laplace equation}\index{Laplace equation} $\Delta u=0$ on a domain $D$ of $\mathbb{R}^d$ with boundary condition $u=g$ on $\partial D$ can be obtained as the minimizer of the functional $E(u)=\int_D \frac{1}{2}\|{\rm grad}~\!u\|^2~\!dx$. Yet, as pointed out by 
%Friedrich Emil Fritz Prym\index{Prym} (1841--1915) in 1871 and 
Karl Theodor Wilhelm Weierstrass\index{Weierstrass, Karl Theodor Wilhelm (1815--1897)} (1815--1897) in 1870, the existence of the minimizer is not obvious. It was in 1899 that Riemann's use of the principle was warranted (Hilbert\index{Hilbert, David (1862--1943)}).

\smallskip

(2) An ``infinitesimal" that shows up as the symbol $dx_i$ in (\ref{eq:infinitesimalform}) is paradoxically thought of as an object which is smaller than any feasible measurement. The symbol $dx$ (called the {\it differential} and exploited by Leibniz\index{Leibniz, Gottfried Wilhelm von (1646--1716)} for the first time) is convenient when expressing tidily formulas that involve differentiation and integration. For example, the expression $\int_cf~\!dx+g~\!dy$ (see (\ref{eq:lineintegral0})) for a line integral along a curve $c(t)=(x(t),y(t))$ $~(a\leq t\leq b)$ is originally $\int_a^b\Big(f(x(t),y(t))\frac{dx}{dt}+g(x(t),y(t))\frac{dy}{dt}\Big)dt$. Because of such a background, calculus had been referred to as the {\it infinitesimal analysis}. Typical are the {\it Analyse des infiniment petits} (1696) by Guillaume Fran\c{c}ois Antoine Marquis de L'H\^opital\index{L'H\^opital, Guillaume Fran\c{c}ois Antoine Marquis de (1661--1704)} (1661--1704), which contributed greatly to the dissemination of calculus in the Continent, and  Euler's\index{Euler, Leonhard (1707--1783)} {\it Introductio in Analysin Infinitorum}, the first systematic exposition of calculus of several variables.

The paradoxical nature of infinitesimal was criticized as incorrect by Berkeley  who described ``$0/0$" in differentiation as the ``ghost of departed quantities,"\footnote{{\it The analyst, or a discourse addressed to an infidel mathematicians} (1734). Berkeley denied the material existence of the world; he says, ``Esse est percipi (To be is to be perceived)."} and provoked controversy among the successors of Leibniz\index{Leibniz, Gottfried Wilhelm von (1646--1716)}. For instance, for Euler\index{Euler, Leonhard (1707--1783)}, infinitely small quantities are actually equal to zero, and differential calculus is simply a heuristic procedure for finding the value of the expression $0/0$" ({\it Institutiones Calculi Differentialis} (\cite{euler}). Yet, with the hindsight of modern geometry, infinitesimals are not altogether absurd and could survive today as ``duals" of certain ``operations" (Sect.~\!\ref{sect:reimenniangeometry}).

It is a significant progress in the history,  however inconspicuous, that differentiation was recognized as an operation that can be manipulated independently of functions to which they are applied. L. F. A. Arbogast was one of the first to make such a conceptual leap ({\it Calcul des d\'{e}rivations}, 1800). It goes without saying that the Leibnizian\index{Leibniz, Gottfried Wilhelm von (1646--1716)} notations played a key part in this comprehension. 
\hfill$\Box$

}
}
\end{rem}

To appreciate Riemann's\index{Riemann, Georg Friedrich Bernhard (1826--1866)} work, we shall step back into Gauss's\index{Gauss, Johann Carl Friedrich (1777--1855)} work, wherein we see the elements of continuity between them as forthrightly expressed in Riemann's \index{Riemann, Georg Friedrich Bernhard (1826--1866)} words, ``In the comprehension of the geometry of surfaces, the data on how a surface sits in space is incorporated with the intrinsic\index{intrinsic} measure-relation in which only the length of curves on the surface is considered" (Part II, \S 3).

We let $(u,v)\mapsto \boldsymbol{S}(u,v)$ be a {\it local parametric representation}\index{local parametric representation} of a surface $S$ in space; that is, $\boldsymbol{S}$ is a smooth one-to-one map from a domain $U$ of $\mathbb{R}^2$ into $S$ such that $\boldsymbol{S}_u(u,v)$ and $\boldsymbol{S}_v(u,v)$ are linearly independent for each $(u,v)\in U$ (thus span the tangent plane of $S$ at $\boldsymbol{S}(u,v)$). Note that the inverse map of $\boldsymbol{S}$ yields a curvilinear coordinate system $(u,v)$ of $S$.

Consider a curve $c(t)$ $(t\in [a,b])$ on $S$ whose coordinates are $\big(u(t), v(t)\big)$, and differentiate both sides of $c(t)=\boldsymbol{S}\big(u(t),v(t)\big)$ to obtain
{\small
$$  
\Big\|\frac{dc}{dt}\Big\|^2=\Big\|\frac{du}{dt}\boldsymbol{S}_u+\frac{dv}{dt}\boldsymbol{S}_v\Big\|^2=\langle \boldsymbol{S}_u,\boldsymbol{S}_u\rangle\left(\frac{du}{dt}\right)^2+2\langle \boldsymbol{S}_u,\boldsymbol{S}_v\rangle\frac{du}{dt}\frac{dv}{dt}+\langle \boldsymbol{S}_v,\boldsymbol{S}_v\rangle\left(\frac{dv}{dt}\right)^2.
$$
}
\noindent\hspace{-0.1cm}so the length of $c$ is given by
{\small 
$$
\int_a^b
\Big\|\frac{dc}{dt}\Big\|~\!dt=\int_a^b\Big\{\langle \boldsymbol{S}_u,\boldsymbol{S}_u\rangle\left(\frac{du}{dt}\right)^2+
2\langle \boldsymbol{S}_u,\boldsymbol{S}_v\rangle\frac{du}{dt}\frac{dv}{dt}+\langle \boldsymbol{S}_v,\boldsymbol{S}_v\rangle\left(\frac{dv}{dt}\right)^2\Big\}^{1/2}~\!dt.
$$}
\noindent \hspace{-0.1cm}Putting 
$E=\langle \boldsymbol{S}_u,\boldsymbol{S}_u\rangle,~\!F=\langle \boldsymbol{S}_u,\boldsymbol{S}_v\rangle,~\!G=\langle \boldsymbol{S}_v,\boldsymbol{S}_v\rangle$ and comparing this formula with (\ref{eq:length}), we see that the first fundamental form\index{first fundamental form} is ${\rm I}:=Edu^2+2Fdudv+Gdv^2$.

Besides, Gauss\index{Gauss, Johann Carl Friedrich (1777--1855)} introduced the coefficients of what we now call {\it the second fundamental form}\index{second fundamental form} ${\rm II}:=Ldu^2+2Mdudv+Ndv^2$ by setting $L=\langle \boldsymbol{S}_{uu},\boldsymbol{n}\rangle, ~M=\langle \boldsymbol{S}_{uv},\boldsymbol{n}\rangle,~N=\langle \boldsymbol{S}_{vv},\boldsymbol{n}\rangle$, where $\boldsymbol{n}$ is the unit normal defined by $\boldsymbol{n}=(\boldsymbol{S}_u\times \boldsymbol{S}_v)/\|\boldsymbol{S}_u\times \boldsymbol{S}_v\|$.  Then $K_S(p)=(LN-M^2)/(EG-F^2)$ ({\it Disquisitiones generales}, Art.~\!10). To show this, we employ the coordinate system $(x,y,z)$ and the function $f(x,y)$ introduced in Sect.~\!\ref{sec:gauss}. We may assume that $\boldsymbol{S}(0,0)=(0,0,0)$ and the direction of the $z$-axis coincides with that of $\boldsymbol{n}$. Let $(x(u,v),$ $y(u,v))$ be the $xy$-coordinate of the orthogonal projection of $\boldsymbol{S}(u,v)$ onto the tangent plane\index{tangent plane}, so that $\boldsymbol{S}(u,v)=\big(x(u,v),$ $y(u,v),$ $f(x(u,v),y(u,v))\big)$.\footnote{For the local parametric representation of $S$ given by the map $(x,y)\mapsto (x,y,f(x,y))$, the second fundamental form\index{second fundamental form} at $(0,0)$ is given by $adx^2+2bdxdy+cdy^2$.} We then have
{\small
\begin{eqnarray*}
\boldsymbol{S}_{uu}\hspace{-0.2cm}&=&\hspace{-0.2cm}\big(
x_{uu},y_{uu}, ax_u{}^2+2bx_uy_u+cy_{u}{}^2\big),~~\boldsymbol{S}_{vv}=\big(
x_{vv},y_{vv}, ax_v{}^2+2bx_vy_v+cy_{v}{}^2\big),\\
\boldsymbol{S}_{uv}\hspace{-0.2cm}&=&\hspace{-0.2cm}\big(
x_{uv},y_{uv}, ax_ux_v+b(x_uy_v+y_ux_v)+cy_{u}y_v\big)\quad ((u,v)=(0,0)).
\end{eqnarray*}
}
\hspace{-0.1cm}From the assumption on the vector $\boldsymbol{n}$, it follows that the coefficients $L,M,N$ are the $z$-coordinates of $\boldsymbol{S}_{uu}$, $\boldsymbol{S}_{uv}$, $\boldsymbol{S}_{vv}$, respectively, so that $L=$ $ax_u{}^2$ $+2bx_uy_u+cy_{u}{}^2,$~$N=ax_v{}^2+2bx_vy_v+cy_{v}{}^2,$~$M=ax_ux_v$ $+b(x_uy_v+y_ux_v)+cy_{u}y_v$. In matrix form, these formulas together may be written as
{\small 
$$
\begin{pmatrix}
L & M\\
M & N
\end{pmatrix}=
\begin{pmatrix}
x_u & y_u\\
x_v & y_v
\end{pmatrix}
\begin{pmatrix}
a & b\\
b & c
\end{pmatrix}
\begin{pmatrix}
x_u & x_v\\
y_u & y_v
\end{pmatrix}.
$$

}
\noindent Therefore, taking the determinant of both sides, we have $LN-M^2=(x_uy_v-y_ux_v)^2(ac-b^2)=(x_uy_v-y_ux_v)^2K_S$, while, noting that $\boldsymbol{S}_u=(x_u,y_u,0)$ and $\boldsymbol{S}_v=(x_v,y_v,0)$, we find that $EG-F^2=\langle \boldsymbol{S}_u,\boldsymbol{S}_u \rangle\langle\boldsymbol{S}_v,\boldsymbol{S}_v \rangle-\langle\boldsymbol{S}_u,\boldsymbol{S}_v \rangle^2=(x_uy_v-y_ux_v)^2
$, from which the claim follows.

We add a few words about the surface element\index{surface element} $d\sigma$ (Sect.~\!\ref{sect:curvnon}). For a point $p\!=\!\boldsymbol{S}(u,v)$, define $\boldsymbol{S}_0\!:\!\mathbb{R}^2\!\rightarrow\!E$ by setting $\boldsymbol{S}_0(x,y)\!=\!p+x\boldsymbol{S}_u\!+\!y\boldsymbol{S}_v$, which is a parametric representation of the tangent plane at $p$. Consider the small parallelogram $P\!=\!\boldsymbol{S}_0\big(\{(x,y)|~\!0\leq x\leq \Delta u,~0\leq y\leq \Delta v\}\big)$. This approximates $\boldsymbol{S}\big(\{(x,y)|~\!0\!\leq x\!\leq\!\Delta u,~\!0\leq y\leq$ $\Delta v\}\big)$. Hence, the surface area of $\boldsymbol{S}\big(\{(x,y)|~0\!\leq\!x\leq\!\Delta u,~\!0\!\leq\!y\!\leq\!\Delta v\}\big)$ is approximated by ${\rm Area}(P)\!=\!\big\{\langle \boldsymbol{S}_u,\boldsymbol{S}_u \rangle\langle\boldsymbol{S}_v,\boldsymbol{S}_v \rangle-\langle\boldsymbol{S}_u,\boldsymbol{S}_v \rangle^2\big\}^{1/2}\Delta u\Delta v\!=\!\sqrt{EG-F^2}\Delta u\Delta v$, so the surface area of $\boldsymbol{S}(U)$ is given by $\iint_U\sqrt{EG-F^2}dudv$, and $d\sigma\!=\!\sqrt{EG-F^2}dudv$ (\!{\it Disquisitiones generales}, Art.~\!17).

\medskip

The reader might suspect that Riemann's\index{Riemann, Georg Friedrich Bernhard (1826--1866)} approach is more or less a direct generalization of Gauss's\index{Gauss, Johann Carl Friedrich (1777--1855)} one. That is not quite correct since Riemann's formulation is entirely intrinsic\index{intrinsic}; therefore, the universe can exist {\it without the outside} (in particular, we may consider surfaces not realized in space, thereby providing a possibility of constructing a model of the non-Euclidean plane).  This should be the point that Gauss\index{Gauss, Johann Carl Friedrich (1777--1855)} spoke of with the greatest praise and rare enthusiasm.

\section{Hyperbolic and projective spaces}\label{sec:after1}
Soon after Riemann\index{Riemann, Georg Friedrich Bernhard (1826--1866)} passed away at the age of 39 on June 28, 1866,  his recondite thesis was published by his friend Dedekind\index{Dedekind, Richard (1831--1916)},\footnote{Abhandlungen der K\"{o}niglichen Gesellschaft der Wissenschaften zu G\"{o}ttingen, {\bf 13} (1868), 133--152.} and was taken up by Beltrami\index{Beltrami, Eugenio (1835--1900)} to put an end to the argument on non-Euclidean geometry. In fact, most people at that time were reluctant to accept the unaccustomed geometry because it might very well be possible to fall in with inconsistency after a more penetrating investigation. So as to convince incredulous people, a {\it model} of non-Euclidean geometry needed to be constructed, resting on well-established geometric ingredients, with which inconsistency can be kept at bay. 

In 1868, the Italian journal {\it Giornale di Mathematiche} carried Beltrami's\index{Beltrami, Eugenio (1835--1900)} article {\it Teoria fondamentale degli spazii de curvatura constante},\footnote{\label{fn:Saggio}This is his second paper on a non-Euclidean model. He did not know Riemann's work when  writing the first paper {\it Saggio di Interpretazione della Geometrica Non-euclidea} (1868), wherein the pseudosphere\index{pseudosphere} (fn.~\!\ref{fn:pseudo}) is exploited in a tricky manner.} in which it was shown that the non-Euclidean space\index{non-Euclidean space} is realized as the 3D open ball $D=\{(x_1,x_2,x_3)|~\! x_1{}^2+x_2{}^2+x_3{}^2<4\}$ with $ds^2=(dx_1{}^2+dx_2{}^2+dx_3{}^2)/\big(1-(x_1{}^2+x_2{}^2+x_3{}^2)/4\big)^2$, where geodesics\index{geodesic} are circular arcs in space crossing perpendicularly $\partial D=\{(x_1,$ $x_2,$ $x_3)|$ $~\! x_1{}^2+x_2{}^2+x_3{}^2=4\}$. The apparent boundary $\partial D$ is at infinity, so that the universe modeled by $D$ has infinite extent.\footnote{In the {\it Teoria fondamentale}, it is indicated that $H^d=\{(x_1,\ldots,x_d)\in \mathbb{R}^d|~x_d>0\}$ with $ds^2=(\sum_{i=1}^ddx_i{}^2)/x_d{}^2$ has constant negative curvature. This manifold, together with other ones, was later to be called the {\it hyperbolic space} (F. Klein\index{Klein, Christian Felix (1849--1925)}, 1871); see \cite{sti}.} In the same vein, the non-Euclidean plane is represented by the 2D open disk $\{(x_1,x_2)|~\! x_1{}^2+x_2{}^2<4\}$ with $ds^2=\displaystyle (dx_1{}^2+dx_2{}^2)/\big(1-(x_1{}^2+x_2{}^2)/4\big)^2$.\footnote{Thinking back over Riemann's\index{Riemann, Georg Friedrich Bernhard (1826--1866)} thesis, we wonder why he did not refer to the non-Euclidean geometry\index{non-Euclidean geometry} since (\ref{eq:constcurve}) for $\alpha<0$ is connotative of Beltrami's model\index{Beltrami, Eugenio (1835--1900)}. In a letter to J.
Ho\"{u}el dated April 4, 1868, Beltrami\index{Beltrami, Eugenio (1835--1900)} says, ``What amazes me is that for all the time I talked with Riemann\index{Riemann, Georg Friedrich Bernhard (1826--1866)} (during the two years he spent in Pisa, shortly before his sad end), he never mentioned these ideas to me, though they must have occupied him for quite a long time, for a fine draft cannot be the work of a single day, even for such a brilliant genius" (see 
\cite{giac}).}

Here, a short comment is in order. Regarding geodesics\index{geodesic} as straight lines, one may logically build up a theory that satisfies all premises of non-Euclidean geometry. What should be emphasized is that straight lines are not necessarily the ones that we visualize on a piece of paper. Actually, as postulated by Hilbert\index{Hilbert, David (1862--1943)} (Remark \ref{rem:list} (2)), ``point" and ``straight line" are {\it undefined objects}, and only the {\it relations} between them such as ``a point $p$ {\it is on} a straight line $\ell$" and ``there is a unique line {\it passing through} given two points" are a vital necessity. Of great significance in this comprehension is that all the non-Euclidean objects are described in the Euclidean terminology. Thus the logical consistency of non-Euclidean geometry reduces to that of Euclidean geometry, and then to that of the real number system  because $\mathbb{R}^3$ is a model of Euclidean space.\footnote{At the ICM 1900, Hilbert\index{Hilbert, David (1862--1943)} proposed to prove the consistency of the real number system within his {\it formalistic} framework. Yet, this turns out to be not feasible as G\"{o}del showed (1931).}

\medskip

In Sect.~\!\ref{sec:proj}, we touched on projective geometry\index{projective geometry} as a branch of classical geometry. For the pedagogical sake, we shall present an alternative perspective that was developed in the same period as Poncelet\index{Poncelet, Jean-Victor (1788--1867)} and Jakob Steiner\index{Steiner, Jakob (1796--1863)} (1796--1863) embarked on it in view of synthetic geometry. We owe the idea in the main to August Ferdinand M\"{o}bius\index{M\"{o}bius, August Ferdinand (1790--1868)} (1790--1868).\footnote{\label{fn:steiner}M\"{o}bius, {\it Der barycentrische Calc\"{u}l}, 1827. A similar approach was taken 
%in the paper {\it Analytisch-geometrische Entwickelungen} (1828 and 1831) 
by J. Pl\"{u}cker. Meanwhile, Steiner stuck to the synthetic method.} It provides another exemplar of the transition to a ``geometry of space itself." 

A ``plane" that one assumes a priori in projective geometry\index{projective geometry} is called today the {\it projective plane}\index{projective plane}. Plane projective geometry is, if we distill its content down to a few essential premises, built on the following simple system of axioms that are stated in terms of relations between points and lines:

\smallskip
(i)  Any two distinct points are contained in one and only one line.

\smallskip
(ii) Any two distinct lines have one and only one point in common.

\smallskip

(iii) there exists four points, no three of which are contained in one line.

\smallskip
 
Although this might cause a sense of discomfort at first, a concrete model of the projective plane\index{projective plane} with ``points" and ``lines", symbolically indicated by $\mathbb{P}$,  is given by the totality of lines in $\mathbb{R}^3$ passing through the origin $O$. Here a ``point" is just an element of $\mathbb{P}$. A ``line" is defined to be an ordinary plane in space passing through $O$. We say that a ``point" $\ell$ is contained in a ``line" $H$ if the line $\ell$ is contained in the plane $H$ in the usual sense. The axioms (i), (ii), (iii) above are easy consequences of the facts that hold for lines and planes.

\begin{figure}[htbp]
%\vspace{-4.5cm}
\vspace{-0.4cm}
\begin{center}
\includegraphics[width=.4\linewidth]
{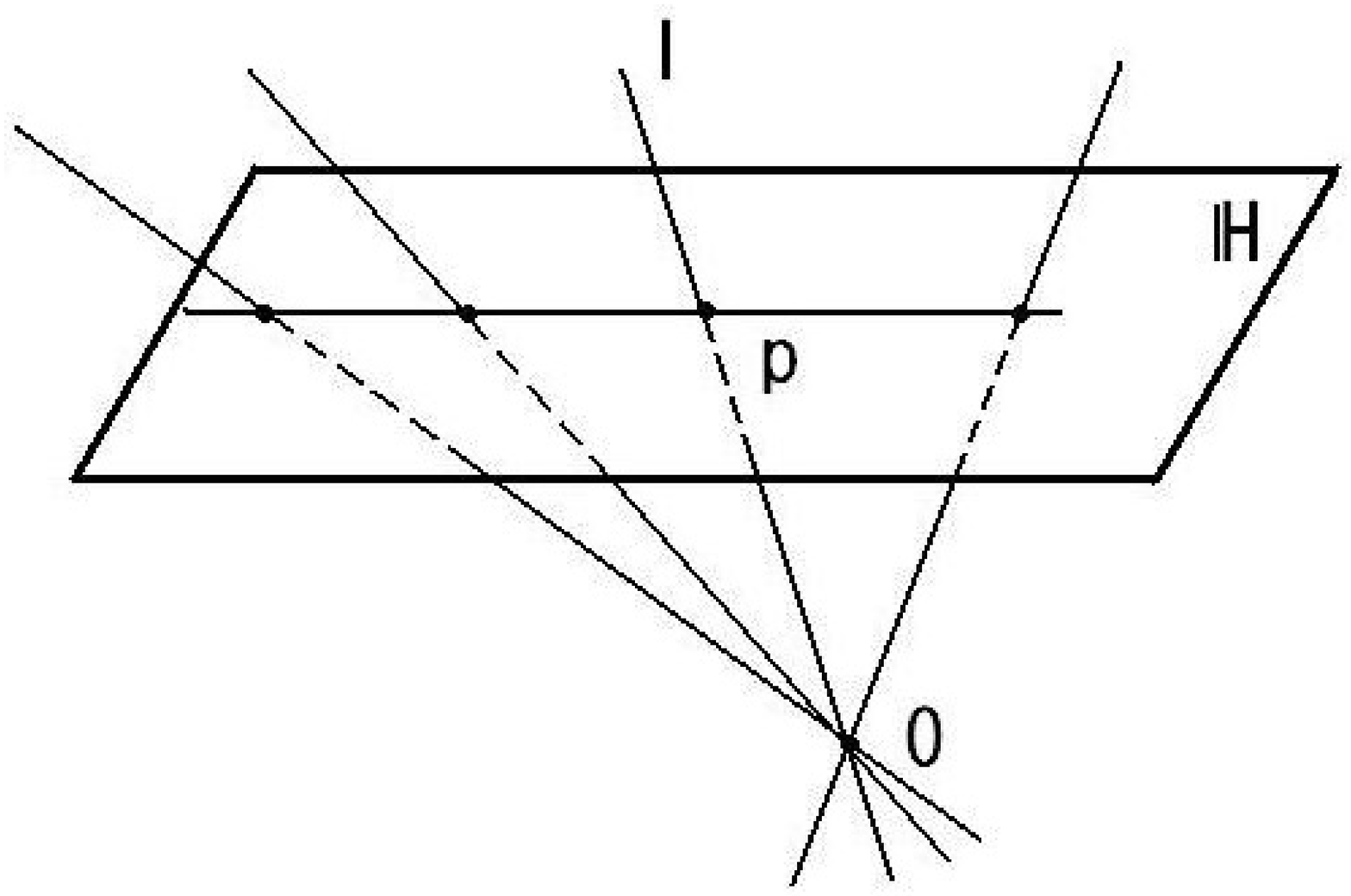}
\end{center}
\vspace{-0.8cm}
\caption{Projective plane}\label{fig:proplane}
\vspace{-0.3cm}
\end{figure}

In Sect.~\!\ref{sec:proj}, we said that, in plane projective geometry\index{projective geometry}, points at infinity are adjoined to the Euclidean plane. To confirm this in our setting, select a plane $\mathbb{H}$ not containing $O$, and to take the line $\ell_p$ containing both $p\in\mathbb{H}$ and $O$. We let $\mathbb{H}(\mathbb{P})$ be the collection of lines $\ell\in \mathbb{P}$ not parallel to $\mathbb{H}$, and let $\mathbb{H}_{\infty}(\mathbb{P})$ be the plane passing through $O$ and parallel to $\mathbb{H}$. Then $p\mapsto \ell_p$ is a one-to-one correspondence\index{one-to-one correspondence} between $\mathbb{H}$ and $\mathbb{H}(\mathbb{P})$, which allows us to identify $\mathbb{H}(\mathbb{P})$ with $\mathbb{H}$. Fig.~\!\ref{fig:proplane} illustrates that, as a point in $\mathbb{H}$ tends to infinity, the line $\ell_p$ converges to a line parallel to $\mathbb{H}$. Accordingly, one may regards $\mathbb{H}_{\infty}(\mathbb{P})$ as the line at infinity.

Now, $\mathbb{P}$ is identified with $P^2(\mathbb{R})$ in a natural manner. We call $(x,y,z)\in \mathbb{R}^3\backslash \{{\bf 0}\}$ the {\it homogeneous coordinate}\index{homogeneous coordinate} of the ``point" $\ell\in \mathbb{P}$ when $\ell$ passes through $(x,y,z)$, and write $\ell=[x,y,z]$ by convention. This concept is immediately generalized to $P^d(\mathbb{R})$, and also to the realm of complex numbers.\footnote{One can consider the projective planes over quaternions\index{quaternion} and octonions as well. It is interesting to point out that Desargues' theorem\index{Desargues' theorem} holds for quaternions\index{quaternion}, but does not for octonions\index{octonion}, and that Pappus' hexagon theorem\index{Pappus' hexagon theorem} does not hold for quaternions.} The simplest case is the {\it complex projective line} $P^1(\mathbb{C})=\{[z,w]|$ $(z,w)\neq (0,0), ~z,w\in \mathbb{C}\}$, which coincides with the {\it Riemann sphere}, the complex plane plus a point at infinity, introduced by Riemann\index{Riemann, Georg Friedrich Bernhard (1826--1866)} (fn.~\!\ref{fn:Abelcschen}). A projective transformation\index{projective transformation} in this case is represented as
$T(z)=(az+b)/(cz+d)~(a,b,c,d\in \mathbb{C},~ad-bc\neq 0)$.

Henri Poincar\'{e}\index{Poincar\'{e}, Henri (1854--1912)} (1854--1912) observed a conflation of $P^1(\mathbb{C})$ and non-Euclidean geometry in his study  of {\it Fuchsian functions} (a sort of ``periodic" functions of one complex variable). Indeed, the above $T$ with $a,b,c,d\in \mathbb{R}$ satisfying $ad-bc\!>\!0$ preserves not only $H\!=\!\{z\!=\!x+yi\!\in\!\mathbb{C}|~y\!>\!0\}$, but also $ds^2=|dz|^2/y^2$ $~(|dz|^2=dx^2+dy^2)$. Since $(H, ds^2)$ is a model of the non-Euclidean plane (fn.~\!\ref{fn:Saggio}), $T$ is a non-Euclidean congruent transformation.\footnote{Poincar\'{e}, {\it Sur les Fonctions Fuchsiennes}, Acta Math., {\bf 1} (1882), 1--62. For this reason, $(H, ds^2)$ is called the {\it Poincar\'{e} half plane}. The line element $ds$ was previously obtained by Liouville\index{Liouville, Joseph (1809--1882)} through transformation of the line element on the pseudosphere\index{pseudosphere} (1850).} More conspicuously, this idea led up to an entirely new phase of {\it analytic number theory} (initiated by Euler\index{Euler, Leonhard (1707--1783)} and definitively established by Dirichlet\index{Dirichlet, Johann Peter Gustav Lejeune (1805--1859)}) through the notion of {\it automorphic form}, a generalization of Fuchsian functions.

Another non-Euclidean model that is cross-fertilized by projective geometry\index{projective geometry} was given by Beltrami\index{Beltrami, Eugenio (1835--1900)} in his {\it Saggio} and independently by F. Klein\index{Klein, Christian Felix (1849--1925)} (1871). This model, called the {\it Klein model}\index{Klein model} and first constructed in 1859 by A. Cayley\index{Cayley, Arthur (1821--1895)} in the context of real projective geometry\index{projective geometry} without reference to non-Euclidean geometry\index{non-Euclidean geometry} (1859), is realized as the open unit disk $D=\{(x,y)\in \mathbb{R}^2|~x^2+y^2<1\}$, in which the ``line segments" are chords in $D$. Worthy of mention is that the cross-ratio\index{cross-ratio} pops up in the distance between two points. To be exact, the distance between two points $P$ and $Q$ in $D$ can be expressed as $d(P,Q)=\frac{1}{2} \log \big|[A, P, Q, B]\big|$, where $A,B$ are the points of intersection on the circle $\partial D$ with the line connecting $P,Q$ (see Remark \ref{rem:spacialrelativity}~\!(5)). This distance is the one associated with $ds^2=\big\{\big(1-(x^2+y^2)\big)(dx^2+dy^2)+(xdx+ydy)^2\big\}/\big(1-(x^2+y^2)\big)^2$.

\section{Absolute differential calculus}\label{sec:after2}

Riemann's\index{Riemann, Georg Friedrich Bernhard (1826--1866)} work was evolved further by Elwin Bruno Christoffel\index{Christoffel, Elwin Bruno (1829--1900)} (1829--1900) and Rudolf Lipschitz\index{Lipschitz, Rudolf (1832--1903)} (1832--1903). The significance of Christoffel's work lies in his introduction of the concept later called {\it covariant differentiation}\index{covariant differentiation} and of the {\it curvature tensor}\index{curvature tensor} which crystallizes the notion of Riemann's curvature (1869).\footnote{{\it \"{U}ber die Transformation der homogenen Differentialausdr\"{u}cke zweiten Grades}, J.~\!f\"{u}r die Reine und Angew.~\!Math., {\bf 70} (1869), 46–-70;~ {\it \"{U}ber ein die Transformation homogener Differentialausdr\"{u}cke zweiten Grades betreffendes Theorem}, ibid., 241--245.} After a while, Gregorio Ricci-Curbastro\index{Ricci, Gregorio Ricci-Curbastro (1853--1925)} (1853--1925) and Tullio Levi-Civita\index{Levi-Civita, Tullio (1873-1941)} (1873-1941) initiated  ``absolute differential calculus\index{absolute differential calculus}"\footnote{{\it M\'{e}thodes de calcul diff\'{e}rentiel absolu et leurs applications}, Math.~\!Ann., {\bf 54} (1900), 125--201.}---calculus, in line with Christoffel's work, dealing with ``invariant" geometric quantities (called {\it tensor fields}\index{tensor field}). Levi-Civita\index{Levi-Civita, Tullio (1873-1941)} further brought in {\it parallel transports}\index{parallel transport} as a partial generalization of parallel translations.\footnote{{\it Nozione di parallelismo in una variet\`{a}qualunque e consequente specificazione geometrica della curvatura Riemanniana}, Rend.~\!Circ.~\!Mat.~\!Palermo, {\bf 42} (1917), 73--205.} Before long, absolute differential calculus\index{absolute differential calculus} was renamed ``tensor calculus"\index{tensor calculus} and used effectively by  Einstein\index{Einstein, Albert (1879--1955)} for his general relativity\index{general relativity}.\footnote{{\it Die Feldgleichungen der Gravitation}. Sitzungsberichte der Preu{\ss}ischen Akademie der Wissenschaften zu Berlin: 844--847 (November 25, 1915).} In this groundbreaking theory, a slight modification of Riemann's theory is required; namely, the matrix $(g_{ij})$ associated with the first fundamental form\index{first fundamental form} is a four-by-four symmetric matrix with three positive and one negative eigenvalue. A manifold with such a fundamental form is called a curved space-time, or a {\it Lorentz manifold}\index{Lorentz manifold} after Hendrik Anton Lorentz\index{Lorentz, Hendrik Anton (1853--1928)} (1853--1928).

\medskip

We shall now plunge ourselves into tensor calculus. Henthforth the matrix $(g_{ij})$ is supposed to be invertible, so that our discussion includes the case of Lorentz manifolds\index{Lorentz manifold}. A tensor field\index{tensor field} (in the classical sense) is roughly a system of multi-indexed functions depending on the choice of a coordinate system and satisfying a certain transformation rule under coordinate transformations so as to yield a {\it coordinate-free} quantity.
Tensor fields thus defined fit in with the {\it principle of relativity} which  Einstein\index{Einstein, Albert (1879--1955)} set up to claim that the laws of physics should have the same form in all admissible coordinate systems of reference. 

An example is the coefficients $\{g_{ij}\}$ of the first fundamental form\index{first fundamental form} (called the {\it metric tensor}). Here we should note that $ds^2=\sum_{i,j=1}^dg_{ij}dx_idx_j$ is a coordinate-free quantity. If we take another coordinate system $(y_1,\ldots,y_d)$, we have a different expression $ds^2=\sum_{h,k=1}^d$ $\overline{g}_{hk}dy_hdy_k$. The chain rule for differentiation applied to the coordinate transformation $x_i=x_i(y_1,\ldots,y_d)$ $~(i=1,\ldots,d)$ is written as $dx_i=\sum_{h=1}^d\frac{\partial x_i}{\partial y_h}dy_h
$. So substituting this for (\ref{eq:infinitesimalform}), we obtain $ds^2=\sum_{h,k=1}^d\sum_{i.j=1}^dg_{ij}\frac{\partial x_i}{\partial y_h}\frac{\partial x_j}{\partial y_k}dy_hdy_k$. Comparing this with $ds^2=\sum_{h,k=1}^d\overline{g}_{hk}dy_hdy_k$, we have the transformation rule $\overline{g}_{hk}=\sum_{i.j=1}^dg_{ij}\frac{\partial x_i}{\partial y_h}\frac{\partial x_j}{\partial y_k}$.

Another example is a {\it contravariant vector}\index{contravariant vector} $\xi^i$ $~(i=1,\ldots,d)$ obeying the rule $\overline{\xi}^i=\sum_{j=1}^d\xi^j\frac{\partial y_i}{\partial x_j}$. For instance, $\xi^i=\frac{dx_i}{dt}$ for a curve $c(t)=\big(x_1(t),\ldots, x_d(t)\big)$ is a contravariant vector because $\frac{dy_i}{dt}=\sum_{j=1}^d\frac{dx_j}{dt}\frac{\partial y_i}{\partial x_j}$ $~(y_i(t)=y_i\big(x_1(t),\ldots,x_d(t)\big))$. Hence, a contravariant vector is reckoned as a vector tangent to the manifold.

A general tensor field\index{tensor field} (of type $(h,k)$) is a system of multi-indexed functions $\xi^{i_1\cdots i_h}_{j_1\cdots j_k}$ $~(i_1,\ldots,i_h,j_1,\ldots,j_k=1,2,\ldots,d)$ satisfying the transformation formula
$$
\overline{\xi}{}^{a_1\cdots a_h}_{b_1\cdots b_k}=\sum_{i_1,\ldots, i_h=1}^d\sum_{j_1,\ldots,j_k=1}^d\xi^{i_1\cdots i_h}_{j_1\cdots j_k}\frac{\partial y_{a_1}}{\partial x_{i_1}}\cdots\frac{\partial y_{a_h}}{\partial x_{i_h}}\frac{\partial x_{j_1}}{\partial y_{b_1}}\cdots\frac{\partial x_{j_k}}{\partial y_{j_k}}.
$$
Hence $g_{ij}$ is a tensor field of type $(0,2)$. For the inverse matrix $(g^{ij})=(g_{ij})^{-1}$, we have  
$
\overline{g}^{hk}=\sum_{i.j=1}^dg^{ij}\frac{\partial y_i}{\partial x_h}\frac{\partial y_j}{\partial x_k}$, so that $g^{ij}$ is a tensor field of type $(2,0)$.

The {\it covariant derivative}\index{covariant derivative} of $\xi^{i_1\cdots i_h}_{j_1\cdots j_k}$ is a tensor field of type $(h,k+1)$ given by
{\small 
\begin{eqnarray*}
\nabla_i\xi^{i_1\cdots i_h}_{j_1\cdots j_k}&:=&\frac{\partial }{\partial x_i}\xi^{i_1\cdots i_h}_{j_1\cdots j_k}+\sum_{a=1}^d\left\{ { {i_1}\atop{a i} } \right\}
\xi^{ai_2\cdots i_h}_{j_1\cdots j_k}+\cdots+\sum_{a=1}^d\left\{ { {i_h}\atop{a i} } \right\}\xi^{i_1\cdots i_{h-1}a}_{j_1\cdots j_k}\\
&&-\sum_{a=1}^d\left\{ { {a}\atop{j_1 i} } \right\}
\xi^{i_1\cdots i_h}_{aj_2\cdots j_k}-\cdots-\sum_{a=1}^d\left\{ { {a}\atop{j_k i} } \right\}\xi^{i_1\cdots i_{h}}_{j_1\cdots j_{k-1}a},
\end{eqnarray*}}
\hspace{-0.1cm}where
{\small \begin{equation}\label{eq:chris}
 \left\{ { {h}\atop{i j} } \right\} = \sum_{a=1}^d \frac{1}{2} g^{h a} \left(\frac{\partial g_{i a} }{\partial x^j } + \frac{\partial g_{j a} }{\partial x^i } - \frac{\partial g_{j i} }{\partial x^a } \right)
\end{equation}}
\hspace{-0.1cm}is what we call the {\it Christoffel symbol}\index{Christoffel symbol}.\footnote{The symbol Christoffel\index{Christoffel, Elwin Bruno (1829--1900)} used is $\left\{ { {i j}\atop{h} } \right\}$. R. Lipschitz\index{Lipschitz, Rudolf (1832--1903)} (reputable for ``Lipschitz continuity") brought in a simiar symbol in his {\it Untersuchungen in Betreff der ganzen homogenen Functionen von n Differentialen}, J.~\!f\"{u}r die Reine und Angew.~\!Math., {\bf 70} (1869), 71-102.} The terms involving the Christoffel symbols\index{Christoffel symbol} is thought of as making the derivation intrinsic\index{intrinsic}. The curvature tensor\index{curvature tensor} $R^i_{jkl}$ is then defined to be a tensor that measures non-commutativity of covariant differentiation; namely, $\nabla_{l}\nabla_{k}\xi^j-\nabla_{k}\nabla_{l}\xi^j=\sum_{i=1}^dR^j_{il k}\xi^i$, where
$$
R^i_{jkl} = \frac{\partial }{\partial x_k}  \Big\{ \displaystyle{ {i}\atop{l j} } \Big\}- \frac{\partial }{\partial x_l}\Big\{ \displaystyle{ {i}\atop{k j} } \Big\} + \sum_{a=1}^d \Big\{ { {i}\atop{k a} } \Big\}\Big\{ { {a}\atop{l j} } \Big\} - \sum_{a=1}^d \displaystyle\Big\{ { {i}\atop{l a} }\displaystyle \Big\}\Big\{ { {a}\atop{k j} } \Big\}. 
$$

Although being implicit, the curvature tensor\index{curvature tensor} turns up in Riemann's paper submitted to the Acad\'{e}mie des Sciences on July 1, 1861 (remained unknown until 1876).\footnote{{\it Commentatio Mathematica, qua respondere tentatur quaestioni ab ${\rm Ill}^{ma}$}, Academia Parisiensi propositae.} In the second half of this paper, which was his answer to the prize question on heat distribution posed by the Academy in 1858, he gave a condition in order that $ds^2=\sum_{i,j=1}^dg_{ij}dx_idx_j$ is flat, or equivalently that there exists a curvilinear coordinate system $(y_1,\ldots,y_d)$ such that $ds^2=dy_1{}^2+\cdots+dy_d{}^2$. His condition can be transliterated into $R^i_{jkl}\equiv 0$ (see \cite{far} for the detail).

We defer the details of how the curvature tensor\index{curvature tensor} is linked with Riemann's sectional curvature\index{sectional curvature} (see Sect.~\!\ref{sect:reimenniangeometry}), and only cite the fact that the Gaussian curvature\index{Gaussian curvature} $K_S$ of a surface $S$ is expressed as $K_S=R^1_{221}/(g_{11}g_{22}-g_{12}{}^2)$, so that $K_S$ is a rational function (not depending on $S$) in $E, F, G$ and their derivatives up to second order as shown by Gauss\index{Gauss, Johann Carl Friedrich (1777--1855)} with ``bare hands" ({\it Disquisitiones generales}, Arts.~\!9-11). More specifically (in the present notation),
{\small 
\begin{eqnarray}
&&\hspace{-1.5cm}K_S=\frac{1}{\sqrt{EG-F^2}}\Big[
\frac{\partial}{\partial u}\Big(\frac{\sqrt{EG-F^2}}{G}\left\{ { {1}\atop{2 2} } \right\} \Big)-\frac{\partial}{\partial v}\Big(\frac{\sqrt{EG-F^2}}{G}\left\{ { {1}\atop{1 2} } \right\} \Big)
\Big],\label{eq:gaussthe}\\
&~~& \nonumber\\
&&\qquad ~~\left\{ { {1}\atop{2 2} } \right\}=\frac{2GF_v-GG_u-FG_v}{2(EG-F^2)},~~
\left\{ { {1}\atop{1 2} } \right\}=\frac{GE_v-FG_u}{2(EG-F^2)}.\nonumber
\end{eqnarray}

}

Now take this for granted. For an isometry $\varPhi: S\rightarrow \overline{S}$, the composition $\overline{\boldsymbol{S}}(u,v)=\varPhi(\boldsymbol{S}(u,v))$ gives a local parametric representation of $\overline{S}$ such that $\overline{E}(u,v)=E(u,v)$, $\overline{F}(u,v)=F(u,v)$, $\overline{G}(u,v)=G(u,v)$ because $\varPhi$ preserves $ds^2$. Therefore the derivatives of $\overline{E},\overline{F},\overline{G}$ of any order coincide with the corresponding derivatives of $E,$ $F,G$, from which Gauss's\index{Gauss, Johann Carl Friedrich (1777--1855)} Theorema Egregium\index{Theorema Egregium} follows.

\begin{rem}{\rm {\small
The following formulas are due to Gauss\index{Gauss, Johann Carl Friedrich (1777--1855)}, K. M. Peterson, G. Mainardi, and D. Codazzi:
\begin{eqnarray}
&&L_v-M_u=\left\{ { {1}\atop{1 2} }\right\}L+\left(\left\{ { {1}\atop{1 1} }\right\}-\left\{ { {2}\atop{1 2} }\right\}\right)M+\left\{ { {2}\atop{1 1} }\right\}N,\label{eq:coda1}\\
&&N_u-M_v=\left\{ { {1}\atop{2 2} }\right\}L+\left(\left\{ { {2}\atop{2 2} }\right\}-\left\{ { {1}\atop{1 2} }\right\}\right)M-\left\{ { {2}\atop{1 2} }\right\}N\label{eq:coda2}.
\end{eqnarray}
In 1867, Bonnet\index{Bonnet, Pierre Ossian (1819--1892)} proved the converse (the {\it fundamental theorem of surface theory}): if functions $E,F,G$ and $L,M,N$ in variables $u$ and $v$ satisfy (\ref{eq:coda1}), (\ref{eq:coda2}), (\ref{eq:gaussthe}), while $EG-F^2>0$, then there exists a unique surface (up to motion) in space which admits $E,F,G$ and $L,M,N$ as the coefficients of the first and second fundamental forms.
\hfill$\Box$

}
}

\end{rem}

In terms of the Christoffel symbol\index{Christoffel symbol}, the equation of a geodesic  $c(t)=\big(x_1(t),$ $\ldots,$ $x_d(t)\big)$ is expressed as 
\begin{equation}\label{eq:geodesic} 
\frac{d^2x_i}{dt^2}+\sum_{j,k=1}^d\Big\{ \displaystyle{ {i}\atop{jk} } \Big\}\frac{dx_j}{dt}\frac{dx_k}{dt}=0\qquad (i=1,\ldots,d).
\end{equation}
This expression in a curved space-time---in alliance with imaginary experiments carried out within a closed chamber in interstellar space to confirm the equivalence of gravity and acceleration\footnote{\label{fn:floating}The freely floating chamber produces a gravity-free state inside, while the accelerated chamber produces a force indistinguishable from gravitation.}---is the starting point of\index{Einstein, Albert (1879--1955)} Einstein's theory\index{general relativity}; that is, we think of (\ref{eq:geodesic}) as an analogue of the Newtonian equation of motion\index{Newtonian equation of motion} $d^2x_i/dt^2=-\partial u/\partial x_i$ under gravitational potential\index{gravitational potential} $u$ (i.e., the equation for ``free fall").\footnote{{\it Grundlage der allgemeinen Relativit\"{a}tstheorie}, Ann.~\!der Physik, {\bf 49} (1916), 769--822.} Pursuing this analogy, he concluded that the first fundamental form\index{first fundamental form} may be loosely regarded as a generalization of the gravitational potential\index{gravitational potential},\footnote{In his thesis (Part III, \S 3), Riemann\index{Riemann, Georg Friedrich Bernhard (1826--1866)} prophetically says, ``We must seek the ground of its metric relations outside it, in {\it binding forces} which act upon it."} and that the fundamental interaction of gravitation as a result of {\it space-time being curved by matter and energy} is manifested by the {\it field equations}\index{Einstein field equations}:
\begin{equation}\label{eq:waveeq}
R_{ij}-\frac{1}{2}Rg_{ij}+\Lambda g_{ij}=\frac{8\pi G}{c^4}T_{ij},
\end{equation}
where $R_{ij}=\sum_{k=1}^d R^{k}_{ikj}$ is the {\it Ricci tensor}\index{Ricci tensor}, $R=\sum_{i,j=1}^dg^{ij}R_{ij}$ is the {\it scalar curvature} (note that $R=2K_S$ for a surface $S$), $\Lambda$ is the {\it cosmological constant},\footnote{This constant is used to explain the observed acceleration of the expansion of the universe. (In 1929, E. Hubble discovered that the universe is expanding.)} $c$ is the speed of light in vacuum, and $T_{ij}$ is the {\it stress--energy tensor} (an attribute of matter, radiation, and non-gravitational force fields).

The left-hand side of (\ref{eq:waveeq}) as an operator acting on $g_{ij}$ is a non-linear generalization of the {\it d'Alembertian} $\displaystyle-\partial^2/\partial t^2+c^2\Delta$ showing up in the description of wave propagation with speed $c$.  This observation brought Einstein\index{Einstein, Albert (1879--1955)} to the prediction of {\it gravitational waves}, which was successfully detected by the Laser Interferometer Gravitational-Wave Observatory (11th February, 2016).

What is remarkable is, as observed by Hilbert\index{Hilbert, David (1862--1943)} (1915),
%\footnote{{\it Die Grundlagen der Physik}, 1915.} 
that the vacuum field equations ($R_{ij}=0$) are the E-L equation\index{Euler-Lagrange equation} associated with the action integral\index{action integral}: $S=\int R\sqrt{-\det(g_{ij})}~\!dx_1dx_2dx_3dx_4$.

\begin{rem}\label{rem:spacialrelativity}{\rm {\small 
(1) D'Alembert\index{d'Alembert, Jean-Baptiste le Rond (1717--1783)} showed that solutions of the wave equation $\frac{\partial^2 y}{\partial t^2}-a^2\frac{\partial^2y}{\partial x^2}=0$ are of the form $f(x+at)$ $+$ $g(x-at)$, where $f$ and $g$ are {\it any} functions (1747). On the other hand, Daniel Bernoulli\index{Bernoulli, Daniel (1700--1782)} (1700--1782), a son of Johann Bernoulli\index{Bernoulli, Johann (1667--1748)} and a friend of Euler\index{Euler, Leonhard (1707--1783)}, 
came up with {\it trigonometric series}\index{trigonometric series} for the first time, when he tried to solve the wave equation (1747). The two ways to express solutions stirred up controversy on the nature of  ``arbitrary" functions (Remark \ref{rem:function}).

\smallskip

(2) 
The Laplacian $\Delta$ named after Laplace\index{Laplace, Pierre-Simon (1749--1827)} appears in the theory of heat conduction by Jean Baptiste Joseph Fourier\index{Fourier, Jean Baptiste Joseph (1768--1830)} (1768--1830). While he was Governor of Grenoble (appointed by Napoleon I), he carried out experiments on the propagation of heat along a metal bar, and used the trigonometric series\index{trigonometric series} to solve the heat equation $\frac{\partial y}{\partial t}=a^2\frac{\partial^2y}{\partial x^2}$ that describes the distribution of heat over a period of time (1807 and 1822); see Remark \ref{rem:function}. An interesting fact is that a special solution is given by the function $(1/\sqrt{4\pi a^2t})\exp(-x^2/4a^2t)$, which shows up as the {\it normal distribution} in mathematical statistics that A. de Moivre, Gauss\index{Gauss, Johann Carl Friedrich (1777--1855)}, and Laplace\index{Laplace, Pierre-Simon (1749--1827)} pioneered.

\smallskip

(3) (\ref{eq:waveeq}) is also regarded as an analogue of {\it Poisson's equation}\index{Poisson's equation} $\Delta u=4\pi G\rho$, where $u$ is the {\it gravitational potential}\index{gravitational potential} associated with a mass density $\rho$.\footnote{The density $\rho$ for a mass distribution $\mu$ is defined by the relation $\rho(\boldsymbol{x})d\boldsymbol{x}=d\mu(\boldsymbol{x})$.} As noticed by Green\index{Green, George (1793--1841)} (1828) and Gauss\index{Gauss, Johann Carl Friedrich (1777--1855)} (1839; {\it Werke} V, 196--242), this equation embodies the law of universal gravitation\index{law of universal gravitation} (and hence $\Delta u=0$ in the vacuum case). Indeed, $u(\boldsymbol{x})=-G\int \|\boldsymbol{x}-\boldsymbol{y}\|^{-1}\rho(\boldsymbol{y})d\boldsymbol{y}$ is a solution of Poisson's equation.

A heuristic (and intrepid) way to derive (\ref{eq:waveeq}) in the vacuum case is to take a look at the tendency of free-falling objects to approach or recede from one another, which is described by the ``variation vector field\index{variation vector field}" $J(t):= \partial \boldsymbol{x}/\partial s\big|_{s=0}$, where $\boldsymbol{x}(t,s)\!=\!(x_i(t,s))$ $(-\epsilon\!<\!s<\!\epsilon)$ is a family of solutions of  $d^2x_i/dt^2\!=\!-\partial u/\partial x_i$. Clearly $J(t)$ satisfies $d^2 J/dt^2+H(J)\!=\!0$, where $H$ is the linear map defined by $H(\boldsymbol{z})\!=\!\sum_{j=1}^d\frac{\partial^2 u}{\partial x_i\partial x_j}z_j$. Notice that ${\rm tr}(H)\!=\!\Delta u\!=\!0$. In turn, the variation vector field $J(t)$ for the geodesic equation (\ref{eq:geodesic}) satisfies the apparently similar equation $\frac{D}{dt}\left(\frac{D}{dt}J\right)+K(J)=0$ (called the {\it Jacobi equation}\index{Jacobi equation}), where $K(\boldsymbol{z})=\sum_{j=1}^dR^{i}_{hjk}\dot{x}_h\dot{x}_kz_j$, and $\frac{DX}{dt}=\frac{dX_i}{dt}+\sum_{h,k=1}^d\left\{ { {i}\atop{h k} } \right\}\frac{dx_h}{dt}X_k$, the covariant derivative of a vector field $X$ along the curve $\boldsymbol{x}(\cdot,0)$ (Sect.~\!\ref{sect:reimenniangeometry}). Since ${\rm tr}(K)=\sum_{h,k=1}^dR_{hk}\dot{x}_h\dot{x}_k$, we may regard the equation $\sum_{h,k=1}^dR_{hk}\dot{x}_h\dot{x}_k=0$  (or $R_{hk}\equiv 0$ since $\dot{\boldsymbol{x}}$ is arbitrary) as an analogue of $\Delta u=0$.

\smallskip

(4) The theory of {\it special relativity}\index{special relativity} created by Einstein\index{Einstein, Albert (1879--1955)} in 1905 is literally a special case of general relativity\index{general relativity}, which has radically changed our comprehension of time and space.\footnote{{\it Zur Elektrodynamik bewegter K\"{o}rper}, Ann.~\!der Physik, {\bf 17} (1905), 891--921.} Roughly speaking, it is a theory of space-time under absence of gravitation. Einstein\index{Einstein, Albert (1879--1955)} assumed that (i) the speed of light $c$ $(=2.99792458 \times 10^8$ m/s) in a vacuum is the same for all observers, regardless of the motion of the light source, and (ii) the physical laws are the same for all non-accelerating observers (this is a special case of the principle of relativity). His setup squares with the theory of electromagnetic waves due to James Clerk Maxwell\index{Maxwell, James Clerk (1831--1879)} (1831--1879); see Sect.~\!\ref{sec:right}. The special relativity was mathematically formulated by Hermann Minkowski\index{Minkowski, Hermann (1864--1909)} (1864--1909).\footnote{Minkowski, {\it Raum und Zeit}, Physikalische Zeitschrift {\bf 10} (1907), 75--88.} His model (called the {\it Minkowski space-time}\index{Minkowski space-time}) is a 4D {\it affine space} with a coordinate system $(x_1,x_2,x_3,t)$ for which the line element is given by $ds^2=dx_1{}^2+dx_2{}^2+dx_3{}^2-c^2dt^2$. Such a coordinate system (stemming from the invariance of the speed of light) is called an {\it inertial system}, which an observer employs to describe events in space-time. If $(y_1,y_2,y_3,s)$ is another inertial system moving in uniform velocity relative to $(x_1,x_2,x_3,t)$, then the relation between them is given by an inhomogeneous {\it Lorentz transformation}\index{Lorentz transformation}: 
\begin{eqnarray}\label{lorenz}
\begin{pmatrix}
\boldsymbol{y}\\
 s
\end{pmatrix}=
\begin{pmatrix}
A & -A\boldsymbol{v} \\
\mp c^{-2}\left(1-\displaystyle\frac{\|\boldsymbol{v}\|^2}{c^2}\right)^{-1/2}~^t\boldsymbol{v} & \pm \left(1-\displaystyle\frac{\|\boldsymbol{v}\|^2}{c^2}\right)^{-1/2}
\end{pmatrix}
\begin{pmatrix}
\boldsymbol{x}  \\
t 
\end{pmatrix}+
\begin{pmatrix}
\boldsymbol{b}  \\
t_0 
\end{pmatrix}
\end{eqnarray}
($\boldsymbol{v}$, $\boldsymbol{b}\in \mathbb{R}^3$,  
and $\boldsymbol{x}={}^t(x_1,x_2,x_3)$, $\boldsymbol{y}={}^t(y_1,y_2,y_3)$). 
The $3\times 3$ matrix $A$ satisfies 
$$
~^t\!AA=I_3+\frac{\|\boldsymbol{v}\|^2}{c^2-\|\boldsymbol{v}\|^2}P_{\boldsymbol{v}},
\quad  |\det A|=\left(1-\frac{\|\boldsymbol{v}\|^2}{c^2}\right)^{-1/2},
$$
where $P_{\boldsymbol{v}}$ is the orthogonal projection of $\mathbb{R}^3$ onto the line $\mathbb{R}\boldsymbol{v}$. Since $\boldsymbol{y}=A(\boldsymbol{x}-t\boldsymbol{v})+\boldsymbol{b}$, the vector $\boldsymbol{v}$ is the {\it relative velocity} of two systems. Equation (\ref{lorenz}) tells us that the relative speed $\|\boldsymbol{v}\|$ cannot exceed $c$. Notice that, formally making $c$ tend to infinity, we obtain a {\it Galilei transformation} compatible with Newtonian mechanics:
\begin{eqnarray*}
\begin{pmatrix}
\boldsymbol{y}\\
 s
\end{pmatrix}=
\begin{pmatrix}
A & -A\boldsymbol{v} \\
0 & \pm 1
\end{pmatrix}
\begin{pmatrix}
\boldsymbol{x}  \\
t 
\end{pmatrix}+
\begin{pmatrix}
\boldsymbol{b}  \\
t_0 
\end{pmatrix}, \qquad ~A\in {\rm O}(3).
\end{eqnarray*}

With (\ref{lorenz}), one may interpret the {\it Lorentz-FitzGerald contraction}  of moving objects. Before special relativity came out, Lorentz\index{Lorentz, Hendrik Anton (1853--1928)} and G. F. FitzGerald posited the idea of contraction to rescue the existence of {\it aether} (or {\it ether}) from a paradox caused by the outcome of Michelson-Morley experiment (1887).\footnote{
%Lorentz, {\it De relatieve beweging van de aarde en den aether}, Zittingsverlag Akad.~\!V.~\!Wet., {\bf 1} (1892), 74--79. 
Aether, whose existence was discarded by Einstein\index{Einstein, Albert (1879--1955)}, had been thought to permit the determination of our absolute motion and also to allow electromagnetic waves to pass through it as elastic waves. Going back into history, ancient thinkers (e.g. Parmenides\index{Parmenides (c. 515--c. 450 BC)}, Empedocles ({\it ca}.~\!490 BCE--{\it ca}.~\! 430 BCE), Plato\index{Plato (c. 428--c. 348 BC)}, Aristotle\index{Aristotle (384--322 BC)}) assumed that aether (\textgreek{A<ij'hr}) filled the celestial regions. Descartes\index{Descartes, Ren\'{e} (1596--1650)} insisted that the force acting between two bodies not touching each other is transmitted through it ({\it Trait\'{e} du monde et de la lumi\`{e}re}, 1629--1633).  Newton\index{Newton, Issac (1642--1726)} suggested its existence  in a paper read to the Royal Society (1675). Riemann\index{Riemann, Georg Friedrich Bernhard (1826--1866)} assumed that both ``gravitational and electrostatic effects" and ``optical and magnetic effects" are caused by aether (1853), while Kelvin\index{Kelvin (William Thomson; 1824--1907)}, holding Newton in reverence, said  that he could not be satisfied with Maxwell's\index{Maxwell, James Clerk (1831--1879)} work until a mechanical model of the aether could be constructed (1884).}

\smallskip
(5) The linear part of the transformation given in (\ref{lorenz}) with $c=1$ belongs to ${\rm O}(3,1)$, where ${\rm O}(n,1)$ is the matrix group (called the {\it Lorentz group}\index{Lorentz group}) that preserves the quadratic form $Q(x_1,\ldots,x_n,x_{n+1})=x_1{}^2+\cdots+x_n{}^2-x_{n+1}{}^2$. Interestingly, the subgroup ${\rm O}^+(2,1)$ of ${\rm O}(2,1)$ that preserves the sign of the last coordinate comes up as the congruence\index{congruence} group for the Klein model\index{Klein model} $D$ (note that the set $P=\{[x_1,x_2,x_3]\in P^2(\mathbb{R})|~\!x_1{}^2+x_2{}^2-x_3{}^2<0\}$ is identified with $D$, and that the projective transformation\index{projective transformation} associated with a matrix in ${\rm O}^+(2,1)$ leaves $P$ (and $D$) invariant and preserves the distance $d(z,w)$ on $D$ because of the projective invariance of cross-ratios\index{cross-ratio}).
\hfill$\Box$

}
}
\end{rem}

\section{Cantor's transfinite set theory}\label{subsec:transfinite}
Mathematical terminology to formulate what Riemann\index{Riemann, Georg Friedrich Bernhard (1826--1866)} wanted to convey had not yet matured enough in his time. Here is an excerpt from Riemann's paper referred to above, which indicates how he took pains to impart the idea of manifolds with immatured terms.

\medskip

{\it If in the case of a notion whose specializations form a continuous manifold, one passes from a certain specialization in a definite way to another, the specializations passed over form a simply extended manifold, whose true character is that in it a continuous progress from a point is possible only in two directions, forwards or backwards. If one now supposes that this manifold in its turn passes over into another entirely different, and again in a definite way, namely so that each point passes over into a definite point of the other, then all the specializations so obtained from a doubly extended manifold. In a similar manner one obtains a triply extended manifold, if one imagines a doubly extended one passing over in a definite way to another entirely different; and it is easy to see how this construction may be continued. If one regards the variable object instead of the determinable notion of it, this construction may be described as a composition of a variability of $n+1$ dimensions out of a variability of $n$ dimensions and a variability of one dimension} (English translation by Clifford).

\medskip

Poincar\'{e}\index{Poincar\'{e}, Henri (1854--1912)} attempted to give a definition of manifold for his own motivation. In the ground-breaking work {\it Analysis Situs} (J.~\!\'{E}cole Polytech, {\bf 1} (1895) 1--121), he offered two definitions, both of which rely on ``constructive procedures." The first is to describe a manifold as the zero set $f^{-1}(0)$  of a smooth function $f:U\longrightarrow \mathbb{R}^k$, where $U$ is an open set of $\mathbb{R}^{d+k}$, and the Jacobian of $f$ is of maximal rank everywhere. In the second, a manifold is produced by a ``patchwork" of a family of open sets of $\mathbb{R}^d$. Although his definitions are considered a precursor to the modern formulation of manifolds, it was not yet satisfactory enough (\cite{sch} and \cite{ohshika}). To remove the ambiguity in the early formulation and to study a global aspect of manifolds, it was necessary to establish the notion of {\it topological space}\index{topological space}, an ultimate  concept of generalized space allowing us to talk about ``finite and infinite" in an entirely intrinsic\index{intrinsic} manner (thus one can say that topological spaces\index{topological space} are ``stark naked" models of the universe).

What should be accentuated in this development, is the invention of ({\it transfinite}) {\it set theory} by Georg Cantor\index{Cantor, Georg (1845--1918)} (1845--1918).
His {\it Grundlagen einer allgemeinen Mannigfaltigkeitslehre} 
(1883) was originally motivated by the work of 
Dirichlet\index{Dirichlet, Johann Peter Gustav Lejeune (1805--1859)} 
and Riemann\index{Riemann, Georg Friedrich Bernhard (1826--1866)} concerning trigonometric series\index{trigonometric series} (Remark \ref{rem:function}),\footnote{\label{fn:realnumber}In {\it \"{U}ber die Ausdehnung eines Satzes aus der Theorie der trigonometrischen Reihen}, Math.~\!Ann., {\bf 5} (1872), 123--132, he defined real numbers, using Cauchy\index{Cauchy, Augustin-Louis (1789--1857)} sequences of rationals.} and made a major paradigmatic shift in the sense that since then mathematicians have accepted publicly ``the actual infinity"\index{actual infinity}, 
as opposed to the traditional attitude towards infinity.\footnote{Bernhard Bolzano\index{Bolzano, Bernhard (1781--1848)} (1781--1848) conceived the notion of set from a philosophical view ({\it Wissenschaftslehre}, 1837).} To grasp the atmosphere before set theory came in, we shall quote Gauss's\index{Gauss, Johann Carl Friedrich (1777--1855)} words. 
In a letter addressed to Schumacher (July 12, 1831), he says, ``I protest against the use of infinite magnitude as something completed, which is never permissible in mathematics. Infinity is merely a way of speaking, the true meaning being a limit which certain ratios approach indefinitely close, while others are permitted to increase without restriction" ({\it Werke}, VIII, p.~\!216). What he claimed is not altogether inapposite as far as ``differentiation" is concerned, but cannot apply to infinity embodied by set theory.

A set (Mengenlehre) is, as Cantor's\index{Cantor, Georg (1845--1918)} declares, a collection of definite, well-distinguishable objects of our intuition or of our thought to be discerned as a whole.\footnote{{\it Beitr\"{a}ge zur Begr\"{u}ndung der transfiniten Mengenlehre}, Math.~\!Ann., {\bf 46} (1895), 481-512. This is Cantor's last major publication.} To tell how set theory is linked to the issue of infinity, let us take a look at natural numbers. We usually recognize natural numbers $1,2,3,\ldots$ as potential infinity\index{potential infinity}; that is, we identify them by a non-terminating process such as adding 1 to the previous number, while Cantor\index{Cantor, Georg (1845--1918)} thinks that one may capture the ``set" of all natural numbers as a completed totality (which we shall denote by $\mathbb{N}$). Once the notion of set as the actual infinity\index{actual infinity} is accepted, one can exploit one-to-one correspondences (OTOC)\index{one-to-one correspondence} to compare two infinite sets as Archimedes\index{Archimedes (c. 287--c. 212 BC)}, Bradwardine\index{Bradwardine, Thomas (c. 1290--c. 1349)}, and Oresme\index{Oresme, Nicole (1320--1382)} did (Sect.~\!\ref{sec:whatinfinity}). In this context, infiniteness of a set is characterized by the property that it has a OTOC\index{one-to-one correspondence} with a proper subset. This characterization (due to Dedekind\index{Dedekind, Richard (1831--1916)}) resolves Galileo's paradox\index{Galileo's paradox} mentioned in Sect.~\!\ref{sec:whatinfinity}, simply by saying that the set of natural numbers is infinite.

The most significant discovery by Cantor\index{Cantor, Georg (1845--1918)} is the existence of different sizes of infinite sets; for instance, the real numbers are more numerous than the natural numbers. Thus $\mathbb{R}$, the set of real numbers, cannot be captured as potential infinity\index{potential infinity} whatsoever. This was mentioned in a letter\footnote{{\it Briefwechsel Cantor-Dedekind}, eds. E. Noether and J. Cavailles, Paris: Hermann.} to Dedekind\index{Dedekind, Richard (1831--1916)} dated December 7, 1873 and published in 1874.
\footnote{
%Cantor, {\it \"{U}ber eine Eigenschaft des Inbegriffes aller reelen algebraischen Zahlen}, J.~\!f\"{u}r die Reine und Angew.~\!Math., {\bf 77} (1874), 123--32. 
His {\it diagonal argument} was given in {\it \"{U}ber eine elementare Frage der Mannigfaltigkeitslehre}. Jahresbericht der Deutschen Mathematiker-Vereinigung, {\bf 1} (1890--1891), 75--78.} Furthermore, he introduced {\it cardinal numbers} as a generalization of natural numbers; that is, a cardinal number is defined to be a name assigned to an arbitrary set, where  two sets have the same name if and only if there is a OTOC\index{one-to-one correspondence} between them.

At any rate, the notion of set {\it per se} is so simple that it seem not necessary to take it up expressly. Indeed, ``set" is nearly an everyday term with synonyms such as ``collection," ``family," ``class," and ``aggregate." In truth, people before Cantor\index{Cantor, Georg (1845--1918)} more or less adapted themselves to ``sets" without noticing that the notion cannot be discussed separately from actual infinity\index{actual infinity}.

Greek thinkers were timid with infinity since it often leads to falsity (remember the erroneous argument by Antiphon\index{Antiphon (c. 480--c. 411 BC)} and Bryson\index{Bryson (born in c. 450 BC)} in ``squaring the circle"\index{squaring the circle}).
Cantor\index{Cantor, Georg (1845--1918)}, too, struck a snag 
when he put forward his theory, if not the same snag the Greeks struck. In fact, the simplicity of the notion of set is deceptive. Unless the concept of set is suitably introduced, we run into a contradiction like Russell's paradox\index{Russell's paradox} discovered in 1901, which had disquieted, for some time, mathematicians of the day who were favorably disposed towards set theory.\footnote{His paradox\index{Russell's paradox}
%---independently found by E. Zermelo---
 results from admitting that $X=\{x|~\!x\not\in x\}$ is a set. In fact, if $X\not\in X$, then $X\in X$ by definition, and if $X\in X$, then $X\not\in X$ again by definition.}

Surely, Cantor's\index{Cantor, Georg (1845--1918)} innovation was unorthodox and daring at that time; he met with belligerent resistance from a few conservative people. For instance, Poincar\'{e}\index{Poincar\'{e}, Henri (1854--1912)} poignantly said that there is no actual infinity\index{actual infinity} since mathematics accepting it has fallen into contradiction (1906),
%\footnote{{\it Les math\'{e}matiques et la logique III}, Rev.~\!M\'{e}taphys.~\!Morale, {\bf 14} (1906), p. 316.} 
and L. J. J. Wittgenstein bluntly dismissed set theory as ``utter nonsense," and ``wrong."\footnote{See 
%V. Rodych (2007), 
``Wittgenstein's Philosophy of Mathematics" in the Stanford Encyclopedia of Philosophy. The frustration and despair caused by the disinterest and cold rejection directed Cantor\index{Cantor, Georg (1845--1918)} to a theological explanation of his theory, which, he believed, would open up a whole new landscape in Christian theology (a letter to C. Hermite, dated January 22, 1894 (\cite{Mesch}).}

Fortunately, various paradoxes derived from a vague understanding were resolved by axiomatizing set theory,\footnote{Axiomatic set theory, established by E. Zermelo in 1908 and A. Fraenkel in 1922, is formulated with a formal logic.} and most mathematicians responded favorably to his theory. Hilbert\index{Hilbert, David (1862--1943)}, who stands squarely against the outcry of conservative people, exclaims, ``No one shall expel us from the Paradise that Cantor\index{Cantor, Georg (1845--1918)} has created" ({\it \"{U}ber das Unendliche}, Math.~\!Ann., {\bf 95} (1926), 161--190).\footnote{In the late 1920s, Hilbert\index{Hilbert, David (1862--1943)} preferred to compare the use of actual infinity\index{actual infinity} to the addition of points at infinity\index{point at infinity} in projective geometry\index{projective geometry}.} 

\section{Topological spaces}
The swift change caused by set theory delimited the history of mathematics. The old thought that mathematics is the science of quantity, or of space and number, has largely disappeared. From then on, mathematics has been built eventually on set theory, indissolubly combined with the concept of {\it mathematical structure}\index{mathematical structure}---in short, any set of objects along with certain relations among those objects---which was put forward in a definitive form by Bourbaki in the campaign for modernization of mathematics. 
An instructive example is $\mathbb{N}$ whose structure is characterized by {\it Peano's axioms}\index{Peano, Giuseppe (1858--1932)}, the setting that epitomizes the transition from potential infinity\index{potential infinity} to actual infinity\index{actual infinity} (1889):
%\footnote{{\it Arithmetices principia, nova methodo exposita}, 1889.}
\smallskip

``The set $\mathbb{N}$ has an element $1$ and an injective map $\varphi:\mathbb{N} \longrightarrow\mathbb{N}$ such that (i) $1\not\in\varphi(\mathbb{N})$,~~(ii) if $S$ is a subset of $\mathbb{N}$ with $1\in$ $S$ and $\varphi(S)\subset S$, then $S=\mathbb{N}$."

\smallskip

Algebraic systems such as {\it groups}, {\it rings}, and {\it vector spaces} are sets with structures as well. Topological spaces are another exemplary, which has, in no small part, evolved out of a long process of understanding the meaning of ``limit" and ``continuity" of functions. 

\begin{rem}\label{rem:function}{\rm {\small 
A primitive form of functions is glimpsed in the {\it Almagest}\index{almagest@Almagest}, Book 1, Chap.~\!11. He made a table of chords of a circle (dating back to Hipparchus\index{Hipparchus (c. 190--c. 120 BC)}),  which, from a modern view, can be thought of as associating the elements of one set of numbers with the elements of another set. After Oresme\index{Oresme, Nicole (1320--1382)} (1350), Galileo\index{Galileo Galilei (1564--1642)} (1638), and Descartes\index{Descartes, Ren\'{e} (1596--1650)} (1637) adumbrated germinal ideas, the great duo of Johann Bernoulli\index{Bernoulli, Johann (1667--1748)} and Leibniz\index{Leibniz, Gottfried Wilhelm von (1646--1716)} adopted the word ``functio" for  ``a quantity formed from indeterminate and constant quantities" (1694). A more lucid formulation was given by Euler\index{Euler, Leonhard (1707--1783)}, who introduced the notation $f(x)$ (1734) and defined a function to be an {\it analytic expression} composed in any way whatsoever of the variable quantity and numbers or constant quantities" ({\it Introductio};  \cite{janke}). But, after a while, he altered this definition as ``When certain quantities depend on others in such a way that they undergo a change when the latter change, then the first are called functions of the second" (\cite{euler}).

Interest in the true nature of functions was renewed when Fourier\index{Fourier, Jean Baptiste Joseph (1768--1830)} claimed that an {\it arbitrary} function $f(x)$ on $[-\pi,\pi]$ can be expressed as a trigonometric series\index{trigonometric series} 
\begin{eqnarray}
&&f(x)=\frac{1}{2}b_0+(a_1\sin x+b_1\cos x)+(a_2\sin 2x+b_2\cos 2x)+\cdots,\label{eq:trigono}\\
&& \qquad a_n=\frac{1}{\pi}\int_{-\pi}^{\pi}f(x)\sin nx~\!dx,\quad b_n=\frac{1}{\pi}\int_{-\pi}^{\pi}f(x)\cos nx~\!dx.\nonumber
\end{eqnarray}
However, he did not turn his eye to the convergence of the series, let alone the meaning of the equality, the unavoidable issues as pointed out by Cauchy\index{Cauchy, Augustin-Louis (1789--1857)} (1820).\footnote{The bud of the issue is seen in the debate among D'Alembert\index{d'Alembert, Jean-Baptiste le Rond (1717--1783)}, Euler\index{Euler, Leonhard (1707--1783)}, and Daniel Bernoulli\index{Bernoulli, Daniel (1700--1782)} about solutions of the wave equation (Remark \ref{rem:spacialrelativity}~\!(1)).} It was in this context that Cauchy\index{Cauchy, Augustin-Louis (1789--1857)} gave a faultless definition of continuous function, using the notion of ``limit" for the first time.\footnote{Euler\index{Euler, Leonhard (1707--1783)} thought of ``$\lim_{n\to \infty}$" (resp. ``$\lim_{x\to 0}$") as ``the value for $n$ infinitely large (resp. the value for $x$ infinitely small"). The modern definition was furnished by Bolzano\index{Bolzano, Bernhard (1781--1848)} in his {\it Der binomische Lehrsatz}, 1816; however, being of little noticed at the time.} Following Cauchy's\index{Cauchy, Augustin-Louis (1789--1857)} idea, Weierstrass\index{Weierstrass, Karl Theodor Wilhelm (1815--1897)} popularized the $\epsilon$-$\delta$ argument\index{epsilon@$\epsilon$-$\delta$ argument} in the 1870's, which made it possible to discuss diverse aspects of convergence. In the meantime, Dirichlet\index{Dirichlet, Johann Peter Gustav Lejeune (1805--1859)} made Fourier's work\index{Fourier, Jean Baptiste Joseph (1768--1830)} rigorous on a clearer understanding of (dis)continuity. With some prodding from Dirichlet,  Riemann\index{Riemann, Georg Friedrich Bernhard (1826--1866)} made further progress by giving a precise meaning to integrability of function.\footnote{\label{fn:firststage}Riemann, {\it \"{U}ber die Darstellbarkeit einer Function durch eine trigonometrische Reihe}, Abhandlungen der K\"{o}niglichen Gesellschaft der Wissenschaften zu G\"{o}ttingen, {\bf 13} (1868), 87--132. This was presented in December of 1853 as his Habilitationsschrift at the first stage.}

The rigor of calculus whose absence was deplored by Abel\index{Abel, Niels Henrik (1802--1829)} was now in place (Remark \ref{rem:analyticalsoc} (4)). Analysis that developed on it (and was accompanied with Cantor's\index{Cantor, Georg (1845--1918)} set theory) contributed to the advent of the concept of topological space\index{topological space}.
\hfill$\Box$

}
}
\end{rem}

A {\it topology}\index{topology} on a set $X$ is a family of subsets (called {\it neighborhoods}) which makes it possible to express the idea of getting closer and closer to a point and also to introduce the notion of continuous maps as a generalization of continuous functions. In his magnum opus {\it Grundz\"{u}ge der Mengenlehre} (1914), Felix Hausdorff\index{Hausdorff, Felix (1869--1942)} (1869--1942) sets the four axioms in terms of neighborhoods, the last being known as the {\it Hausdorff separation axiom}; thus his topological spaces\index{topological space} are what we now call {\it Hausdorff spaces}\index{Hausdorff space}, a proper generalization of {\it metric spaces}.\footnote{The earlier results on topology\index{topology}---ascribed to F. Riesz, M. R. Fr\'{e}chet, and others---fitted naturally into the framework set up by Hausdorff\index{Hausdorff, Felix (1869--1942}. Currently, there are several ways to define topological spaces; e.g. by the axioms on the closure operation (1922) due to K. Kuratowski, and by the axioms on open sets (1934) due to W. F. Sierpi\'{n}ski. The notion of metric space was first suggested by Fr\'{e}chet in his dissertation paper (1906).}

Now, a manifold with the first fundamental form\index{first fundamental form} is a metric space\index{metric space}, so we are tempted to discuss {\it finiteness} of the universe in the topological framework. Finiteness of a figure in $\mathbb{R}^N$ is interpreted by the terms {\it closed} and {\it bounded}. Thus the issue boils down to the unequivocal question of how we define such terms in an intrinsic\index{intrinsic} manner, without referring to ``outside." The answer is afforded by the notion of {\it compactness} (Fr\'{e}chet, 1906), which expresses ``closed and bounded in itself" so to speak, whose origin is the Bolzano-Weierstrass theorem\index{Bolzano, Bernhard (1781--1848)}\index{Weierstrass, Karl Theodor Wilhelm (1815--1897)} and the Heine-Borel theorem in analysis.

Armed with the impregnable language of finiteness, one may call into question, without assessing relevance to metaphysical or theological meaning, whether our universe is finite. Yet, this overwhelming question is beyond the scope of this essay and belongs entirely to astronomy. The recent observation (by the {\it Planck}, the space observatory operated by the European Space Agency) lends support to an ``almost" flat model of the universe in the large. However this does not entail its infiniteness since ``flatness" does not imply ``non-compactness." In turn, if there are only finitely many particles in the whole universe, then the universe should be finite (here it is assumed that the mass would be uniformly distributed in the large).\footnote{In the {\it Letter to Herodotus}, Epicurus (341 BCE--270 BCE), an adherent of atomism and advocate of ``billiard ball universe," states that the number of atoms is infinite.} This seems not to be a foolhardy hypothesis altogether. The {\it cosmological principle} in modern cosmology claims that the distribution of matter in the universe is homogeneous\index{homogeneous} and isotropic\index{isotropic} when viewed on a large enough scale. Einstein\index{Einstein, Albert (1879--1955)} says, ``I must not fail to mention that a theoretical argument can be adduced in favor of the hypothesis of a finite universe."\footnote{{\it Geometry and Experience}, an address to the Berlin Academy on January 27th, 1921.} Truth be told, however, nobody could definitively pronounce on the issue at the moment. On this point, our intellectual adventure cannot finish without frustration since the issue of finiteness is the most fundamental in cosmology.

\section{Towards modern differential geometry}\label{sect:reimenniangeometry}

After the notion of topological space\index{topological space} took hold, mathematicians got ready to move forward. What had to be done initially was to link the concept of curvilinear coordinate system with topological spaces\index{topological space}. An incipient attempt was made by Hilbert\index{Hilbert, David (1862--1943)}, who, using a system of neighborhoods, tried to characterize the plane.\footnote{Appendix IV (1902)
to {\it Grundlagen der Geometrie}.} 
Subsequently, Hilbert's\index{Hilbert, David (1862--1943)} former student Hermann Klaus Hugo Weyl (1885--1955)\index{Weyl, Hermann Klaus Hugo (1885--1955)} published  
{\it Die Idee der Riemannschen Fl\"{a}che} 
(1913), a classical treatise that
lays a solid foundation for complex analysis initiated by Riemann\index{Riemann, Georg Friedrich Bernhard (1826--1866)} (fn.~\!\ref{fn:Abelcschen}). 
His perspicuous exposition, though restricted to the 2D case, opened up the modern synoptic view of geometry and analysis on manifolds.
Stimulated with this, 
O. Veblen and his student J. H. C. Whitehead gave the first general 
definition of 
manifolds in their book {\it The foundations of differential geometry}, 1933. At around the same time, Hassler Whitney\index{Whitney, Hassler (1907--1989)} (1907--1989) and others clarified the foundational aspects of manifolds during the 1930s. One of significant outcomes 
is that any smooth manifold (with an additional mild property) can be {\it embedded} in a higher-dimensional Euclidean space\index{Euclidean space} (1936), thus giving a justification to one of Poincar\'{e}'s\index{Poincar\'{e}, Henri (1854--1912)} definitions of manifold (Sect.~\!\ref{subsec:transfinite}).
%\footnote{Whitney, {\it Differentiable manifolds}, Ann.~\!of Math., {\bf 37} (1936), 645--680.} 

Meanwhile, differential geometry of surfaces was expanded by Jean-Gaston Darboux\index{Darboux, Jean-Gaston (1842--1917)} (1842--1917).
He
published, between 1887 and 1896, four enormous volumes entitled {\it Le\c cons sur la th\'{e}orie g\'{e}n\'{e}rale des surfaces et les applications g\'{e}om\'{e}triques du calcul infinit\'{e}simal}, including most of his earlier work and touching on global aspect of surfaces to some extent. Darboux's\index{Darboux, Jean-Gaston (1842--1917)} spirit was then inherited by \'{E}lie Cartan\index{Cartan, \'{E}lie (1869--1951)} (1869--1951).
Most influential is his meticulous study
of {\it symmetric space} (1926), a class of 
manifolds with an analogue of point symmetry at each point,
%\footnote{{\it Sur une classe remarquable d'espaces de Riemann, I, II}, Bulletin de la Soci\'{e}t\'{e} Math\'{e}matique de France, {\bf 54} (1926), 214--216; {\bf 55} (1927), 114--134.} 
including $\mathbb{R}^d$, $H^d$,  $S^d$, and $P^d(\mathbb{R})$; 
thus unifying various geometries in view of symmetry\index{symmetry} and also following the Greek tradition in seeking symmetric shapes.
Remarkable is that Cartan's\index{Cartan, \'{E}lie (1869--1951)} list includes symmetric spaces described in terms of quaternions\index{quaternion} and octonions\index{octonion}.
Futhermore Cartan\index{Cartan, \'{E}lie (1869--1951)} enlarged Klein's\index{Klein, Christian Felix (1849--1925)} ``Erlangen Program"\index{Erlangen Program} (fn.~\!\ref{fn:erlangen}) so as to encompass 
general geometries.
%\footnote{{\it Espaces \`{a} connexion affine, projective et conforme}, Acta Math., {\bf 48} (1926), 1--42.} 

A decade from 1930 was the period when geometry began to be intensively studied from a global point of view, in concord with a new field launched by Riemann\index{Riemann, Georg Friedrich Bernhard (1826--1866)} and Poincar\'{e}\index{Poincar\'{e}, Henri (1854--1912)} (see Sect.~\!\ref{subsec:topodes}). What is more, the fundamental concepts have been made more transparent by means of a {\it coordinate-free setup}, which is, in a sense, conform to the principle of relativity in a perfect way.\footnote{Grassmann\index{Grassmann, Hermann (1809--1877)} was the first to realize the importance of the coordinate-free concepts.}

\medskip

In what follows, just to make our story complete, we shall quickly review, at the cost of some repetition,  how some of the ingredients introduced by forerunners are reformulated in modern terms. The first is the definition of smooth manifolds 
as a {\it terminus ad quem}
of our journey to ``generalized curved spaces."

\smallskip

\noindent{\bf Definition} (1) 
A Hausdorff space\index{Hausdorff space} $M$ is said to be a {\it $d$-dimensional topological manifold} (or $d$-{\it manifold}) if 
each point of $M$ has an open neighborhood homeomorphic to an open set in $\mathbb{R}^d$; therefore we have an {\it atlas of local charts}; that is, a family of curvilinear coordinate systems that covers $M$.

\smallskip

(2) 
A topological manifold $M$ is called smooth if there is an atlas 
such that every coordinate transformation is smooth (such an atlas is said to be smooth). 

\smallskip

Using a smooth atlas, we may discuss ``smoothness" of various objects attached to smooth manifolds; say, {\it smooth functions} and {\it smooth maps}.
What comes next to mind is the question whether it is possible to define {\it tangent space}\index{tangent space} as the generalization of tangent planes\index{tangent plane} of a surface. 
The idea to conceptualize tangent space without reference to a coordinate system and ambient spaces dates back to the work of C. Bourlet who, paying special attention to the product rule, gave an algebraic characterization of differentiation.\footnote{{\it Sur les op\'{e}rations en g\'{e}n\'{e}ral et les \'{e}quations diff\'{e}rentielles lin\'{e}aires d'ordre infini}, Ann.~\!Ec.~\!Normale,
{\bf 14} (1897), 133--190.} 

To explain the idea in depth, we shall first observe that a tangent plane\index{tangent plane} of a surface can be intrinsically defined.
Let $T_pS$ be the space of vectors tangent to a surface $S$ at $p$, and let $C^{\infty}_p(S)$ be the set of smooth functions defined on a neighborhood of $p$.
Define the action of $\boldsymbol{\xi}
\in T_pS$ on $C^{\infty}_p(S)$ by setting 
$
\boldsymbol{\xi}(f)$ $=\frac{d}{dt}f\big(c(t)\big)\big|_{t=0}
$, 
where $c:(-\epsilon,\epsilon)\longrightarrow S$ is a curve with $c(0)=p$ and $\dot{c}(0)=\boldsymbol{\xi}$. 
The correspondence $f\in C^{\infty}_{p}(S)\mapsto \boldsymbol{\xi}(f)\in \mathbb{R}$ has the following properties:

\medskip

{\rm (1) (Linearity)}\quad $\boldsymbol{\xi}(af+bg)=a\boldsymbol{\xi}(f)+b\boldsymbol{\xi}(g)$ $~~a,b\in \mathbb{R},~f,g\in C^{\infty}_{p}(S)$.

\smallskip

{\rm (2) (Product rule)}\quad $\boldsymbol{\xi}(fg)=f(p)\boldsymbol{\xi}(g)+g(p)\boldsymbol{\xi}(f)$.

\smallskip

{\rm (3)} $(a\boldsymbol{\xi}+b\boldsymbol{\eta})(f)=a\boldsymbol{\xi}(f)+b\boldsymbol{\eta}(f)$ $~~a,b\in \mathbb{R},~\boldsymbol{\xi},\boldsymbol{\eta}\in T_{p}S$, $~f\in C^{\infty}_{p}(S)$.

\smallskip

{\rm (4)} If $\boldsymbol{\xi}(f)=0$ for any $f\in C^{\infty}_{p}(S)$, then $\boldsymbol{\xi}=\boldsymbol{0}$. 

\smallskip
We let $\tau_{p}S$ be the vector space consisting of maps $\boldsymbol{\omega}:C^{\infty}_{p}(S)\rightarrow \mathbb{R}$ with the properties $\boldsymbol{\omega}(af +bg)=a\boldsymbol{\omega}(f)+b\boldsymbol{\omega}(g)$ and $\boldsymbol{\omega}(fg)=f(p)\boldsymbol{\omega}(g)+g(p)\boldsymbol{\omega}(f)$. 
In view of (1), (2), $\boldsymbol{\xi}$ as a map of $C^{\infty}_{p}(S)$ into $\mathbb{R}$
belongs to $\tau_pS$. Hence we obtain a map $\iota$ of $T_pS$ into $\tau_pS$, which turns out to be a linear isomorphism. Identifying $T_pS$ with $\tau_pS$ via $\iota$, we have an intrinsic\index{intrinsic} description of $T_pS$.

The foregoing discussion suggests how to define the tangent spaces\index{tangent space} 
of a smooth manifold $M$.
Let $C^{\infty}_p(M)$ be the set of smooth functions defined on neighborhoods of $p\!\in\!M$, and define $T_pM$ to be the vector space consisting of operations $\omega\!:\!C^{\infty}_p(M)\rightarrow \mathbb{R}$ satisfying linearity and the product rule.\footnote{This definition is given in Claude Chevalley's {\it Theory of Lie Groups} (1946).} 
For a curvilinear coordinate system $(x_1,\ldots,x_d)$ around $p$, partial differentiations 
$\partial_i|_p\!:=\!\displaystyle \partial/\partial x_i\big|_p$~($i\!=\!1,\ldots,d$) form a basis of $T_pM$. Moreover, the {\it velocity vector} $\dot{c}(t)\in T_{c(t)}M$ for a curve $c$ in $M$ is defined by putting $\dot{c}(t)f=\frac{d}{dt}f(c(t))$.

The notion of tangent space\index{tangent space} allows us to import various concepts in calculus into the theory of manifolds. An example is 
the {\it differential} $\varPhi_{*p}:T_pM\longrightarrow T_{\varPhi(p)}N$ of a smooth map $\varPhi:M\longrightarrow N$ defined by
$
\varPhi_{*p}(\boldsymbol{\xi})(f)= \boldsymbol{\xi}(f\circ\varPhi),~\boldsymbol{\xi}\in T_pM,~f\in C^{\infty}_{\varPhi(p)}(N)
$, 
which is an intrinsic\index{intrinsic} version of the {\it Jacobian matrix}.

A {\it Riemannian manifold}\index{Riemannian manifold} is a manifold with a ({\it Riemannian}) {\it metric}\index{Riemannian metric}, a smooth family of inner products $\langle \cdot,\cdot\rangle_p$ on $T_pM$~($p\in M$). The relation between the metric and 
the first fundamental form\index{first fundamental form}
is given by $g_{ij}(p)=\langle \partial_i|_p,\partial_j|_p\rangle_p$. 
A smooth map $\varPhi:M\longrightarrow N$ between Riemannian manifolds is called an {\it isometry}\index{isometry} if 
$\langle \varPhi_{*p}(\boldsymbol{\xi}),\varPhi_{*p}(\boldsymbol{\eta})\rangle_{\varPhi(p)}=
\langle\boldsymbol{\xi},\boldsymbol{\eta} \rangle_p~~(\boldsymbol{\xi},\boldsymbol{\eta}\in T_pM)$.
The differential $\varPhi_{*p}$ of an isometry $\varPhi$ is obviously injective for every $p$. Conversely, for a smooth map $\varPhi$ of $M$ into a Riemannian manifold $N$ 
such that $\varPhi_{*p}$ is injective for every $p$, one can equip a metric on $M$ (called the {\it induced metric}) which makes $\varPhi$ an isometry.
The first fundamental form\index{first fundamental form} 
on a surface $S$ is nothing but the metric induced by the inclusion map of $S$ into $\mathbb{R}^3$.\footnote{In 1956, J. F. Nash showed that every Riemannian manifold has an isometric inclusion into some Euclidean space\index{Euclidean space} $\mathbb{R}^N$.} 
Furthermore, the metric (\ref{eq:constcurve}) given by Riemann\index{Riemann, Georg Friedrich Bernhard (1826--1866)} coincides essentially with the standard metric on the sphere $S^{d}(R)= \{(x_1,\ldots,x_{d+1})|~\!x_1{}^2+\cdots+x_{d+1}{}^2=R^2\}$ 
that is read in the curvilinear coordinate system derived by the stereographic projection\index{stereographic projection} $\varphi:S^d(R)\backslash\{(0,\ldots,0,R)\}\longrightarrow\mathbb{R}^d$: 
$
\varphi(x_1,\ldots,x_{d+1})=\big(Rx_1/(R-x_{d+1}),\ldots,Rx_d/(R-x_{d+1})
\big)$.

No less important than tangent spaces\index{tangent space} is 
the dual space $T^*_pM$ of $T_pM$, which is distinguished from $T_pM$ though 
they are isomorphic as vector spaces.
\footnote{The distinction between a vector space and its dual was definitively established following the work of S. Banach and his school, though it had been already reflected for a long time by the distinction between {\it covariant} and {\it contravariant}, and also between {\it cogredient} and {\it contragredient}, the concepts in group-representation theory.} 
We denote by $\{dx_1|_p,\ldots,dx_d|_p\}$ the {\it dual basis} of $\{\partial_1|_p,\ldots,\partial_d|_p\}
$; i.e., $dx_i|_p(\partial_j|_p)=\delta_{ij}$. The symbol $dx_i|_p$ 
defined in this way 
justifies the notion of differential 
(Remark \ref{rem:infinixx}~\!(2)). 
A {\it tensor of type}\index{tensor} $(h,k)$ 
is then a multilinear functional 
$
T:\underbrace{T_p^*M\times\cdots\times T_p^*M}_h\times \underbrace{T_pM\times \cdots\times T_pM}_k\longrightarrow \mathbb{R}
$.
which is related to a classical tensor by $T^{i_1\cdots i_h}_{j_1\cdots j_k}(p):=T(dx_{i_1}|_p,\ldots,dx_{i_h}|_p,\partial_{j_1}|_p,$ $\ldots,\partial_{j_k}|_p)
$.
A {\it tensor field of type}\index{tensor field} $(h,k)$ is a smooth assignment to each $p\in M$ of a tensor $T(p)$ of type $(h,k)$. 
A {\it vector field} is a tensor field of type $(1,0)$. The vector space consisting of all vector fields on $M$ is denoted by $\mathfrak{X}(M)$.

Essential to global analysis are {\it differential forms}\index{differential form} and their {\it exterior derivative} introduced by Poincar\'{e}\index{Poincar\'{e}, Henri (1854--1912)} (1886) 
%(without much definition) 
and Cartan\index{Cartan, \'{E}lie (1869--1951)} (1899).
%\footnote{Poincar\'{e}, {\it Sur les residus des integrales doubles}, Acta Math., {\bf 9} (1886/7) 321--380. Cartan\index{Cartan, \'{E}lie (1869--1951)}, {\it Sur certaines expressions diff\'{e}rentielles et le probl\`{e}me de Pfaff}, Ann.~\!Ec.~\!Normale, {\bf 16} (1899), 239--332.} 
To give a current definition of differential forms\index{differential form}, we denote by $\wedge^kT_p^*M$ the totality of tensors $\omega$ of type $(0,k)$ satisfying $\omega(\boldsymbol{\xi}_{\sigma(1)},\ldots,$ $\boldsymbol{\xi}_{\sigma(k)})={\rm sgn}(\sigma)\omega(\boldsymbol{\xi}_{1},\ldots,\boldsymbol{\xi}_{k})$ $~(\boldsymbol{\xi}_i\in T_pM)$ for every permutation $\sigma$ of the set $\{1,\ldots,k\}$. 
A smooth assignment to each $p\!\in\!M$ of $\omega(p)\!\in\!\wedge^kT_p^*M$ is called a ({\it differential}) {\it $k$-form}\index{differential form}, all of which constitute a vector space, denoted by $A^k(M)$. A $k$-form is informally expressed as
$$
\omega=
\displaystyle\sum_{1\leq i_1<\cdots<i_k\leq d}f_{i_1\cdots i_k}(x_1,\ldots,x_d)dx_{i_1}\wedge\cdots \wedge dx_{i_k}, 
$$
where the ``wedge product" $\wedge$\index{wedge product} is supposed to satisfy  $(dx_h\wedge dx_i)\wedge dx_j=dx_h\wedge (dx_i\wedge dx_j),$ $~dx_i\wedge dx_j=-dx_j\wedge dx_i$.\footnote{The totality of $\wedge^kT_p^*M$ $(k=0,1,\ldots,d)$
, equipped with the wedge product, constitutes the {\it exterior algebra}\index{exterior algebra} introduced by Grassmann\index{Grassmann, Hermann (1809--1877)} (Sect.~\!\ref{sec:dimension}).} With this notation, the {\it exterior derivation}\index{exterior derivation} $d_k:A^k(M)\longrightarrow A^{k+1}(M)$ is defined by
$$
d_k\omega(=d\omega):=\sum_{i=1}^d\sum_{1\leq i_1<\cdots<i_k\leq d}\frac{\partial f_{i_1\cdots i_k}}{\partial x_i}dx_i\wedge dx_{i_1}\wedge\cdots \wedge dx_{i_k}.
$$
Note here that $df=f_{x_1}dx_1+\cdots+f_{x_d}dx_d$ for a smooth function $f\in A^0(M)$. A significant fact derived from the equality $f_{xy}=f_{yx}$ is the identity $d_k\circ d_{k-1}=0$. Thus we obtain the vector space $H^k_{\rm dR}(M)$ $:={\rm Ker}~\!d_k/{\rm Image}~\!d_{k-1}$, which is called the {\it de Rham cohomology group}; Remark \ref{rm:form} (3).

Differential forms give an incentive to justify the formal expression $dx_1\cdots dx_d$ in multiple integrals, because the basic rules for the wedge product make the appearance of {\it Jacobian determinants} automatic: $dy_1\wedge\cdots\wedge dy_d=\frac{\partial(y_1,\ldots,y_d)}{\partial(x_1,\ldots,x_d)}dx_1\wedge\cdots \wedge dx_d$, where $y_i=y_i(x_1,\ldots, y_d)$ is a coordinate transformation. This is in conformity with a change of variables in multiple integrals\!:\footnote{This formula was first proposed by Euler\index{Euler, Leonhard (1707--1783)} for double integrals (1769), then generalized to triple ones by Lagrange\index{Lagrange, Joseph-Louis (1736--1813)} (1773). Ostrogradski extended it to general multiple integrals (1836).} $dy_1\cdots dy_d=\Big|\frac{\partial(y_1,\ldots,y_d)}{\partial(x_1,\ldots,x_d)}\Big|dx_1\cdots dx_d$. Thus we are motivated to make the following definition\!: ``$M$ is said to be {\it orientable}\index{orientable} if it has a smooth atlas  such that the Jacobian determinant of every coordinate transformation is positive" (Poincar\'{e}\index{Poincar\'{e}, Henri (1854--1912)}, {\it Analysis Situs}, 1895; see Sect.~\!\ref{sec:right}). Then, what springs to mind is to define the integral over an orientable\index{orientable} $d$-manifold $M$ of a $d$-form $\omega=f(x_1,\ldots,x_d)~\!dx_1\wedge\cdots \wedge dx_d$ by setting $\int_M\omega=\int f(x_1,\ldots,x_d)~\!dx_1\cdots dx_d$.

Now, given $\omega\in A^k(N)$ and a map $\varPhi$ of an orientable\index{orientable} $k$-manifold $M$ (possibly with boundary) into $N$, we define the integral of $\omega$ along $\varPhi$ by $\int_{\varPhi} \omega:=\int_M\varPhi^*\omega$, where $(\varPhi^*\omega)(p)\in \wedge^kT_p^*M$ is defined by setting $(\varPhi^*\omega)(\boldsymbol{\xi}_{1},\ldots,\boldsymbol{\xi}_{k})=\omega(\varPhi_*(\boldsymbol{\xi}_{1}),\ldots,\varPhi_*(\boldsymbol{\xi}_{k}))$. This is a generalization of line integrals, with which we may unify Stokes' formula\index{Stokes' formula}, Gauss's\index{Gauss, Johann Carl Friedrich (1777--1855)} divergence theorem\index{divergence theorem}, and Green's\index{Green, George (1793--1841)} theorem\index{Green's theorem}, as the {\it generalized Stokes' formula}\!:
$$
\int_Md\omega=\int_{\partial M}\omega\left(=\int_{\iota}~\omega\right),
$$
where $\omega$ is a $(d-1)$-form on an orientable\index{orientable} $d$-manifold $M$ with boundary $\partial M$ (appropriately oriented), and $\iota:\partial M\longrightarrow M$ is the inclusion map.
%\footnote{Poincar\'{e}\index{Poincar\'{e}}, {\it Les M\'{e}thodes Nouvelles de la M\'{e}caniques C\'{e}leste}, 1899.} 

The notion of vector field is quite old. In his {\it Principes g\'{e}n\'{e}raux du mouvement des fluides} (1755), Euler\index{Euler, Leonhard (1707--1783)} used it to represent a fluid's velocity. Related to vector fields is the work  {\it Theorie der Transformationsgruppen} (1888--1893) by Sophus Lie\index{Lie, Sophus (1842--1899)} (1842--1899). Under the influence of\index{Klein, Christian Felix (1849--1925)} Klein's Erlangen Program\index{Erlangen Program}, he made the study of {\it infinitesimal group actions}, which eventually evolved into the theory of {\it Lie groups}. As he observed, giving such an action is equivalent to giving a subspace $\mathfrak{g}$ of $\mathfrak{X}(M)$ satisfying ``$X,Y\in \mathfrak{g}$ $\Rightarrow$ $[X,Y]\in \mathfrak{g}$," where $[X,Y](p)f=X(p)(Yf)-Y(p)(Xf)$. Here the binary operation $[X,Y]$ (the {\it Lie bracket}) satisfies  the {\it Jacobi identity} $[[X,Y],Z]+[[Y,Z],X]+[[Z,X],Y]=0$. An algebraic system with a similar operation was to be called a {\it Lie algebra}.\footnote{Lie algebra (``infinitesimal group" in Lie's term) was independently invented by W. Killing in the 1880s with quite a different purpose.}

\medskip

 A {\it covariant differentiation}\index{covariant differentiation} acting on tensor fields is an intrinsic\index{intrinsic} generalization of directional differentiation acting on vector-valued functions, which is also designated a {\it connection}\index{connection}, a term introduced by H. Weyl\index{Weyl, Hermann Klaus Hugo (1885--1955)} (1918).
%\footnote{{\it Reine infinitesimalgeometrie}, Math.~\!Z., {\bf 2} (1918), 384--411.} 
In the case of vector fields, it is a bilinear map $\nabla:(\boldsymbol{\xi},X)\in T_pM\times \mathfrak{X}(M)\mapsto \nabla_{\boldsymbol{\xi}}X\in T_pM$ satisfying
$\nabla_{\boldsymbol{\xi}}(fX)=\boldsymbol{\xi}(f)X+f\nabla_{\boldsymbol{\xi}}X~(a,b \in \mathbb{R},~f\in C^{\infty}(M))$ (J.-L. Koszul).
Define the functions $\left\{ { {k}\atop{i j} } \right\}$ by the relation
$\nabla_{\partial_i}\left(\partial_j\right)=\sum_{k=1}^d\left\{ { {k}\atop{i j} } \right\}\partial_k$.
Then 
{\small 
$$
\nabla_{\boldsymbol{\xi}}X=\sum_{i=1}^d\xi_i\Big(
\frac{\partial f_k}{\partial x_i}+\sum_{j=1}^d\left\{ { {k}\atop{i j} } \right\}f_j
\Big)\partial_k\quad (\boldsymbol{\xi}=\sum_{i=1}^d\xi_i\partial_i,~X=\sum_{j=1}^df_j\partial_j).
$$

}

The {\it Levi-Civita connection}\index{Levi-Civita connection} is a unique connection satisfying
\begin{equation}\label{eq:1a}
\nabla_XY-\nabla_YX=[X,Y],~~\boldsymbol{\xi}\langle X,Y\rangle=\langle\nabla_{\boldsymbol{\xi}}X,Y\rangle+\langle X,\nabla_{\boldsymbol{\xi}}Y\rangle.
\end{equation}
The first equality is equivalent to $\left\{ { {k}\atop{i j} } \right\}=\left\{ { {k}\atop{j i} } \right\} ~(i,j,k=1,\ldots,d)$, while the second one tells us that $\left\{ { {k}\atop{i j} } \right\}$ coincides with the Christoffel symbol\index{Christoffel symbol} (\ref{eq:chris}).

Let $c(t)=(c_1(t),\ldots,c_d(t))$ $~(a\leq t\leq b)$ be a smooth curve in $M$. If we write $f_i(t)$ for $f_i(c(t))$ for brevity, then 
{\small 
$$
\nabla_{\dot{c}(t)}X=\sum_{i,k=1}^d\frac{dc_j}{dt}\Big(
\frac{\partial f_k}{\partial x_i}+\sum_{j=1}^d\left\{ { {k}\atop{i j} } \right\}f_j
\Big)\partial_k=\sum_{i,k=1}^d\Big(\frac{df_k}{dt}+\sum_{j=1}^d\left\{ { {k}\atop{i j} } \right\}\frac{dc_i}{dt}f_j\Big)\partial_k.
$$}
\hspace{-0.1cm}With this formula in mind, we define the covariant derivative of a vector field $X(t)=\displaystyle f_1(t)\partial_1|_{c(t)}+\cdots+f_d(t)\partial_d|_{c(t)}$ along the curve $c$ by setting 
{\small
$$
\frac{D}{dt}X=\sum_{k=1}^d\Big(\frac{df_k}{dt}+\sum_{i,j=1}^d\left\{ { {k}\atop{i j} } \right\}\frac{dc_i}{dt}f_j\Big)\partial_k.
$$}
\hspace{-0.1cm}If $DX/dt\equiv 0$, then $X$ is said to be {\it parallel}. 
Since $DX/dt\equiv 0$ is a system of linear equations of first order, there exists a parallel vector field $X$ along $c$ satisfying $X(c(a))=\boldsymbol{\xi} \in T_{c(a)}M$. The {\it parallel transport} $P_c:$ $T_{c(a)}M\longrightarrow T_{c(b)}M$ is then defined by $P_c(\boldsymbol{\xi})=X(c(b))$; thus a connection literally yields a ``bridge"\index{connection} between two tangent spaces\index{tangent space}, originally unrelated to each other. For the Levi-Civita connection\index{Levi-Civita connection}, we have $\langle P_c(\boldsymbol{\xi}),P_c(\boldsymbol{\eta})\rangle_{c(b)}$ $=\langle \boldsymbol{\xi},\boldsymbol{\eta}\rangle_{c(a)}$ in view of (\ref{eq:1a}). Therefore $P_c$ preserves the angle between two tangent vectors.

A geodesic\index{geodesic} is defined to be a curve $c$ whose velocity vector $\dot{c}(t)$ is parallel along $c$, i.e., 
{\small $\frac{D}{dt}\frac{dc}{dt}=0
$}, which, 
in the case of Levi-Civita connection\index{Levi-Civita connection}, turns out to be
the E-L equation\index{Euler-Lagrange equation} associated with the functional $E(c)=\int_a^b\|\dot{c}(t)\|^2dt$.

\medskip

Finally we shall give a modern formulation of curvature. A multilinear map $R(\cdot,\cdot)\cdot:T_pM\times T_pM\times T_pM\longrightarrow T_pM$ may be defined so as to satisfy
\begin{equation*}
R(X,Y)Z=\nabla_X\nabla_YZ-\nabla_Y\nabla_XZ-\nabla_{[X,Y]}Z~~(X,Y,Z\in \mathfrak{X}(M)).
\end{equation*}
Putting $R(\omega,\boldsymbol{\xi},\boldsymbol{\eta},\boldsymbol{\zeta}):=\omega(R(\boldsymbol{\xi},\boldsymbol{\eta})\boldsymbol{\zeta})~(\omega\in T_p^*M,~\boldsymbol{\xi},\boldsymbol{\eta},\boldsymbol{\zeta}\in T_pM)$, we obtain  the tensor field $R$ of type $(1,3)$ such that $R(dx_l,\partial_i,\partial_j,\partial_k)\!=\!R^{l}_{kij}$. Moreover, taking the subspace $H$ of $T_pM$ spanned by linearly independent vectors $\boldsymbol{\xi}, \boldsymbol{\eta} \in T_pM$, we obtain the sectional curvature
$
K(H)\!=\!\langle R(\boldsymbol{\xi},\boldsymbol{\eta})\boldsymbol{\eta}, \boldsymbol{\xi}\rangle\big/\{\|\boldsymbol{\xi}\|^2\|\boldsymbol{\eta}\|^2-\langle\boldsymbol{\xi},\boldsymbol{\eta}\rangle^2\}^{1/2}
$. 
The {\it Ricci curvature}\index{Ricci curvature} is the quadratic form $R(\boldsymbol{\xi},\boldsymbol{\eta})\!:=\!\text{tr}(\boldsymbol{\zeta}\mapsto R(\boldsymbol{\zeta},\boldsymbol{\eta})\boldsymbol{\xi})$ on $T_pM$, to which 
the Ricci tensor\index{Ricci tensor} is linked by $R_{ij}=R(\partial_i,\partial_j)$. 

\begin{rem}\label{rem:paratra}{\rm {\small
(1) Parallel transports on a surface $S$ 
is related to classical parallelism. Let $X$ is a vector field along a curve $c:[0,a]\rightarrow S$, and let $H$ be the plane tangent to $S$ at $c(0)$. Imagine that $c$ and $X$ are freshly painted and still wet. Roll $S$ on $H$ in such a way that $c(t)$ $~(0\leq t\leq a)$ is a point where the rolled surface at time $t$ is tangent to $H$ ({\it Cartan\index{Cartan, \'{E}lie (1869--1951)} rolling}; Fig.~\!\ref{fig:trans}). Then the curve $c$ is transferred to a curve $c'$ in $H$, $c(t)$ to a point $c'(t)$, and $X(t)$ to a vector $X'(t)$ at $c'(t)$. With a little nudge, we can show that $X$ is parallel if and only if $X'$ along $c'$ is parallel in the classical sense.

\begin{figure}[htbp]
\vspace{-0.7cm}
\begin{center}
\includegraphics[width=.6\linewidth]
{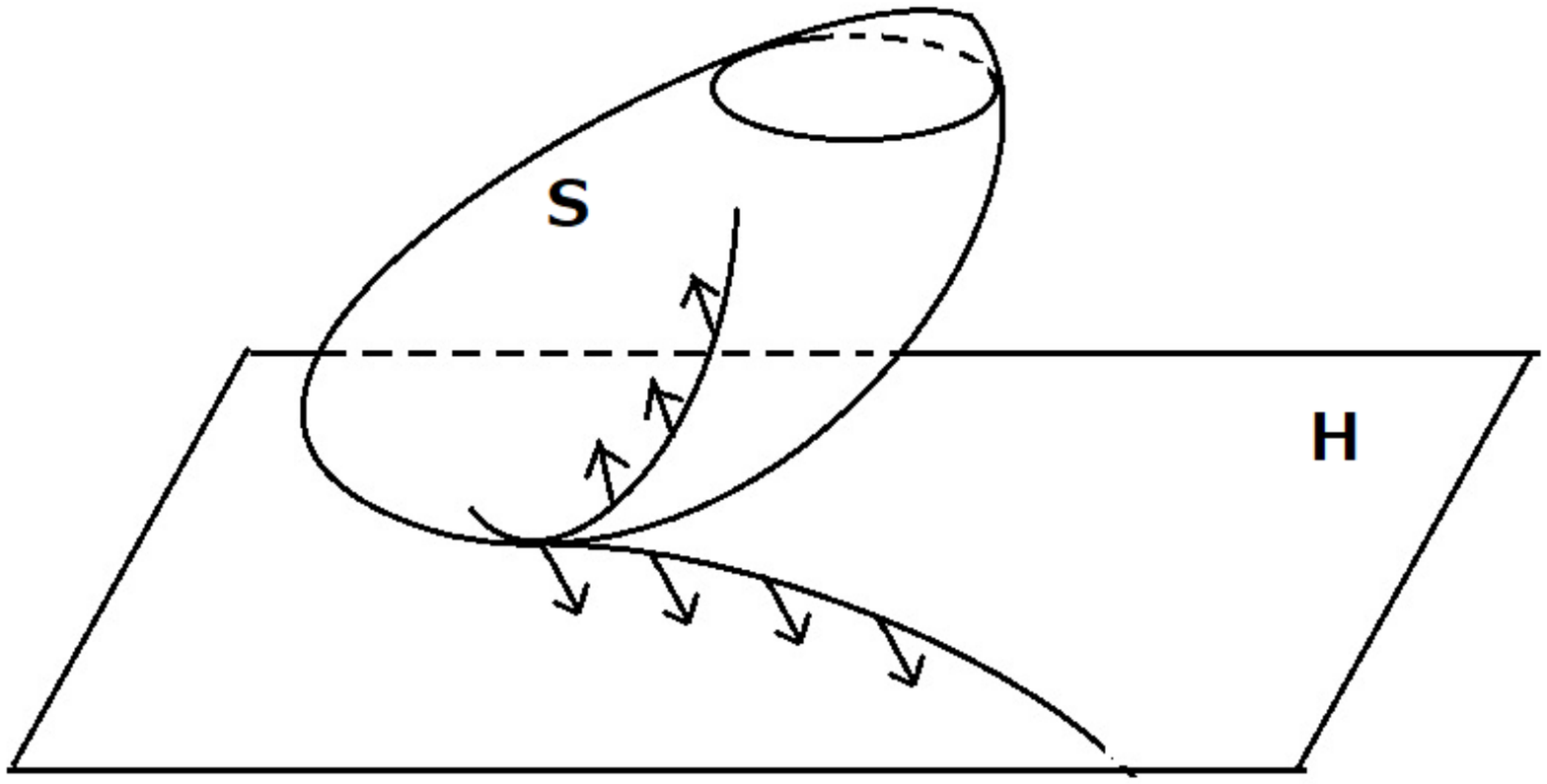}
\end{center}
\vspace{-2.4cm}
\caption{Parallel transport}\label{fig:trans}
\vspace{-0.2cm}
\end{figure}

Pushing further this idea, we may construct a one-to-one correspondence between smooth curves $c$ in a manifold $M$ with $c(0)\!=\!p$ and smooth curves $c_0$ in $T_pM$ with $c_0(0)\!=\!{\bf 0}$. Moreover this correpondence extends to the one for continuous curves by means of {\it stochastic integrals}, with which we may define {\it Brownian motion} on $M$,\footnote{``Brownian motion" originally means the random motion of small particles suspended in fluids. It was named for the botanist R. Brown, the first to study such phenomena (1827). In 1905, Einstein\index{Einstein, Albert (1879--1955)} made its statistical analysis, and observed that the probabilistic behavior of particles is described by the heat equation (Remark \ref{rem:spacialrelativity} (2)).} 
%This fact is 
the starting point of ``stochastic differential geometry" initiated by K. Ito and P. Malliavin.

\smallskip

(2)  The parallel transport\index{parallel transport} around closed loops informs us how a manifold is curved. We shall see this by looking at the parallel transport $P_{\triangle ABC}: T_AS\longrightarrow T_AS$ along the perimeter $A\to$ $B\to$ $C\to A$ of a geodesic triangle\index{geodesic triangle} $\triangle ABC$ on a surface $S$.  It coincides with the rotation of vectors in $T_AS$ through the angle $\theta_{\triangle ABC}:=|\pi-(\angle A+\angle B+\angle C)|$ since the parallel transport\index{parallel transport} preserves the angle of two tangent vectors. By virtue of Gauss's\index{Gauss, Johann Carl Friedrich (1777--1855)} formula (\ref{eq:gaussformula}), we find $\theta_{\triangle ABC}=\left|\iint_{\triangle ABC}K_S~\!d\sigma\right|$.\footnote{The idea explained here is generalized in terms of {\it holonomy groups}, the concept introduced by \'{E}. Cartan\index{Cartan, \'{E}lie (1869--1951)} (1926).}

\smallskip

(3) The idea of connection makes good sense for a {\it vector bundle}, a family of vector spaces parameterized by a manifold which locally looks like a direct product, and are used to formulate {\it gauge theory} in modern physics that provides a unified framework to describe fundamental forces of nature.
%\footnote{Gauge theory has its origin in Weyl's work\index{Weyl, Hermann Klaus Hugo (1885--1955)} in 1918.}
%{\it Gravitation und Elektrizit\"{a}t}, Akademie der Wissenschaften, Berlin, 465-480 (1918).}
\hfill$\Box$

}
}

\end{rem}

\section{Topology of the universe}\label{subsec:topodes}
It cannot be completely denied that the universe may have a complicated structure. Hence it is irresistible to ponder the question ``how can we describe the complexity of the universe?" In pursuing this, it will be appropriate to use topological terms as in the case of the question about finiteness.

In passing, the term ``topology"\index{topology} used to describe a structure on an abstract set also indicates the branch of geometry which principally put a premium on the properties preserved under homeomorphisms\index{homeomorphism}. What, in addition, is important is the concept of {\it homotopy}\index{homotopy}, coined by Dehn and P. Heegaard in 1907, and employed by L. E. J. Brouwer in the current meaning of the word (1911).

\begin{figure}[htbp]
%\vspace{-9.2cm}
\vspace{-0.3cm}
\begin{center}
\includegraphics[width=.58\linewidth]
{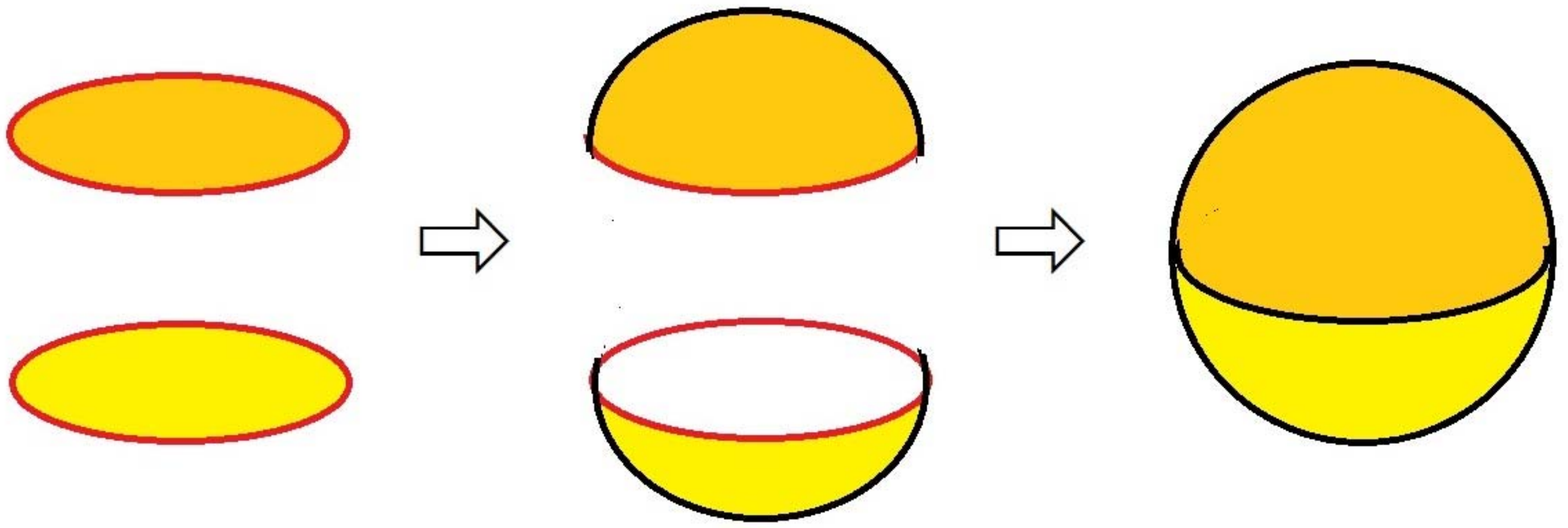}
\end{center}
\vspace{-2.8cm}
\caption{From two disks to a sphere}\label{fig:glue}
\vspace{-0.2cm}
\end{figure}

 A distinct advantage of topology\index{topology} compared with classical geometry is that we can employ lenient operations such as ``gluing (pasting)" figures together. For instance, we obtain a sphere by gluing boundaries of two disks (Fig.~\!\ref{fig:glue}). 

Likewise one can glue boundaries (spheres) of two balls to get the $3$-sphere $S^3=\{(x,y,z,w)\in \mathbb{R}^4|~x^2+y^2+z^2+w^2=1\}$ where one of two balls is the mirror image of another. But this operation cannot be illustrated with drawings, so we shall proceed the other way around; that is, we start with the resulting figure $S^3$, and put
$
S^3_+=\{(x,y,z,w)\in S^3|~w\geq 0\},~~ S^3_{-}=\{(x,y,z,w)\in S^3|~w\leq 0\}$. Clearly $S^3=S^3_+\cup S^3_{-}$ and $S^3_+\cap S^3_{-}=\{(x,y,z)|~x^2+y^2+z^2=1\}=S^2$. Furthermore, $(x,y,z,w)\mapsto (x,y,z)$ yields homeomorphisms of $S^3_{\pm}$ onto the ball $D=\{(x,y,z)\in \mathbb{R}^3|~x^2+y^2+z^2\leq 1\}$. This signifies that $S^3$ is obtained by gluing two balls along their boundaries.

This construction of $S^3$ reminds us of an allegorical vision of Christian afterlife in  Dante's\index{Dante, Alighieri (1265--1321)} poem {\it Divina Commedia}. In Canto XXVIII of {\it Paradiso} (line 4--9), he ``constructs" the empyrean as a {\it mirror image} of the Aristotelian universe. Thus, if it were taken at face value (though not a little farfetched), Dante\index{Dante, Alighieri (1265--1321)} would have obtained $S^3$ as a model of the ``universe" (\cite{Peter}).

Apart from the Dante's\index{Dante, Alighieri (1265--1321)} fantastic imagination, the 3-sphere is surely the easiest conceivable model of a (homogeneous\index{homogeneous} and isotropic\index{isotropic}) finite universe, if not possible to look out it all at once like $S^2$. Einstein\index{Einstein, Albert (1879--1955)} once said
that the universe (at any instant of time) can be viewed globally as $S^3$ (see Sect.~\!\ref{subsec:transfinite}). He wrote, ``Now this is the place where the reader's imagination boggles. `Nobody can imagine this thing,' he cries indignantly. `It can be said, but cannot be thought. I can imagine a spherical surface well enough, but nothing analogous to it in three dimensions.' We must try to surmount this barrier in mind, and the patient reader will see that it is by no means a particularly difficult task" (\cite{ein1}).

J.-P. Luminet and his colleagues put forward an alternative model of a finite universe with constant positive curvature, which, called the {\it Poincar\'{e} homology sphere}\index{Poincar\'{e}, Henri (1854--1912)}, is the quotient of $S^3$ by the {\it binary icosahedral group}, a finite subgroup of the spin group\index{spin group} ${\rm Spin}(3)$.\footnote{\label{fn:icosa}J.-P. Luminet, {\it et al.}, {\it Dodecahedral space topology as an explanation for weak wide-angle temperature correlations in the cosmic microwave background}, Nature, {\bf 425} (2003), 593--595.} Their model---pertinent to the Poincar\'{e} conjecture\index{Poincar\'{e} conjecture} on which we shall comment in due course---explains an apparent periodicity in the {\it cosmic microwave background}, electromagnetic radiation left over from an early stage of the universe in Big Bang cosmology.

\begin{rem}\label{spinxxx}{\rm {\small
The binary icosahedral group is the preimage of the icosahedral group under the homomorphism of ${\rm Spin}(3)$ onto ${\rm SO}(3)$ (Sect.~\!\ref{sec:gauss}), and is explicitly given as the union of 24 quaternions\index{quaternion} $\{\pm 1, \pm i,\pm j, \pm k, \frac{1}{2}(\pm 1\pm i\pm j\pm k)\}$ with 96 quaternions obtained from $\frac{1}{2}(0\pm i\pm \phi^{-1}j\pm \phi k)$ by an even permutation of all the four coordinates $0, 1, \phi^{-1}, \phi$, and with all possible sign combinations ($\phi$ being the golden ratio\index{golden ratio}).
\hfill $\Box$

}
}
\end{rem}

Solely for the convenience of the reader, let us provide a shorthand historical account of topology\index{topology}. The etymology of ``topology"\index{topology} is the German word ``Topologie," coined by J. B. Listing, in his treatise ``Vorstudien zur Topologie" (1847) as a synonym for the ``geometry of position." He learned this discipline from Gauss\index{Gauss, Johann Carl Friedrich (1777--1855)}, and launched the study of several geometric figures like screws, knots, and links from the topological perspective. It was Leibniz\index{Leibniz, Gottfried Wilhelm von (1646--1716)}, if traced back to the provenance of topology\index{topology}, who offered a first rung on the ladder. In a letter to Huygens\index{Huygens, Christiaan (1629--1695)} dated September 8, 1679, Leibniz\index{Leibniz, Gottfried Wilhelm von (1646--1716)} communicated that he was not satisfied with the algebraic methods in geometry and felt for a different type of calculation leading to a new geometry to be called {\it analysis situs}.
%\footnote{He expressed an outlook of the same sort in the manuscript {\it Charcteristica geometrica}.} 
Leibniz\index{Leibniz, Gottfried Wilhelm von (1646--1716)} himself did not put his plan into practice, but his intuition was not leading him astray, was embodied by Euler\index{Euler, Leonhard (1707--1783)} in the solution of the celebrated problem ``Seven Bridges of K\"{o}nigsberg"\index{Seven Bridges of K\"{o}nigsberg} (1736),
%\footnote{{\it Solutio problematis ad geometriam situs pertinentis}, Commentarii Academiae Scientiarum Petropolitanae, {\bf 8} (1736), 128--140.} 
and also in his polyhedron formula $v-e+f=2$ for a convex polyhedron with $v$ vertices, $e$ edges, and $f$ faces, which, stated in a letter to his friend C. Goldbach (November 14, 1750), is recognized as giving a topological characteristic of the sphere (Sect.~\!\ref{sect:curvnon}),\footnote{A copy of Descartes'\index{Descartes, Ren\'{e} (1596--1650)} work around 1630, which was taken by Leibniz\index{Leibniz, Gottfried Wilhelm von (1646--1716)} on one of his trips to Paris, reveals that he obtained an expression for the sum of the angles of all faces of a polyhedron, from which Euler's polyhedron formula\index{Euler's polyhedron formula} can be deduced.} and contains some seeds of {\it combinatorial topology} pioneered by Poincar\'{e}\index{Poincar\'{e}, Henri (1854--1912)} that deals with geometric figures based on their decomposition into combinations of elementary ones.

Gauss\index{Gauss, Johann Carl Friedrich (1777--1855)} took part in the formation of topology\index{topology}, though he opted not to publish any work on topology\index{topology} as usually happened with him. On January 22, 1833 (still in his prime), he noted down a summary of his consideration over the past few months,\footnote{{\it Nachlass zur Electrodynamik} in Gauss {\it Werke}, V, translated by Ricca and Nipoti \cite{ricca}).} which was to usher us into the theory of {\it knots} and {\it links}.

\medskip
{\it 
Of the {\it geometria situs}, which was foreseen by Leibniz\index{Leibniz, Gottfried Wilhelm von (1646--1716)}, and into which only a pair of geometers (Euler\index{Euler, Leonhard (1707--1783)} and Vandermonde
%\footnote{A.-T. Vandermonde, {\it Remarques sur les probl\`{e}mes de situation}, M\'{e}m.~\!Acad.~\!R.~\!Sci.~\!Paris, (1771), 566--574.}) 
were granted a bare glimpse, we know and have, after a century and a half, little more than nothing.

A principal problem at the {\it interface} of {\it geometria situs} and {\it geometria of magnitudinis} will be to count the intertwining of two closed or endless curves.

Let $x,y,z$ be the coordinates of an undetermined point on the first curve; $x',y',z'$ those of a point on the second and let
{\small 
\begin{equation}\label{eq:link1}
\iint \frac{(x'-x)(dydz'-dzdy')+(y'-y)(dzdx'-dxdz')+(z'-z)(dxdy'-dydx')}
{\big((x'-x)^2+(y'-y)^2+(z'-z)^2\big)^{3/2}}=V
\end{equation}
}
\hspace{-0.15cm}\noindent then this integral taken along both curves is $
4m\pi
$, 
$m$ being the number of intertwinings {\rm [called the {\it linking number}\index{linking number} today]}. The value is reciprocal, i.e., it remains the same if the curves are interchanged. 
}

\medskip

His pretty observation has something to do with electromagnetism\index{electromagnetism}. Indeed, (\ref{eq:link1}) reminds us of the {\it Biot--Savart law}\index{Biot--Savart law} (1820), an equation describing the magnetic field $\boldsymbol{B}$ generated by an electric current $\boldsymbol{i}$, which is named after J.-B. Biot and F. Savart. It is stated as 
\begin{equation}\label{eq:biot}
\boldsymbol{B}(\boldsymbol{x})=\frac{\mu_0}{4\pi}\int_{\mathbb{R}^3}\frac{\boldsymbol{i}(\boldsymbol{y})\times (\boldsymbol{x}-\boldsymbol{y})}{\|\boldsymbol{x}-\boldsymbol{y}\|^3}d\boldsymbol{y},
\end{equation}
where $\mu_0$ is the magnetic permeability of vacuum (fn.~\!\ref{diel}).
%\footnote{{\it Note sur le Magn\'{e}tisme de la pile de Volta}, Annales Chim. Phys., {\bf 15} (1820), 222--223.}

Another (more pedagogical) example
in topology\index{topology} to which Gauss\index{Gauss, Johann Carl Friedrich (1777--1855)} slightly contributed  is the notion of {\it winding number}\index{winding number}. It is defined analytically by 
\begin{equation}\label{eq:winding}
W(c,p_0)=
\frac{1}{2\pi}\int_c\frac{-(y-y_0)dx+(x-x_0)dy}{(x-x_0)^2+(y-y_0)^2},
\end{equation}
where $c$ is a closed (not necessarily simple) plane curve that does not pass through $p_0=(x_0,y_0)$. This is an integer because, if we parameterize the curve $c$ as $c(t)=(x_0+r(t)\cos\theta(t), y_0+r(t)\sin\theta(t))$ $~(0\leq t\leq 1)$ with continuous functions $r(t)$ and $\theta(t)$, then, without costing much effort, we see that the integral in (\ref{eq:winding}) is transformed into $\frac{1}{2\pi}\int_0^1~\!d\theta(t)=\frac{1}{2\pi}[\theta(1)-\theta(0)]$, which, because $c(0)=c(1)$, represents the total (net) number of times that $c$ travels around $p_0$. 

From the analytic expression (\ref{eq:winding}) and the fact that an integral-valued continuous function is constant, it follows that, if $c_1$ and $c_2$ are homotopic as maps of the circle $S^1$ into $\mathbb{R}^2\backslash\{p_0\}$, then $W(c_1,p_0)=W(c_2,p_0)$. Moreover, $W(c,p)$ does not depend on $p$ so far as $p$ is in the connected component $D_0$ of $\mathbb{R}^2\backslash c$ containing $p_0$; hence one may put $W(c,D_0):=W(c,p_0)$.
%\footnote{The winding number is a special case of the {\it degree} of a map between two manifolds.} 

The homotopy invariance of winding numbers\index{winding number} stands us in good stead in proving the fundamental theorem of algebra\index{fundamental theorem of algebra} (cf.~\!Gauss\index{Gauss, Johann Carl Friedrich (1777--1855)} {\it Werke} III, 31--56). Given a polynomial $f(z)\!=\!z^n+a_1z^{n-1}+\cdots+a_{n-1}z+a_n$ $(a_n\neq 0)$, one can take $R$ $>$ $0$ such that $c_s(t)\!:=\!f_s(Re^{2\pi\sqrt{-1}t})\neq 0$ $~(0\leq s, t\leq 1)$, where $f_s(z)=z^n+$ $s(a_1z^{n-1}+\cdots+a_{n-1}z +a_n)$, so that $W(c_1,0)=W(c_0, 0)=n(\neq 0)$. On the other hand, if $f(z)$ $\neq 0$ for any $z\in \mathbb{C}$, then $c^r(t)\!:=\!f(re^{2\pi\sqrt{-1}t})\neq 0$ $(0$ $\leq r\leq R,~0\leq t\leq 1)$, so that $W(c_1,0)=W(c^R,0)=W(c^0,0)=0$; thereby a contradiction, and hence $f(z)$ $=0$ must have a solution. 

What deserves a good deal of attention is that the integral (\ref{eq:winding}) is neatly expressed as the complex line integral:
$
\frac{1}{2\pi\sqrt{-1}}$ $\int_c\frac{1}{z-z_0}~\!dz$ $~(z_0=x_0+y_0\sqrt{-1})$.
Thus as 
Gauss\index{Gauss, Johann Carl Friedrich (1777--1855)} briefly noticed in the letter to Bessel in 1811 (Sect.~\!\ref{sec:gauss}), the winding number\index{winding number} reveals the nature of the complex logarithm\index{complex logarithm} (see Remark \ref{fn:complexline}).

 Winding numbers appeared in the study of ``signed areas"\index{signed area} by Albrecht Ludwig Friedrich Meister\index{Meister, Albrecht Ludwig Friedrich (1724--1788)} (1724--1788).\footnote{{\it Generalia de genesifigurarum planarum, et inde pendentibus earum affectionibus}, Novi Commentarii Soc.~\!Reg.~\!Scient.~\!Gott., {\bf 1} (1769/1770), 144--180.} 
Here, the signed area {\it surrounded} by a general closed curve $c$ is defined to be $\frac{1}{2}\int_c(-ydx+xdy)$. If $c$ is simple and counterclockwise, this is the ``genuine" area ${\rm Area}(D)$ of the domain $D$ surrounded by $c$ in virtue of Green's theorem\index{Green's theorem} (\ref{eq:lineintegral0}). For a general $c$ and the bounded connected components $D_1,\ldots,D_N$ of $\mathbb{R}^2\backslash c$, we have  
$$
\frac{1}{2}\int_c(-ydx+xdy)=\sum_{k=1}^NW(c,D_k){\rm Area}(D_k).
$$
In fact, this is a consequence of the formula $\int_c\omega=\sum_{k=1}^NW(c,D_k)\int_{\partial D_k}\omega$ holding for any 1-form $\omega$ on $\mathbb{R}^2$, where the boundary $\partial D_k$ is assumed to have counterclockwise orientation. This suggests a formal framework allowing us to write ``$c=\sum_{k=1}^NW(c,D_k)\partial D_k$." Indeed, this identity can be justified by the notion of ``chain" in {\it homology theory}\index{homology} (Remark \ref{rm:form} (2)).

\medskip
Moving away from the small treasure trove of Gauss's\index{Gauss, Johann Carl Friedrich (1777--1855)} ``toys", we shall see how topology\index{topology} evolved thereafter. 

The period from the mid-19th century to the early 20th century was the infancy of topology\index{topology}. It coincides with the era when mathematics began to develop autonomously on a rigid base. Prompted by such a state of affairs, geometers created a wealth of novel ideas. One of the marked achievements in this period is the discovery of {\it one-sided} surfaces\index{one-sided surface} that differ from ordinary surfaces in the way of their disposition in the space.
Remember that a smooth surface $S$ is two-sided if and only if it admits a global continuous (unit) normal vector field, or equivalently $S$ is orientable\index{orientable}. (This equivalence is deduced from the fact that, if $\boldsymbol{S}_1:U_1\longrightarrow S$ and $\boldsymbol{S}_2:$ $U_2$ $\longrightarrow$ $S$ are two local parametric representations\index{local parametric representation} of a surface $S$ with $\boldsymbol{S}_1(U_1)\cap\boldsymbol{S}_2(U_2)\neq \emptyset$, then $(\boldsymbol{S}_1)_u\times (\boldsymbol{S}_1)_v=\frac{\partial(z,w)}{\partial(u,v)}(\boldsymbol{S}_2)_z\times (\boldsymbol{S}_2)_w$,  where $\boldsymbol{S}_2{}^{-1}\circ\boldsymbol{S}_1(u,v)$ $=(z(u,v),w(u,v))$ on $U_1\cap \boldsymbol{S}_1{}^{-1}(\boldsymbol{S}_2(U_2)).)$ The foremost examples of one-side surfaces are the {\it M\"{o}bius band}\index{M\"{o}bius band} and the {\it Klein bottle}\index{Klein bottle}. The former was discovered by M\"{o}bius\index{M\"{o}bius, August Ferdinand (1790--1868)} (1865) and Listing (1862).
%\footnote{M\"{o}bius, {\it \"{U}ber die Bestimmung des Inhaltes eines Poly\"{e}ders}, Ber.~\!Verh.~\!S\"{a}chs.~\!Ges.~\!Wiss., {\bf 17} (1865), 31--68. Listing, {\it Der Census R\"{a}umlicher Complexe oder Verallgemeinerung des Eulerschen Satzes von der Polyedern}, 1862.} 
The latter was constructed by Klein\index{Klein, Christian Felix (1849--1925)} (1882).
%\footnote{{\it \"{U}ber Riemanns Theorie der Algebraischen Functionen und Ihre Integrale}, 1882.} 
The projective plane\index{projective plane} $P^2(\mathbb{R})$ is another example\index{one-sided surface}.

\begin{figure}[htbp]
%\vspace{-11.5cm}
\vspace{-0.5cm}
\begin{center}
\includegraphics[width=.72\linewidth]
{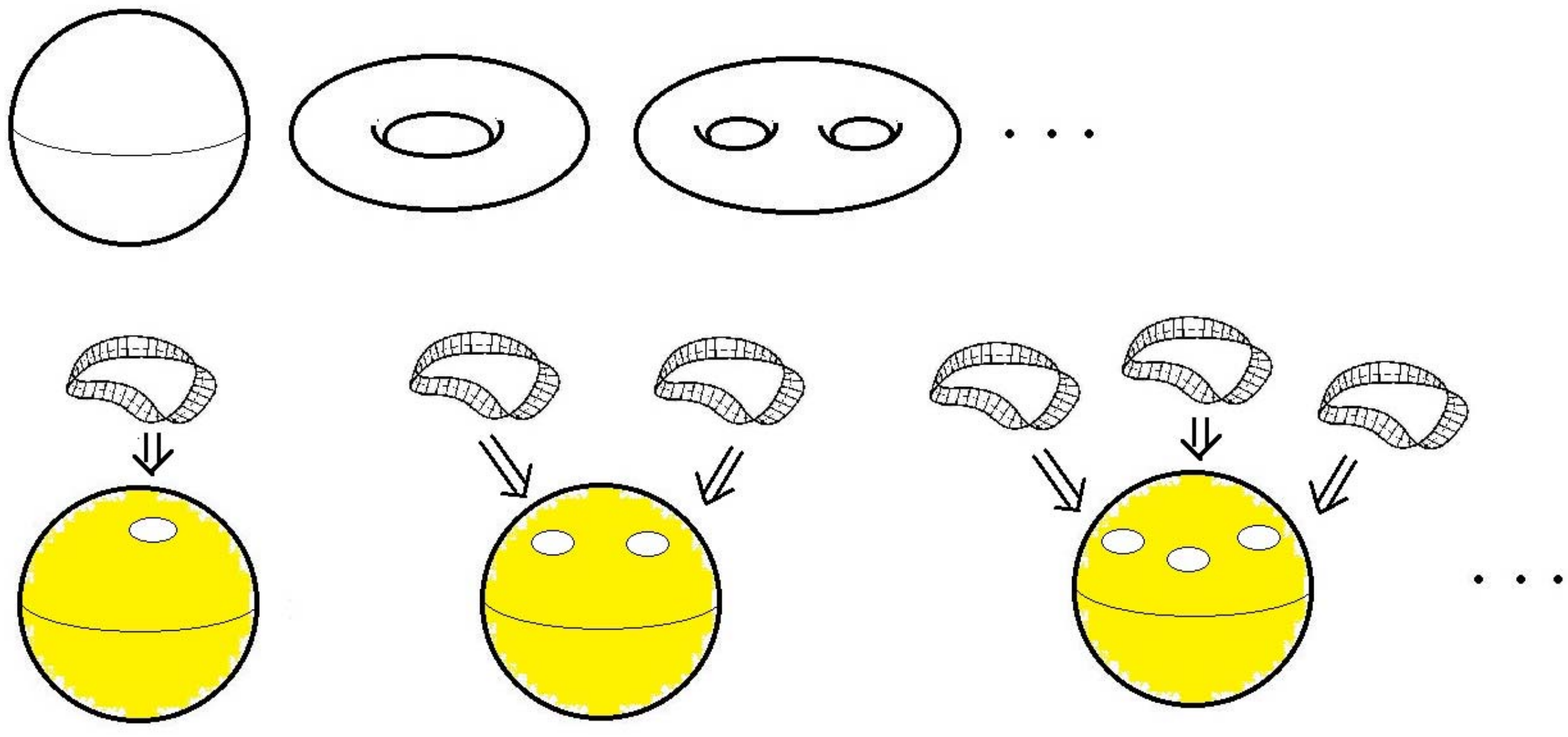}
\end{center}
\vspace{-2.6cm}
\caption{Classification of closed surfaces}\label{fig:surfaceclass}
\vspace{-0.2cm}
\end{figure}

The classification of closed surfaces is the problem attacked in the late 19th century. The case of two-sided surfaces was independently treated by M\"{o}bius\index{M\"{o}bius, August Ferdinand (1790--1868)} (1863) and M. E. C. Jordan (1866).
%\footnote{M\"{o}bius, {\it Theorie der elementaren Verwandschaften}, Abhandlung S\"{a}chsische Gesellschaft der Wissenschaften, 1863. Jordan, {\it Sur la d\'{e}formation des surfaces}, Journal de Math\'{e}matique, {\bf 11} (1866), 105--109.}
One-sided closed surfaces were classified by von Dyck\index{von Dyck, Walther Franz Anton (1856--1934)} (1888).
%\footnote{{\it Beitr\"{a}ge zur Analysis situs} I, Math.~\!Ann., {\bf 32} (1888), 457--512.} 
As depicted in Fig.~\!\ref{fig:surfaceclass}, every closed one-sided surface\index{one-sided surface} can be constructed from $S^2$ by attaching a finite number of M\"{o}bius bands\index{M\"{o}bius band} (the leftmost is the projective plane\index{projective plane}, and the next is the Klein bottle\index{Klein bottle}). Thus the ``shapes" of 2D models of the {\it finite} universe are completely comprehended.

Now, one may ask, ``what about the higher-dimensional case?" If restricted to 3-  or 4-manifolds, this is a problem concerning a possible structure of our universe, and hence has a much higher profile. More specifically, one may pose the question: Does our universe have a structure like a surface with many holes? Here, the meaning of ``hole" is well elucidated in terms of {\it algebraic topology}
which borrows tools from abstract algebra crystallized  in the early 20th century.

Riemann\index{Riemann, Georg Friedrich Bernhard (1826--1866)} had a clear perception of what we now call the {\it Betti numbers}\index{Betti number}, significant topological invariants related to the Euler characteristic\index{Euler characteristic} and defined today as the rank of the {\it homology group} (see the remark below). His grandiose program---showing again that he was far ahead of the times---was published posthumously as the {\it Fragment aus der Analysis Situs} \cite{re}. Later, Enrico Betti\index{Betti, Enrico (1823--1892)} (1823--1892), an Italian friend of Riemann, made clearer what Riemann thought.\footnote{{\it Sopra gli spazi di un numero qualunque di dimensioni}, Ann.~\!Mat.~\!Pura Appl., {\bf 4} (1871), 140-158. A letter dated October 6, 1863 from Betti to his colleague tells us that he got an accurate conception about the connectivity of spaces through a conversation with Riemann.} Finally, their idea came to fruition as homology theory by Poincar\'{e}\index{Poincar\'{e}, Henri (1854--1912)} (1899),
%\footnote{{\it Compl\'{e}ment \`{a} l'Analysis Situs}, Rend.~\!Circ.~\!Mat.~\!Palermo, {\bf 13} (1899), 285--343.} 
who also introduced the {\it fundamental group}, another algebraic tool and a modern version of ``Greek geometric algebra"\index{geometric algebra} in a sense because the algebraic system in question are directly constructed from geometric figures. 

In the meantime, G. de Rham related homology to differential forms\index{differential form}, proving what we now call {\it de Rham's theorem}\index{de Rham's theorem}, one of the earliest outcomes in global analysis (1931). This line---coupled with Chern's generalization of the Gauss-Bonnet formula\index{Gauss-Bonnet formula} (1944) and the theory of {\it harmonic integrals} developed by Weyl\index{Weyl, Hermann Klaus Hugo (1885--1955)}, W. V. D. Hodge and K. Kodaira---culminated in the {\it Atiyah-Singer theory}, one of the most exhilarating adventures of the 20th century. Additionally, the work of M. Gromov on large-scale aspects of manifolds can be considered as far-reaching generalizations of a batch of results on relations between curvature\index{curvature} and topology\index{topology} obtained by geometers on and after the 1950s.

The success of the taxonomy of closed surfaces impelled topologists to attack the case of closed 3-manifolds, but the matters turned out to be more complicated than expected. In 1982, W. P. Thurston offered the {\it geometrization conjecture} as an initial template, which daringly says that all 3-manifolds admit a certain kind of decomposition involving the {\it eight geometries} (one of them is hyperbolic geometry\index{hyperbolic geometry}). This spectacular supposition was proved by G. Perelman in 2003. His arguments include the proof of the long-standing {\it Poincar\'{e} conjecture}\index{Poincar\'{e} conjecture} which claims, ``a simply connected closed 3-manifold is homeomorphic to the 3-sphere."\footnote{Poincar\'{e}, {\it Cinqui\`{e}me compl\'{e}ment \`{a} l'analysis situs}, Rend.~\!Circ.~\!Mat.~\!Palermo, {\bf 18} (1904), 45--110). He claimed, at first, that homology is sufficient to tell if a closed 3-manifold is homeomorphic to $S^3$ (1900), but 4 years later, he found that what we now call the {\it Poincar\'{e} homology sphere} gives a counterexample.} A conspicuity of his proof is that a non-linear evolutional equation involving Ricci curvature\index{Ricci curvature} 
%($\partial_t(g_{ij})=-2R_{ij}$) 
is exploited in an ingenious way; thereby displaying a miraculous trinity of differential geometry, topology, and analysis.

The situation gets out of hand when trying to classify closed manifolds $M$ with ${\rm dim}~\!M\geq 4$.
%\footnote{Andrey Andreyevich Markov Jr., Proceedings of the ICM, Edinburgh 1958.} 
This is because an arbitrary {\it finitely presented group} (a group defined by a finite number of generators, and a finite number of defining relations) can be the fundamental group of a closed 4-manifold, and there is no algorithm\index{algorithm} to decide whether two finitely presented groups are isomorphic.\footnote{A. Turing and A. Church rigorously mathematized the concept of algorithm\index{algorithm} by analyzing the meaning of computation. They independently treated the 
%{\it Entscheidungsproblem} (
{\it decision problem} challenged by Hilbert\index{Hilbert, David (1862--1943)} in 1928, and showed that a general solution is impossible (1936).}

\medskip

We close this section with an instructive example in {\it differential topology} which has been intensively studied in the latter half of the last century.

\begin{figure}[htbp]
%\vspace{-3.8cm}
\vspace{-0.6cm}
\begin{center}
\includegraphics[width=.92\linewidth]
{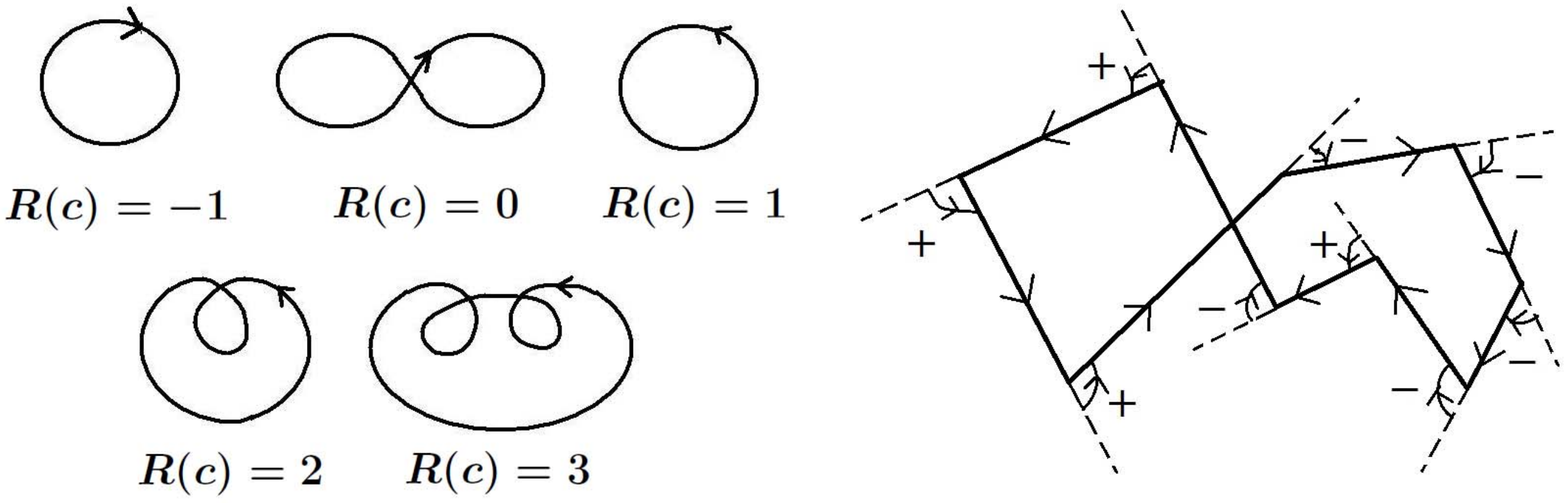}
\end{center}
\vspace{-4.5cm}
\caption{Rotation number and exterior angles}\label{fig:polygonal1}
\vspace{-0.2cm}
\end{figure}

We shall say that a smooth curve $c$ in $\mathbb{R}^2$ is {\it regular}\index{regular curve} if $\dot{c}(t)\neq \boldsymbol{0}$ for every $t$. Given a closed regular curve $c$, we define the {\it rotation number}\index{rotation number} $R(c)$ by setting $R(c)=W(\dot{c},{\bf 0})$. In plain language, $R(c)$ is the total number of times that a person walking once around the curve turns counterclockwise.\footnote{As known, at least informally, by Meister\index{Meister, Albrecht Ludwig Friedrich (1724--1788)} ({\it Generalia}, 1770), the rotation number\index{rotation number} is defined for a closed oriented polygonal curve. It is the sum of all (signed) exterior angles divided by $2\pi$ (the right of Fig.~\!\ref{fig:polygonal1}). In the {\it Geometria Speculativa} ({\it ca}.~\!1320),  Bradwardine\index{Bradwardine, Thomas (c. 1290--c. 1349)} (Sect.~\!\ref{sec:whatinfinity}) studied a class of general polygonal curves, which seems to have a relevance with rotation number\index{rotation number} (Jeff Erickson, {\it Generic and Regular Curves}, PDF, 2013).} This notion came up in Gauss's\index{Gauss, Johann Carl Friedrich (1777--1855)} work of experimental nature; he observed, among others ({\it Werke}, VIII, 271--281), that, if we denote by $n(c)$ the number of self-intersection points of $c$, then we have the sharp inequality: $n(c)\geq ||R(c)|-1|$.

In 1937, H. Whitney\index{Whitney, Hassler (1907--1989)} proved that a closed regular curve can be deformed to another one through a ``smooth" family of closed regular curves (this being the case, two curves are said to be {\it regularly homotopic}) if and only if two curves have the same rotation number\index{rotation number}.
%\footnote{{\it On regular closed curves in the plane}, Compositio Math., {\bf 4} (1937), 276--284.} 
In particular, a clockwise oriented circle ($R=-1$) is not regularly homotopic to a counterclockwise oriented circle ($R=1$). This fact is rephrased as ``one cannot turn a circle {\it inside-out}," because giving an orientation to a closed regular curve $c$ is equivalent to choosing one of two normal unit vector fields along $c$, and the inside and outside of a circle correspond to the inner normal and outer normal vector fields along it, respectively. Of particular note is the unexpected result by S. Smale that one can turn the sphere $S^2$ inside-out, quite by contrast to the case of the circle (1959). 

\begin{rem}\label{rm:form}{\small{\rm 
(1) The rotation number\index{rotation number} is linked with the signed curvature\index{signed curvature}  (Remark \ref{rem:signedcurvature}) via the formula   $\int_c\kappa(c(s))ds$ $=2\pi R(c)$. Indeed, 
$$
2\pi R(c)=\int_c\frac{\ddot{y}(s)\dot{x}(s)-\ddot{x}(s)\dot{y}(s)}{\dot{x}(s)^2+\dot{y}(s)^2}~\!ds=\int_c[\ddot{y}(s)\dot{x}(s)-\ddot{x}(s)\dot{y}(s)]~\!ds=\int_c\kappa(c(s))~\!ds.
$$
We thus have a fusion of topology\index{topology} and differential geometry at a very basic level.

\smallskip

(2) In the {\it Analysis situs}, Poincar\'{e}\index{Poincar\'{e}, Henri (1854--1912)} defined ``homology classes" in a somewhat fuzzy way. The group structure on homology classes, which he did not explicitly indicate, was studied by E. Noether and others in the period 1925--28. Later on, S. Eilenberg  developed the {\it singular homology theory} that allows to define the homology group for a general topological space\index{topological space} (1944).

A central role in homology theory is played by a {\it chain complex}, a series of homomorphisms 
$
\cdots\stackrel{\partial_{k+1}}{\longrightarrow} C_k(X,\mathbb{Z})\stackrel{\partial_k}{\longrightarrow}\cdots\stackrel{\partial_2}{\longrightarrow} C_1(X.\mathbb{Z})\stackrel{\partial_1}{\longrightarrow} C_0(X,\mathbb{Z})\stackrel{\partial_0}{\longrightarrow}0$ satisfying $\partial_k\circ \partial_{k+1}=0$, where, in the singular case, $C_k(X,\mathbb{Z})$ is the free abelian group generated by ``singular $k$-simplices" (continuous maps $\sigma_k$ from the $k$-simplex into a topological space $X$). The $k$-simplex here is a generalization of a segment ($k\!=\!1$), a triangle ($k\!=\!2$), and a tetrahedron ($k\!=\!3$). The homomorphism $\partial_k$ brings $\sigma_k$ to the sum of the singular $(k-1)$-simplices represented by the restriction of $\sigma_k$ to the faces of the $k$-simplex, with an alternating sign to take orientation into account. The factor group $H_k(X,\mathbb{Z})$ $=$ ${\rm Ker}~\partial_k/{\rm Image}~\partial_{k+1}$ is what we call the $k$-th homology group of $X$. The $k$-th Betti number\index{Betti number} $b_k(X)$ is the rank of $H_k(X,\mathbb{Z})$ (if finite). When $b_k(X)=0$ for $k>d$, the Euler characteristic\index{Euler characteristic} of $X$ is defined to be $\chi(X)=\sum_{k=0}^{d}(-1)^kb_k(X)$, which, for a surface, agrees with the definition relying on a triangulation (Sect.~\!\ref{sect:curvnon}).
 
\smallskip

(3) For a $k$-form $\omega$ on a convex domain $D\subset\mathbb{R}^d$ satisfying $d\omega=0$, there is a $(k-1)$-form $\eta$ with $d\eta=\omega$; that is, $H^k_{\rm dR}(D)=\{0\}$. This (the {\it Poincar\'{e} lemma}\index{Poincar\'{e} lemma}),\footnote{\label{fn:poincare}Poincar\'{e}, {\it  Les M\'{e}thodes nouvelles da la M\'{e}canique c\'{e}leste}, vol. 3, Gauthier-Villars, Paris, 1899, pp.9--15. 
This fact, which, in the case of 1-forms, dates back to Euler\index{Euler, Leonhard (1707--1783)} (1724/1725) and Clairaut, is essentially ascribed to Vito Volterra (1860--1940).} in tandem with the generalized Stokes' formula\index{generalized Stokes' formula}, is crucial in the proof of de Rham's theorem\index{de Rham's theorem} asserting that $b_k(M)={\rm dim}~\!H^k_{\rm dR}(M)$ for a closed manifold $M$.

\smallskip

(4) Define the $d$-form $\Omega$ on a closed even-dimensional orientable\index{orientable} manifold $M$ by 
$$
\hspace{-0.3cm}\Omega = \frac{(-1)^n}{2^{2n}\pi^nn!}\sum_{\sigma=(i_1,\ldots,i_d)}{\rm sgn}(\sigma)
g^{i_1i_2}g^{i_3i_4}\cdots g^{i_{2n-1}i_{2n}}
\Omega_{i_1i_2}
\wedge \Omega_{i_3i_4}\wedge\cdots\wedge \Omega_{i_{2n-1}i_{2n}},
$$

\vspace{-0.5cm}
$$
\Omega_{ij}=\sum_{m,k,l=1}^dg_{im}R^m_{jkl}dx_k\wedge dx_l \quad ({\rm dim}~\!M=d=2n).
$$
Then Chern's generalization of the Gauss-Bonnet formula\index{Gauss-Bonnet formula} is expressed as 
$$
\chi(M)=\int_M \Omega. 
$$ 

(5) Hodge's theorem says that the space of {\it harmonic forms} $\mathcal{H}^k(M)=\{\omega\in A^k(M)|~\!d_k\omega=d_{k-1}^*\omega=0\}$ is isomorphic to $H^k_{\rm dR}(M)$, where $d_{k-1}^*$ is the adjoint of $d_{k-1}$ (with respect to natural inner products on $A^*(M)$). Using this fact, we obtain 
\begin{equation}\label{eq:indexthm}
\dim {\rm Ker}~\!D-\dim {\rm Coker}~\!D=\chi(M),
\end{equation}
 where $D=d+d^*$, an operator from the space of even forms to the space of odd forms. A special feature of $D$ is that it is {\it elliptic}; i.e., its {\it principal symbol} is invertible.\footnote{A differential operator $P$ of order $m$ on $M$ (acting on vector-bundle valued functions) is locally expressed as $P\!=\!\sum_{|\alpha|\leq m}A_{\alpha}(\boldsymbol{x})D^{\alpha}$, where $\alpha\!=\!(\alpha_1,\ldots,\alpha_d)$ denotes a $d$-tuple of non-negative integers, $|\alpha|=\alpha_1+\cdots+\alpha_d$, $D^{\alpha}\!=\!(\partial^{\alpha_1}/\partial x_1^{\alpha_1})$ $\cdots (\partial^{\alpha_d}/\partial x_d^{\alpha_d})$, and $\{A_{\alpha}(\boldsymbol{x})\}$ are  matrix-valued functions. The principal symbol of $P$ is defined to be $\sigma_P(\boldsymbol{x},\boldsymbol{\xi})\!:$ $=$ $\sum_{|\alpha|=m}A_{\alpha}(\boldsymbol{x})\xi_1^{\alpha_1}\cdots\xi_d^{\alpha_d}$ ($\boldsymbol{\xi}\neq {\bf 0}\in T^*M$).} 

The left-hand side of (\ref{eq:indexthm}) is called the {\it analytic index} of $D$. The Atiyah-Singer theorem\index{Atiyah-Singer theorem} asserts that the analytic index (${\rm Ind}~\!P$) of an elliptic operator $P$ equals the {\it topological index} defined in terms of topological data, and that there is a $d$-form $\Omega$ involving both curvature tensor\index{curvature tensor} and principal symbol of $P$ such that ${\rm Ind}~\!P=\int_M\Omega$.
\hfill$\Box$

}
}
\end{rem}

\section{Right and left in the universe}\label{sec:right}

\begin{figure}[htbp]
%\vspace{-3.8cm}
\vspace{-0.5cm}
\begin{center}
\includegraphics[width=.35\linewidth]
{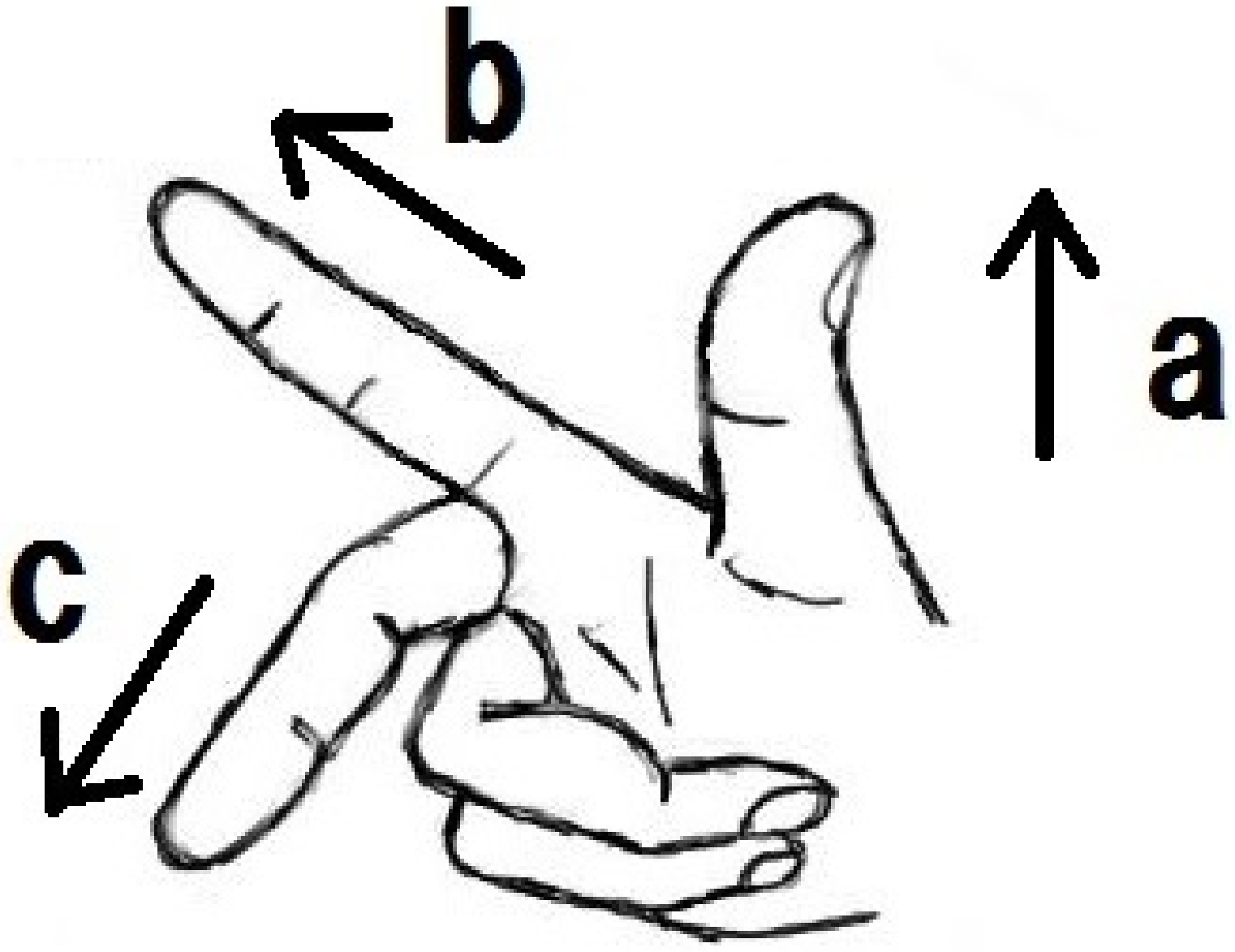}
\end{center}
\vspace{-3.6cm}
\caption{The frame corresponding to the gun-shaped right hand}\label{fig:righthand}
\vspace{-0.2cm}
\end{figure}

An intriguing question about our universe is whether it has ``one side" or ``two sides." In discussing sideness, it seems necessary, if looking back to the case of surfaces (Sect.~\!\ref{subsec:topodes}), to assume the existence of ``outer side" of the universe. The fact is that it is genuinely intrinsic\index{intrinsic} in character as in the case of finiteness. Why is it so? This time, we bring in the idea of ``right and left," the terms used not only in everyday life, but also in mathematics, which are usually designated by the human hands. No less obvious perhaps to our eyes, but no less essential, is the recognition that the human hands are represented by {\it frames}\index{frame} (ordered basis), the notion obtained by whittling away the extraneous details of hands; see Fig.~\!\ref{fig:righthand} illustrating the {\it right-handed frame} $(\boldsymbol{a}, \boldsymbol{b},\boldsymbol{c})$.\footnote{The fact that there are just two kinds of frames in $\mathbb{R}^d$ is equivalent to that the general linear group ${\rm GL}_d(\mathbb{R})$ has two connected components.}

To clarify the nature of right and left, we shall refer to the notorious {\it mirror paradox}\index{mirror paradox} that is best stated as a question: Why do flat mirrors reverse left and right, but not top and bottom? {\it Prima facie}, this sounds puzzling,  but the mirror paradox is simply a paradox of ``a red herring" or a word play to confuse or startle people by making them believe that the wording ``left and right" is the same sort of ``top and bottom." Indeed, the latter is not described by a frame, but by a single vector; say, the vector represented by the directed segment joining bottom and top.\footnote{Shin-itiro Tomonaga argues that this paradox cannot be resolved by neither geometric optics nor mathematics (1963), while Martin Gardner says that the mirror does not reverse right and left, but does reverse front and back (\cite{gard}).}

Now, the claim ``the universe has two sides" is rephrased as ``the right-handed frame\index{right-handed} and the left-handed\index{left-handed} frame\index{frame} are distinct wherever they are" in the sense that these two frames cannot be superposed no matter how one is moved to another in the universe. Otherwise expressed, if the universe has ``only one side," then the distinction between the right-handed\index{right-handed} frame and the left-handed\index{left-handed} frame does not make sense; namely there is only one kind of frame. What is more, as in the case of surfaces, two-sidedness\index{two-sidedness} and orientability\index{orientability} of the universe are equivalent when we regard the universe as a smooth manifold. Thus the distinguishability between right and left depends on the global intrinsic structure of the universe.

Even if our universe is Euclidean, there remains a subtle problem on ``right and left" that dates back to the dispute between Newton\index{Newton, Issac (1642--1726)} and Leibniz\index{Leibniz, Gottfried Wilhelm von (1646--1716)} and have fascinated philosophers and cosmologists since then (\cite{gross}). Newton\index{Newton, Issac (1642--1726)} insisted that there is a `left' and `right' in the universe, while Leibniz opposed this view and argued that left and right are in no way different from each other.\footnote{Leibniz's critique of Newton\index{Newton, Issac (1642--1726)} is unfolded in his third letter to S. Clarke (fn.~\!\ref{fn:clarke}).} Puzzled by {\it enantiomorphs} for decades, Kant\index{Kant, Immanuel (1724--1804)} mused about this issue in connection with Leibniz's\index{Leibniz, Gottfried Wilhelm von (1646--1716)} claim that all spatial facts should be reduced to facts about relative distances between material bodies
%\footnote{
({\it Prolegomena}, 1738).
% zu einer jeden k\"{u}nftigen Metaphysik}, 1783.} 
He surmised in his prolix discussion that the difference between similar but not congruent things cannot be made intelligible by understanding about any concept, and are known only through sensuous intuition. He even set up an extreme thought-experiment to adduce conclusive evidence of his assertion.\footnote{{\it Von dem ersten Grunde des Unterschieds der Gegenden im Raume}, 1768. In it, he considers a marble hand broken off a statue that is supposed to be the only object in the universe, and asks whether it makes sense to say that it is still either a right hand or a left.} According to him, this suggests, after all, that space is not independent of the mind that perceives it.

The purport of what Newton\index{Newton, Issac (1642--1726)}, Leibniz\index{Leibniz, Gottfried Wilhelm von (1646--1716)}, and Kant\index{Kant, Immanuel (1724--1804)} thought up is clarified if we replace Kant's claim by the statement, ``one cannot furnish any mathematical characterization of right-handedness of a frame." This might perplex the reader because, in several places of this essay, we used {\it vector product}\index{vector product}, whose (casual) definition relies on the right-handed frame; thus seemingly incomplete as a mathematical concept. As a matter of fact, we may select any frame as a reference frame; if we choose another kind of frame, then the resulting vector has the opposite direction, and there is no trouble caused at all.

\begin{figure}[htbp]
%\vspace{-5.4cm}
\vspace{-0.3cm}
\begin{center}
\includegraphics[width=.39\linewidth]
{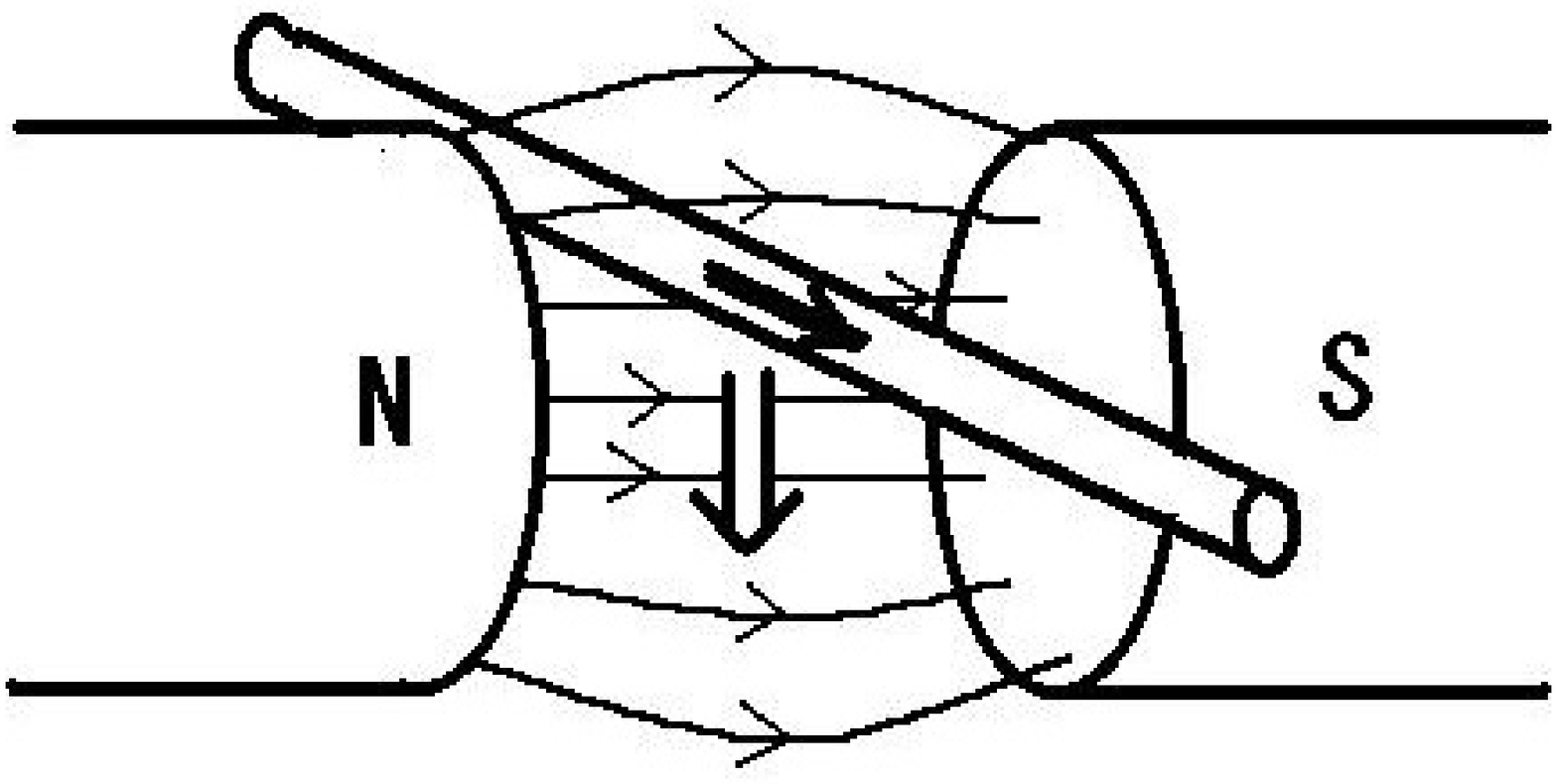}
\end{center}
\vspace{-1.7cm}
\caption{Fleming's left-hand rule}\label{fig:magnet}
\vspace{-0.2cm}
\end{figure}

How about the issue of right and left in physics?  Ernst Mach\index{Mach, Ernst (1838--1916)} (1838--1916)\footnote{\label{fn:Mach}Mach is best-known for the sensation-based theory of reality. In {\it Die Mechanik in ihre Entwicklung historisch-kritisch dargestellt} (1883), he criticized Newton's\index{Newton, Issac (1642--1726)} conclusion on absolute motion based on the bucket experiment (fn.~\!\ref{fn:bucket}), saying that centrifugal force may act on water in the stationary bucket if the universe rotates. According to him, the inertia and acceleration of a body are determined by all of the matter of the universe. This view was called the {\it Mach principle}\index{Mach principle} by Einstein\index{Einstein, Albert (1879--1955)} (1918).} cogently argued that only {\it asymmetric laws}\index{asymmetric law} can distinguish left and right (1886). Here, a law is said to be asymmetric (or {\it to violate parity}) if it is not invariant under transformations interchanging left and right (specifically, under the {\it parity transformation} $(x,y,z)\mapsto (-x,-y,-z)$). In this connection, the neophyte of physics might recollect {\it Fleming's left-hand rule}\index{Fleming's left-hand rule} about the direction of force acting on a current-flowing wire under an external magnetic field (Fig.~\!\ref{fig:magnet}). As the name suggests, the rule is described by using human hands, and hence seem to imply that the electromagnetic phenomenon is asymmetric. But that is not actually the case because of the nature of magnetic fields. In order to put it clearly, let us recall the {\it Lorentz force} $e\left[\boldsymbol{E} + (\boldsymbol{v} \times \boldsymbol{B})\right]$ experienced by a particle of charge $e$ moving with velocity $\boldsymbol{v}$ in the presence of an electric field $\boldsymbol{E}$ and a magnetic field $\boldsymbol{B}$. Here $\boldsymbol{B}$ is not an ordinary ({\it polar}) vector, but an {\it axial} vector gaining an additional sign flip under a reflection.\footnote{The distinction between ``polar" and ``axial" vectors was made by W. Voigt in 1896. Note that $polar \times polar=axial$ and $polar \times axial= polar$.} 

The nature of magnetic fields is more well explained by identifying $\boldsymbol{B}=(B_1,$ $B_2,B_3)$ with the differential 2-form $\omega=B_3dx_1\wedge dx_2+B_1dx_2\wedge dx_3+B_2dx_3\wedge dx_1$ since $\omega$ is ``handedness-invariant."  This identification is not altogether artificial, though, at first sight, differential forms\index{differential form} might look an elaborate mathematical fabrication, having nothing to do with the real world. To convince ourselves of this, let us employ differential forms\index{differential form} to express Maxwell's equations\index{Maxwell's equations}. 

The original form of Maxwell's equations\index{Maxwell's equations} in a vacuum is:
\begin{eqnarray}\label{eq:maxwell}
\begin{cases}
\epsilon_0 ~{\rm div}~\boldsymbol{E}(t,{\bf x}) = \rho(t,\boldsymbol{x})~~(\text{Gauss's\index{Gauss, Johann Carl Friedrich (1777--1855)} flux law}), & \\  \mu_0^{-1}{\rm rot}~\boldsymbol{B}(t,{\bf x})-\epsilon_0 \displaystyle\frac{\partial \boldsymbol{E}}{\partial t}  = \boldsymbol{i}(t,\boldsymbol{x}) ~~(\text{Amp\`{e}re-Maxwell law}), & \\
{\rm div}~\boldsymbol{B}(t,\boldsymbol{x})=0 ~~(\text{Absence of magnetic monopoles}), &\\ 
 {\rm rot}~\boldsymbol{E}(t,\boldsymbol{x}) +\displaystyle\frac{\partial \boldsymbol{B}}{\partial t}  =  0 ~~(\text{Faraday's law}), &
\end{cases}
\end{eqnarray}
where $\epsilon_0$ stands for the dielectric constant of a vacuum,\footnote{\label{diel}$\epsilon_0 =8.854187817\ldots\times 10^{-12}\ \text{F/m}$, $~\mu_0 =1.2566370614 \cdots \times 10^{-6} ~ \mathrm{H/m}$. Strictly speaking, these equations are Heaviside's version (1888/1889). Maxwell's\index{Maxwell, James Clerk (1831--1879)} original equations (consisting of 20 equations in 20 variables) are quite complicated.
It follows from (\ref{eq:maxwell}) that both $\boldsymbol{E}$ and $\boldsymbol{B}$ satisfiy the wave equation $\displaystyle\mu_0\epsilon_0\partial^2 \boldsymbol{f}/\partial t^2-\Delta \boldsymbol{f}=\boldsymbol{0}$ under the absence of $\rho$ and $\boldsymbol{i}$. Hence electromagnetic waves propagate at the speed $1/\sqrt{\mu_0\epsilon_0}=0.299792458\times 10^8~\!{\rm m}/{\rm s}$, which coincides with the speed of light. Maxwell\index{Maxwell, James Clerk (1831--1879)} thus predicted that light is an electromagnetic wave (1865). This was confirmed experimentally by H. Hertz in 1887.} and $\rho$ (resp. $\boldsymbol{i}$) is the density of electric charge (resp. the electric current density). Writing $\boldsymbol{E}$ $=\!(E_1,E_2,E_3)$ and lumping together the magnetic field and the electric field, define a 2-form $\Omega$ by setting $\Omega\!=\!B_3dx_1\wedge dx_2\!+\!B_1dx_2\wedge dx_3\!+\!B_2dx_3\wedge dx_1\!+\!E_1dx_1\wedge dt\!+\!E_2dx_2\wedge dt\!+\!E_3dx_3\wedge dt$. Writing $\boldsymbol{i}\!=\!(i_1,i_2,i_3)$ likewise, we put $\eta$ $=\!\mu_0 (i_1dx_1\!+\!i_2dx_2\!+\!i_3dx_3)\!-\!\epsilon_0^{-1}\rho dt$. Then (\ref{eq:maxwell}) are rewritten as $d\Omega=0,~d^{*}\Omega\!=\!\eta$, which is manifestly invariant under Lorentz transformations\index{Lorentz transformation}. Here $d^*$ is the {\it adjoint} of the exterior derivation $d$ with respect to the Minkowski metric\index{Minkowski, Hermann (1864--1909)}; specifically, 
{\small
\begin{eqnarray*}
d^{*}\omega&\!=\!&\Bigl(\frac{\partial f_3}{\partial x_2}-\frac{\partial f_2}{\partial x_3}-c^{-2}\frac{\partial g_1}{\partial t}\Bigr)dx_1+\Bigl(\frac{\partial f_1}{\partial x_3}-\frac{\partial f_3}{\partial x_1}-c^{-2}\frac{\partial g_2}{\partial t}\Bigr)dx_2\\
& &+\Bigl(\frac{\partial f_2}{\partial x_1}-\frac{\partial f_1}{\partial x_2}-c^{-2}\frac{\partial g_3}{\partial t}\Bigr)dx_3-\Bigl(\frac{\partial g_1}{\partial x_1}+\frac{\partial g_2}{\partial x_2}+\frac{\partial g_3}{\partial x_3}\Bigr)dt
\end{eqnarray*}}
\hspace{-0.1cm}for $\omega\!=\!f_3dx_1\wedge dx_2\!+\!f_1dx_2\wedge dx_3\!+\!f_2dx_3\wedge dx_1\!+\!g_1dx_1\wedge dt\!+\!g_2dx_2\wedge dt\!+\!g_3dx_3\wedge dt$.

\medskip

It was unthinkable that anyone should question the validity of symmetry under a mirror reflection until a genuine asymmetric law\index{asymmetric law} was predicted by T.-D. Lee and C.-N. Yang in 1956 as the parity violation of the {\it weak interaction} (a force that governs all matter in the universe). Right after their announcement, the mind-boggling prediction was confirmed by C.-S. Wu and her collaborators through an experiment monitoring the beta decay of cobalt-60 atoms. Expressed in words, this implies that the nature at a very fundamental level can tell the characteristic difference between left-\index{left-handed} and right-handed\index{right-handed} (see \cite{gard}). 

In addition to the weak interaction, there are the three other fundamental interactions in nature; say, the electromagnetic force, the strong interaction, and gravitation. The strong interaction is, loosely put, the mechanism responsible for the strong nuclear force and is known to be ``P-symmetric", i.e., invariant under the parity transformation like the electromagnetic force. In this connection, a few more words about symmetry in quantum physics will not be out of place; 
besides the P-symmetry, we have the ``C-symmetry" under the {\it charge conjugation} that reverses the electric charge and all the internal quantum numbers, and the ``T-symmetry" under the {\it time reversal} replacing $t$ by $-t$. Each of these symmetries can be violated individually. Theoretically, however, there exists no physical phenomenon that violates the ``CPT-symmetry," the combination of all three symmetries. This fact is called the ``CPT theorem."

\section{Conclusion}

In retrospect, the theorem of angle-sum\index{theorem of angle-sum} and the Pythagorean Theorem\index{Pythagorean Theorem} in classical antiquity were the fresh impetus in the far-reaching feats by Gauss\index{Gauss, Johann Carl Friedrich (1777--1855)} and Riemann\index{Riemann, Georg Friedrich Bernhard (1826--1866)}. The triangle inequality\index{triangle inequality} as well was a fundamental source of the theory of topological spaces\index{topological space}, with which we can get a better understanding of intrinsicness of curved spaces. Moreover, the effort to refine the ancient approach to ``infinitesimals" came to fruition as calculus by Newton\index{Newton, Issac (1642--1726)} and Leibniz\index{Leibniz, Gottfried Wilhelm von (1646--1716)}, which, combined with the idea of coordinate system (no little indebted to Descartes' method\index{Descartes, Ren\'{e} (1596--1650)}), provided us with a powerful tool to investigate our space. What we ought to remember in particular is Cantor's\index{Cantor, Georg (1845--1918)} set theory that not only encompasses all the necessary stuff for mathematizing cosmology, but also spurred us on to give a probing interpretation to infinity whose nature had been a perennial controversial issue passed from antiquity.

Sullivan \cite{sul} says felicitously in a general context, ``A history of mathematics is largely a history of discoveries which no longer exist as separate items, but are merged into some more modern generalization, these discoveries have not been forgotten or made valueless. They are not dead, but transmuted."

%Our universe is now recognized as a 3D Riemannian manifold at any instant of time (except for the early stages of the universe where we must take account of quantum effects).

%\begin{quote}
%{\it The most incomprehensible thing about the universe is that it is comprehensible.}\index{Einstein, Albert (1879--1955)} ({\it Physics and Reality}, 1936) 

%\end{quote}

\selectlanguage{english}

{\small

%\printindex

}

\end{document}